\documentclass[11pt]{amsart}

\usepackage{amsmath, amssymb, amscd, cancel, graphicx, soul, stmaryrd, chngcntr, comment, tikz}
\usepackage[dvipsnames,svgnames]{xcolor}
\usepackage{mathdots, enumitem}
\usetikzlibrary{decorations.markings}

\usepackage[all]{xy}
\usepackage{multirow, longtable, array}

\newtheorem{thm}{Theorem}[section]
\newtheorem{conj}[thm]{Conjecture}
\newtheorem{cor}[thm]{Corollary}
\newtheorem{lem}[thm]{Lemma}

\newtheorem{defin}[thm]{Definition}
\newtheorem{clm}[thm]{Claim}
\newtheorem{remark}[thm]{Remark}

\newtheorem{ctr}[thm]{Criterion}

\newtheorem{addendum}[thm]{Addendum}

\def\A{{\mathcal A}}
\def\cB{{\mathcal B}}

\def\N{{\mathcal N}}

\def\Gam{\Pi}

\def\cN{{\mathcal N}}

\def\bZ{{\mathbb Z}}

\def\del{{\partial}}

\def\mod{{\textup{mod} \;}}

\usepackage{mathtools}
\usepackage{hyperref} % Comment out this line if we don't want hyperref
\providecommand{\texorpdfstring}[2]{#1}

\newcommand{\Z}{\mathbb{Z}}
\newcommand{\id}{\mathbf{1}}
\newcommand{\sm}{\!\setminus\!}
\newcommand{\sse}{\subseteq}

\newcommand{\Mod}[1]{\ (\mathrm{mod}\ #1)}
\newcommand{\wt}[1]{{\bk{\bf #1}}}
\def\mm{{\text{min}}}
\def\MM{{\text{max}}}
\newcommand{\pdfmath}[2]{\texorpdfstring{$#1$}{#2}}

\counterwithout{table}{section}
\counterwithout{table}{subsection}
\counterwithout{table}{subsubsection}

\DeclarePairedDelimiter{\pth}{(}{)}
\DeclarePairedDelimiter{\bk}{[}{]}
\DeclarePairedDelimiter{\bc}{\{}{\}}
\DeclarePairedDelimiter{\abs}{\lvert}{\rvert}

\DeclarePairedDelimiter{\norm}{\lVert}{\rVert}

\begin{document}

\title[Genus zero simple knots]
{Genus zero simple knots}

\author[Joshua Evan Greene]{Joshua Evan Greene}

\address{Department of Mathematics, Boston College\\ Chestnut Hill, MA 02467}

\email{joshua.greene@bc.edu}

\author[Shiyu Liang]{Shiyu Liang}

\address{Department of Mathematics, University of Texas at Austin\\ Austin, TX 78712}

\email{shiyu.liang@utexas.edu}

\author[John Luecke]{John Luecke}

\address{Department of Mathematics, University of Texas at Austin\\ Austin, TX 78712}

\email{luecke@math.utexas.edu}

%\author[Joshu Evan Greene,
%Shiyu Liang, John Luecke]{Joshua Greene, %Shiyu Liang, John Luecke}

%\address{Department of Mathematics, UT Austin\\ Austin, Texas 78712}

%\email{joshua.greene@bc.edu, 
%shiyu.liang@math.utexas.edu, luecke@math.utexas.edu}

\maketitle

%\[
%\begin{array}{rl}
%
%\text{email:}& {\tt joshua.greene@bc.edu} \cr
%
%\text{phone:}& 410-570-2782 \\

%\text{address:}& \text{Department of Mathematics} \\
%& \text{301 Carney Hall} \\
%& \text{Boston College} \\
%& \text{Chestnut Hill, MA 02467-3806}

%\end{array}
%\]

% \medskip

\noindent {\bf Abstract.} We classify the simple knots in lens spaces with planar rational Seifert surfaces. 
As a consequence, we classify the homology classes of all knots in lens spaces that admit planar rational Seifert surfaces. 
% not sure this really goes into the abstract
%For a knot in a lens space with irreducible exterior, having a planar rational Seifert surface is equivalent to the knot admitting a surgery containing a non-separating $2$-sphere, to the exterior of the knot being a Heegaard Floer homology solid torus, and to the exterior of the knot being a generalized solid torus.

% \noindent {\bf MSC classes:}  05C21, 05C50, 11H55, 57M25, 57M27

\section{Introduction}\label{sec:intro}

Simple knots form an important and well-studied class of knots in lens spaces.
Consider a standard genus-1 Heegaard splitting of a lens space $L(p,m)$ ($p > 1$) and take one meridian disk in each Heegaard solid torus, so that the boundaries of these two disks meet in $p$ points.
By definition, a knot is simple if it is isotopic to the union of a pair of arcs, one properly embedded in each meridian disk.
%they are obtained by taking the standard genus-1 Heegaard splitting of the lens space, with a meridian disk in each H
%the knots obtained by taking the union of a pair of arcs with the same endpoints, one properly embedded in a meridian disk in the standard genus-1 Heegaard 
%Let $L(p,m)$ denote a lens space.
%By convention, we assume that $p > 1$; in particular, $S^3$ and $S^1 \times S^2$ are not lens spaces.
%Choose meridian disks $D_0$, $D_1$ for the solid tori in a Heegaard splitting of $L(p,m)$ whose boundaries intersect in $p$ points.
%A {\bf simple knot} $K$ in $L(p,m)$ is the union of a pair of properly embedded arcs $\delta_0 \subset D_0$ and $\delta_1 \subset D_1$. See Figure~\ref{fig:simple} and section~\ref{sec:setup}.
%\begin{defin} Let $D_{\alpha}, D_{\beta}$ be meridian disks in a standard genus-one Heegaard splitting of a lens space $L(p,m)$. A {\bf simple knot} $K$ in $L(p,m)$ is the union of a pair of arcs, $\delta_{\alpha}, \delta_{\beta}$, properly in $D_{\alpha}, D_{\beta}$. See Figure~\ref{fig:simple} and section~\ref{sec:setup}.
%\end{defin}
Each non-trivial class of $H_1(L(p,m))$ contains a unique simple knot.
%For the purposes of this paper, $S^3$ and $S^1 \times S^2$ are not lens spaces. For the lens space $L(p,m)$, we take the convention that $p>1$.
Simple knots rose to prominence in connection with the Berge conjecture in Dehn surgery and in Floer homology, and they possess many pleasant properties \cite{berge2018surgeriesLensSpaces,greene2022fiberedSimpleKnots,niwu2014correctionTerms}.

%$\spadesuit$ {\em (I have been taking the unknot in the lens space as not being simple, e.g. simple knots have irreducible exteriors. Maybe this is not standard.) }

In this paper, we classify the simple knots that admit genus $0$ (i.e.~planar) rational Seifert surfaces. 
The property of possessing a planar Seifert surface is known to be equivalent to some properties pertaining to Dehn surgery and Floer homology, collected below in Theorem 
\ref{thm:generalizedsolidtorus}.
(See the beginning of section \ref{sec:setup} for the meanings of basic terms like Seifert surface, slope, distance, and Dehn surgery.)

%These knots appear in several interesting contexts as described by the following.

%\newcounter{equivalence}
%\setcounter{equivalence}{\value{thm}}
%\begin{thm}\label{thm:generalizedsolidtorus} Let $K$ be a knot in a lens space. Assume the exterior of $K$, $X_K$, is irreducible. The following are equivalent
%\begin{enumerate}
    %\item $K$ admits a surgery to $(S^1 \times S^2) \# N$ where $N$ is a rational homology sphere.
    %\item $X_K$ is a generalized solid torus in the sense of \cite{rasmussenrasmussen2017floerSimpleManifolds}.
    %\item $K$ has a planar rational Seifert surface.
    %\item $X_K$ is a Heegaard Floer solid torus in the sense of Watson.
%\end{enumerate}
%\end{thm}

%That is, classifying simple knots in any of these situations classifies those simple knots in any other.

Assume that $Q$ is a planar Seifert surface for $K$, and let $\lambda$ denote the slope of the boundary components of $Q$ on the boundary of the knot exterior $X_K$.
Note that $K(\lambda)$, the result of Dehn surgery along $K$ with slope $\lambda$, contains a non-separating $2$-sphere gotten by capping off the components of $\partial Q$ by meridian disks of the surgery solid torus.
Thus, $K(\lambda)$ is diffeomorphic to  $(S^1 \times S^2) \# N$, where $N$ is a rational homology sphere (possibly $S^3$).

The classification of simple knots with planar Seifert surfaces falls into three cases: when the slope $\lambda$ is non-longitudinal (i.e.~at distance greater than one from the knot meridian $\mu_K$ in $\del X_K$); when the slope $\lambda$ is longitudinal, and the surgery $K(\lambda)$ is not prime; and when the slope $\lambda$ is longitudinal, and the surgery $K(\lambda)$ is $S^1 \times S^2$.
These three cases are treated respectively in Theorems \ref{thm:main:nonlongitudinal}, \ref{thm:mainNextended}, and \ref{thm:Gr-Li-Lu}.

In the first case, the Cyclic Surgery Theorem \cite{cullergordonlueckeshalen1987dehnSurgery} leads to the following classification (e.g.~see \cite{darcysumners2000rationalTangleDistances}).
To state it, let $\mu \in H_1(L(p,m))$ denote the class of the core of one of the Heegaard solid tori of $L(p,m)$ (compare Definition \ref{def:mu} below); it should not be confused with the knot meridian $\mu_K \subset \del X_K$.

\newcounter{maindelta>1}
\setcounter{maindelta>1}{\value{thm}}
\newcounter{secdelta>1}
\setcounter{secdelta>1}{\value{section}}

\begin{thm}
[Non-longitudinal theorem]
\label{thm:main:nonlongitudinal}
Let $K$ be a knot in $L(p,m)$ that admits a planar Seifert surface with $q>1$ boundary components and non-longitudinal slope $\lambda$.
%, with $\Delta(\lambda, \mu_K)>1$.
Then 
\begin{enumerate}
    \item $K(\lambda)=S^1 \times S^2$;
    \item $L(p,m)$ is diffeomorphic, possibly orientation-reversing, to $L(k q^2, 1+kdq)$, where $\gcd (d,q)=1$, $q>d>0$, and $k=\Delta(\lambda, \mu_K)>1$;
    \item $K$ is the $(q,q)$-torus knot in $L(k q^2, 1+kdq)$;
    %(i.e. $K$ lies on the Heegaard torus for $L(k q^2, 1+kdq)$ with winding number $q$ on each side);
    \item $\pm [K]=q \mu \in H_1(L(k q^2, 1+kdq))$; and
    \item $K$ is simple.
    \end{enumerate}
\end{thm}

%$\spadesuit$ {\em (Should we add:
%Remark: When q=1, the simple knot K, which has irreducible exterior, is a core of L(p,m). The core curves of a lens space always admit surgeries to $S^1 \times S^2$ along slope $\lambda$ where $\Delta(\lambda, \mu_K)=p$.)
%}

By the $(q,q)$-torus knot in $L(k q^2, 1+kdq)$, we mean that $K$ lies on the Heegaard torus for $L(k q^2, 1+kdq)$ with winding number $q$ on each side. This property determines $K$. Note that Theorem~\ref{thm:main:nonlongitudinal} does not assume $K$ is simple but deduces that it must be.
Theorem~\ref{thm:realizationNlongitudinal} of section~\ref{sec:Appendix} gives a converse: the knots described in Theorem~\ref{thm:main:nonlongitudinal} do indeed have planar Seifert surfaces of non-longitudinal slope. 

The focus of the paper is when $Q$ has longitudinal slope and $K$ is simple. The following theorem captures the second case, when $K(\lambda)$ is not prime.

\newcounter{mainr>1}
\newcounter{secr>1}
\setcounter{mainr>1}{\value{thm}}
\setcounter{secr>1}{\value{section}}
\begin{thm}
[Longitudinal, non-prime theorem]
\label{thm:mainNextended}
Let $K \subset L(p,m)$ be a simple knot with a planar Seifert surface with $q>1$ boundary components and longitudinal slope $\lambda$, and suppose that $K(\lambda)$ is not prime.
Then there are integers $k>1$ and $0 < d < q$ such that up to diffeomorphism of $L(p,m)$, possibly orientation-reversing, $p=kq^2$, $m = kdq + 1$, and  $K(\lambda)=(S^1 \times S^2) \# L(k,1)$.
Moreover, either
\begin{enumerate}
\item
$\pm [K]=kdq \mu \in H_1(L(p,m))$ with $\gcd(q,d)=1$; or
\item
$\pm [K]=kq \mu \in H_1(L(p,m))$, and
$K$ is a $k$-fold cable on the $(q,q)$-torus knot in $L(p,m)$.
%a companion $L$ that lies on the Heegaard torus for $L(p,m)$ with winding number $q$ on each side. 
%a companion, $L$, that lies on the Heegaard torus of $L(p,m)$; furthermore $L$ has winding number $q$ in both Heegaard solid tori.
%In this case, .
\end{enumerate}     
%Let $K \subset L(p,m)$ be a simple knot with a planar Seifert surface with $q$ boundary components and longitudinal slope $\lambda$. Then either
%\begin{enumerate}
%    \item $K(\lambda)=S^1 \times S^2$; or
%    \item  up to diffeomorphism, possibly orientation-reversing, there are integers $k>1$ and $0 < d < q$ such that $K(\lambda)=(S^1 \times S^2) \# L(k,1)$, $p=kq^2$, and $m = kdq + 1$.
%    In this case, either
%    \begin{enumerate}
%    \item $\pm [K]=kdq \mu \in H_1(L(p,m))$ with %$\gcd(q,d)=1$; or
%    \item $K$ is an $k$-fold thm:mainNextended on a companion, $L$, that lies on the Heegaard torus of $L(p,m)$; furthermore $L$ has winding number $q$ in both Heegaard solid tori. In this case, $\pm [K]=kq \mu \in H_1(L(p,m))$.
%    \end{enumerate}
%\end{enumerate}
\end{thm}

Once again,
Theorem~\ref{thm:realizationNlongitudinal} shows that the knots described in Theorem~\ref{thm:mainNextended} do indeed have planar Seifert surfaces of longitudinal slope and Dehn surgeries along that slope as described. We note that Theorem~\ref{thm:mainNextended} resolves Question 16 of \cite{bakercabling}, classifying the doubly primitive knots in $(S^1 \times S^2) \# N$ when $N \neq S^3$.

More than half of this paper is then involved in the third case, which is to classify those simple knots $K \subset L(p,m)$ with planar Seifert surfaces for which $K(\lambda)=S^1 \times S^2$. Note that attaching a 2-handle to $L(p,m)$ along $K$, a cancelling $3$-handle, and a 4-handle shows that $L(p,m)$ bounds a rational homology ball. Lisca classified all such lens spaces, thereby giving restrictions on the lens spaces $L(p,m)$ containing genus zero knots \cite[Corollary 1.3]{lisca2007lensSpacesRibbonConjecture}.
On the other hand, 
%Berge and Gabai \cite{berge1991knotsSolidTorus, gabai1990bridgeBraids} classified the knots in a solid torus upon which longitudinal surgery gives back a solid torus. Let $K'$ be such a knot, a {\em Berge-Gabai} knot, and identify the solid torus in which it sits with a genus-one Heegaard torus of $S^1 \times S^2$. Then $K'$ becomes a knot in $S^1 \times S^2$ upon which the corresponding longitudinal surgery produces a lens space. 
J. Rasmussen discovered that the lens spaces on  Lisca's list all have knots $K$ that admit longitudinal surgeries to $S^1 \times S^2$ via embeddings of Berge-Gabai knots (see \cite[Section 1.5]{greene2013lensSpaceRealization}). Together these give  
\begin{thm}\label{thm:lisca-rasm}(\cite{lisca2007lensSpacesRibbonConjecture}, Rassmussen via \cite{greene2013lensSpaceRealization}, see also \cite[Theorem 1.2]{baker2016knotsS1xS2})
Let $M$ be a lens space. $M$ contains a knot that admits longitudinal surgery to $S^1 \times S^2$ if and only if up to diffeomorphism, possibly orientation-reversing, there are $d,q \in \Z$ such that $M$ is $L(q^2,m)$ where
\begin{itemize}
\item $m=dq+1$ with $\gcd(d,q)=1,2$, or
\item $m=d(q-1)$ where $d$ is odd, and $d$ divides $q-1$ or $2q+1$.
\end{itemize}
\end{thm}

%$\spadesuit$ {\em (if we go with this, need to change the statement of 6.1 correspondingly.)}

 The correspondence between the knots described by Rasmussen and the lens spaces determined by Lisca, raised the question of whether these might be the only knots in lens spaces surgering longitudinally to $S^1 \times S^2$. Baker, Buck, and Lecuona, however, found other examples. In \cite{baker2016knotsS1xS2}, they list all the knots currently known to surger longitudinally to $S^1 \times S^2$. All these knots are simple, so each is determined by the corresponding homology class. 

 \begin{thm}[\cite{baker2016knotsS1xS2}]\label{thm:BBL}
Let $M=L(q^2,m)$ be a lens space of Theorem~\ref{thm:lisca-rasm}, and $\mu$ be the canonical generator of $H_1(M)$.
Then $S^1 \times S^2$ can be obtained by longitudinal surgery on a simple knot $K$ in $M$ with $\pm [K]=\kappa \mu \in H_1(L(q^2,m))$ where $\kappa$ is as listed below: 
\begin{enumerate}
\item If $m=dq+1$ such that $\gcd(d, q) = 1$, then $\kappa =q,dq$; %(bgi, gofk);
\item If $m=dq+1$ such that $\gcd(d, q) = 2$, then $\kappa =q$;
%? \kappa = ±q so that K is the dual to a bgii knot;
\item If $m=d(q-1)$ such that $d$ divides $2q+1$, then $\kappa =q,dq$;
or
\item If $m=d(q-1)$ such that $d$ is odd and divides $q-1$, then $\kappa =q, 2q, dq$. Furthermore, the case with $\kappa=2q$ only occurs when $q=1-2d$. 
% \kappa=±dm so that K is the dual to a bgiii knot,
%? \kappa = ±m so that K is the dual to a bgv knot, or
%? \kappa=±2m=±4m so that K is the dual to a spor knot if m=1?2d; 
\end{enumerate}
\end{thm}

If a knot in a lens space surgers to $S^1 \times S^2$, then it must have a planar Seifert surface (see Theorem~\ref{cor:Korder}, or Theorem~\ref{thm:generalizedsolidtorus} below). Let $M$ be a lens space as in Theorem~\ref{thm:lisca-rasm}. By classifying those simple knots in $M$ with a planar Seifert surface of longitudinal slope, we show that this list is complete.

\newcounter{mainthm}
\setcounter{mainthm}{\value{thm}}
\newcounter{mainsec}
\setcounter{mainsec}{\value{section}}
\begin{thm}
[Spherical Theorem]
\label{thm:Gr-Li-Lu}
Let $K$ be a simple knot in a lens space that admits a longitudinal surgery to $S^1 \times S^2$. Then $K$ is listed in Theorem~\ref{thm:BBL}. 
\end{thm}

\begin{remark}\label{rem:cebanuwork}
Cebanu, in his thesis \cite[Theorem 4.8.5]{cebanu2012propertyR}, proved Theorem~\ref{thm:Gr-Li-Lu} for simple knots in the class of lens spaces $L(q^2, dq+1)$ (as in $(1),(2)$ of Theorem~\ref{thm:BBL}). In personal communication (see \cite{baker2016knotsS1xS2}), he claimed to have subsequently proven  Theorem~\ref{thm:Gr-Li-Lu} for the remaining classes, where $m=d(q-1)$ (as in $(3),(4)$  of Theorem~\ref{thm:BBL}). To our knowledge, this extension to the remaining classes has never appeared. The proof given here uses different techniques (e.g. Morse-theoretic) which may be of independent interest.
\end{remark}

The cores of the solid tori in a genus-1 splitting of a lens space (which are unique up to isotopy) are exactly those knots in the lens space with irreducible exterior and with a disk for a planar Seifert surface. These disks will have non-longitudinal slope (as $S^3$, $S^1 \times S^2$ are excluded) and these knots are simple. The surgery along the slope of the disk yields $S^1 \times S^2$. The knots with planar Seifert surfaces which are not disks and the lens spaces in which they sit are much rarer.
Theorems~\ref{thm:main:nonlongitudinal}, \ref{thm:mainNextended}, and \ref{thm:Gr-Li-Lu} give a complete list of simple knots in lens spaces with planar Seifert surfaces that are not disks ($q>1$). We point out that the knots that surger to $S^1 \times S^2$ fiber over the circle (with fiber the planar Seifert surface), see Theorem~\ref{thm:cebanufiber} or Theorem 1.1 of \cite{nivafaee2019nullSurgery}. Those that surger to $(S^1 \times S^2) \# N$, where $N \ne S^3$, do not fiber over the circle.

%The classification of all surgeries from lens spaces to $S^1 \times S^2$ (and vice versa) would follow from the Theorems~\ref{thm:main:nonlongitudinal} and \ref{thm:Gr-Li-Lu} and following:

%\begin{conj}[\cite{baker2016knotsS1xS2}, \cite{greene2013lensSpaceRealization}]
%If a knot in a lens space admits a longitudinal surgery to $S^1 \times S^2$, then the knot is simple.
%\end{conj}

For a knot in a lens space, the existence of a planar Seifert surface is equivalent to conditions pertaining to Dehn surgery and Floer homology.

\newcounter{equivalence}
\setcounter{equivalence}{\value{thm}}
\begin{thm}\label{thm:generalizedsolidtorus} Let $K$ be a knot in a lens space with irreducible exterior $X_K$.
The following are equivalent:
\begin{enumerate}
    \item $K$ admits a surgery to $(S^1 \times S^2) \# N$, where $N$ is a rational homology sphere.
    \item $X_K$ is a generalized solid torus \cite[Definition 7.2] {rasmussenrasmussen2017floerSimpleManifolds}.
    \item $K$ has a planar Seifert surface.
    \item $X_K$ is a Heegaard Floer homology solid torus \cite[Definition 26]{hanselmanrasmussenwatson2022torusBoundary}.
\end{enumerate}
\end{thm}

Together, the results above give the following.

%Watson defines the notion of a {\it Heegaard Floer homology solid torus} for a rational homology solid torus, $M$, which Alishahi-Lidman-Lipshitz \cite{alisha2025detectingHeegaardFloerSolidTori} show is equivalent to a Dehn filling of $M$ containing a non-separating 2-sphere. Thus for the exteriors of knots in lens spaces, the notion of a generalized solid torus and a Floer homology solid torus are the same.  We obtain the following characterization 

\begin{cor}\label{cor:classify1-4} Let $K$ be a simple knot in a lens space with irreducible exterior. $K$ satisfies any of the equivalent conditions $(1)$--$(4)$ of Theorem~\ref{thm:generalizedsolidtorus} if and only if $K$ is either
\begin{itemize}
    \item the core of a Heegaard solid torus; or
    \item one of the knots listed in Theorems~\ref{thm:main:nonlongitudinal}, \ref{thm:mainNextended}, or \ref{thm:BBL}.
\end{itemize}
\end{cor}

 Theorem 27 of \cite{hanselmanrasmussenwatson2022torusBoundary} shows that if one filling of a rational homology solid torus
$M$ is an L-space, then $M$ is a Heegaard Floer homology solid torus if and only if all fillings except along the rational
longitude are L-spaces. Thus the knots listed in Theorems~\ref{thm:main:nonlongitudinal}, \ref{thm:mainNextended}, and \ref{thm:BBL} are exactly the simple knots in lens spaces upon which all Dehn surgeries are L-spaces except for that along the rational longitude. 

A simple knot realizes the smallest genus Seifert surface among all knots in its homology class \cite{niwu2014correctionTerms}. Thus if any knot in a lens space has a planar Seifert surface, then the simple knot in its class must also. 

\newcounter{KmeansK'also}
\setcounter{KmeansK'also}{\value{thm}}
\newcounter{secKmeansK'also}
\setcounter{secKmeansK'also}{\value{section}}

\begin{thm}\label{thm:r>0reduciblemeansreductiblesimple}
Let $K$ be a knot in $L(p,m)$ with irreducible exterior. Assume that, for some slope $\lambda$, $K(\lambda)=(S^1 \times S^2) \# N$, where $N$ is a rational homology sphere. Then $K$ is a Floer simple knot and has a planar Seifert surface $Q$ of slope $\lambda$. Furthermore, $p=q^2 |H_1(N)| \Delta(\lambda,\mu_K)$, where $q=|\partial Q|$.

Let $K'$ be the simple knot in the same homology class as $K$. 
Then the following hold:
\begin{enumerate}
    \item $K'$ has a planar Seifert surface, $Q'$ of some slope $\lambda'$;
    \item $|\partial Q'|=q$ and $\Delta(\lambda, \mu_{K})=\Delta(\lambda', \mu_{K'})$; and
    \item $K(\lambda')=(S^1 \times S^2) \# N'$ for a lens space $N'$ with $|H_1(N')|=|H_1(N)|$.
\end{enumerate}
\end{thm}

%$\spadesuit$ {\em (Probably drop the following Corollary.)} In the case when $K$ is simple, we have the following corollary.

%\begin{cor}\label{cor:redsurganygenus}
  %  Let $K\subset L(p, m)$ be a simple knot. Assume that, for some $\lambda$, $K(\lambda)$ is reducible. Then either
   % \begin{enumerate}
       % \item $K(\lambda) = S^1\times S^2$;
       % \item $K(\lambda) = (S^1\times S^2)\# L_1$, for some lens space $L_1$; or
       % \item $K(\lambda) = L_2\# L_3$, for some lens spaces $L_2, L_3$.
    %\end{enumerate}
  %  In the cases of $(1)$ and $(2)$, $K$ has a planar Seifert surface.
%\end{cor}

%\begin{proof}
   % As $K$ is simple, it has tunnel number $1$, hence any filling of $X_K$, in particular $K(\lambda)$, has Heegaard genus at most $2$. By the additivity of Heegaard genus under connect sum and the reducibility of $K(\lambda)$, it has at most two summands, and each has genus $1$. As $\rk(H_1(X_K)) = 1$, $K(\lambda)\neq (S^1\times S^2)\# (S^1\times S^2)$. Thus, one of $(1)$--$(3)$ holds. If $(1)$ or $(2)$ holds, then Theorem~\ref{thm:r>0reduciblemeansreductiblesimple} implies $K$ has a planar Seifert surface.
%\end{proof}

A knot $K$ in a rational homology sphere $M$ is called {\em Floer simple} if its knot Floer homology is the simplest possible, that is $\widehat{HFK}(K) = \mathbb{Z}^{|H_1(M)|}$. For example, a simple knot in a lens space is Floer simple. 

Theorem~\ref{thm:r>0reduciblemeansreductiblesimple} is proven in section~\ref{sec:homologydetection}. We now use it to show that the classification of simple knots with planar Seifert surfaces also specifies the homology classes (and the lens spaces in which these classes lie) of {\it all} knots in a lens space with planar Seifert surfaces.

\begin{cor}\label{cor:homologyclasses} Let $K \subset L(p,m)$ have irreducible exterior. If $K$ satisfies any of the equivalent conditions $(1)$--$(4)$ of Theorem~\ref{thm:generalizedsolidtorus} then one of the following holds
\begin{itemize}
    \item $K$ is the core of a Heegaard solid torus.
    \item $L(p,m)$ and $[K] \in H_1(L(p,m))$ are as in Theorems~\ref{thm:BBL}, \ref{thm:mainNextended}, or \ref{thm:main:nonlongitudinal}.
\end{itemize}
In particular, let $\mu_K$ be the meridian of $K$. If $\lambda$ is a slope such that $K(\lambda)=(S^1 \times S^2) \# N$ 
where $N$ is a rational homology sphere 
then one of the following holds
\begin{enumerate}[label=(\roman*)] 
    \item $K$ is the core of a Heegaard solid torus, and $K(\lambda)=S^1 \times S^2$.
    \item $L(p,m)$ and $[K]$ are as in Theorem~\ref{thm:BBL}, $\Delta(\lambda, \mu_K)=1$, and $K(\lambda)=S^1 \times S^2$. 
    \item $L(p,m)$ and $[K]$ are as in Theorem~\ref{thm:mainNextended}, $\Delta(\lambda, \mu_K)=1$, where $|H_1(N)|>1$. 
    \item $L(p,m)$ and $K$ are as in Theorem~\ref{thm:main:nonlongitudinal}, $\Delta(\lambda, \mu_K)>1$.
\end{enumerate}
\end{cor}

%\begin{thm}\label{thm:homologyclassification}
%Let $K$ be a knot in a lens space that admits a longitudinal Dehn surgery to $S^1 \times S^2$. Let $K^\ast$ be the simple knot in the same homology class as $K$. Then $K^\ast$ admits a longitudinal Dehn surgery to $S^1 \times S^2$. In particular, 
%the homology class of $K$ must be as listed in Theorem~\ref{thm:BBL}.
%\end{thm}

\begin{proof}
Let $K \subset L(p,m)$ have irreducible exterior. We show that if there is a slope $\lambda$ such that $K(\lambda)=(S^1 \times S^2) \# N$ then $(i)-(iv)$ holds. The first equivalence stated then follows from this and Theorem~\ref{thm:generalizedsolidtorus}.

%Clearly, if $(i)-(iv)$, holds then $K(\lambda)=(S^1 \times S^2) \# N$ for slope $\lambda$.

Assume $K(\lambda)=(S^1 \times S^2) \# N$. 
%As $K$ is a knot in a lens space with irreducible exterior, $N$ is necessarily a rational homology sphere. 
By Theorem~\ref{thm:r>0reduciblemeansreductiblesimple}, $K$ has a planar Seifert surface, $Q$, of slope $\lambda$. Let $q=|\partial Q|$.

If $q=1$, then by the irreducibility of $X_K$ we see that $X_K$ is a solid torus. This implies that $K$ is isotopic to a core of a Heegaard solid torus of $L(p,m)$. $(i)$ now follows. We hereafter assume $q>1$.

If $\Delta(\lambda,\mu_K)>1$, Theorem~\ref{thm:main:nonlongitudinal} says that $(iv)$ holds.

We then assume $\Delta(\lambda,\mu_K)=1$ ($\lambda$ is longitudinal). Then $[K]$ has non-trivial torsion $q$ in $L(p,m)$; consequently, $[K] \neq 0 \in H_1(L(p,m))$.
Let $K'$ be the simple knot with $[K']=[K]$, and let $\mu_{K'}$ be its meridian. By Theorem~\ref{thm:r>0reduciblemeansreductiblesimple}, $K'$, has a planar Seifert surface $Q'$ of slope $\lambda'$ and with $q'$ boundary components. As $\Delta(\lambda', \mu_{K'})=\Delta(\lambda, \mu_{K})=1$, $\lambda'$ is longitudinal. Theorem~\ref{thm:r>0reduciblemeansreductiblesimple}$(3)$ says that $K'(\lambda')=(S^1 \times S^2) \# N'$ where $N'$ is a lens space of order $|H_1(N)|$. It also then shows that both $K$ and $K'$ are Floer simple.  

If $N'$ is $S^3$, then Theorem~\ref{thm:cebanufiber} (see also Theorem 1.1 of \cite{nivafaee2019nullSurgery}) says that $K'$ fibers over the circle with fiber $Q'$ (the fiber is unique up to isotopy). Corollary 5.3 of \cite{niwu2014correctionTerms} then implies that $K$, in the same homology class as $K'$ and both Floer simple, fibers over the circle with $Q$ as its fiber. Thus $K(\lambda)=S^1 \times S^2$ as claimed. Theorem~\ref{thm:Gr-Li-Lu} applied to $K'$ says that $(ii)$ holds ($[K]=[K'])$, as claimed. 

If $N' \ne S^3$, then Theorem~\ref{thm:main:nonlongitudinal} says that $(iii)$ holds, noting that $|H_1(N)|=|H_1(N')|>1$. 
\end{proof}

A simple knot in a lens space is Floer simple. A conjecture of Baker-Grigsby-Hedden says that the converse is true.

\begin{conj}\label{conj:floersimple}(\cite{baker2008gridDiagrams}, Conjecture 1.5). If $K$ is a knot in a lens space that is Floer simple, then $K$ is simple.
\end{conj}

%$\spadesuit$ {\em (Need to check that K Floer simple is equivalent to the exterior being Floer simple in this context (\cite{rasmussenrasmussen2017floerSimpleManifolds}) then either express this in terms of the exterior being simple, or define Floer simple knot.) }

Note that if a knot $K$ satisfies $(1)$--$(4)$ of Theorem~\ref{thm:generalizedsolidtorus}, then Theorem~\ref{thm:r>0reduciblemeansreductiblesimple} says that $K$ is Floer simple.
Thus Conjecture~\ref{conj:floersimple} would imply that $K$ is simple and hence listed in Corollary~\ref{cor:classify1-4}.

In the direction of Conjecture~\ref{conj:floersimple}, we note the following result of Hedden.

\begin{thm}\label{thm:hedden}(\cite{hedden2011floerBerge}, Proposition 3.3)
A $1$-bridge knot in a lens space is Floer simple if and only if it is simple.
\end{thm}

We also point out the following result of Li, Moriah, and Pinsky. 

\begin{thm}\label{thm:lmp}
(\cite{lpm2025}, Theorem 7) 
Let $K \subset M$ be a tunnel number one knot with irreducible knot exterior, where M is either $S^3$ or $(S^1 \times S^2)\#L(r, s)$ (where $L(r, s)$ is any lens space). If a non-trivial Dehn surgery on $K$ yields a lens space, then $K$ has a doubly primitive presentation.
\end{thm}

Ken Baker shows (private communication) that in Theorem~\ref{thm:lmp} the dual knot in the lens space is simple.

Towards the classification of genus zero knots in lens spaces, we note that Karatzas \cite{karatzas2026} has recently shown that if a knot $K$ has a planar rational Seifert surface with at most four boundary components and lies in a lens space without Klein bottles, then K is not hyperbolic.  

\subsection{Outline of the paper.}\label{subsec:pfandoutline}

Here is an outline of the paper. 

%In sections~\ref{sec:setup} and \ref{sec:atleast3},  we assume $K$ is a simple knot in a lens space with a longitudinal, planar, rational Seifert surface. In sections~\ref{sec:estimates}, \ref{sec:dq+1}, and \ref{sec:d(q-ep)}, we specialize to the case where $K$ surgers to $S^1 \times S^2$ and verify the supporting material in the proof of Theorem~\ref{thm:Gr-Li-Lu}.

$K \subset L(p,m)$ is a simple knot with a planar, rational Seifert surface of slope $\lambda$. Except for sections~\ref{sec:homologydetection} --\ref{sec:Appendix}, $\lambda$ has longitudinal slope. Note that sections~\ref{sec:alphabetawts}--\ref{sec:alexwidth} focus on the case where $K(\lambda)=S^1 \times S^2$ (and the proof of Theorem~\ref{thm:Gr-Li-Lu}). 

Section~\ref{sec:setup} sets up a Morse presentation of the rational Seifert surface, $Q$, of $K$ within an interval of (twice punctured) Heegaard $2$-tori, $T$. We pick transverse levels $T \times \{z\}$ between each of the $2q-2$ saddle singularities of $Q$. We denote by $\Gamma_i$, $i=0, \dots, 2q-2$, these level surfaces along with their intersections with $Q$, and the union of these is called $\Gamma$. The parameters $y,x=q-y, y', x'=q-y'$ are defined here and several properties of the $\Gamma_i$ established. 

Section~\ref{sec:atleast3} shows that if at any level $\Gamma_i$, except $i=0,2q-2$, there are at most two isotopy classes of arcs, then $K$ is determined---and among those listed in Theorem~\ref{thm:BBL} and Theorem~\ref{thm:mainNextended}. This allows us to assume in the proof of Theorem~\ref{thm:Gr-Li-Lu} that there are always at least three isotopy classes of arcs in $\Gamma_i$ when $i \neq 0,2q-2$.

In section~\ref{sec:r>1}, using the main result of section~\ref{sec:atleast3}, we classify those simple knots $K$ having a planar rational Seifert surface with longitudinal slope $\lambda$ but for which $K(\lambda) \neq S^1 \times S^2$, Theorem~\ref{thm:mainNextended}.

Section~\ref{sec:alphabetawts} defines the $\alpha,\,\beta$-weights of a level arc, and how a level arc with these weights decomposes the homology class of $K$ (Lemma~\ref{lem:decomposeK}). This is the structure that will be used in sections~\ref{sec:alphabetawts} through \ref{sec:alexwidth} to prove Theorem~\ref{thm:Gr-Li-Lu}, as outlined at the beginning of section~\ref{sec:alphabetawts}.
We endow the arcs of intersection of $\Gamma_i$, $i=0, \dots, 2q-2$, with these $\alpha$-weights, and define the total ($\alpha$-) weight of a $\Gamma_i$. Note that the indexing of $\Gamma$ is from from $V_{\alpha}$ to $V_{\beta}$. We reindex the level surfaces of $\Gamma$ from $V_{\beta}$ to $V_{\alpha}$, defining $\Gamma_i'=\Gamma_{2q-2-i}$, $i=0, \dots, 2q-2$. $\Gamma'$ is the union of the $\Gamma_i'$. We endow the arcs of intersection of $\Gamma_i'$ with their $\beta$-weights, and define the total ($\beta$-) weight of each $\Gamma_i'$. Note that $\Gamma'$ is basically $\Gamma$ turned upside down (with the roles of $V_{\alpha}$ and $V_{\beta}$ switched), so conclusions about $\Gamma$ will have analogs in $\Gamma'$.
%We call the level surfaces with their arcs of intersection $\Gamma_i$, $i=0, \dots, 2q-2$, where the constituent arcs are weighted by their $\alpha$-weights  (described by the total weight of $\Gamma_i$) and the $\Gamma_i$ are indexed going from $V_{\alpha}$ to $V_{\beta}$ (so $\Gamma_0$ belongs to $T \times \{0\}$, $\Gamma_{2q-2}$ belongs to $T \times \{1\}$). The collection of the $\Gamma_i$ is referred to as $\Gamma$.  We call these same level surfaces and their arcs $\Gamma_i'$, $i=0, \dots, 2q-2$, when the constituent arcs are weighted by their $\beta$-weights and the $\Gamma_i'$ are indexed going from $V_{\beta}$ to $V_{\alpha}$ (so $\Gamma_i'=\Gamma_{2q-2-i}$ with $\alpha$-weights changed to $\beta$-weights). 
Each arc in $\Gamma_i$ runs between distinct boundary components, so there are at most four distinct isotopy classes of arcs in $\Gamma_i$, and by section~\ref{sec:atleast3} there are at least three. Theorem~\ref{lem:b'>b} describes how the $\alpha$-weights change as one moves between levels of $\Gamma$ with three isotopy classes of arcs (analogously for $\Gamma'$). Lemma~\ref{thm:m:b<x} and Theorem~\ref{thm:=x}  describe the weights of the early level surfaces ($\Gamma_i$ when $i \leq x$). The analogous results hold for $\Gamma'$. Theorem~\ref{thm:=x} introduces the useful constants $r,\, s,\, e,\, r',\, s',\, e'$.

Section~\ref{subsec:specialize} gives the proof of the Spherical Theorem, Theorem~\ref{thm:Gr-Li-Lu}, using the results of  sections~\ref{sec:estimates}, \ref{sec:dq+1}, \ref{sec:d(q-ep)}, and \ref{sec:alexwidth}. Here we set the context for these next four sections, including Theorem~\ref{thm:lisca} which is a reformulation of Theorem~\ref{thm:lisca-rasm}. 

Section~\ref{sec:estimates} provides estimates on the weights of the arcs of $\Gamma_b$, as $b$ increases. The main results are the Initial Growth Theorem~\ref{thm:M:b<x+1} and the Engine Theorem~\ref{thm:q+y,m_b=1}. In subsection~\ref{subsec:b<q}, we prove Theorem~\ref{thm:M:b<x+1} bounding the maximal weight at levels $\Gamma_{b}$ when $b$ is still less than half-way through the $2q-2$ levels. In subsection~\ref{subsec:engine}, the main result is Theorem~\ref{thm:q+y,m_b=1} which gives estimates on the weights for levels more than half-way through the $2q-2$ levels. For a first pass, we note that the proof of Theorem~\ref{thm:q<3y,b>q+2e} will give the reader a good sense of the type of arguments used to prove the Engine Theorem~\ref{thm:q+y,m_b=1}. By symmetry, the estimates of section~\ref{sec:estimates} apply as well to $\Gamma'$  and its $\beta$-weighting.

We consider the three different families of lens spaces given by Theorem~\ref{thm:lisca}. In section~\ref{sec:dq+1}, we prove that $\pm [K]=q \mu, dq \mu$ when $K$ is in a lens space of the form $L(q^2,m)$ where $m=dq+1$. In section~\ref{sec:d(q-ep)}, we consider the lens spaces of the form $L(q^2,m)$ where $m=d(q \pm 1)$. When $d \, | \, 2q \mp 1$, we show that $\pm [K]=q \mu, dq \mu$. When $d \, | \, q \pm 1$, we show that $\pm [K]=q \mu, dq \mu, 2q \mu$. In both sections~\ref{sec:dq+1} and \ref{sec:d(q-ep)}, we apply Corollary~\ref{cor:Mmdecom} in conjunction with Theorems~\ref{thm:q+y,m_b=1} and \ref{thm:M:b<x+1}
to write $\pm[K]=c \mu + f \mu'$ where $c,f$ are small, then apply the numerics of the lens space to conclude that $K$ is in the target homology classes. We note that the shorter section~\ref{sec:dq+1}, will give the reader a good sense of how the argument works in both sections~\ref{sec:dq+1} and \ref{sec:d(q-ep)}.

Section~\ref{sec:alexwidth} uses Alexander polynomial calculations to prove that when $m=d(q \pm 1)$ and $d \, | \, q \pm 1$, then $\pm [K]=2q \mu$ only when $q=1-2d$. 

In section~\ref{sec:homologydetection}, we prove Theorem~\ref{thm:generalizedsolidtorus}, which gives several equivalent notions to a knot in a lens space having a planar rational Seifert surface. We then prove Theorem~\ref{thm:r>0reduciblemeansreductiblesimple} showing that if a knot satisfies any of these notions, then so does the simple knot in the same homology class.

In section~\ref{sec:nonlongitudinal}, we classify the knots admitting non-longitudinal planar Seifert surfaces, Theorem~\ref{thm:main:nonlongitudinal}.

Section~\ref{sec:Appendix} realizes the knots described in  
Theorems~\ref{thm:main:nonlongitudinal} and \ref{thm:mainNextended}
as knots with planar, rational Seifert surfaces and with the corresponding surgeries along the their rational longitude.
\bigskip

\section*{Acknowledgments.}
We thank Cameron Gordon and Mark Brittenham for helpful conversations. We thank Tye Lidman and Jake Rasmussen  for educating us about Heegaard Floer homology solid tori, and in particular Theorem~\ref{thm:generalizedsolidtorus}.
We thank Ken Baker for many helpful conversations during the course of this work. In particular, we thank him for showing us the examples of Theorem~\ref{thm:realizationNlongitudinal}, that inspired Theorem~\ref{thm:mainNextended}, and for the many discussions about them. 

\section*{AI Methodology.}
No AI was used in this preparation of this paper.
As we were making final edits, we asked AI to prove Theorem 1.5 using techniques similar to those in Section 10.
It generated a proof which had similar length and complexity to ours, but without any geometric intuition.
SL intends to clarify and extend this argument in a follow-up paper.

\newpage
\tableofcontents
\newpage

\section{Setup and definitions}\label{sec:setup}% Definition of $\Gamma,\, \Gamma'$, $\alpha, \,\beta$-weights, $r, \,s, \,e$.

%Let $K$ be a knot in the lens space $L(p,m)$, and let $X_K$ denote its exterior.
%A {\em (rational) Seifert surface} for $K$ is a connected, orientable surface $Q$ properly embedded in $X_K$ whose boundary $\partial Q$ consists of a non-empty collection of coherently oriented curves on $\partial X_K$ (under an orientation from $Q$).
%\begin{defin}\label{def:distance}
%and is referred to as the rational longitude of $K$. 
%\end{defin}

%When $\lambda$ is a slope on $\partial X_K$, we denote by $K(\lambda)$ the Dehn surgery of $K$ along $\lambda$.

We begin with a few conventions.
For a knot $K$ in a manifold $\mathcal{M}$, let $X_K$ denote the knot exterior, i.e.~ the complement in $\mathcal{M}$ of an open tubular neighborhood of $K$.
A {\it slope} on $X_K$ is an isotopy class of essential simple closed curve on $\partial X_K$.
The {\it distance} $\Delta(\alpha,\beta)$ between two slopes $\alpha, \beta$ is defined to be their geometric intersection on $\partial X$.
A slope $\lambda$ on $\partial X_K$ is {\em longitudinal} if $\Delta(\lambda,\mu_K)=1$, where $\mu_K$ is the meridian of $K$.
Denote by $K(\alpha)$ the Dehn surgery on $K$ with slope $\alpha$.
A {\it (rational) Seifert surface} for $K$ is a properly embedded, connected surface $Q \subset X_K$, such that an orientation on the surface induces coherent orientations of $\partial Q$ on $\partial X_K$. 
The {\em slope} of a Seifert surface $Q$ of $K$ is the slope of any component of $\partial Q$.

When $\mathcal{M}$ is a solid torus or $S^1 \times S^2$, the {\em winding of number} of $K$ in $\mathcal{M}$, $\omega$, is the positive integer such that the image under inclusion $H_1(K) \to H_1(\mathcal{M})$ is $\omega \mathbb{Z}$.
Let $\mathcal{M}$ be a solid torus with meridian $\underline{m}$ and a choice of longitude $\underline{l}$. Let $c_{p,q}$ be the curve on $\partial \mathcal{M}$ given by the class $q \underline{m}+p \underline{l} \in H_1(\partial \mathcal{M})$. 
%Pushing $c_{p,q}$ into the interior of $\mathcal{M}$, we call $K_{p,q}$ the {\em $(p,q)$-cable knot} in $\mathcal{M}$. Note that $\abs{q}$ is the winding number of $K_{p,q}$ in $\mathcal{M}$. 
If $L$ is a framed knot in a manifold $\mathcal{R}$, then we can identify the solid torus $\mathcal{M}$ with a neighborhood of $L$ in such a way that $\underline{l}$ is identified with the framing curve of $L$. Under this identification, $c_{p,q}$ becomes a knot, $L_{p,q}$, in $\mathcal{R}$ called the {\em $(p,q)$-cable of L}. In this context $L$ is called the {\em companion} of $L_{p,q}$ in $\mathcal{R}$. Suppressing $p$, we refer to the $(p,q)$-cable of $L$ as a {\em $|p|$-fold cable} of $L$. 
%In the lens space $L(p,m)$ with Heegaard splitting $V_\alpha \cup_{\hat{T}} V_\beta$ (see below), the {\em $(r,s)$-torus knot} is the $(r,s)$-cable of the core of $V_\alpha$ with framing from $V_\alpha=D^2 \times S^1$.
%curve on the Heegaard torus with respect to its standard coordinates.  as the $(p,q)When $L$ is the unknot in $S^3$, $L_{p,q}$ is the $(p,q)$-torus knot. 

\bigskip

Beginning in this section and up to and including section~\ref{sec:alexwidth}, we assume $K$ is a simple knot in the lens space $L(p,m)$, $p>m>0$ (working up to diffeomorphism), with a planar rational Seifert surface $Q$ of longitudinal slope $\lambda$ and with $q>1$ boundary components.
%Then $Q$ is incompressible and $\partial$-incompressible in $X_K$, the exterior of $K$. 
Let $\hat{Q} \subset K(\lambda)$ be the non-separating $2$-sphere gotten by capping off $Q$ in $K(\lambda)$ with disks. Let $K' \subset K(\lambda)$ be the dual knot to the surgery. 

\begin{lem}\label{lem:p,N}
The homology class $[K] \in H_1(L(p,m))$ has order $q$. There exists an integer $k>0$ such that $p=kq^2$. Furthermore, $K(\lambda) = (S^1 \times S^2) \# N$, where $N$ is a lens space or $S^3$, and $|H_1(N)|=k$.
\end{lem}

\begin{proof}
The fact that $[K]$ has order $q$ follows from \cite[Lemma 2.13]{calegarigordon2013smallRationalGenus} and the fact that $Q$ has $q$ longitudinal boundary components. 
We have $K(\lambda)=(S^1 \times S^2) \# N$, where $N$ is the capped-off complement of a product neighborhood of $\hat{Q}$. The summand $N$ is a rational homology sphere, because $b_1(K(\lambda)) \le b_1(L(p,m))+1 = 1$. 
Because $K$ is simple, it is a $(1,1)$-knot, so $K(\lambda)$ has Heegaard genus 2.
Because Heegaard genus is additive under connected sum, it follows that $N$ has Heegaard genus 1.
Hence $N$ is either a lens space or $S^3$. Finally, Theorem~\ref{thm:r>0reduciblemeansreductiblesimple} (or an elementary argument involving Smith normal form) implies $p=kq^2$. 
\end{proof}

We follow the notation of sections~2.1 and 2.2 of \cite{greene2022fiberedSimpleKnots} in describing the simple knot $K$ and take $Q$ to be the norm-minimizing Seifert surface of $K$ described there (which must also be planar with $q$ boundary components of slope $\lambda$), carried by a branched surface, $\mathcal{B}$. Please refer to that to fill in details. Note that $Q$ is incompressible and $\partial$-incompressible in $X_K$ and that an orientation on $Q$ induces an orientation on the components of $\partial Q$ which are in the same direction (coherently oriented) in $\partial X_K$.

\begin{figure}[ht]
    \centering
    \includegraphics[width=5in]{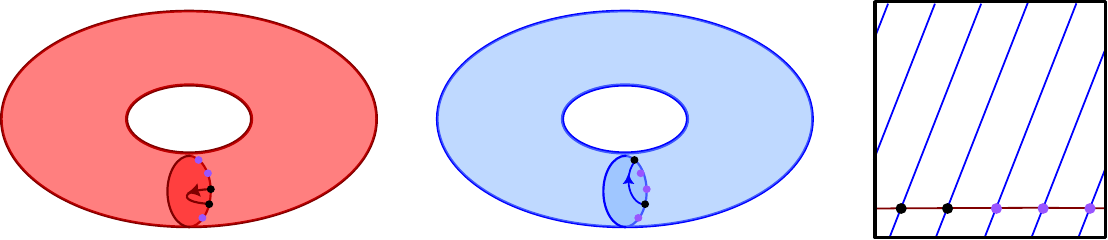}
    \put(-305, 80){$V_\alpha$}
    \put(-165, 80){$V_\beta$}
    \put(-40, 82){$\hat{T}$}
    \put(-310, 15){\small$\delta_{\alpha}$}
    \put(-164, 11){\small$\delta_{\beta}$}
    \put(-49, 12){\small$x_{\kappa}$}
    \put(-65, 12){\small$x_{0}$}
    \put(-292, 8){\tiny$x_{0}$}
    \put(-292, 15){\tiny$x_{\kappa}$}
    \put(-148, 8){\tiny$x_{0}$}
    \put(-157, 30){\tiny$x_{\kappa}$}
    \caption{A simple knot $\delta_\alpha \cup \delta_\beta$ in a doubly-pointed Heegaard diagram.}\label{fig:simple}
\end{figure}

Following sections~2.1 and 2.2 of \cite{greene2022fiberedSimpleKnots}, take a genus-1 Heegaard decomposition of the lens space $L(p,m) = V_\alpha \cup_{\hat T} V_\beta$ with Heegaard torus $\hat T$ and solid tori $V_\alpha, V_\beta$. Choose meridian disks $D_\alpha \subset V_\alpha, \,D_\beta \subset V_\beta$ so that $\alpha = \del D_\alpha$ and $\beta = \del D_\beta$ meet transversely in $p$ points of intersection.
Orient $\alpha, \beta$ so that $\alpha \cdot \beta =p$. Let $x_0,\dots,x_{p-1}$ denote a cyclic order along the oriented $\alpha$ of its points of intersection with $\beta$.
Orienting $K$, there is a $0 < \kappa < p$ such that the oriented $K$ is the simple knot  $K(p,m,\kappa) \subset L(p,m)$ gotten by taking the union of properly embedded, oriented arcs $\delta_\alpha \subset D_\alpha, \,\delta_\beta \subset D_\beta$ with $\del \delta_\alpha = - \del \delta_\beta = x_0 - x_\kappa$. See Figure~\ref{fig:simple}.

We will be describing the homology classes of knots in $H_1(L(p,m))$ in terms of canonical generators $\mu,\, \mu'$:

\begin{defin}
\label{def:mu}
Define  $\mu' \in H_1(L(p,m))$ to be the class of the core of $V_{\alpha}$ which links $\alpha$ positively in $V_\alpha$, and $\mu \in H_1(L(p,m))$ the class of the core of $V_{\beta}$ which links $\beta$ negatively in $V_\beta$. 
\end{defin}
One checks then that
\[\mu'=m \mu,\quad [K]=[K(p,m,\kappa)] = \kappa \mu.\]

\begin{figure}[ht]
    \centering
    \includegraphics[width=5in]{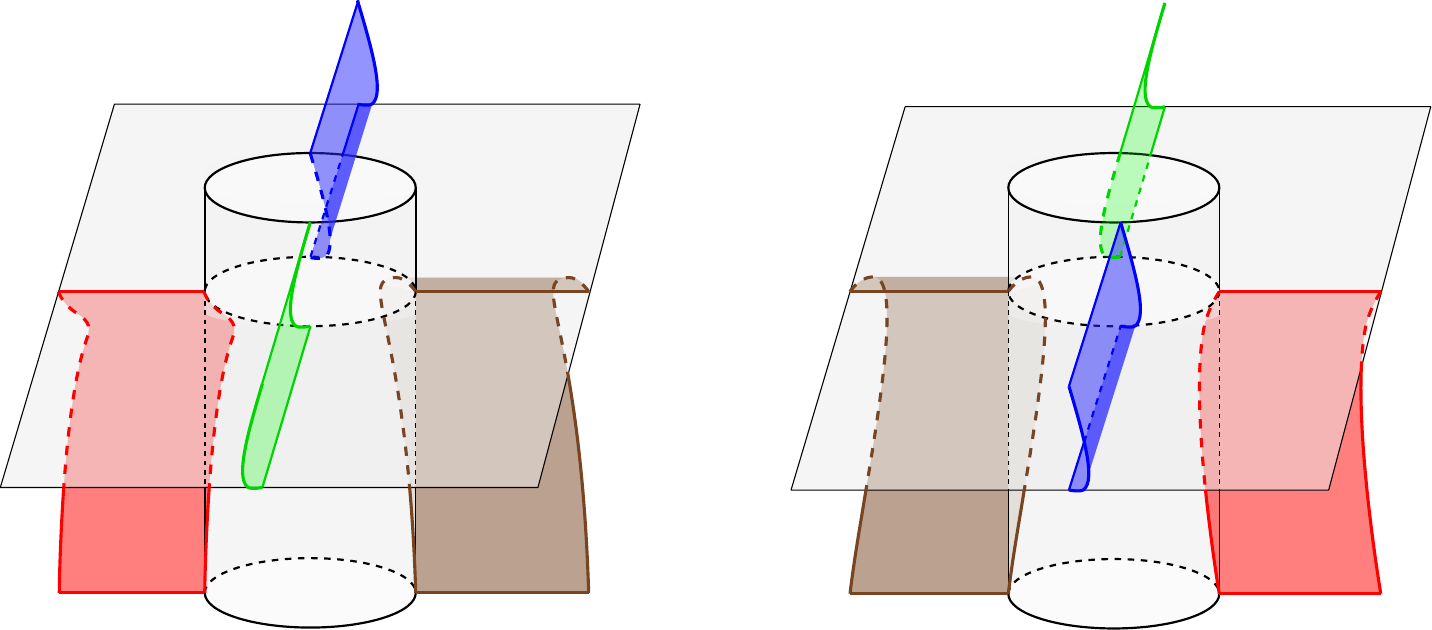}
\put(-10,5){$D_{\alpha}^1$}
\put(-330,92){$\alpha_1$}
\put(-245,92){$\alpha_2$}
\put(-40,92){$\alpha_1$}
\put(-273,110){$\beta_1$}
\put(-285,55){$\beta_2$}
\put(-70,110){$\beta_2$}
\put(-80,55){$\beta_1$}
\put(-133,92){$\alpha_2$}
\put(-365,5){$D_{\alpha}^1$}
\put(-210,5){$D_{\alpha}^2$}
\put(-165,5){$D_{\alpha}^2$}
\put(-65,160){$D_{\beta}^2$}
\put(-267,160){$D_{\beta}^1$}
\put(-90, 25){$D_{\beta}^1$}
\put(-297,25){$D_{\beta}^2$}
\put(-300,150){\textcolor{magenta}{$q-t$}}
\put(-330,123){\textcolor{magenta}{$v+t'$}}
\put(-265,123){\textcolor{magenta}{$v+t+t'-q$}}
\put(-355,40){\textcolor{magenta}{$v$}}
\put(-250,40){\small \textcolor{magenta}{$v+t$}}
\put(-342,12){\textcolor{white}{$t'$}}
\put(-242,12){\textcolor{magenta}{$q - t'$}}
\put(-300,65){\textcolor{magenta}{$t$}}
\put(-130,123){\small \textcolor{magenta}{$u+t'-t$}}
\put(-30,123){\textcolor{magenta}{$u+t'$}}
\put(-80,150){\textcolor{magenta}{$t$}}
\put(-25,12){\textcolor{white}{$t'$}}
\put(-35,40){\textcolor{magenta}{$u$}}
\put(-155,40){\small \textcolor{magenta}{$u+q-t$}}
\put(-140,12){\small \textcolor{magenta}{$q-t'$}}
\put(-115,65){\textcolor{magenta}{$q-t$}}
    \caption{ 
    The branched surface $\cB$ in $X_{K}$  near $K \cap \hat{T}$, where $\hat{T}$ is the Heegaard surface for $L(p,m)$. The left is near $x_0$, the right near $x_{\kappa}$. $Q$ is carried by $\cB$, with the pictured weights where $u,v \geq 0$ and $t>0$ are integers. 
    %on the Heegaard torus in the neighborhood of $x_0$ and $x_\kappa$.
    }
    \label{fig: branch1}
\end{figure}

Let $T$ be the twice-punctured torus gotten by removing open neighborhoods of  $x_0,\, x_\kappa$ from the Heegaard torus $\hat T$, and label the respective boundary components of $T$ as $A$, $B$.
Note that $K \cap V_{\alpha}=\delta_{\alpha}$ is oriented from $x_\kappa$ to $x_0$, and $K \cap V_{\beta}=\delta_{\beta}$ is oriented from $x_0$ to $x_\kappa$. 
%Recall that the exterior of $K$ fibers over the circle with fiber a $q$-punctured sphere $Q$. Then $Q$ is the unique mininiml-genus Seifert surface for $K$ in $L(p,m)$. $\spadesuit$ {\it (Need to change this to apply it to non-simple case.)}
Recall that we are taking $Q$ to be the Seifert surface carried by the branched surface $\mathcal B$ of $X_K$ described in Section~2.2 of \cite{greene2022fiberedSimpleKnots}.
The branched surface, $\mathcal B$, is formed as the union of $T$ along with four half-disks $D_{\alpha}^1,\, D_{\alpha}^2,\, D_{\beta}^1,\, D_{\beta}^2$ smoothed according to co-orientations on the sectors of $\mathcal{B}$ (where $D_{\alpha}^1 \cup  D_{\alpha}^2 = D_\alpha \cap X_K$ and  $D_{\beta}^1 \cup D_{\beta}^2 =D_\beta \cap X_K$). See Figure~\ref{fig: branch1}.
Choose a core curve $c_{\alpha}$ of $V_{\alpha}$ so that it intersects $K$ in a single point, and choose a small tubular neighborhood of $c_{\alpha}$, $C_{\alpha}$, that intersects $D_{\alpha}$ in a meridional disk, $D_\alpha^C$, of $C_\alpha$.
Similarly choose a core $c_{\beta}$ for $V_{\beta}$ and tubular neighborhood $C_\beta$ with $D_\beta \cap C_\beta = D_\beta^C$, a meridian of $C_\beta$.
Removing the interior of the neighborhoods of these cores from $L(p,m)$ leaves $\hat{T} \times \bk{0,1}$, with $\hat{T} \times \{0\} = \partial C_\alpha \subset V_{\alpha}$, $\hat{T} \times \{1\}= \partial C_\beta \subset V_{\beta}$, and $\hat{T} \times [0,1]$ meeting $K$ in two arcs $(\{x_0\} \cup \{x_\kappa\}) \times [0,1]$.
Removing a neighborhood of $K$ from this product leaves $T \times [0,1]$.

$Q$ is carried by $\mathcal B$.
Let $t',\, q-t'$ be the weights of $Q$ in the sectors of $\mathcal{B}$ corresponding, respectively, to $D_{\alpha}^1,\, D_{\alpha}^2$ as in Figure~\ref{fig: branch1}. Then $Q \cap (C_\alpha \cap X_K)$ is a collection of $t'$ half-disks which are copies of $D_\alpha^1 \cap (C_\alpha \cap X_K)$ along with $q-t'$ half-disks which are copies of $D_\alpha^2 \cap (C_\alpha \cap X_K)$. Thus $Q \cap (T \times \{0\})=Q \cap (\partial C_\alpha \cap X_K)$ is a collection $t'$ arcs in one isotopy class of $T \times \{0\}$ and $q-t'$ arcs in another isotopy class, each connecting different components of $\partial T \times \{0\}$. Similarly, let $q-t,\, t$ be the weights of $Q$ in the sectors of $\mathcal{B}$ given, respectively, by $D_{\beta}^1,\, D_{\beta}^2$. Then $Q \cap (C_\beta \cap X_K)$ is a collection of $q-t$ half-disks which are copies of $D_\beta^1 \cap (C_\beta \cap X_K)$ and $t$ half-disks which are copies of $D_\beta^2 \cap (C_\beta \cap X_K)$. Thus $Q \cap (T \times \{1\})=Q \cap (\partial C_\beta \cap X_K)$ is a collection $q-t$ arcs in one isotopy class of $T \times \{1\}$ and $t$ arcs in another isotopy class, each connecting different components of $\partial T \times \{1\}$.

\begin{defin}\label{def:alpha1}
Let $\bc{x,y} = \bc{t', q-t'}$ with $y \leq x$, and $\bc{x', y'} = \bc{t, q-t}$ with $y' \leq x'$.
\begin{itemize}
\item Let $\hat{\alpha}$ be the isotopy class in $T$ of the arcs $Q \cap (T \times \{0\})$ with $y$ members. 
\item Let $\hat{\alpha}'$ be the isotopy class in $T$ of the arcs $Q \cap (T \times \{0\})$ with $x$ members. 
%Then rel endpoints $\hat{\alpha}$ is isotopic in $V_{\alpha}$ to the oriented $K \cap V_{\alpha}$.
 
\item Let $\hat{\beta}$ be the isotopy class in $T$ of the arcs $Q \cap (T \times \{1\})$ with $y'$ members. 

\item Let $\hat{\beta}'$ be the isotopy class in $T$ of the arcs $Q \cap (T \times \{1\})$ with $x'$ members.
%Then rel endpoints $\hat{\beta}$ is isotopic to the oriented $K \cap V_{\beta}$. 
\end{itemize}
\end{defin}

Note that $\{\hat{\alpha}, \hat{\alpha}'\}=\{\alpha_1, \alpha_2\}$ and $\{\hat{\beta}, \hat{\beta}'\}=\{\beta_1, \beta_2\}$  where $\alpha_1, \alpha_2, \beta_1, \beta_2$ are the isotopy classes of the arcs depicted in Figure~\ref{fig: branch1}. 

Capping off $\partial T$ with disks to get $\hat{T}$, we may consider $\hat{\alpha}, \hat{\alpha}'$ as arcs running between $x_0$ and $x_\kappa$. Orienting from $x_\kappa$ to $x_0$, $\hat{\alpha}, \hat{\alpha}'$ is then isotopic rel endpoints to the oriented $K \cap V_\alpha$ in $V_\alpha$. Similarly, we may consider $\hat{\beta}, \hat{\beta}'$ as arcs in $\hat{T}$ running between $x_0$ and $x_\kappa$. Orienting them from $x_0$ to $x_\kappa$, each is isotopic rel endpoints to the oriented $K \cap V_\beta$ in $V_\beta$.

In~\cite{greene2022fiberedSimpleKnots} it is shown that
$\kappa=t'(p/q)=t'kq$ and that $l=tkq$ where $l<p$ is the least positive number such that $lm \equiv \kappa \pmod{p}$ and where $p=kq^2$.
Rewriting this in terms of the parameters $y,y'$, we have
%Thus $[K]=\kappa \mu=t'kq\mu$, $[K]=l \mu'=tkq \mu'$. As $y \in \bc{t', q-t'}$, in $H_1(L(p, m))$, up to sign, $[K]=t'kq\mu$ is $ykq \mu$. Similarly, up to sign $[K]$ is $y'kq \mu'$. We hereafter write this as
\begin{equation}\label{eq:yhomology}
\pm [K]=(y kq) \mu, \quad \pm [K]=(y' kq) \mu'
\end{equation}
whereby we mean that $[K]$ is determined up to sign by these classes.

Note that, as $[K]$ has order $q$ in $H_1(L(kq^2,m))$, both $y$ and $y'$ are relatively prime to $kq$. Thus $q>x>y$ and $q>x'>y'$.

%We use the product structure to identify isotopy classes of arcs in $T$ with those of $T \times \{i\}$ for $i \in [0,1]$.

%\begin{defin}\label{def:alpha1'}
%Let $\bc{x,y} = \bc{t', q-t'}$ with $y \leq x$.
%Let $\hat{\alpha} \in \{\alpha_1, \alpha_2\}$ be the isotopy class in $T \times \{0\}$ of the level arcs with $y$ members, oriented from $B$ to $A$.
%Then rel endpoints $\hat{\alpha}$ is isotopic in $V_{\alpha}$ to the oriented $K \cap V_{\alpha}$.
%Let $\bc{x', y'} = \bc{t, q-t}$ with $y' \leq x'$. Let $\hat{\beta} \in \{\beta_1, \beta_2\} $ be the isotopy class in $T \times \{1\}$ of the level arcs with $y'$ members, oriented from $A$ to $B$.
%Then rel endpoints $\hat{\beta}$ is isotopic to the oriented $K \cap V_{\beta}$. 
%\end{defin}

%Removing $\hat{T} \times (0,1)$ from $L(p,m)$ produces a pair of solid tori $V'_\alpha \subset V_\alpha$ and $V'_\beta \subset V_\beta$.

%As $y \in \bc{t', q-t'}$, in $H_1(L(p, m))$, up to sign, $[K]=t'kq\mu$ is $ykq \mu$. Similarly, up to sign $[K]$ is $y'kq \mu'$. We hereafter write this as
%\begin{equation}\label{eq:yhomology}
%\pm [K]=y kq \mu, \quad \pm [K]=y' kq \mu'
%\end{equation}
%Note that, as $[K]$ has order $q$ in $H_1(L(kq^2,m))$, both $y$ and $y'$ are relatively prime to $kq$. Thus $q>x>y$ and $q>x'>y'$. 

\begin{lem}
\label{lem:morseposition}
The intersection $Q \cap (T \times [0, 1])$ may be isotoped fixing $Q \cap (T \times \{0\} \cup T \times \{1\})$ so that there are only saddle and circle tangencies with respect to projection to $[0,1]$. Under this embedding, we have the following.
\begin{enumerate}[label=(\alph*)]
    \item There are exactly $2q-2$ saddle singularities.
    \item Let $C$ be a simple closed curve of intersection of $Q$ with $T \times \{i\}$, either transverse or circle singularity. Then $C$ is essential in each of $Q$, $T \times \{i\}$, and $X_K$. 
    \item Let $T \times \{i\}$ be transverse to $Q$. Then any arc of intersection with $Q$ must connect different components of $\partial (T \times \{i\})$.
\end{enumerate}
\end{lem}

The argument follows the techniques of section~4 of Thurston's thesis, {\it Foliations of Three-Manifolds Which Are Circle Bundles} found in   \cite{thurston2022collectedWorksFoliations}, in particular the Evolution Lemma and the proof of Theorem 3; or see Theorem 4 of \cite{thurston1986homologyNorm}. 

\begin{proof}
As $Q$ is transverse to $T \times \{0\} \cup T \times \{1\}$, we may isotope $Q$ in the interior of $T \times [0,1]$ to be Morse with respect to projection to $[0,1]$. The level curves emanating above/below a center tangency bound disks in the corresponding level surfaces $T \times \{i\}$. Eventually, following this disk through the levels, the limiting disk meets a saddle tangency. If the boundary of this disk meets the saddle in a single point, then one can use the trace of these disks to cancel these two tangencies. If it meets the saddle tangency in two points producing a pinched disk, then one can isotop to reduce the number of center tangencies. If it meets the saddle tangency in two points producing a pinched annulus, then one can isotop $Q$ to have a center tangency along the core of the annulus, thereby replacing the saddle and center tangencies by a circle tangency. In this way we can in sequence get rid of all center tangencies, first doing all possible cancellations of the first two types, then doing all possible circle tangency creations.

Under the above embedding of $Q$, $T \times [0,1]$ pulls back to a singular foliation, $\mathcal{F}_Q$, of the subsurface $Q \cap (T \times [0,1])$ with only saddle (index one) singularities.
%Each boundary component of $Q$ has a neighborhood in $Q$ which meets this singular foliation as follows.
Furthermore, $\mathcal{F}_Q$ is transverse to $\del Q$ in this subsurface of $Q$.
Each boundary component of $Q$ meets two  arcs of $\mathcal{F}_Q$, one from $T \times \{0\}$ and the other from $T \times \{1\}$, each of which is boundary parallel in $Q$.
The boundary parallelisms of these arcs into $\del Q$ on $Q$ is the union of half-disks which is the complement of $\mathcal{F}_Q$ in $Q$. 

We can then extend $\mathcal{F}_Q$ to a singular foliation of the capped-off sphere $\hat{Q}$ by adding $2q$ center (index zero) singularities, two for each component of $\partial Q$.
By the index formula, there then must be $2q-2$ saddle singularities.

Let $C$ be a simple closed curve of intersection of $Q$ and $T \times \{i\}$. As $C$ is a leaf in $\mathcal{F}_Q$, which has no center singularities, $C$ cannot bound a disk in $Q$.  Thus $C$ must be essential in $Q$. As $Q$ is incompressible in the exterior of $K$, $C$ must also be essential in $X_K$; hence, in $T \times \{i\}$.

Finally, assume $T \times \{i\}$ is transverse to $Q$. An orientation on $Q$ induces orientations on the boundary components of $Q$ which are parallel, as oriented curves, on $\partial X_K$---i.e.~the component of $\partial Q$ are {\it coherently oriented} on $\partial X_K$. As $X_K$ is orientable, any arc component of $Q \cap (T \times \{i\})$ must then connect different components of $\partial (T \times \{i\})$, as claimed.
%Note that in the spirit of Theorem 4 of \cite{thurston1986homologyNorm}, we could have gotten rid of circle tangencies as follows. Note as above that they can only occur when there are $q$ arcs in the level of tangency. Find a circle tangency whose ascending/descending annulus touches a saddle. 
\end{proof}

\begin{remark}
%In the spirit of Theorem 4 of \cite{thurston1986homologyNorm}, we can get rid of circle tangencies as follows. Argue that a circle tangency can only occur when the $q$ arcs in the level of tangency lie in exactly two isotopy classes of $T$.  Now apply the argument in Theorem 3 of \cite{thurston2022collectedWorksFoliations}, Case A2(b), to conclude that there must be a circle tangency whose ascending/descending annulus hits first a saddle singularity, $s$, coming from a level arc of $Q$. Isotop $Q$ to interchange the heights of $s$ and the saddle corresponding to the circle tangency. Then use the center tangency corresponding to the circle tangency to cancel $s$. We thereby replace the circle tangency with its constituent saddle tangency. Continue in this manner until all circle tangencies are replaced.
%However, in the current course of the argument we will naturally arrive at the point where, as far as Theorem~\ref{thm:Gr-Li-Lu} is concerned, we can assume no center tangencies occur. See Corollary ~\ref{cor:nocircletangciesr}.
%%We argue that indeed circle tangency leads to a transverse level $T \times \{i\}$ where there are only two isotopy classes of level arcs and no simple closed curves of intersection. The argument of Theorem~\ref{thm:annulus} then shows that we may assume that there is no such level.
In a similar vein, one should be able to get rid of circle tangencies.
However, these will not affect our arguments.
Furthermore, in the course of the argument we will naturally arrive at the point where  we can assume no center tangencies occur. See Corollary ~\ref{cor:nocircletangciesr}.
\end{remark}

%$\spadesuit$ {\em (Since we get rid of center tangencies basically at the end of the next section as now written, maybe let's just change this remark to something like:

We think of the saddle singularities in the above Morse (level) presentation of $Q$ as given by bandings of level curves, $Q \cap T \times \{i\}$, as $i$ changes. 

\begin{lem}\label{lem:bandingdiffclasses}
Assume that $T \times \{i-\varepsilon\}$, $T \times \{i+\varepsilon\}$ ($\epsilon$ positive or negative) are transverse to $Q$ and that a single saddle singularity occurs for $Q$ between these levels. This saddle corresponds to a banding $b$ between two level curves of $T \times \{i-\varepsilon\}$ to get the level curves of $T \times \{i+\varepsilon\}$. 
\begin{enumerate}
\item If $b$ connects a level arc to itself in $T \times \{i-\varepsilon\}$, then it produces a level arc and an essential level closed curve in $T \times \{i+\varepsilon\}$. In particular, if there are at least three isotopy classes of level arcs in $T \times \{i-\varepsilon\}$, then there is no banding of an arc to itself.
\item If $b$ connects two different level arcs of $T \times \{i-\varepsilon\}$, then these level arcs must be in different isotopy classes in $T \times \{i-\varepsilon\}$. 
\item Assume that $T \times \{i-\varepsilon\}$ has at least two isotopy classes of level arcs.  
Each of the isotopy classes of arcs in $T \times \{i-\varepsilon\}$ either gains or loses exactly one arc in $T \times \{i+\varepsilon\}$.
Further, assume $b$ connects two different level arcs of $T \times \{i-\varepsilon\}$. Then corresponding classes (which must be different) each lose an arc, and two different isotopy classes of arcs are created that are different from those banded.
\item Each isotopy class of arcs in $T \times \{i-\varepsilon\}$ loses at most one arc in going to $T \times \{i+\varepsilon\}$.
\end{enumerate}
\end{lem}

\begin{proof}
%Let $W_{-}$,$W_{+}$  be the sides of $T \times \{i-\varepsilon\}$, $T \times \{i+\varepsilon\}$ (resp.) in $X_K$ intersecting $V_{\alpha}$. Then $Q \cap W_{+}$ is obtained from $ Q \cap W_{-}$ by adding a 1-handle (corresponding to $b$). Thus if $b$ is a banding between of an arc of $T \times \{i-\varepsilon\}$ to itself, the banding yields and arc and simple closed curve component of $T \times \{i+\varepsilon\}$, as claimed. If $T \times \{i-\varepsilon\}$ contains at least three isotopy Note that $b$ cannot connect two level arcs of $T \times \{i-\varepsilon\}$ as every level arc of $T \times \{i+\varepsilon\}$ must connect $A$ to $B$. 

Note that we may assume no circle tangencies, as they only involve level simple closed curves and will not alter the statements below.

(1) Assume that $b$ bands an arc of $T \times \{i-\varepsilon\}$ to itself.
Then this gives rise to banding of a leaf to itself in the pullback foliation of $T \times [0,1]$ in $Q \cap T \times [0,1]$, which is transversely orientable.
The transverse orientation implies that the band must band the same side of the level arc to itself, thereby giving an arc and a simple closed curve, which is essential in $T$ by Lemma~\ref{lem:morseposition}(b).
If $T \times \{i-\varepsilon\}$ contains at least three isotopy classes, then the banding must occur in a disk.
Therefore, the banding cannot band an arc to itself, since it would produce an inessential curve.

(2) If $b$ was between two arcs in the same isotopy class in $T$ this would result in a level arc connecting a boundary component to itself in $T \times \{i+\varepsilon\}$, contradicting Lemma~\ref{lem:morseposition}$(c)$. 

(3) Assume $T \times \{i-\varepsilon\}$ contains at least two isotopy classes of level arcs.
If there is a level simple closed curve in $T \times \{i-\varepsilon\}$ (which must be essential in $T \times \{i-\varepsilon\}$), then $T \times \{i-\varepsilon\}$ has exactly two isotopy classes of arcs.
Now $b$ cannot band two simple closed curves as this would give a trivial simple closed curve in $T \times \{i+\varepsilon\}$, contradicting Lemma~\ref{lem:morseposition}(b). Thus $b$ bands an arc to a simple closed curve or to itself. In this case, the banded arc class loses an arc and the unbanded class gains one, as claimed. We hereafter assume there are no level simple closed curves. If $b$ bands an arc to itself, then again there must be exactly two isotopy classes of arcs and as above the banded arc class loses an arc and the unbanded class gains one. So we assume $b$ bands two different arcs.
Then by conclusion $(2)$ above they must connect different isotopy classes. Thinking of these classes as lying in a punctured annulus in $T \times \{i-\varepsilon\}$, one sees that $b$ must go from one side of the annulus to the other (as the resulting arcs must go from $A$ to $B$ in $T \times \{i+\varepsilon\}$). Thus there is a simple closed curve that intersects each of these new arcs once and is disjoint from each of the original arcs. Thus these arcs are in different isotopy classes in $T \times \{i+\varepsilon\}$ than the original arcs; furthermore, they are in different isotopy classes from each other by the argument of conclusion $(2)$. As $T \times \{i+\varepsilon\}$ has at most four different isotopy classes, these arcs will add an arc to any isotopy class in $T \times \{i-\varepsilon\}$ different than those banded. Thus the isotopy classes banded each lose an arc, and any others in $T \times \{i-\varepsilon\}$ gain one. 

(4) Note that only level curves involved in the banding are lost in the transition. Two arcs vanishing in the same isotopy class would contradict conclusion $(2)$. 
\end{proof}

Define $T_0=T \times \bc{0}$ (the $V_{\alpha}$-end of the height function) and  $T_{2q-2}=T \times \bc{1}$ (the $V_{\beta}$-end). Recall that there are $2q-2$ saddle points of $Q$ in $T \times [0,1]$. 
For each $i=1, \dots, 2q-3$, we pick a level surface $T \times \{z\}$ transverse to $Q$ between the $i^{th}$ and $(i+1)^{th}$ saddle-point critical levels and label this level surface $T_i$.
%, taking that labelled $0$ to be $T \times \bc{0}$ (the $V_{\alpha}$-end of the height function) and that labelled $2q-2$ to be $T \times \bc{1}$ (the $V_{\beta}$-end). 
At each $T_i$ %there are no simple closed curves of intersection---that is, 
we see $q$ properly embedded arcs of the level surface coming from its intersection with $Q$.

\begin{figure}
    \includegraphics[width=2.5in]{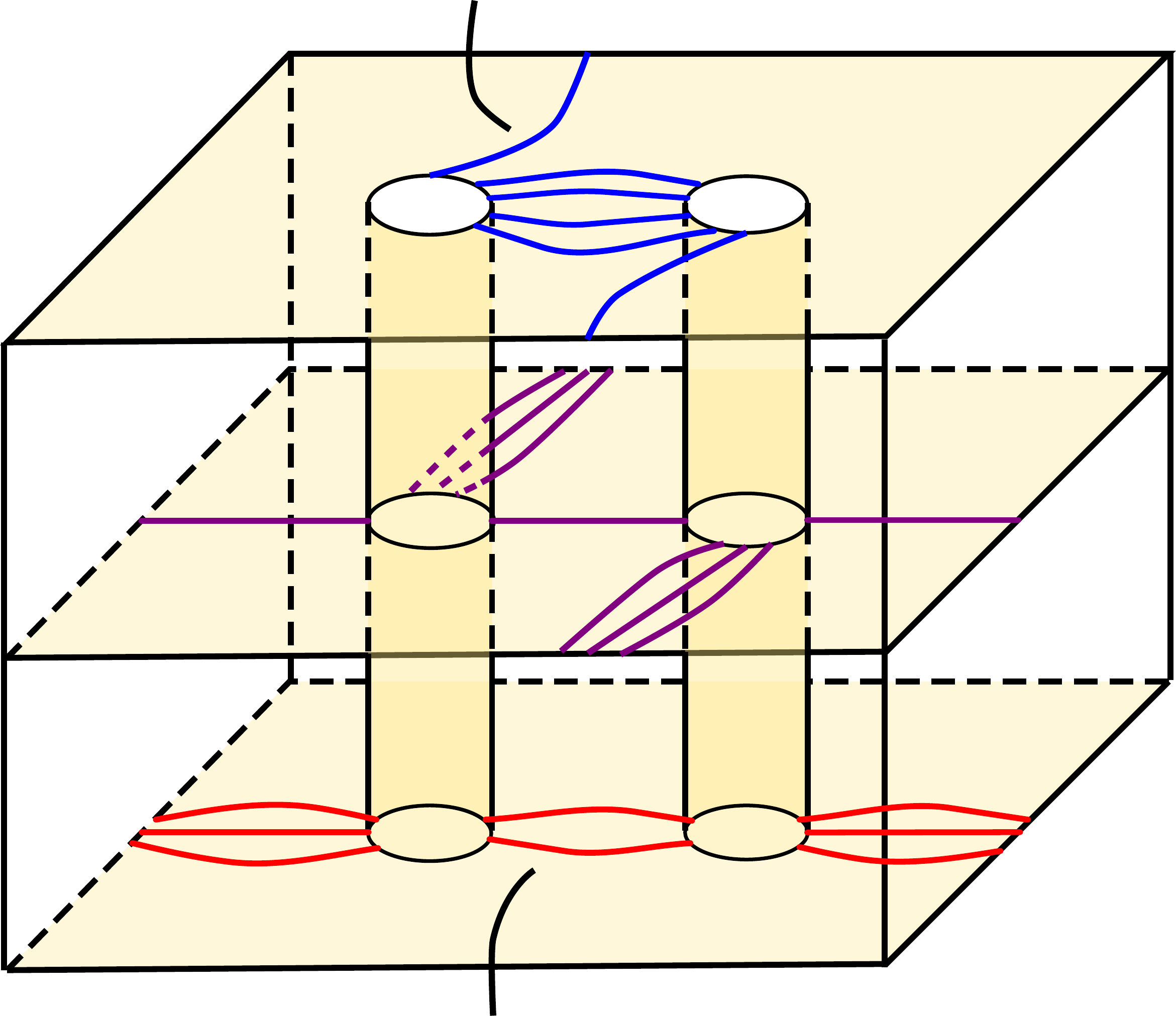}
    \put(-120, 160){$y'$ arcs}
    \put(5, 130){$\Gamma_{2q-2}$}
    \put(10, 80){$\Gamma_{i}$}
    \put(10, 30){$\Gamma_{0}$}
    \put(-115, -10){$y$ arcs}
    \caption{$T \times [0, 1]$ in $X(K)$; $Q \cap (T \times \{0\})=\Gamma_0$ is $y$ arcs in class $\hat{\alpha}$ and $q - y=x$ in class $\hat{\alpha}'$. $Q \cap (T \times \{1\})=\Gamma_{2q-2}$ is $y'$ arcs in class $\hat{\beta}$, and $q - y'=x'$ in class $\hat{\beta}'$.}
    \label{fig:level}
\end{figure}

\begin{defin} Let $\Gamma_i,\, i = 0, \dots, 2q-2$, be $T_i$ (above) along with $Q \cap T_i$ ($q$ arcs of intersection along with any simple closed curves of intersection). Define $\Gamma = \Gamma_0 \cup \dots \cup \Gamma_{2q-2}$. See Figure~\ref{fig:level}.

Re-index $\Gamma_i$ from $V_{\beta}$ to $V_{\alpha}$ by setting $\Gamma_i'=\Gamma_{(2q-2)-i}$. Define $\Gamma'=\Gamma_0' \cup \dots \cup \Gamma_{2q-2}'$.
\end{defin}

The union of surface graphs $\Gamma$ records how the arcs of $Q \cap (T \times \{i\})$ change as we move from $i=0$ to $i=1$ (i.e. from $V_{\alpha}$ to $V_{\beta}$).
Likewise, $\Gamma'$ records this change in going from $i=1$ to $i=0$ (from 
$V_{\beta}$ to $V_{\alpha})$.

Recall that we label the components of $\partial T$ as $A,B$. The components of $\del T \times [0,1]$ are therefore $A \times [0,1]$ and $B \times [0,1]$.
Number the components of $\partial Q$ from $1$ to $q$.
Each component of $\del Q$ has one arc $\omega_A$ in $A \times [0,1]$ and one arc $\omega_B$ in $B \times [0,1]$. 
The endpoints of these arcs demarcate the intersection of this component with $T \times \{0\}$ and $T \times \{1\}$. 
Also, any arc in a level surface of $\Gamma$ has one endpoint at $A$ and one endpoint at $B$ (Lemma~\ref{lem:morseposition}(c)), and it has a label from $\{1,\dots,q\}$ at each endpoint (the corresponding component of $\partial Q$).
Under the projection $T \times [0,1]$ to the first component, an arc $\tau \subset \Gamma_b$ projects to an arc in $T$, and the projection extends by radial arcs in the capping disks to an arc in $\hat{T}$ with endpoints $x_0$ and $x_\kappa$.
We refer to this arc as the projection of $\tau$ to $\widehat{T}$.

\begin{lem}\label{lem:sameisotopyclasses}
Let $\tau \subset \Gamma_b$ and $\tau' \subset \Gamma_{b'}$ be arcs with $b<b'$.
Suppose that $\tau$ and $\tau'$ have the same labels at $A$ and at $B$, and suppose that the projection of $\tau \cup \tau'$ to $\hat{T}$ is trivial in $H_1(L(p,m))$.
Then $\tau$ is parallel to $\tau'$ in $Q$.
Moreover, no critical point of $Q$ occurs between $\tau$ and $\tau'$ in $Q$, and the level (Morse) presentation of $Q$ identifies these arcs between levels $b,\, b'$.
\end{lem}

\noindent
{\em Remark.}
By symmetry, the same holds if $\Gamma$ is replaced by $\Gamma'$.

\begin{proof}
Suppose that in $\Gamma$, the arcs $\tau$ and $\tau'$ have label $i$ at endpoint $A$ and label $j$ at endpoint $B$.
Then on the surface $Q$, $\tau$ and $\tau'$ are a pair of arcs from boundary component $i$ to boundary component $j$.
At boundary component $i$, both have an endpoint in the $\omega_A$ arc, and at boundary component $j$, both have an endpoint in the $\omega_B$ arc.
The union of $\tau$, $\tau'$, the portion of the $\omega_A$ arc between their endpoints on component $i$, and the portion of the $\omega_B$ arc between their endpoints on component $j$ forms a closed curve $\gamma \subset Q$.
It cobounds with $\del Q$ a $t$-holed disk $P$ for some $0\le t \le q-1$.
The holes are coherently oriented curves, each isotopic to $K$ in $L(p,m)$.
It follows that $P$ is a homology from $[\gamma]$ to $t [K] \in H_1(L(p,m))$.
On the other hand, $\gamma$ is homologous to the projection of $\tau \cup \tau' \subset T \times [0,1]$ to $\hat{T}$ (the portions of $\gamma$ composed of arcs in $\del Q$ get projected to points).
By hypothesis, this curve is nullhomologous.
We infer that $t [K] = 0 \in H_1(L(p,m))$.
However, $[K]$ has order $q$ in $H_1(L(p,m))$, and $0 \le t \le q-2$.
It follows that $t = 0$.
Thus, $P$ is a parallelism between $\tau$ and $\tau'$ on $Q$.
It follows that the singular foliation of $Q$ induced by the Morse function is nonsingular on $P$.
The last part of the Lemma now follows.
\end{proof}

%\begin{cor}
%If $b < y$, then in going from $\Gamma_b$ to $\Gamma_{b+1}$, there cannot be an increase in the number of arcs in class $\hat{\alpha}$. If $b < x$, then in going from $\Gamma_b$ to $\Gamma_{b+1}$ there there cannot be an increase in the number of arcs in class $\hat{\alpha}'$, where $\{\hat{\alpha},\hat{\alpha}'\}=\{\alpha_1, \alpha_2\}$.

%Similarly, if $b < y'$, then in going from $\Gamma_b'$ to $\Gamma_{b+1}'$, there cannot be an increase in the number of arcs in class $\hat{\beta}$. If $b < x'$, then in going from $\Gamma_b'$ to $\Gamma_{b+1}'$ there there cannot be an increase in the number of arcs in class $\hat{\beta}'$, where $\{\hat{\beta},\hat{\beta}'\}=\{\beta_1, \beta_2\}$. 
%\end{cor}

%\begin{proof}
%\end{proof}

\begin{cor}\label{cor:xdisappears}
If $b \leq y$, then there are $y-b$ arcs in class $\hat{\alpha}$ in $\Gamma_b$. If $b \leq x$, then there are $x-b$ arcs in class $\hat{\alpha}'$ in $\Gamma_b$.
Similarly, if $b \leq y'$, then there are $y'-b$ arcs in class $\hat{\beta}$ in $\Gamma_b'$. If $b \leq x'$, then there are $x'-b$ arcs in class $\hat{\beta}'$ in $\Gamma_b'$.

\end{cor}

\begin{proof}
We give the proof for $\Gamma$. The same argument works for $\Gamma'$.

In $\Gamma_0$, there are two isotopy classes of arcs, $\hat{\alpha}$ with $y$ arcs and $\hat{\alpha}'$ with $x$ arcs.
%Each arc has $\alpha$-weight zero.
%Furthermore, from the branched surface picture Figure~\ref{fig: branch1} we can see that e
Each arc in $\Gamma_0$ begins and ends at the same component of $\partial Q$, cobounding a disk in $Q$ with that component of $\partial Q$ (a half-disk of $Q \cap (C_\alpha \cap X_K)$ discussed earlier).
%Labelling the endpoints as set out before Lemma~\ref{lem:sameisotopyclasses}, w
We see that for each positive integer $i\leq q$ there is exactly one arc in $\Gamma_0$ with label $i$ at both ends. 
 By Lemma~\ref{lem:bandingdiffclasses}$(4)$, if $b<y$ ($b<x$), then there must be arcs in class $\hat{\alpha}$ ($\hat{\alpha}'$, resp.) of $\Gamma_b$. Furthermore, if any arc in an isotopy class of $\Gamma_b$ has both labels the same, then all do. In particular, if $b<y$ ($b<x$), all the arcs in class $\hat{\alpha}$ ($\hat{\alpha}'$, resp.) have the same label.

If $b<y$, then the number of arcs in class $\hat{\alpha}$, which is non-zero, cannot increase from $\Gamma_b$ to $\Gamma_{b+1}$. Otherwise, there is an arc $\tau$ in class $\hat{\alpha}$ of $\Gamma_{b+1}$  created by the saddle singularity between levels $b$ and $b+1$. Because it is in class $\hat{\alpha}$, and class $\hat{\alpha}$ is non-empty in $\Gamma_b$, $\tau$ must have the same label, $i$, on its endpoints. On the other hand, there is an arc $\tau'$ in $\Gamma_0$ with label $i$ on its endpoints in class $\hat{\alpha}$ or $\hat{\alpha}'$. By Lemma~\ref{lem:sameisotopyclasses} ($\tau \cup \tau'$ projects into a neighborhood of a meridian of $\hat{T}= \del V_\alpha$), $\tau,\, \tau'$ are isotopic in $Q$. This contradicts that $\tau$ arose from a saddle singularity between levels $b, b+1$. 

Similarly, if $b<x$, then the number of arcs in class $\hat{\alpha}'$, which is non-zero, cannot increase from $\Gamma_b$ to $\Gamma_{b+1}$. 

For any $b<x$, there are arcs in class $\hat{\alpha}'$ of $\Gamma_b$ by the first paragraph. If all $q$ arcs were in class $\hat{\alpha}'$ then at some preceding level the number of arcs in $\hat{\alpha}'$ must have increased, contradicting the second paragraph. Thus for any $b<x$ there are at least two isotopy classes of arcs in $\Gamma_b$. By Lemma~\ref{lem:bandingdiffclasses}$(3)$ and the fact that $\hat{\alpha}$ and $\hat{\alpha}'$ cannot gain arcs, we conclude that $\hat{\alpha}$ ($\hat{\alpha}'$) loses exactly one arc in going from $\Gamma_b$ to $\Gamma_{b+1}$ when $b<y$ ($b<x$, resp.), proving the Corollary.
\end{proof}

\section{Two isotopy classes at a level graph \pdfmath{\Gamma_b}{γb}}\label{sec:atleast3}

As in section~\ref{sec:setup}, let $K$ be a simple knot in $L(p,m)$ with a planar Seifert surface $Q$ of longitudinal slope $\lambda$ and with $q$ boundary components. 
By Lemma~\ref{lem:p,N}, there is a value $k>0$ with $p=kq^2$.
Let $\Gamma$ be the collection of level surfaces for $K$ described in section~\ref{sec:setup}.
The goal of this section is to prove the following result.
It is used to prove Theorem~\ref{thm:mainNextended} in the next section.
%, which will allow us to assume that each $\Gamma_b$ has at least three isotopy classes of arcs as long as $b \neq 0, 2q-2$. 

In the lens space $L(kq^2, kdq \pm 1)$, with $k>0,\, q>d>0$ such that $(d,q)=1$, there is a unique knot lying on the Heegaard torus with winding number $q$ in each Heegaard solid torus. We refer to this knot as the {\em $(q,q)$-torus knot} in $L(kq^2, kdq \pm 1)$. Writing $L(kq^2, kdq \pm 1)$ as the standard Heegaard splitting $V_\alpha \cup_{\hat{T}} V_\beta$ (section~\ref{sec:setup}), we have that, up to orientation on $L$, $[L]=d [\alpha] + q [l] \in H_1(\partial V_\alpha)$, where $\alpha, l$ are meridian, longitude and $L$ is the $(q,q)$-torus knot.

\begin{thm}
\label{thm:annulus}
Suppose that there exists an index $0 < b < 2q-2$ such that the arcs of $\Gamma_b$ belong to at most two isotopy classes.
Then up to diffeomorphism of $L(p,m)$, there exist integers $k>0,\, q>d>0$ such that $(d,q)=1$, $p=kq^2$, and $m=kdq + 1$.
Furthermore, either
\begin{enumerate}
\item $\pm [K]=kdq \mu \in H_1(L(p,m))$; or
\item $K$ is an $k$-fold cable on the $(q,q)$-torus knot, so $\pm [K]=kq \mu \in H_1(L(p,m))$. 
\end{enumerate}
\end{thm}

%\begin{thm}
%\label{thm:annulus}
%$K \subset L(p,m)$ for $p=kq^2$ and $k,m>0$. Let $b \in \{0, \dots, 2q-2 \}$,
%\begin{enumerate}
%\item $\Gamma_b$ must have at least two different isotopy classes among its arcs. 
%\item Assume $\Gamma_b$ has no simple closed level curves of intersection. If $\Gamma_b$ has exactly two distinct isotopy classes among its arcs, then either
%\begin{itemize}
%\item $b=0$ or $b=2q-2$; or
%\item up to diffeomorphism of $L(p,m)$, there is a $q>d>0$ with $m=kdq + 1$. Furthermore, one of the following holds
%\begin{enumerate}
%\item $\pm [K]=kdq \mu \in H_1(L(p,m))$ and $(d,q)=1$; or
%\item $K$ is an $k$-fold cable on a companion, $L$, that lies on the Heegaard torus for $L(p,m)$ with winding number $q$ on each side. In particular, $\pm [K]=kq \mu \in H_1(L(p,m))$. 
%\end{enumerate}
%\end{itemize}
%\end{enumerate}
%The same hold for any $\Gamma_{b}'$ (since $\Gamma_b'=\Gamma_{(2q-2-b)}$).
% \end{thm}

We build up to the proof of Theorem~\ref{thm:annulus} through a sequence of lemmas.
We first reduce to the case in which the graph $\Gamma_b$ does not contain any simple closed curves.

\begin{lem}
\label{lem:alphaweightr}
If there exists an index $0 < b < 2q-2$ such that the arcs of $\Gamma_b$ lie in at most two isotopy classes, then we can choose $b$ such that the arcs of $\Gamma_b$ lie in exactly two isotopy classes, and $\Gamma_b$ contains no simple closed curves.
\end{lem}

\begin{proof}
Pick $0 < b < 2q-2$ minimal with the property that the arcs of $\Gamma_b$ lie in at most two isotopy classes.
We claim that $b$ is the desired index.

To begin, the arcs of $\Gamma_0$ are in two isotopy classes, and $\Gamma_0$ has no simple closed curves.
These two isotopy classes are contained in an annulus on $\hat{T} \times \{0\}$ whose slope is meridional in the solid torus with boundary $\hat{T} \times \{0\}$ in $V_\alpha$.
The saddle between $\Gamma_0$ and $\Gamma_1$ corresponds to a banding between arcs of intersection.
It cannot be a banding of an arc with itself, as this would result in a meridional simple closed curve in $\Gamma_1$:
this curve is nullhomotopic in the knot exterior $X_K$, contradicting Lemma~\ref{lem:morseposition}(b).
It follows by Lemma~\ref{lem:bandingdiffclasses}(c), using the fact that $x > y \ge 1$, that $\Gamma_1$ has arcs in at least three different isotopy classes.
Thus, $1<b<2q-2$ and $\Gamma_{b-1}$ has at least three isotopy classes of arcs.
The complementary regions to these arcs in $T_0$ are simply-connected, so $\Gamma_{b-1}$ has no simple closed curves of intersection.
By Lemma~\ref{lem:bandingdiffclasses}(1), the saddle singularity between $\Gamma_{b-1}$ and $\Gamma_{b}$ must be between two different isotopy classes of arc on $\Gamma_{b-1}$.
Thus, $\Gamma_b$ has no simple closed curves of intersection and has arcs in at least two different isotopy classes, as desired.
\end{proof}

Fix the index $b$ fulfilling the conclusion of Lemma~\ref{lem:alphaweightr} for the rest of the section.
%For any value $b \in \{1, \dots, 2q-3 \}$ and l
Let $T=T_b$ be the twice-punctured torus corresponding to $\Gamma_b$, and let $\hat{T}$ be the Heegaard torus at level $b$ obtained from capping $\partial T$.
%Recall that $Q$ is a planar surface with $q$ boundary components and that we denote the exterior of $K$ by $X_K$. 
%By Lemma~\ref{lem:morseposition}, any simple closed curve of $T \cap Q$ is essential in both $Q$ and $T$.
Because there are two isotopy classes of arcs of $\Gamma_b$, they lie within an annulus $\mathcal{A} \subset \hat{T}$.

Label the components of $\partial Q$ from $1$ to $q$ in sequence along $\partial X_K$, and label the components of $\partial T$ as $A,B$. 

\begin{defin}
\label{def:graphs}
{\textup{(Graphs of intersection $G_Q, G_T$.)}}

Define the {\em labelled graph of intersection} $G_Q$ to be $Q$ along with the arcs of $Q \cap T$, where we label as $A$ or $B$ the endpoints of the arcs according to the corresponding boundary components of $T$. We think of $G_Q$ as a graph in the sphere, where the components of $\partial Q$ are considered the (fat) vertices and the arcs considered the (labelled) edges of $G_Q$.

Similarly, define the {\em labelled graph of intersection} $G_T$ to be $T$ along with the arcs of intersection of $Q \cap T$, where the endpoints of the arcs are labelled $1, \dots, q$ by the corresponding components of $\partial Q$. We think of $G_T$ as a graph in the torus $\hat{T}$. 

%A {\em Scharlemann cycle}, or simply {\em Schycle}, in $G_Q$ is a cycle of edges which bounds a disk in $G_Q$. An innermost cycle of edges in $G_Q$ is a Schycle. We pick one such, $S$. Let $s$ be the number of edges in $S$.
\end{defin}

Note that $G_T$ is just $\Gamma_b$, and each arc of intersection runs from component $A$ to component $B$ of $\partial T$.

Because each component of $\partial Q$ intersects each component of $\partial T$ exactly once, at each vertex of $G_Q$ we see exactly two edges incident: one with label $A$, the other with label $B$. Thus, the edges and vertices of $G_Q$ group in cycles of edges. Because the arcs of $Q \cap T$ run from component $A$ of $\partial T$ to component $B$, all edges of $G_Q$ have one end with label $A$ and one with label $B$. In a cycle of edges, orienting each edge from label $A$ to $B$ creates an oriented cycle. In the graph $G_Q$ we must have at least two disk faces. 

\medskip

\begin{lem}
    \label{lem:diskface}
Suppose that $S$ is a disk face of $G_Q$ with $s$ sides.
Then $s_1 \ge 1$ sides lie in one isotopy class of arc on $\Gamma_b$ and the other $s_2 \ge 1$ sides lie in the other isotopy class.
\end{lem}

The notation $S$ is chosen to evoke ``Scharlemann cycle".

\begin{proof}
Let $U$ be the Heegaard solid torus on the side of $\hat{T}$ containing $S$.  As all edges of $G_Q$ lie in $\mathcal{A}$ on $\hat{T}$, $S$ involves at most two edge classes of $T$.
Let $s_1$, $s_2$ be the number of edges of $S$ in each class, so $s=s_1+s_2$.

We first we claim that $s>1$.
Otherwise, let $e$ be the single edge of $S$.
Without loss of generality, suppose that $U \subset V_{\alpha}$.
Then $S$ contains an arc $e'$ of $\Gamma_0$ in its interior, and $S$ gives a parallelism between $e,e'$.
That is, there is no saddle singularity of $Q$ between $e \subset \Gamma_0$ and $e' \subset \Gamma_b$. 
%which must be in class $S$ gives an isotopy rel. endpoints of $e$ to $h_U$ which in turn is isotopic to $\hat{\alpha}$ (recall Definition~\ref{def:alphaweights}). There is an arc, $e'$ of $\Gamma_0$ with the same endpoint labels as $e$ and in class $\hat{\alpha}$ or $\hat{\alpha}'$. By Lemma~\ref{lem:sameisotopyclasses}, there can be no saddle between $e$, $e'$ on $Q$. Thus they are in the same class, either $\hat{\alpha}$ or $\hat{\alpha}'$ on $T$. 
By Corollary~\ref{cor:xdisappears}, $b<x$.
On the other hand, by Lemma~\ref{thm:m:b<x}, if $b<x$ and the edges of $\Gamma_b$ lie in at most two classes, then $b=0$, a contradiction. 

Next, we claim that both $s_1, s_2 \ge 1$.
Suppose for a contradiction that $s_2=0$.
Let $e$ be an edge of $S$.
Then a neighborhood of $e$ on $\hat{T}$ will be a disk, $D$, that contains vertices $A,\, B$ and all edges of $S$.
Let $h_U$ denote the arc $K \cap U$. 
Then the regular neighborhood $\N (D \cup h_U \cup S)$ will be a punctured lens space of order $s$ in the solid torus $U$, a contradiction.
\end{proof}

We use the notation $\N (Z)$ going forward to denote a regular neighborhood of a subset $Z$ of a manifold. 
Let $U$, $V$ be the Heegaard tori on either side of $\hat{T}$, with $U$ being the side containing $S$.
Let $h_U=K \cap U$ and $h_V=K \cap V$ (thinking of ``h" for ``handle").
% Note that if all arcs in $\Gamma_b$ were of the same isotopy classes, then one of $s_1$, $s_2$ would have been $0$. On the other hand, because $y,y' \neq 0$, $\Gamma_{0}, \Gamma_{2q-2}$ must contain two isotopy classes of arcs.
The space $H = \N(\mathcal{A} \cup h_U) \cap U$ is a genus-$2$ handlebody.
We distinguish a set of meridian disks in $H$.
First, $D_0$ is the co-core of the handle $\N(h_U)$.
Second, choose a co-core $d_i$ of the annulus $\mathcal{A}$ disjoint from the feet of $h_U$ and that only hits the edges in the parallelism class containing $s_i$ edges.
Then $D_i$ is the product of $d_i$ with $I$ in the identification $\mathcal{A} \times I = \N(\mathcal{A}) \subset H$.
The boundary $\del S$ within $\del H$ has intersection number $\pm s$ with $\del D_0$ and intersection number $\pm s_i$ with $\del D_i$ on $\del H$.
On this basis we see that $D_0$,$D_i$ are a pair of cutting disks for $H$ for $i=1,2$, the oriented cycle $\partial S$ intersects each of $\partial D_0, \partial D_1, \partial D_2$ in the same direction, and $\del S$ is a non-separating curve on $\del H$.
Hence $U' = H \cup \N(S)$ has torus boundary, it is contained in $U$, and $U' \cap \hat{T} = \mathcal{A}$.  
Hence $\mathcal{A}' = \del U' - int(\mathcal{A})$ is a properly embedded annulus in $U$. Let $U'' = U - int(U')$. %Note that 
%smallskip
%\noindent {The core of $A$ is not a meridian of $U$.}
%smallskip
%Otherwise, there is a meridian disk of $U$ that intersects $U'$ in a boundary component of $A$. A neighborhood of $U'$ in union with a disk, would give punctured lens space of order $s>1$ in $U$, a contradiction.

%$A'$ is isotopic rel boundary into the boundary of $U$.
%Therefore, $U'$ and its complement $U'' = U - int(U')$ are solid tori; one of the two guides the isotopy of $A'$ into $T$, and the other is a deformation retract of $U$.
%We will now determine which.

\begin{lem}
\label{l: solid torus}
Let $a$ be a component of $\partial \mathcal{A}$.
Then $[a]$ represents $s$ times a generator in $H_1(U)$ and $H_1(U')$, and it represents a generator in $H_1(U'';\bZ)$.
Both $U'$ and $U''$ are solid tori.
Moreover, $U$ deformation retracts to $U'$, and $U''$ guides an isotopy of $\mathcal{A}'$ into $\hat{T}$.
\end{lem}

\begin{proof}
First, suppose by way of contradiction that $a$ bounds a disk $\Delta$ in $U$.
%Suppose first that $\Delta \subset V$.
%Attaching the 2-handle $\N(\Delta)$ to $H$ along $a$ forms a solid torus with meridian $D_0$.
%Further attaching $\N(S)$ to it gives a punctured lens space $M$ of order $s \leq q$ embedded in $L(kq^2,m)$, a contradiction.
%If instead $\Delta \subset U$, then we work towards the same contradiction, as follows.
Such a disk $\Delta$ is a meridian disk of $U$.
Since $h_U$ is contained in a meridian disk of $U$, we can arrange by an isotopy of $\Delta$ that $\Delta \cap h_U = \emptyset$.
It follows that $\Delta \cap S$ consists of simple closed curves of intersection.
By successively compressing $\Delta$ along innermost loops of intersection on $S$, we arrive at a disk $\Delta$ such that $\Delta \cap S = \emptyset$.
Now, $H \cup \N(\Delta)$ is a solid torus.
It has meridian disk $D_0$.
The disk $S$ (less a collar neighborhood) meets $H \cup N(\Delta)$ solely in its boundary, and it runs $s$ times over the boundary of $D_0$.
Therefore, $H \cup \N(\Delta) \cup \N(S)$ is a copy of a punctured lens space of order $s \le q$ in $L(kq^2,m)$, a contradiction.

Since $a$ does not bound a disk in $U$, it follows that $\mathcal{A}'$ is a properly embedded, incompressible annulus in $U$.
Every properly embedded, incompressible annulus in a solid torus is boundary parallel, and it cuts off a solid torus to either side: one is the boundary parallelism, and the other is a deformation retract of the original solid torus \cite[Lemma 1.10]{hatcher3mfds}.
We examine which is the case for $U'$ and $U''$.
As $M$ is gotten by attaching a 2-handle along $U'$, we conclude that 
$[a]$ represents $s$ times a generator of $H_1(U')$. As $s>1$ and $U,U''$ are solid tori,  $a$ must represent the core of $U''$. Thus $U''$ guides an isotopy of $\mathcal{A}'$ into $\hat{T}$. In particular, $U'$ deformation retracts to $U$, and $[a]$ represents $s$ times the generator of $H_1(U)$.
%Attach $\N(S)$ and a 2-handle along $a$ to $\del H$ to form a 3-manifold $M$.
%The matrix of intersection pairings in $\del H$ of the pair $(\del D_0, \del D_i)$ with $(\del S, a)$ is, up to sign, $I = \left( \begin{matrix} s & 0 \\ s_i & 1\end{matrix} \right)$.
%Therefore, $H_1(M) \approx \mathrm{coker}(I) \approx \bZ / s \bZ$.
%Note that $M$ is the same as attaching a 2-handle along $a$ to the solid torus $U'$, so it is a punctured lens space of order $s$.
%It follows that $[a]$ represents $s$ times a generator of $H_1(U')$. As $s>1$, we conclude that $a$ (isotopic to the core of $A'$) must represent a generator in the solid torus $U''$ and $U''$ guides an isotopy of $A'$ into $hat{T}$. Thus $U'$ deformation retracts to $U$, and $[a]$ represents $s$ times the generator of $H_1(U)$.
%On the other hand, whichever of $U'$ or $U''$ guides the isotopy of $A'$ into $T$ deform retracts onto $A'$ and in turn to $a$, so $[a]$ generates its homology.
%Since $[a]$ does not generate $H_1(U')$, it follows that $[a]$ generates the homology of $H_1(U'')$ and $U''$ guides the isotopy of $A'$ into $T$.
%Further, $U$ deform retracts to $U'$, so $[a]$ represents $s$ times the generator of $H_1(U)$.
\end{proof}

The handle $\N(h_U)$ meets $\hat{T}$ in the fat vertices of $G_T$ labeled $A,B$.
    The complement of their interiors in $\del (\N(h_U))$ is an annulus $\mathcal{A}_{AB}$ properly embedded in $U$.
$S$ intersects $\mathcal{A}_{AB}$ in $s$ parallel arcs (cocores).
$S$ decomposes these fat vertices into $s$ arcs apiece called {\em corners}.
$S$ therefore decomposes $\mathcal{A}_{AB}$ into $s$ rectangles, where each rectangle consists of two arcs of $S$, a corner of fat vertex $A$, and a corner of fat vertex $B$.
A corner of fat vertex $A$ and a corner of fat vertex $B$ are {\em contiguous} if they occur in the boundary of the same rectangle.
A corner of fat vertex $A$  or fat vertex $b$  is an {\em outer corner} if it is in the same component of $(\mathcal{A} \cap X_K) - S$ as a component of $\del \mathcal{A}$.
An outer corner of fat vertex $A$ is {\em opposite} to an outer corner of fat vertex $B$ if they are contained in different components of $(\mathcal{A} \cap X_K) - S$.
Thus, the outer corners of fat vertex $A$ and fat vertex $B$ come in two pairs of opposites.

\begin{lem}
\label{l: opposite corners}
Either $\min \{s_1,s_2\} = 1$, or else there exists a pair of contiguous opposite corners.
%of fat vertex $1$ and fat vertex $-1$.
\end{lem}

\begin{figure}
\includegraphics[width=4in]{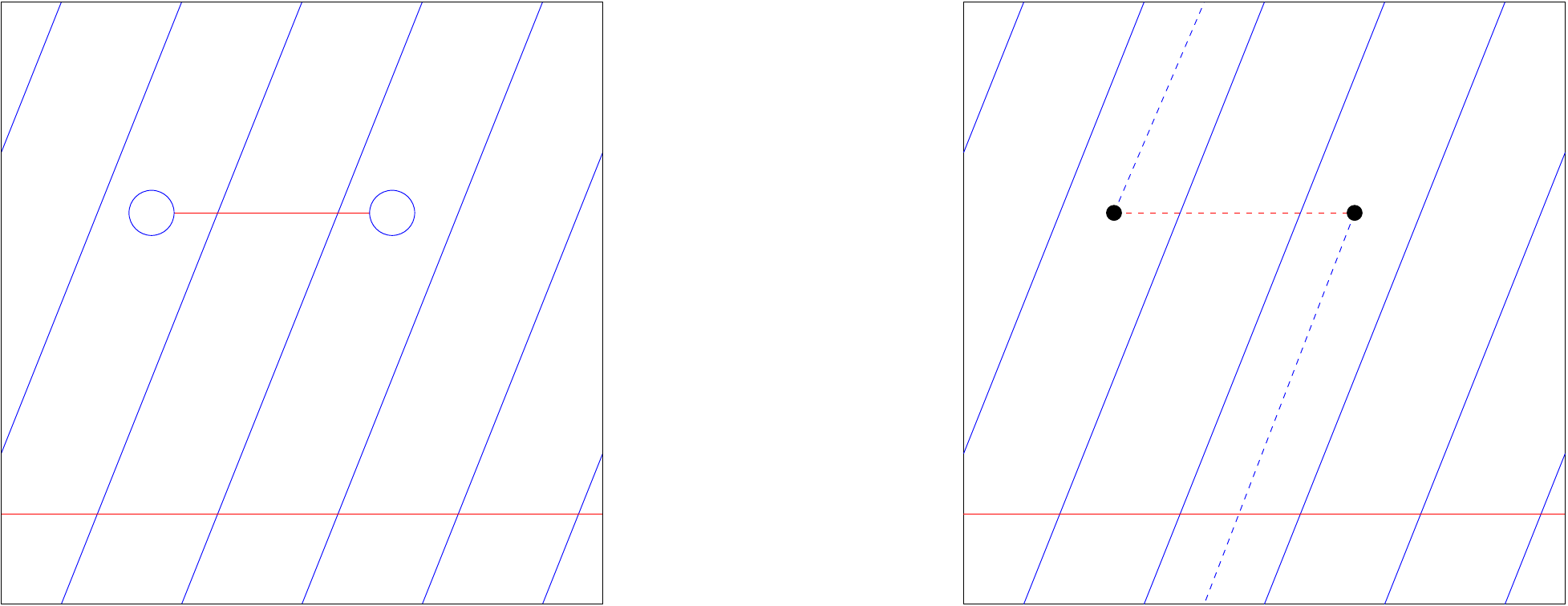}
\caption{At left, the genus-2 Heegaard diagram $(\Sigma, \{ \del D_0, \del D_1 \}, \{\del S, a \})$ of $L'$.
The surface $\Sigma$ results from gluing opposite sides of the square and attaching an annulus to the two circles.
$\del D_0$ is represented by the long horizontal arc and $\del D_1$ by the short horizontal arc.
$\del S$ is represented by the diagonal arcs and $a \simeq a'$ are the two circles.
At right, the induced genus-1 doubly-pointed Heegaard diagram $(\Sigma_1,\del D_0, \del S, z, w)$ of $K' \subset L'$, along with horizontal dashed arc $a_1$ and diagonal dashed arc $b_1$.}
\label{f: simple}
\end{figure}

\begin{proof}
We will embed the solid torus $U'$ in a lens space $L'$ and argue that it is the exterior of a simple knot $K'$.
By the uniqueness of genus-1 Heegaard splittings of lens spaces, it will follow that $K'$ is isotopic to the core of a Heegaard solid torus of $L'$.
By calculating the homology classes of $K'$ and the two cores, we will see that the two outcomes of the Lemma correspond to which of the two cores $K'$ is isotopic.

We keep the notation used in the proof of Lemma \ref{l: solid torus}.
Let $\Sigma = \del H$ and fill $\del M \approx S^2$ with a ball $B$ to form a lens space $L'$.
Set $i=1$ or $2$.
Then $(\Sigma, \{ \del D_0, \del D_i \}, \{\del S, a \})$ is a Heegaard diagram of $L'$.
See Figure \ref{f: simple}.
By construction, $a$ meets $\del D_i$ in a single point and is disjoint from $\del D_0$.
It follows that $(\Sigma, \{ \del D_0, \del D_i \}, \{\del S \})$ is a Heegaard diagram for the exterior of a knot $K' \subset L'$, where $\N(K')$ is the union of $B$ and the $2$-handle attached along $a$.

We produce a doubly-pointed Heegaard diagram for $K' \subset L'$ in the following standard way.
Destabilize the Heegaard diagram by removing the pair $(\del D_i,a)$ and compressing $\Sigma$ along $a$.
The result, for either value of $i$, is a standard genus-1 Heegaard diagram $(\Sigma_1, \del D_0, \del S)$ of $L'$.
$\Sigma_1$ inherits a pair of disks from the compression.
Place basepoints $z$ and $w$, one in each disk.
The result is a doubly-pointed Heegaard diagram $(\Sigma_1, \del D_0, \del S,z,w)$ of $K' \subset L'$.
See Figure \ref{f: simple}.
%We see at once that $K'$ is a simple knot in $L'$ (reference).

We obtain the homology class $[K'] \in H_1(L')$ up to sign, as $K'$ is unoriented, in the following standard way (compare \cite{greenelevine2016strongLSpaces}).
Connect $z$ and $w$ by an arc $a_1$ in the complement of $\del D_0$ that meets $\del S$ in the fewest number of points.
This value equals $\min \{s_1,s_2\} \le s/2$.
Similarly, connect $z$ and $w$ by an arc $b_1$ in the complement of $\del S$ that meets $\del D_0$ in the fewest number of points.
Again, this value is bounded above by $s/2$.
We have $|a_1 \cdot \del S| = |a_1 \cap S|$ and $|b_1 \cdot \del D_0| = |b_1 \cap \del D_0|$.
Write $L' = U_1 \cup_{\Sigma_1} U_2$, where $U_1$ denotes the handlebody containing $D_0$ and $U_2$ denotes the handlebody containing $S$.
Let $c_i$ denote the core circle of $U_i$, $i=1,2$.
Then $[K'] = |a_1 \cap \del S| \cdot [c_2]$ and $[K'] = |b_1 \cap \del D_0| \cdot [c_1]$ up to sign in $H_1(L') \approx \bZ / s \bZ$.

%The meridian of $K'$ is represented by $a \simeq a'$, the two components of $\del A$.

On the other hand, $(\Sigma, \{ \del D_0, \del D_i \}, \{\del S \})$ is a Heegaard diagram for $U'$.
Thus, the exterior of $K'$ in $L'$ is a solid torus.
A genus-1 Heegaard splitting of $L'$ is unique up to isotopy (\cite{bonahonotal1983scindements}).
It follows that $K'$ is isotopic and hence homologous to one of the cores $c_1$ or $c_2$.
Therefore, either $|a_1 \cap \del S| \equiv \pm1 \pmod s$ or $|b_1 \cdot \del D_0| \equiv \pm1 \pmod s$.
As $|a_1 \cap \del S|, |b_1 \cap \del S_2| \le s/2$, it follows that either $\min \{s_1, s_2 \} = |a_1 \cap \del S| = 1$ or else $|b_1 \cap \del D_0|=1$.
The second possibility implies the existence of an arc that runs from one boundary component of $\mathcal{A}$, through an outer corner of $\del_1 T$, across a rectangle of $\mathcal{A}_{\pm1}$, and out of the opposite corner of $\del_{-1} T$ to the other boundary component of $\mathcal{A}$.
Thus, the second possibility leads to the second conclusion of the Lemma.
% In the first case, $1 \in \{s_1,s_2\}$, and in the second case, the other conclusion of the Lemma follows.
\end{proof}

%$\spadesuit$ {\em (An alternate proof of Lemma~\ref{l: opposite corners} is to apply Lemma 4.2 of cite[Gordon-Luecke Non-integral toroidal surgeries] using Theorem 1.2 of cite[Osborne-Zieschang] to show that $S$ cannot be a bad binary face (If $S$ were bad then it would give rise to a non-primitive relator in the two generator group $\pi_1(H)$ and $U'$ would not be a solid torus.)}

\begin{lem}\label{l:schycleweights}
Let $e_i$ denote the isotopy class in $T$ containing $s_i$ edges of $S$ for $i=1,2$. Then $h_U \cup e_i$ has winding number $s_j$ in $H_1(U)$, where $\{i,j\}=\{1,2\}$. 
\end{lem}

\begin{proof}
Orient $h_U$ from $B$ to $A$.
Orient $e_i \subset T$ from $A$ to $B$, and extend it to an oriented arc $f_i$ in the capped-off torus $\hat{T}$ with the same endpoints as $h_U$.
Let $x_i = [h_U \cup f_i] \in H_1(U)$, $i=1,2$.
The disk $S$ gives the nullhomology $s_1 x_1 + s_2 x_2 = 0$.
Hence $x_1 = s_2 t \mu$ and $x_2 = -s_1 t \mu$ for some integer $t$.
On the other hand, $x_1-x_2 = [a] \in H_1(U)$, for the implied orientation on $a$, and $[a] = \pm s \mu$ by Lemma \ref{l: solid torus}.
it follows that $(s_1+s_2)t \mu = \pm s \mu$.
Hence $t=\pm1$ and we obtain the desired conclusion.
\end{proof}

\begin{lem}\label{l:vertices} The vertices of $G_Q$ come in groups of $s$, each lying in a cycle of $s$ edges of $G_Q$.  That is, $G_Q$ is union of such {\em s-cycles}.
\end{lem}

\begin{proof}
Consider a cycle of $G_Q$ consisting of $j$ vertices and edges of $G_Q$. Its push-off on the $U$-side can be isotoped to lie on the boundary of the genus-two handlebody $\N(\mathcal{A} \cup h_U)$, running $j$ times over the $1$-handle $\N(h_U)$ and disjoint from $\mathcal{A}$. Likewise, the boundary of $S$ is a curve that can be taken to lie on the boundary of $\N(\mathcal{A} \cup h_U)$, running $s$ times over the $1$-handle $\N(h_U)$ and disjoint from $\mathcal{A}$. That is, if we attach a 2-handle to $\N(\mathcal{A} \cup h_U)$ along $a$ we would see both curves on the boundary of solid torus winding $j$, $s$ times respectively around the torus. As these curves are disjoint, we conclude $j=s$.
\end{proof}

As $T$ has no simple closed curves of intersection with $Q$, each face $f$ of $G_Q$ is a planar surface where,
%We say that $f$ {is of type} {\bf n} if it is a planar surface with $n$ boundary components. 
as described in the lemma above, each boundary component of $f$ is an $s$-cycle of edges and `half' vertices. Recall that $U,V$ are the Heegaard solid tori on either side of $\hat{T}$, with $U$ being the side containing $S$. Define a {\em $U$-face} (resp.~{\em $V$-face}) of $G_Q$ to be one on the $U$-side (resp.~$V$-side) of $\hat{T}$. A $U$-face (resp.~$V$-face) $f$ of $G_Q$ can be viewed as a properly embedded planar surface in $U-\N(K)$ (resp.~$V-\N(K)$); furthermore, because an orientation on $Q$ induces coherent orientations of the components of $\partial Q$ on $\partial \N(K)$, an orientation on $f$ induces orientations on $\partial f$ that run coherently over $\N(K) \cap U$ (resp.~$\N(K) \cap V$).

\begin{lem}\label{lem:cu2s}
Every $U$-face of $G_Q$ is a disk. There is a single face $V$-face of $G_Q$, and it has $\upsilon$ boundary components, where $q=\upsilon s$. 
%The winding number of $a$ in $U$ is $s$.
The winding number of $a$ in $V$ is a positive multiple of $\upsilon^2 s$. When $\upsilon=1$, the winding number of $a$ in $V$ is exactly $s$.
\end{lem}

\begin{proof}
%Let $H,S$ be as above and $a$ a component of $\partial \mathcal{A}$. By Lemma~\ref{l: solid torus}, $a$ has winding number $s$ in $U$, as claimed.
The curve $a$ is not meridional in $V$, by a similar argument to the first paragraph of the proof of Lemma~\ref{l: solid torus}: we spell it out, and it is even simpler in this case.
Otherwise, there is a disk $\Delta$ in $L(p,m)$ whose intersection with $H \cup S$ is $a$. Then $\N(\Delta) \cup H \cup \N(S)$ is a punctured lens space in $L(p,m)$  of order $s \leq q<p=kq^2$ in homology, a contradiction.

Attach a solid torus $\mathcal{V}$ to $U$ so that $a$ bounds a meridian disk $D$ of $\mathcal{V}$. Then $N=U \cup \mathcal{V}$ is  is a lens space of order $s$ in homology.  
Let $L$ be a core of the solid torus $D \cup H$ in $N$, and let $f$ be a $U$-face of $G_Q$ with $\upsilon$ boundary components. Then $f$ gives a connected Seifert surface for $L$ in $N$. As $f$ has $\upsilon$ boundary components, each of which winds $s$ times around $L$, Lemma 2.13 of \cite{calegarigordon2013smallRationalGenus} says that $L$ has order $\upsilon s$ in $N$. Thus $\upsilon=1$. That is, all $U$-faces of $G_Q$ are disks.

Thus there is a single $V$-face, $f'$, which has $\upsilon$ boundary components, where $\upsilon s=q$. Note that if $\upsilon=1$, the same argument as for $S$ in $U$ shows that $a$ has winding number $s$ in $V$, as claimed. We finish by showing for $\upsilon>1$, the winding number of $a$ in $V$ is a positive multiple of $\upsilon^2s$. Let $h_V=K \cap V$ and $H' = \N(\mathcal{A} \cup h_V) \cap V$. Attach a solid torus $\mathcal{U}$ to $V$ so that $a$ bounds a meridian disk, $D'$, of $\mathcal{U}$. Because $a$ is not a meridian in $V$, the space $N'=V \cup \mathcal{U}$ is either $S^3$ or a lens space, and $|H_1(N')|$ equals the winding number of $a$ in $V$. Let $L'$ be a core of the solid torus $D' \cup H'$ in $N'$, and let $X_{L'}$ be the exterior of $L'$ in $N'$. Then $f$ gives a Seifert surface for $L'$ in $N'$. This Seifert surface has $\upsilon$ boundary component of slope $\lambda$ on $L'$, where $\lambda$ winds $s$ times around $L'$. Now $[\lambda] \in H_1(X_{L'})$ has order which is a factor of $\upsilon$. By Lemma 2.13 of \cite{calegarigordon2013smallRationalGenus} (via Dehn filling $X_{L'}$ along a slope distance one from $\lambda$), it has order exactly $\upsilon$. Note then that order of the torsion summand of $H_1(X_{L'})$ is a multiple of $\upsilon$. Apply Lemma 3.2 of \cite{watson2012khovanovObstructions}, where $M=X_{L'}, \,\lambda_M=\lambda, \,|l|=\upsilon$, and $\alpha$ is the meridian of $L'$ in $N'=X_{L'}(\alpha)$. We conclude that  $|H_1(N')|=c_M s$ where $c_M$ is a product of $\upsilon$ and the order of the torsion subgroup of $H_1(X_{L'})$. That is, $|H_1(N')|$ is a multiple of $\upsilon^2s$; consequently so is the winding number of $a$ in $V$.
\end{proof}

\begin{lem}
\label{lem:annulus1}
If $p=kq^2$, then $s=q$, and there is a value $0 < d < q, \, (d,q)=1$, such that $m=kdq \pm 1$. 
%The component $a$ of $\partial \mathcal{A}$ represents the $(q, \pm d)$-torus knot in $L(p,m)=L(kq^2, kdq \pm 1)$ (resp.).
\end{lem}

\begin{proof}
Let $\upsilon=q/s$ be the number of boundary components of the $V$-face of $G_Q$. We first show that $\upsilon=1$, implying that $s=q$, as claimed.

Let $\mu_U$ (resp.~$\mu_V$) be the canonical generator of $H_1(U)$ (resp.~$H_1(V))$ in $L(p, m)$. Then $\{\mu_U, \mu_V\}=\{\mu, \mu'\}$ in the language of section~1, and $\mu_V=m^\delta \mu_U$ in $H_1(L(p,m))$ for some $\delta \in \{-1,1\}$. By Lemmas~\ref{l: solid torus} and ~\ref{lem:cu2s}, $\pm [a]=s\mu_U \in H_1(U)$ and $\pm [a]=hs\upsilon^2\mu_V \in H_1(V)$ for some positive integer $h$; furthermore, $h=1$ when $\upsilon=1$. Thus from $H_1(L(kq^2,m))$ we get the congruence
\[s \equiv \pm hs \upsilon^2 \cdot m^\delta \pmod{kq^2}.\]

As $q=\upsilon s$, this becomes the congruence \[ 1\equiv \pm h \upsilon^2 \cdot m^\delta \pmod{k\upsilon q}.\]

This implies that $\upsilon=1$, else the left side of the congruence is coprime to the modulus but the right side is not. Thus we conclude $s=q$ as claimed. Furthermore, we conclude that $h=1$ and the above congruence now becomes
\[ m\equiv \pm 1 \pmod{kq}.\]
As $p>m>0$, there is a $q>d>0$ such that $m=kdq \pm 1$. 

Now $a$ is a curve on the Heegaard surface of $L(p,m)=L(kq^2, kdq \pm 1)$ with winding number $q$ in each Heegaard torus. This implies that $q$ and $d$ are relatively prime.
\end{proof}

We are now prepared to prove Theorem~\ref{thm:annulus}, the purpose of this section.

\begin{proof}
[Proof of Theorem~\ref{thm:annulus}]
First, apply Lemma~\ref{lem:alphaweightr}.
Then $\Gamma_b$ contains two isotopy classes of arcs, which lie in an annulus $\mathcal{A}$ on $\hat{T}$, and $\Gamma_b$ contains no simple closed curves.
By Lemma~\ref{lem:annulus1}, for $0 < d < q, \, (d,q)=1$, $m=kdq \pm 1$ and $s=q$.
In addition, the edges of $G_Q$ form a single cycle of length $s$, that bounds two disk faces, the $U$-face $S$ and a $V$-face $S'$.

Let $e_1, e_2$ be the two isotopy classes of arcs in $\Gamma_b$ with $s_1,s_2$ edges of $S$, respectively. Let $h_U=K \cap U$, $U'=\N(\mathcal{A} \cup h_U \cup S) \cap U$ and $U''=U-int(U')$. Similarly, define $h_V=K \cap V$, $V'=\N(\mathcal{A} \cup h_V \cup S') \cap V$ and $V''=V- int(V')$. By Lemma~\ref{l: solid torus} and its $S'$-analog, we conclude that a component of $\partial \mathcal{A}$ has winding number $q$ in both $U'$ and $V'$ and has winding number $1$ in $U''$ and $V''$. Let $W'=U' \cup V'$ and $W''=U'' \cup V''$. Note that $W''$ is a solid torus. 

%Each of $S,S'$ is a schycle of order $s=q$. Let $s_1, s_2$ be the number of edges in the two isotopy classes of edges. 

First assume $\min\{s_1,s_2\}=1$.
Without loss of generality, we take $s_1=1$. By Lemma~\ref{l:schycleweights} 
and its $S'$-analog,  $h_U \cup e_2$ (similarly $h_V \cup e_2$) has winding number $1$ in $U$ (resp. $V$).
Without loss of generality, assume $U=V_\alpha, V=V_\beta$.
%Then the $\alpha$-weight of $e_2$ is $\pm 1$, and the $\beta$-weight of $e_2$ is $\pm 1$. By Lemma~\ref{lem:decomposeK}, we have that $\pm [K]=\mu \pm \mu'=[1 \pm (kdq\pm1)] \mu$. 
Because $[K] = [h_U \cup e_2] + [h_v \cup e_2]$, it follows that $\pm [K] = \mu \pm \mu' \in H_1(L(p,m))$ (signs independent of one another).
Thus, $\pm [K] = [1 \pm (kdq \pm 1)] \mu$.
As $[K]$ has order $q$ in $H_1(L(p,m))$, we conclude that $\pm [K]=kdq \mu$ as claimed.
%Moreover, since $p=kq^2$ and $[K]$ has order $q$ in $H_1(L(p,m))$, we conclude that $(d,q)=1$.
That is, $L(p,m)=L(kq^2,kdq \pm 1)$, $(d,q)=1$, and $\pm [K]=kdq \mu$, which gives the conclusion Theorem~\ref{thm:annulus}(1).

We hereafter assume $s_1,s_2>1$. By Lemma~\ref{l: opposite corners},
each of the faces $S,S'$ has opposite corners that are contiguous. Indeed, we may take these opposite corners to be the same for $S$ and for $S'$.  Pick points $p_1,p_2$ on different components of $\partial \mathcal{A}$. Let $H=\N(\mathcal{A} \cup h_U) \cap U$. Because $S$ has contiguous opposite corners, there is an arc $\sigma_U$ that runs from $p_1$ to $p_2$ in $\partial H - int(\mathcal{A})$ and is disjoint from $S$ (running over the rectangle connecting the contiguous corners). In particular, $\sigma_U$ runs from $p_1$ to $p_2$ on $\partial U'-\mathcal{A}$. Similarly, the contiguous corners of $S'$ allow an arc $\sigma_V$ that runs from $p_1$ to $p_2$ along $\partial V' - \mathcal{A}$.
Let $\sigma=\sigma_U \cup \sigma_V$.
Then $\sigma$ is a knot in $L(p,m)$ that lies on $\partial W'$.
At the same time, $\sigma$ is isotopic in $L(p,m)$ to $K$; furthermore, under this isotopy the framing of $K$ given by a component of $\partial Q$ matches the framing on $\sigma$ induced by $\partial W'$, since the opposite corners are the same for $S$ and $S'$.
%why was this remarked --
%In particular, Dehn surgery on $\sigma$ along the framing by $\partial W'$ is $S^1 \times S^2$. 
%---and why could it not be a connected sum?
The curve $\sigma \subset \del W' = \del W''$ does not bound a disk in $W''$, for then $K = \sigma$ would be an unknot, in contradiction to the assumption that $q > 1$.
%Then $\sigma$ lies on $\partial W''=\partial W'$.
%As $K$ does not bound a disk in $L(p,m)$, $\sigma$ is not a meridian of the solid torus $W''$.
Thus, the winding number of $\sigma$ in $W''$ is $w > 0$.
Let $L$ be the core of $W''$.
Then $K=\sigma$ is the $w$-fold cable of $L$ in $L(p,m)$.
As $a$ has winding number $1$ in $W''$, we may isotop $L$ in $W'' \subset L(p,m)$ to a component of $\partial \mathcal{A}$.
That is, $L$ is isotopic to a curve on the Heegaard torus of $L(p,m)$ with winding number $q$ in $U, V$ by Lemma~\ref{lem:cu2s}.
In particular, $X_L$ is Seifert fibered over the disk with two singular fibers of order $q$. 
%Note that by Lemma~\ref{lem:annulus1}
%$L$ represents the $(q, \pm d)$-torus knot in $L(p,m)=L(kq^2, kdq \pm 1)$.

By \cite[Lemma VI.7]{jaco1980lecturesThreeManifoldTopology}, an orientable Seifert fibered space is irreducible unless it is $S^1 \times S^2$ or $P^3 \# P^3$. If $w=1$, then $K=L$. We may assume $k>1$ (else Theorem~\ref{thm:annulus}(2) holds). As $K(\lambda)$ is reducible, it cannot be Seifert fibered, so $\lambda$ must be the slope of the Seifert fiber of $X_K=X_L$ and $K(\lambda)$ is the connected sum of lens spaces of order $q$, contradicting the assumption that $K(\lambda)$ has a non-separating sphere. Thus $w>1$, and $K$ is a $w$-fold cable on $L$. Again using the fact that $K(\lambda)$ cannot be Seifert fibered or the connected sum of two lens spaces of order $q$, we see $\lambda$ is the slope on $K$ of the cabling annulus. Surgery along this slope produces a lens space summand of order $w$. As $K(\lambda) = (S^1 \times S^2) \# N$ where $N$ is a lens space of order $k$, we conclude that $w=k$. Thus $K$ is the $k$-fold cable on the $(q,q)$-torus knot in $L(p,m)=L(kq^2, kdq \pm 1)$. As $L$ has winding number $q$ in each Heegaard solid torus,  $[K]=\pm k[L]=\pm kq \mu$, as claimed.

 Noting that $L(kq^2,kdq-1)$ is diffeomorphic to $L(kq^2,kdq+1)$ preserving the stated homology classes up to sign and the $(q,q)$-torus knot (see the proof of Theorem~\ref{thm:lisca}), we get the statement of Theorem~\ref{thm:annulus}(2).
%Again, by Lemma~\ref{l: opposite corners}, each of the faces $S,S'$ has opposite corners that are contiguous. This implies that $s_1 \neq s_2$. WLOG we assume that $s_1 > s_2$. Let $e_1$ be the isotopy class in $\Gam_b$ corresponding to the $s_1$ edges of $S$. By Lemma 3.5 of [GLL], $h_U \cup e_1$ has winding number $s_2$ in $U$, and $h_V \cup e_1$ has number $\pm s_2$ in $V$. By Lemma~\ref{lem:decomposeK} of [GLL], we have that $\pm [K]= s_2\mu \pm s_2 \mu'=s_2 [1 \pm (rdq\pm1)] \mu$. As $[K]$ has order $q$ and $r>1$, we conclude that $\pm [K]=s_2drq \mu$. Thus, with the above, we conclude $rq \equiv \pm s_2drq \pmod{rq^2}$. That is, $1 \equiv s_2d \pmod{q}$. 
%Furthermore, an $(r,t)$ surgery on $L$ produces $S^1 \times S^2$, as claimed. ($K$ is a copy of the slope of the Seifert surface for $L$ and wraps $r$ times around $L$. Can this be simple?)
%As this lives as a submanifold in $S^1 \times S^2$, we conclude that $\sigma$ has winding number $1$ in $W''$. As $\partial \mathcal{A}$ also has winding number $1$ in $W''$, we may isotop $K=\sigma$ in $W'' \subset L(q^2,m)$ to a component of $\partial \mathcal{A}$. Then $K$ has winding number $q$ in $U$. In particular, $[K]=\pm q \mu$, as claimed.
\end{proof}

\begin{cor}\label{cor:nocircletangciesr}
Under the assumptions set out at the beginning of the section, either the conclusion of Theorem~\ref{thm:annulus} holds, or else no circle tangencies occur in the level presentation of $Q$.
\end{cor}

%\begin{cor}\label{cor:nocircletangciesr}
%Either
%\begin{enumerate}[label=(\alph*)]
%\item $p=kq^2, m=kdq + 1$, and $\pm [K]=kq^2 \mu, kdq \mu$; or
%\item There are no circle tangencies in the level presentation of $Q$.
%\end{enumerate}
%\end{cor}

\begin{proof} 
Suppose that the conclusion of Theorem~\ref{thm:annulus} does not hold.
Then the hypothesis of Theorem~\ref{thm:annulus} does not hold.
In tandem with Lemma~\ref{lem:alphaweightr}, this implies that every graph $\Gamma_b$ with $0 < b < 2q-2$ contains at least three isotopy classes of arcs.
It follows that no level surface between levels $1$ and $(2q-2)-1$ can contain a circle tangency, by Lemma~\ref{lem:morseposition}$(b)$. %Lemma~\ref{lem:alphaweightr} then verifies Corollary~\ref{cor:nocircletangciesr} for any level between $\Gamma_1$ and $\Gamma_{2q-3}$.
For any level surface between $\Gamma_0$ and $\Gamma_1$, note that any circle tangency must be trivial in $T$ or parallel to the meridian $\alpha$, hence trivial in $X_K$.
Again, Lemma~\ref{lem:morseposition}$(b)$ gives a contradiction. Similarly, there can be no circle tangency between $\Gamma_{2q-3}$ and $\Gamma_{2q-2}$.
\end{proof}

%Similarly, following the proof of Lemma~\ref{lem:alphaweightr} using Lemma~\ref{lem:betaweightTxi} in place of Lemma~\ref{lem:alphaweightTxi}, we get the following

%\begin{lem}\label{lem:betaweight}
%Either 
%\begin{enumerate}[label=(\alph*)]
%\item $p=kq^2, m=kdq \pm 1$, and $\pm [K]=kq \mu, kdq \mu$; or
%\item For every $i \neq 0, 2q-2$, $\Gamma_i'$ has at least three isotopy classes of level arcs. For every $i=0, \dots, 2q-2$ there are no simple closed level curves in $\Gamma_i'$, and the sum of the $\beta$-weights of the arcs of $\Gamma_i'$ is $0$ in ${\mathbb Z}$.
%\end{enumerate}
%\end{lem}

\section{When \pdfmath{K(\lambda)}{K(λ)=(S1xS2)\#N} is not prime}\label{sec:r>1}

%\section{When \pdfmath{K(\lambda)=(S^1 \times S^2) \# N, N \neq S^3}{K(λ)=(S1xS2)\#N}}\label{sec:r>1}

The goal of this section is to prove Theorem \ref{thm:mainNextended}.
The hypotheses of Theorem~\ref{thm:mainNextended} are that $K$ is a simple knot in $L(p,m)$, $Q$ is a planar Seifert surface for $K$ of longitudinal slope $\lambda$ and with $q > 1$ boundary components, and $K(\lambda)$ is not prime.
We assume these hypotheses for the duration of this section.
We seek to apply Theorem~\ref{thm:annulus}, whose conclusion matches that of Theorem~\ref{thm:mainNextended}.
Thus, we must argue that its hypothesis is met: there exists an index $0 < b < 2q-2$ such that the arcs of $\Gamma_b$ belong to at most two isotopy classes.
The idea is that the presence of the assumed separating, essential two-sphere places significant constraints on $\Gamma$, leading to this conclusion.

%\newcounter{savedthm4}
%\newcounter{savedsec1}
%\setcounter{savedthm4}{\value{thm}}
%\setcounter{savedsec1}{\value{section}}
%\setcounter{thm}{\value{mainr>1}}
%\setcounter{section}{\value{secr>1}}
%\begin{thm}
%Let $K \subset L(p,m)$ be a simple knot with a planar Seifert surface with $q$ boundary components and longitudinal slope $\lambda$. Then either
%\begin{enumerate}
%    \item $K(\lambda)=S^1 \times S^2$; or
%    \item  Up to diffeomorphism, possibly orientation-reversing, there are integers $k>1$, $0<d<q$, such that $K(\lambda)=(S^1 \times S^2) \# L(k,1)$, $p=kq^2$, and $m=kdq + 1$.
%    In this case, either
%    \begin{enumerate}
%    \item $\pm [K]=kdq \mu \in H_1(L(p,m))$ with $(q,d)=1$; or
%    \item $K$ is an $k$-fold cable on a companion, $L$, that lies on the Heegaard torus of $L(p,m)$; furthermore $L$ has winding number $q$ in both Heegaard solid tori. $\pm [K]=kq \mu \in H_1(L(p,m))$.
%    \end{enumerate}
%\end{enumerate}     
%\end{thm}

%\setcounter{thm}{\value{savedthm4}}
%\setcounter{section}{\value{savedsec1}}

%\begin{remark} The assumption that $\lambda$ is longitudinal should not be necessary, but the conclusions will differ. The same argument should pretty much apply.
%\end{remark}

%Note $X_K$ is irreducible.
%Up to orientation-reversing diffeomorphism we may assume $p,m>0$.
By Lemma~\ref{lem:p,N}, there is a value $k>0$ such that $p=kq^2$ and $K(\lambda)=(S^1 \times S^2) \# N$, where $N$ is a lens space with $|H_1(N)|=k$. 
Under the assumption that $K(\lambda)$ is not prime, we therefore have $k > 1$.
The planar surface $Q$ is incompressible and $\partial$-incompressible in the knot exterior $X_K$. Let $\hat{Q} \subset K(\lambda)$ be the non-separating $2$-sphere gotten by capping off $Q$ in $K(\lambda)$ with disks. Let $K' \subset K(\lambda)$ be the dual knot to the surgery. Then $K'$ intersects $\hat{Q}$ in $q$ points. 

%\noindent {\bf We hereafter assume $N$ in Lemma~\ref{lem:p,N} is not $S^3$. Note that this implies that {\bf $K$ is not fibered} and $k>1$.}

\begin{lem}\label{lem:separatinginc} There is a connected, separating, incompressible, $\partial$-incompressible planar surface $P$ in $X_K$, disjoint from $Q$, with boundary slope $\lambda$.
\end{lem}

\begin{proof} 
Let $\hat P$ be the unique (up to isotopy) separating, essential sphere in $K(\lambda)$ disjoint from $\hat{Q}$. Isotop $K'$ so that it intersects $\hat P$ minimally, and let $n=|K' \cap \hat{P}|$.

Let $\tau$ be an arc of $K' - \cN(\hat{Q})$. Then $R=\partial (\cN(\hat{Q} \cup \tau))$ is an essential separating $2$-sphere disjoint from $\hat{Q}$. Thus $R$ is isotopic to $\hat{P}$ and $|K' \cap R|= 2q-2$. Thus $n \leq 2q-2$.

Let $\overline{P}=\hat{P} \cap X_K$. Note that $\partial \overline{P}$ has $n$ components. 
If $D$ is a (boundary-)compressing disk of $\overline{P}$, then, since $Q$ is irreducible and boundary-irreducible, we may isotope $D$ disjoint from $Q$.
In this way, we maximally compress $\overline{P}$ in $X_K$, resulting in a surface $P'$ which is disjoint from $Q$. $P'$ can contain no boundary parallel annuli, as such would supply an isotopy of $K$ intersecting $\hat P$ fewer times.
We are done unless each component of $P'$ is non-separating or is a $2$-sphere which bounds a 3-ball in $X_K$ (by the irreducibility of $X_K$).
As $\hat{P}$ is separating, $P'$ must have at least two non-separating components.
Thus there is a non-separating component $P''$ of $P'$ with fewer than $q$ boundary components.
After tubing oppositely oriented boundary components of $P''$ along $\partial \cN(K')$ we get a Seifert surface $P'''$ for $K$ with fewer than $q$  boundary components (note that the number of boundary components is non-zero as there is no non-separating surface in $L(p,m)$). But Lemma~2.13 of \cite{calegarigordon2013smallRationalGenus} applied to $Q$ says that $[K]$ has order $q$ in $H_1(L(p,m))$ and applied to $P'''$ says that the order of $[K]$ is less than $q$.

We conclude that some component of $P'$ is the desired separating, incompressible planar surface. 
\end{proof}

%Let $K' \subset (S^1 \times S^2) \# N$ be the knot dual to the surgery and let $\hat{Q}$ be the non-separating $2$-sphere gotten by capping off $Q$ in $K(\lambda)$. 
%Let $\hat{P}$ be a separating $2$-sphere in $K(\lambda)$ disjoint from $\hat{Q}$ that intersects $K'$ the fewest times. 
Let $P$ be the separating, incompressible, planar surface in $X_K$ disjoint from $Q$ given by Lemma~\ref{lem:separatinginc}.
 Let $\rho$ be the number of components of $\partial P$.
 Note that an orientation on $Q$ induces coherent orientations of $\partial Q$ on $\partial X_K$; whereas, as $P$ is separating, an orientation on $P$ induces orientations on $\partial P$ that alternate in direction along $\partial X_K$. The $\rho/2$ components of $\partial P$ in one direction we call {\em positive} and the remaining $\rho/2$ in the other direction are {\em negative}.

We follow the set-up of sections~\ref{sec:setup} and \ref{sec:atleast3}.
%[Josh] I don't understand this remark:
%In particular, we take $Q$ to be that coming from \cite{greene2022fiberedSimpleKnots}.
Let $V_{\alpha} \cup_{\hat{T}} V_{\beta}$ be the canonical Heegaard splitting of $L(p,m)$, where $\alpha, \beta$ is a meridian of $V_{\alpha}, V_{\beta}$ (resp.). As $K$ is simple in $L(p,m)$ we choose cores $c_{\alpha}, c_{\beta}$ of $V_{\alpha}, V_{\beta}$ that each intersect $K$ in a single point. Let $C_{\alpha}, C_{\beta}$ be small neighborhoods of $c_{\alpha},c_{\beta}$ such that 
$L(p,m)- (C_{\alpha} \cup C_{\beta}) \cong \hat{T} \times [0,1]$. Let $D_{\alpha}$ ($D_{\beta}$) be a meridian disk of $C_{\alpha}$ ($C_{\beta}$) containing $K \cap C_{\alpha}$ ($K \cap C_{\beta}$). As in section~\ref{sec:setup}, $Q \cap C_{\alpha}$ ($Q \cap C_{\beta}$) is a collection of half-disks lying in $\N(D_{\alpha})-\N(K)$ ($\N(D_{\beta})-\N(K))$). We may assume the same for $P$ except that $P$ may also intersect $C_{\alpha} \cup C_{\beta}$ in meridional disks of $C_{\alpha}, C_{\beta}$ away from $K$. We isotop $Q \cup P$ in $X_K$ to minimize the number of components (i.e. disks and half-disks) of intersection with $C_\alpha \cup C_\beta$. 

Recall that $T=\hat{T}-\N(K)$ is a twice-punctured torus, and $X_K \cap (\hat{T} \times [0,1]) \cong T \times [0,1]$.
The surface $Q$ is transverse to $T \times \{0\}$ ($T \times \{1\}$), with $q$ arcs of intersection running from one boundary component of $T$ to the other, all lying in an annulus in $\hat{T} \times \{0\}$ ($\hat{T} \times \{1\}$), and representing two different isotopy classes.
The surface $P$ is transverse to $T \times \{0\}$ ($T \times \{1\}$), with $\rho$ arcs of intersection running from one boundary component of $T$ to the other and all lying in an annulus (same annulus as for $Q$) in $\hat{T} \times \{0\}$ ($\hat{T} \times \{1\}$),
% ($\blacklozenge$ \textbf{May say by Lemma~\ref{lem:separatinginc}})$\spadesuit$ {\em (I thought this followed just from the fact that $P$ intersects $C_{\alpha}$ in half-disks and meridional disks. We would want incompressibility to say no trivial simple closed curves of intersection at other levels. Is that what you meant?)}
along with some simple closed curves of intersection all bounding meridians in $C_{\alpha} \cup C_{\beta}$. %Let $P'$ be $P$ with any meridians of $C_{\alpha} \cup C_{\beta}$
We shall argue momentarily that no such closed curves occur.
%Lemma \ref{lem:nomeridians} 

We apply the argument of Lemma~\ref{lem:morseposition}, noting that the closed curves of intersection of $P$ with $\hat{T} \times \{0,1\}$ are essential in $\hat{T} \times \{0,1\}$, to show we can isotop $Q \cup P$ in $T \times [0,1]$, fixing $T \times \{0,1\}$, so that only saddle and circle tangencies appear. Furthermore, the conclusions of Lemma~\ref{lem:morseposition} still hold for $Q$. We similarly note that no curve, $c$, of $P \cap (T \times \{i\})$ can be trivial $T \times \{i\}$. Otherwise, by the incompressibility of $P$, an innermost $c$ must bound a disk $D_P$ in $P$.  As there are no center singularities of $P \cap (T \times [0,1])$, $D_P$ must intersect $C_{\alpha} \cup C_{\beta}$. The irreducibility of $X_K$ allows us to replace $D_P$ with the disk that $c$ bounds in $T \times \{i\}$,  thereby reducing the components of $P \cap (C_{\alpha} \cup C_{\beta})$. After perturbing $P$ to introduce a center singularity then applying 
the argument of Lemma~\ref{lem:morseposition} to eliminate it in $T \times [0,1]$, without increasing  $|P \cap  (C_{\alpha} \cup C_{\beta})|$, we 
obtain a contradiction to the minimality of $P$.

Note that $\partial P \cap (T \times [0,1])$ consists of $\rho$
transverse arcs on each component of $\partial T \times [0,1]$.

\begin{lem}\label{lem:nomeridians}

\begin{enumerate}
\item $T \times \{0, 1\}$ has no simple closed curves of intersection with $Q \cup P$.
\item If $T \times \{j\}$ is transverse to $Q \cup P$, then no simple closed curve of $T \times \{j\} \cap (Q \cup P)$ can be trivial (bound a disk) in $X_K$.
\end{enumerate}

\end{lem}

\begin{figure}
    {
    \includegraphics[width=2in]{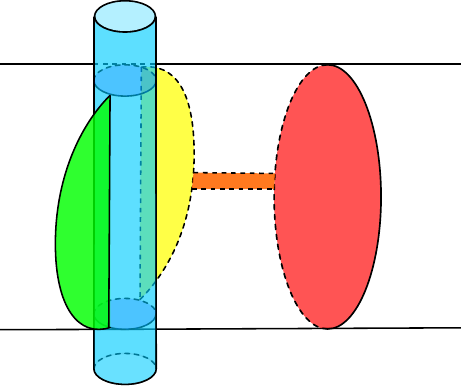}
    \put(-50,55){$D_2$}
    \put(-125,50){$D_3$}
    \put(-97,75){$D_1$}
    \put(-30,110){$T\times \bc{i+\varepsilon}$}
    \put(-76,50){$B$}
    }
    %\quad
    %\raisebox{12mm}{
    %\includegraphics[width=3.5in]{image/res_meridian_on_P.pdf}
    %\put(-20,16){$D_3$}
    %\put(-216,16){$D_1$}
    %\put(-165,16){$D_2$}
    %\put(-40,60){On $P$}
   % }
    \caption{Getting rid of meridians on $P$}
    \label{fig:gettingridofmeridians}
\end{figure}

\begin{proof}

Suppose for a contradiction that $T \times \{0\}$ has a simple closed curve of intersection with $Q \cup P$. Any such comes from an intersection $P \cap c_\alpha$.
Assume the first saddle singularity (starting from $0$) of $Q \cup P$ occurs at level $T \times \{i\}$. 
Pick $\varepsilon > 0$ so that $T \times [i-\varepsilon,i+\varepsilon]$ is transverse to $Q \cup P$ except at $T \times \{i\}$.

If the saddle singularity occurs on $Q$, then there remains a meridional curve of intersection between $P$ and $T \times \{i + \varepsilon \}$.
This forces an arc of $Q$ to band to itself through the singularity and produce a simple closed curve of intersection between $Q$ and $T \times \{i + \varepsilon\}$, parallel to one of $P$.
This curve, bounding a disk in $X_K$, contradicts Lemma~\ref{lem:morseposition}.
Therefore, the first saddle singularity must occur on $P$.

Suppose first that this saddle singularity connects an arc, $\tau_{i-\varepsilon}$, with a simple closed curve, $c_{i-\varepsilon}$, of $P \cap (T \times \{i-\varepsilon\})$ creating an arc, $\tau_{i+\varepsilon}$, of  $P \cap (T \times \{i+\varepsilon\})$. As $P$ is disjoint from $Q$, the non-trivial $c_{i-\varepsilon}$ must be a meridian in the solid torus $(\hat{T} \times [0,i-\varepsilon]) \cup C_{\alpha}$. Then $c_{i-\varepsilon}$ bounds a disk $D_2$ such that $P \cap D_2 = c_{i-\varepsilon}$ and where $D_2$ is a vertical annulus in  $T \times [0,i-\varepsilon]$ in union with a meridian of $C_{\alpha}$.
The arc $\tau_{i-\varepsilon}$ bounds a half-disk $D_1$ in $P$. By the incompressibility of $P$, $c_{i-\varepsilon}$ must also bound a disk $D_P$ in $P$. Note that $D_P$ must intersect $C_{\alpha} \cup C_{\beta}$ in at least one meridian disk (as $P \cap (T \times [0,1]$) has no center singularities). The irreducibility of $X_K$ allows us to isotop $D_P$ to $D_2$, fixing $P - \text{Int}(D_P)$ and $Q$. After the saddle at $T \times \{i\}$, we then see that $\tau_{i+\varepsilon}$ bounds a half-disk in $P$; furthermore, we may isotope this half-disk, to standard form---a half-disk, $D_3$, intersecting $C_{\alpha}$ in $\N(D_{\alpha})-\N(K)$. See Figure~\ref{fig:gettingridofmeridians} where the banding $B$ indicates the result of the saddle.
We have thereby isotoped $P$ into a standard position in which it intersects $C_\alpha$ in fewer components (half-disks and meridian disks) than originally, contradicting minimality.

Suppose second that the saddle singularity on $P$ connects two simple closed curves, $c_1,c_2$ of intersection of $P \cap (T \times \{i-\varepsilon\})$ in moving to $T \times \{i\}$. 
The result would be a trivial simple closed curve in $T \times \{i+\varepsilon\}$---which we have seen cannot happen.

Suppose third that the saddle singularity on $P$ connects distinct arcs of $P$.
Then we would see a boundary-compression of $P$ following the saddle, a contradiction.

Suppose fourth and finally that this first saddle in $T \times \{i\}$, connects a single arc, $\tau_{i-\varepsilon}$, of $P \cap (T \times \{i-\varepsilon\})$ to itself, producing an essential, simple closed curve, $c_{i+\varepsilon}$, and an arc, $\tau_{i+\varepsilon}$, of $P \cap (T \times \{i+\varepsilon\})$. In $P$, $c_{i+\varepsilon} \cup \tau_{i+\varepsilon}$ bounds, along with $\partial P$, an annulus gotten by banding a half-disk (bounded by $\tau_{i-\varepsilon}$) to itself.  On the other hand $c_{i+\varepsilon}$ bounds a disk $D$ such that $P \cap D=c_{i+\varepsilon}$ and where $D$ is a vertical annulus in  $T \times [0,i+\varepsilon]$ in union with a meridian of $C_{\alpha}$. By the incompressibility of $P$, $c$ must also bound a disk $D_P$ in $P$. $D_P$ must intersect $C_{\alpha} \cup C_{\beta}$ in at least one meridian disk (as $P \cap (T \times [0,1]$) has no center singularities). The irreducibility of $X_K$ allows us to isotop $D_P$ to $D$, fixing $P - \text{Int}(D_P)$ and $Q$. The result is a half-disk bounded by $\tau_{i+\varepsilon}$; furthermore, we may isotope this half-disk to standard form---a half-disk intersecting $C_{\alpha}$ in $\N(D_{\alpha})-\N(K)$.  But this reduces the number of components of $P \cap (C_{\alpha} \cup C_{\beta})$---a contradiction. 

%%%%%%%%%The next paragraph is an alternate argument for the fourth case above %%%%%%%%%
%As observed before, $c_{i+\varepsilon}$ must be essential in the level torus, hence a meridian of $V_{\alpha}$. This $c$ then bounds a meridian disk $D$ of $V_{\alpha}$ such that $P \cap D=c$ and such that $D \cap (T \times [0,i+\varepsilon])$ is a vertical annulus. By the incompressibility of $P$, $c$ must also bound a disk $D_P$ in $P$. The irreducibility of $X_K$ allows us to isotop $D_P$ to $D$, fixing $P - \text{Int}(D_P)$ and $Q$. $P$ will have the same number of intersections with $c_{\alpha} \cup c_{\beta}$ before and after this isotopy (it cannot have fewer afterwards). Now we see a circle tangency of $P$ at $c$ in $T \times \{i+\varepsilon\}$. Perturb the circle tangency to be Morse with a center and saddle tangency, calling this saddle $s'$. Isotop $P$ in $T \times [0, i+\varepsilon]$ so that $s'$ occurs before $s$, say at level $T \times \{i-\varepsilon/2\}$, and cancel $s$ with the center tangency. After this isotopy, $s$ becomes the first saddle of $P \cup Q$ and connects an arc with a simple closed curve of $P \cap T \times \{i-\varepsilon\}$. We  apply the argument for this situation given above to reduce the number of intersections of $P$ with $C_{\alpha}$, a contradiction.

Thus $T \times \{0\} \cap (Q \cup P)$ contains no simple closed curves of intersection.  The same argument applies to prove $(1)$ for $T \times \{1\}$.

For $(2)$, assume $T \times \{j\} \cap (Q \cup P)$ contains a simple closed curve $c$ which bounds a disk in $X_K$. Because $Q \cup P$ is incompressible in $X_K$, we conclude that $c$ bounds a disk in $Q \cup P$. As $(Q \cup P) \cap 
(T \times \{0, 1\})$ contains no simple closed curves of intersection, this disk lies in $T \times [0, 1]$. But then $(Q \cup P) \cap (T \times [0, 1])$ has a center singularity in its Morse presentation, a contradiction.
\end{proof}

By Lemma~\ref{lem:nomeridians}$(1)$, $(Q \cup P) \cap C_\alpha$ is a collection of $q$ half-disks in $Q$, one for each component of $\partial Q$, and $\rho$ half-disks in $P$, one for each component of $\partial P$. Similarly for $(Q \cup P) \cap C_\beta$. As $T \times [0,1]$ induces a foliation of $Q \cup P$ without center singularities, we see that $Q$ admits $2q-2$ saddle singularities and $P$ admits $2\rho -2$.

Starting with $T \times \{0\}$, and ending with $T \times \{1\}$, pick transverse levels $\Gam_b$, $b = 0, \dots 2q-2+2\rho-2$,  between the saddle singularities of $P \cup Q$ and call the union of these $\Gam$. $\Gam_b$ is the corresponding level surface along with the arcs and simple closed curves of intersection of $Q \cup P$ with this level surface. 

We collect some facts that follow as in sections~\ref{sec:setup} and \ref{sec:atleast3}. Label the components of $\partial T$ as $A,B$.

\begin{enumerate} 
\item For each $j=0, \dots, 2(q+\rho)-4$, each of the $q$ arcs of $Q \cap \Gam_j$ runs from component $A$ to $B$. 
\item In $\Gam_0$ the $q$ arcs of $\Gam_0 \cap Q$ lie in two non-empty isotopy classes of $\Gam_0$, $\hat{\alpha}, \hat{\alpha}'$ (Definition~\ref{def:alpha1}), whose projection to $\hat{T}$  forms the class of the meridian of $V_{\alpha}$.
The same pertains for the $q$ arcs of $\Gam_{2(q+\rho)-4} \cap Q$ with respect to $C_\beta$.
\item For each $j$, no arc of $\Gam_j$ is incident to the same component of $\partial \Gam_j$ and isotopic rel endpoints into that component (as $Q,P$ are $\partial$-incompressible). 
\item For each $j$, the arcs of $Q \cap \Gam_j$ lie in at least two isotopy classes on $\Gam_j$. (Lemma~\ref{lem:diskface})
\item In $\Gam_0$ the $\rho$ arcs of $\Gam_0 \cap P$ also lie in isotopy classes $\hat{\alpha}, \hat{\alpha}'$, though possibly all in one class. The same pertains for the $\rho$ arcs of $\Gam_{2(q+\rho)-4} \cap Q$ with respect to $C_\beta$. %There may be simple closed curves of $P \cap \Gam_0$ which are meridians of $V_{\alpha}$. Similarly for $\Gam_{2(q+\rho)-4}$.
\item $[K]$ has order $q$ in $H_1(L(p,m))$. 
%\item Any simple closed curves of intersection in $\Gam_b$ are essential in $\Gam_b$ and in $P$ or $Q$.
\end{enumerate}  

\begin{defin}\label{def:G_P} For each level $\Gam_j$ we consider labelled graphs of intersection as in section~\ref{sec:atleast3}: 
\begin{itemize}
    \item $G_P^j$ is $P \cap \Gam_j$ on $P$. We think of $G_P^j$ as a graph in the 2-sphere $\widehat{P}$: the disks filling components of $\partial P$ referred to as (fat) {\em vertices}, the arc components as {\em edges}. The endpoint of each edge is labelled $A$ or $B$ according to the component of $\partial \Gam_j$ in which it lies. Recall that half the components of $\partial P$ are positive and half negative.
We call a vertex of $G_P^j$ positive (negative) if the corresponding component of $\partial P$ is positive (negative). 
\item $G_Q^j$ is $Q \cap \Gam_j$ on $Q$. We think of $G_Q^j$ as a graph in $\widehat{Q}$, the 2-sphere obtained by capping the components of $\partial Q$ with disks: the disks filling components of $\partial Q$ referred to as (fat) {\em vertices}, the arc components as {\em edges}. The endpoint of each edge is labelled $A$ or $B$ according to the component of $\partial \Gam_j$ in which it lies.
\item We can think of $\Gam_j$ as a labelled graph in the same way, with fat vertices $A$ and $B$ and the edges labelled by the components of $\partial Q, \partial P$. Note that fact $(3)$ above says that $\Gam_j$, viewed as a graph, has no 1-sided faces.
\item Any simple closed curves of intersection with $\Gam_j$ are included in $G_Q^j, G_P^j$.
 \end{itemize}
\end{defin}

For each $j$, the level surface $\Gam_j$ has an $\alpha$-side and a $\beta$-side. The faces of $G_P^j$ and $G_Q^j$ are thereby divided into {\em $\alpha$-faces} and {\em $\beta$-faces}. Note that the $\alpha$ faces of $G_P^j$ are contained in the $\alpha$ faces of $G_P^{j+1}$, and the $\beta$ faces of $G_P^{j+1}$ are contained in the $\beta$ faces of $G_P^j$. 

The following uses the fact that $Q,P,\Gam, X_K$ are orientable. 

\noindent {\bf Parity Rule}: For every $j$,
\begin{itemize}
    \item An edge of $G_P^j$  connects vertices of the same sign iff it has label $A$ at one end and label $B$ at the other.
    \item an edge of $G_Q^j$ has label $A$ at one end and label $B$ at the other.
\end{itemize}

\begin{lem}
\label{lem:AtoA}
There is a value $b$ such that $G_P^b$ has either
\begin{itemize} 
\item a simple closed curve of intersection, or
\item two edges connecting vertices of opposite sign, one with label $A$ at both ends and the other with label $B$ at both ends. Both arise from the saddle in $P$ between levels $\Gam_{b-1}$ and $\Gam_b$.
\end{itemize}

\end{lem}

\begin{proof}
We assume $G_P^b$ has no simple closed curves of intersection for any $b$.
$G_P^0$ and $G_Q^0$ has a collection of $1$-sided $\alpha$-faces, one based at each vertex. That is, for each (fat) vertex, there is an edge incident at both ends which is parallel into that vertex. The complement of these $1$-sided faces is a single $\beta$-face. Similarly, $G_P^{2(q+\rho)-4}$ and $G_Q^{2(q+\rho)-4}$ has a collection of $1$-sided $\beta$-faces, one based at each vertex, and one $\alpha$-face. In particular, in $G_P^{2(q+\rho)-4}$ the $\alpha$-face is incident to vertices of both signs. On the other hand, in $G_P^0$, there is no $\alpha$-face incident to vertices of both signs. Let $b$ be the smallest value such that $G_P^{b}$ has an $\alpha$-face, $f$, which is incident to vertices of both signs. Then going from $G_P^{b}$ to $G_P^{b-1}$ must involve a saddle of $P$ occurring within $f$. If there is an edge of $G_P^{b}$ connecting vertices of opposite sign, we take $f$ to be the face on the $\alpha$-side of this edge. By the Parity Rule, this edge must have the same label on both ends, say label $A$. Then along this same boundary component of $f$ will be an edge with label $B$ at both ends and connecting vertices of opposite sign. By definition of $b$, the saddle between level $b-1$ and $b$ must involve these edges, and these edges are as claimed. Thus we assume all edges of $G_P^{b}$ connect vertices of the same sign. 
%By the Parity Rule, each edge of $G_P^{b}$ has label $A$ at one end and label $B$ at the other. 
The saddle must connect two edges $e_1, e_2$ of $f$. If $e_1$ and $e_2$ both connect, say, positive vertices, then the resulting $\alpha$-face in $G_P^{b-1}$ will still have vertices of both signs incident, contradicting the choice of $b$. Thus $e_1$ connects, say, positive vertices, and $e_2$ connects negative vertices. The resulting edges $e_1',e_2'$ in $G_P^{b-1}$ will each connect a positive vertex to a negative one. But the $\alpha$-face on one side of, say, $e_1'$ will contradict the definition of $b$.
\end{proof}

\begin{defin}
Let $b_Q$ be the smallest value such that there is a saddle of $Q$ between the levels $\Gam_{b_Q-1}$ and $\Gam_{b_Q}$. Similarly, let $b'_Q$ be the largest value such that there is a saddle of $Q$ between the levels $\Gam_{b'_Q-1}$ and $\Gam_{b'_Q}$. Note that $\Gam_{b_Q-1} \cap Q$ is isomorphic as a labelled graph to $\Gam_0 \cap Q$. %and $\Gam_b \cap Q$ is isomorphic to $\Gamma_1 \cap Q$. 
Similarly, 
%$\Gam_{b'_Q-1} \cap Q$ is isomorphic as a labelled graph to $\Gamma_{2q-3} \cap Q$, and
$\Gam_{b_Q'} \cap Q$ is isomorphic to $\Gam_{2(q+\rho)-4} \cap Q$. 

\end{defin}

\begin{lem}\label{lem:annulusr>1}
Let $b$ be as in Lemma~\ref{lem:AtoA}. Then
\begin{enumerate}
    \item The edges of $\Gam_b \cap Q$ lie in two isotopy classes on $\Gam_b$.
    %\item All the edges of $\Gam_{b-1} \cap (Q \cup P)$ lie in the same two isotopy classes in $\Gam_{b-1}$. 
    \item $b'_Q>b > b_Q$. 
\end{enumerate}
\end{lem}

%$\spadesuit$ {\em (Don't think I used that $b$ was the smallest, just that there is one. Need to check that. I also need to add a couple of figures to illustrate case II.)}

\begin{proof} 
All the edges of $\Gam_b \cap Q$ run from vertex $A$ to vertex $B$.
As argued in section~\ref{sec:atleast3}, these edges must lie in at least two isotopy classes, and can lie in at most four. 

We apply Lemma~\ref{lem:AtoA}.
First suppose that $G_P^b$ contains a simple closed curve $c$. In $\Gam_b$, $c$ must be essential and there must be exactly two isotopy classes of edges in $\Gam_b \cap Q$, verifying conclusion $(1)$. 
%As $b$ is smallest, $(2)$ then follows as the saddle from $\Gam_b$ to $\Gam_{b-1}$ must involve $c$, verifying that $(2)$ holds. 
If $b \leq b_Q$, then $b < b_Q$ (because the saddle below $\Gam_{b_Q}$ is in $Q$). The edges of $\Gam_b \cap Q$ are isotopic to those of $\Gam_0 \cap Q$, forcing $c$ to be a meridian of $C_{\alpha}$, and contradicting Lemma~\ref{lem:nomeridians}$(2)$.
Thus $b>b_Q$.
Similarly, if $b_Q' \leq b$, then $c$ is a meridian of $C_{\beta}$, another contradiction. % As $Q$ is incompressible in $X_K$, $c$ bounds a disk in $Q$. As $Q \cap T \times [0,1]$ is incompressible in $T \times [0,1]$, we conclude that $Q \cap \Gam_0$ contains a simple closed curve intersection, a contradiction. $\spadesuit$ {\em (Maybe just check $(7)$ and quote that here)}
Hence conclusion (2) holds as well in this case.

We now suppose the second outcome of Lemma~\ref{lem:AtoA} holds, so that there are edges $e_A^b,e_B^b$ of $G_P^b$, with the ends of $e_A^b$ labelled $A$ and those of $e_B^b$ labelled $B$. That is, as edges of $\Gam_b$, $e_A^b$ connects vertex $A$ to itself and $e_B^b$ connects vertex $B$ to itself. We rule out the possibilities that there are three or four isotopy classes in $\Gam_b$ of the edges in $\Gam_b \cap Q$.

\bigskip

\noindent {\bf I. $\Gam_b \cap Q$ has four isotopy classes of edges}. 

There are two faces of $\Gam_b \cap Q$ on $\Gam_b$ which are 4-gons, and all others are bigons. As each edge of $\Gam_b \cap Q$ connects $A$ to $B$, we see that each face has alternating corners from $A$ and $B$ (corners coming from the fat vertices $A,B$ of $\Pi_b$). Then $e_A^b$ must lie in one of these 4-gons and the other in $e_B^b$, connecting diagonal corners. But this contradicts that $e_A$ and $e_B$ arise from the saddle between $\Gam_{b-1}$ and $\Gam_b$.

\begin{figure}
\includegraphics[width=4in]{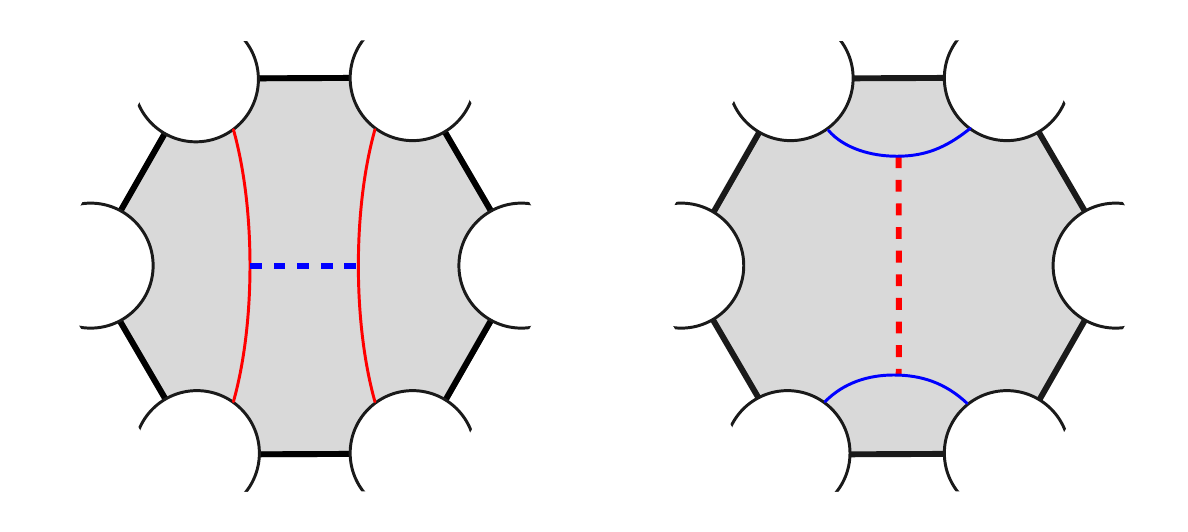}
    \put(-245,105){$A$}
    \put(-195,105){$B$}
    \put(-170,62){$A$}
    \put(-270,62){$B$}
    \put(-245,19){$A$}
    \put(-195,19){$B$}
    \put(-100,105){$A$}
    \put(-50,105){$B$}
    \put(-25,62){$A$}
    \put(-125,62){$B$}
    \put(-100,19){$A$}
    \put(-50,19){$B$}
    \put(-243,62){$e_A^b$}
    \put(-198,62){$e_B^b$}
    \put(-78,96){$e_1^{b-1}$}
    \put(-78,24){$e_2^{b-1}$}
    \put(-230,-5){$f\subset \Pi_b$}
    \put(-90,-5){$f\subset \Pi_{b-1}$}
    \caption{In the left half of the figure, $e_A^b, e_B^b$ of $\Pi_b \cap P$ (in red) lie within a hexagonal face, $f$, of $\Pi_b \cap Q$. In the right half, $e_1^{b-1}, e_2^{b-1}$ of $\Pi_{b-1}$ (in blue) arise from a saddle of $P$, indicated by the dashed lines, in going between levels $\Pi_{b}$ and $\Pi_{b-1}$. }

    \label{fig:eb}
\end{figure}

\noindent {\bf II. $\Gam_b \cap Q$ has exactly three isotopy classes of edges}. 

There is one $6$-gon face $f$ of $\Gam_b \cap Q$ on $\Gam_b$, and all other faces are bigons. The corners of $f$ are alternating $A$- and $B$-corners. 
Recall that $\Gam_b$ has no $1$-sided faces, as any such would 
imply a boundary compression of $P \cup Q$.
Thus,
$e_A^b, e_B^b$ must lie in $f$ and be diagonal edges in $f$ (Figure~\ref{fig:eb}).
As $\Gam_{b-1}$ has no $1$-sided faces and $e_A^b, e_B^b$ arise from the saddle between $\Gam_{b-1}$ and $\Gam_b$, we see that $e_A^b$ and $e_B^b$ arise by a banding of edges $e_1^{b-1},e_2^{b-1}$ of $\Pi_{b-1} \cap P$ in $\Pi_{b-1}$ that are parallel to diametrically opposite edges within $f$, each connecting $A$ to $B$.
All bandings take place in $f$, since all other faces are bigons, and the graphs do not contain 1-sided faces.
Under examination of Figure~\ref{fig:eb}, the next saddle singularity must either (a) band arcs of $P$, one in the isotopy class of $e_A^b$ and the other in the isotopy class of $e_B^b$, or else (b) band arcs of $(P \cup Q) \cap \Gam_b$, one from the isotopy class of $e_1^{b-1}$ and the other from the isotopy class of $e_2^{b-1}$: any other putative banding would result in a 1-sided face, which cannot occur.
Moreover, in case (b), the arcs banded together must both come from $P$, as otherwise $Q \cap \Gam_{b+1}$ would contain an arc with both endpoints labeled the same, violating the Parity Rule.
Note that if case (b) occurs, then it is not the final banding, since $Q \cap \Gam_{b+1}$ still contains three isotopy classes of arcs.
Iterate this observation.
We ultimately arrive to a level $b' \ge b$ where, staying within $f$: 

\begin{enumerate}
    \item In going from $\Gam_{b'-1}$ to $\Gam_{b'}$, we see edges $e_1^{b'-1}, e_2^{b'-1}$ of $\Gam_{b'-1}$ in the same relative isotopy class as $e_1^{b-1},e_2^{b-1}$ that are banded to produce edges $e_A^{b'}, e_B^{b'}$ of $\Gam_{b'}$ that are in the same relative isotopy class of $f$ as $e_A^b,e_B^b$.
    \item In going from $\Gam_{b'}$ to $\Gam_{b'+1}$, we see the edges $e_A^{b'}, e_B^{b'}$ of $\Gam_{b'}$ banded to produce edges $e_1^{b'+1}, e_2^{b'+1}$ of $\Gam_{b'+1}$ that are in the same relative isotopy class as $e_1^{b'-1},e_2^{b'-1}$. 
\end{enumerate}
See Figure~\ref{fig:eb'}.

\begin{figure}
\includegraphics[width=6in]{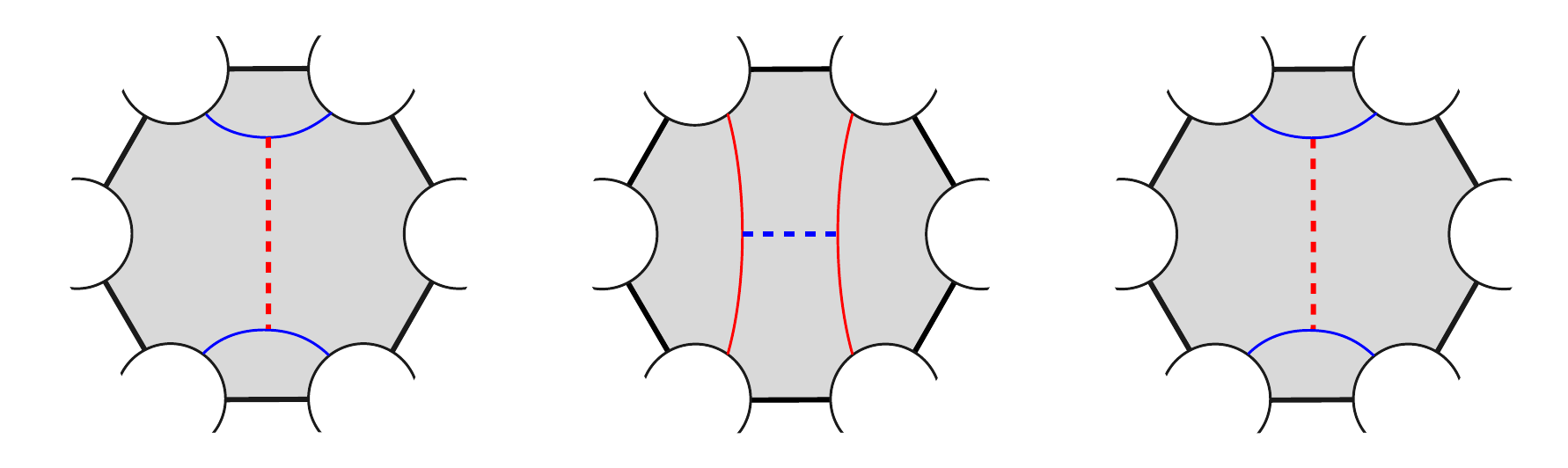} 
    % left labels
    \put(-392,110){$A$}
    \put(-338,110){$B$}
    \put(-310,62){$A$}
    \put(-420,62){$B$}
    \put(-392,15){$A$}
    \put(-338,15){$B$}
    % mid labels
    \put(-247,110){$A$}
    \put(-193,110){$B$}
    \put(-165,62){$A$}
    \put(-275,62){$B$}
    \put(-247,15){$A$}
    \put(-193,15){$B$}
    % right labels
    \put(-102,110){$A$}
    \put(-48,110){$B$}
    \put(-20,62){$A$}
    \put(-130,62){$B$}
    \put(-102,15){$A$}
    \put(-48,15){$B$}
    \put(-371,96){$e_1^{b'-1}$}
    \put(-371,24){$e_2^{b'-1}$}
    \put(-243,62){$e_A^{b'}$}
    \put(-198,62){$e_B^{b'}$}
    \put(-82,96){$e_1^{b'+1}$}
    \put(-82,24){$e_2^{b'+1}$}
    \put(-389,100){\small $\gamma_A$}
    \put(-342,100){\small $\gamma_B$}
    \put(-244,100){\small $\gamma_A$}
    \put(-198,100){\small $\gamma_B$}
    \put(-99,100){\small $\gamma_A$}
    \put(-53,100){\small $\gamma_B$}
    \put(-380,-5){$f\subset \Pi_{b'-1}$}
    \put(-230,-5){$f\subset \Pi_{b'}$}
    \put(-90,-5){$f\subset \Pi_{b'+1}$}
    \caption{Saddles from $P$ between levels $\Gam_{b'-1}, \Gam_{b}$ and $\Gam_{b'+1}$ are indicated by dashed lines. Edges in black are from $Q$, in blue or red from $P$. Note that $e_1^{b'-1}$ ($e_2^{b'-1}$) is in the same relative isotopy class as $e_2^{b'+1}$ ($e_2^{b'+1}$).}
    \label{fig:eb'}
\end{figure}

% \begin{figure}
% \includegraphics[width=5in]{fig.eb'.jpg}
   
%     \caption{Saddles from $P$ between levels $\Gam_{b'-1}, \Gam_{b}$ and $\Gam_{b'+1}$ are indicated by dashed lines.  Note that $e_1^{b'-1}$ ($e_2^{b'-1}$) is in the same relative isotopy class as $e_2^{b'+1}$ ($e_2^{b'+1}$).}

%     \label{fig:eb'}
% \end{figure}

%To see the above, recall that $\Gam$ has no 1-sided faces and edges of $\Gam_j \cap Q$ always connect $A$ to $B$. So until the above happens, all saddles are from $P$; furthermore, starting with $j=b-1$, in going from $\Gam_j$ to $\Gam_{j+1}$, there are edges of $\Gam_j$ in the same isotopy class of $f$ as $e_1^{b-1}, e_2^{b-1}$ that are banded to get edges in the same isotopy class of $f$ as edges $e_A^b, e_B^b$, all staying within $f$. 

Thus we see edges $e_1^{b'-1}, e_2^{b'-1}$ in $\Gam_{b'-1}$ that are in the same relative isotopy class of $f$ as $e_1^{b'+1}, e_2^{b'+1}$. Let $\gamma_A$ ($\gamma_B$) be the label of $e_1^{b'-1}$ at vertex $A$ ($B$) of $\Pi_{b'-1}$. Then $e_1^{b'+1}$ in $\Pi_{b'+1}$ has labels $\gamma_A, \gamma_B$ at the respective vertices.

%Let $\gamma_{1,A}$ ($\gamma_{1,B}$) be the label of $e_1^{b'-1}$ at vertex $A$ ($B$) of $\Pi_{b'-1}$. Let $\gamma_{2,B}$ ($\gamma_{2,B}$) be the label of $e_2^{b'-1}$ at vertex $A$ ($B$) of $\Pi_{b'-1}$. By the Parity Rule, $\gamma_{1,A}, \gamma_{1,B}$ correspond to vertices of the same sign in $G_P^{b'-1}$ (possibly the same vertex) and similarly $\gamma_{2,A}, \gamma_{2,B}$ are the same sign. The labels of $e_A^{b'}$ ($e_B^{b'}$) in $\Pi_{b'}$ are $\gamma_{1,A},\gamma_{2,A}$ ($\gamma_{1,B},\gamma_{2,B}$). By the Parity Rule then $\gamma_{1,A},\gamma_{2,A}$ correspond to vertices of opposite sign in $G_P^{b'}$. Thus $\{\gamma_{1,A}, \gamma_{1,B} \} \cap \{\gamma_{2,A}, \gamma_{2,B} \} = \emptyset$. Finally, note that
%$e_1^{b'-1}$ in $\Pi_{b'-1}$ has the same labels at its endpoints as $e_1^{b'+1}$ in $\Pi_{b'+1}$ at the respective vertices. 

\begin{figure}
\includegraphics[width=2in]{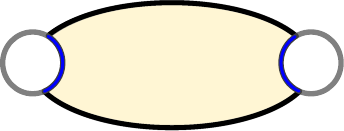}
    \put(-82,60){$e_1^{b'-1}$}
    \put(-82,-13){$e_1^{b'+1}$}
    \put(-136,26){$\gamma_A$}
    \put(-18,26){$\gamma_B$}
    \put(-116,26){\textcolor{Blue}{$\tau_A$}}
    \put(-38,26){\textcolor{Blue}{$\tau_B$}}
    \caption{
    The edges $e_1^{b'-1}, e_1^{b'+1}$ are pictured in $P$. The curve $c=e_1^{b'-1} \cup  e_1^{b'+1} \cup \tau_A \cup \tau_B$ is homotopically trivial in $X_K$, hence bounds a disk in $P$. But this disk contains a saddle in $P$, contradicting that $\tau_A, \tau_B$ are transverse to the height function on $X_K$. }

    \label{fig:c}
\end{figure}

%Similarly, $e_B^{b'-1}$ has the same labels at its endpoints as $e_B^{b'+1}$ at the respective vertices. 
On $P$, we see the arcs $e_1^{b'-1}, e_1^{b'+1}$ connecting vertices $\gamma_{A}, \gamma_{B}$.
%the same pair of components (possibly the same component) of $\partial P$, $\gamma_1, \gamma_2$. 
Let $\tau_A, \tau_B$ be arcs on $\gamma_{A},\gamma_{B}$ that lie between levels $\Gam_{b'-1}$ and $\Gam_{b'+1}$. See Figure~\ref{fig:c}. 
The curve $c=e_1^{b'-1} \cup e_1^{b'+1} \cup \tau_A \cup \tau_B$ is a closed curve in $P$ which lies in the ball $f \times I \subset X_K$ between the levels $\Gam_{b'-1}$ and $\Gam_{b'+1}$; consequently $c$ is homotopically trivial in $X_K$. As $P$ is incompressible in $X_K$ ($\pi_1$-injective), $c$ must bound a disk in $P$ on the $\alpha$-side of $e_1^{b'+1}$. This disk must contain the saddle between $\Gam_{b'}$ and $\Gam_{b'+1}$, and $e_2^{b'+1}$ also lies in this disk with endpoints in $\tau_A \cup \tau_B$. But $\tau_A, \tau_B$ are transverse to the height function in $X_K$, and intersect $\Gam_{b'+1}$ only at $e_1^{b'+1}$, a contradiction.  (This argument was inspired by Lemma 3.3.1 of \cite{culler1982}, see Figure 16 there.) 
%This contradicts that $\Gam_{b'}$ contains no $1$-sided faces or trivial closed curves. 
%But this contradicts that in $P$ we see that $c$ has components of $\partial P$ on both sides. In particular, $\gamma_{1,A}, \gamma_{1,B}$ are on one side and $\gamma_{2,A}, \gamma_{2,B}$ on the other. 
%a saddle of $P$ between $e_1^{b'-1}$ and $e_1^{b'+1}$, contradicting that $P$ has no center singularities in $\hat{T} \times [0,1]$.

\bigskip  

We conclude that $\Gam_b \cap Q$ has exactly two isotopy classes of edges, and conclusion $(1)$ of the Lemma holds.
We now verify $(2)$.
For contradiction, assume $b \leq b_Q$. Then $b < b_Q$ (the saddle preceding $\Gam_{b_Q}$ is in $Q$). Then $\Gam_b \cap Q$ is isomorphic to $\Gam_0 \cap Q$. All the edges of $\Gam_b \cap Q$ lie in classes $\hat{\alpha}, \hat{\alpha}'$; furthermore, each such edge has the same label on both of its endpoints. Let $\A$ be an annulus on the level torus corresponding to $\Gam_b$ containing the edges of $\Gam_b \cap Q$.

\begin{figure}
\includegraphics[width=3in]{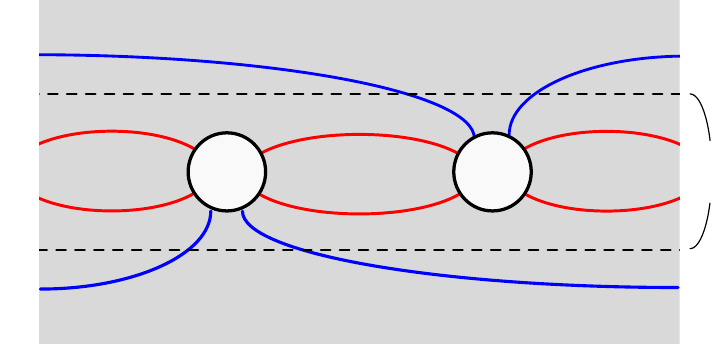}
    \put(-135,46){\tiny{$1$}}
    \put(-85,46){\tiny{$1$}}
    \put(-166,46){\tiny{$q$}}
    \put(-55,46){\tiny{$q$}}
    \put(-7,48){$\mathcal{A}$}
    \put(-113,90){$\Pi_b$}
    \put(-153,48){$A$}
    \put(-73,48){$B$}
    \put(-192,6){\textcolor{Blue}{$e_A^{b}$}}
    \put(-192,92){\textcolor{Blue}
    {$e_B^{b}$}}
    \caption{Edges of $\Gam_b \cap Q$ are in red, edges of $\Gam_b \cap P$ are in blue. As $\Gam_b$ has no $1$-sided faces, the number of labels from $P$ between labels $1$ and $q$ (from $Q$) must be different on vertex $A$ than vertex $B$, a contradiction. }

    \label{fig:oppositesidesofA}
\end{figure}

Recall from Definition~\ref{def:G_P} that $\Gam_b$ has no $1$-sided faces. Note that $e_A^b, e_B^b$ must connect opposite sides of $\A$. Otherwise, $e_A^b$ connects one side of $\A$ to itself and $e_B^b$ connects the other side to itself, as in Figure~\ref{fig:oppositesidesofA}. We would find a pair of labels of $\Gam_b \cap Q$ where the number of edges of $\Gam_b \cap P$ incident between these labels on $A$ is different than the number incident on $B$.
That is a contradiction, because the whole set of cyclically ordered labels around $A$ and $B$ match one another in reverse.

\begin{figure}
\includegraphics[width=4in]{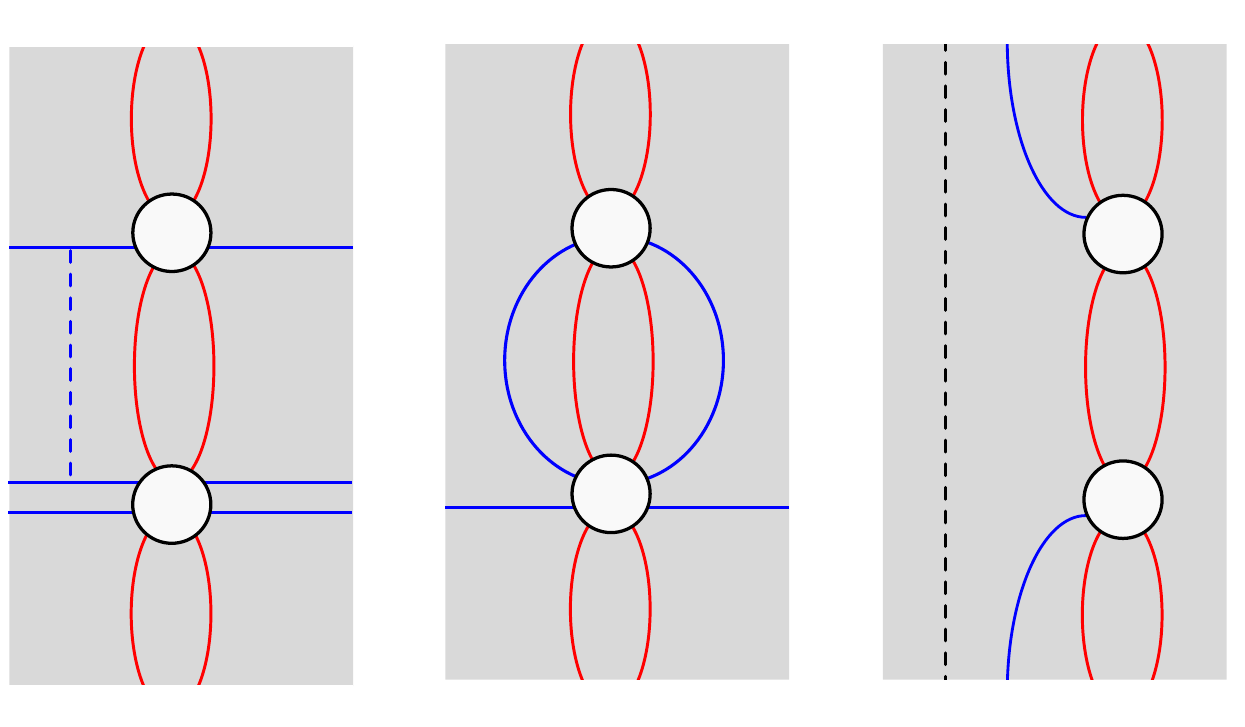}
    \put(-254,63){\tiny{$1$}}
    \put(-152,65){\tiny{$1$}}
    \put(-33,64){\tiny{$1$}}
    \put(-254,99){\tiny{$1$}}
    \put(-152,100){\tiny{$1$}}
    \put(-33,99){\tiny{$1$}}
    \put(-254,38){\tiny{$q$}}
    \put(-152,40){\tiny{$q$}}
    \put(-34,39){\tiny{$q$}}
    \put(-254,127){\tiny{$q$}}
    \put(-152,128){\tiny{$q$}}
    \put(-34,126){\tiny{$q$}}
    \put(-162,55){\tiny{$x$}}
    \put(-162,112){\tiny{$x$}}
    \put(-43,51){\tiny{$x$}}
    \put(-43,113){\tiny{$x$}}
    \put(-65,80){$c$}
    \put(-260,-2){$\Pi_{b'-1}$}
    \put(-155,-2){$\Pi_{b'}$}
    \put(-35,-2){$\Pi_0$}
    \put(-254,48){$B$}
    \put(-152,50){$B$}
    \put(-33,49){$B$}
    \put(-254,111){$A$}
    \put(-152,112){$A$}
    \put(-33,110){$A$}
    \put(-280,38){\textcolor{Blue}{$e_B^{b}$}}
    \put(-280,117){\textcolor{Blue}
    {$e_A^{b}$}}
    \put(-183,85){\textcolor{Blue}
    {$e^{b'}$}}
    \put(-60,120){\textcolor{Blue}
    {$e^{0}$}}
    \caption{A saddle of $P$ creates from $e_A^b$ an edge, $e^{b'}$, that when projected to $\Gam_0$ is parallel to an edge, $e^0$, of $P \cap \Gam_0$.  $e^{b'}$ and $e^0$ form a curve $c$ in $P$ that bounds a disk in $X_K$ ($c$ is isotoped to a meridian of $C_\alpha$ in the Figure). We conclude that $e^0, e^{b'}$ are parallel in $P$, contradicting the intervening saddle.}

    \label{fig:b<bQ}
\end{figure}

There must be a $b_Q>b'>b$ where $(i)$ $e_A^b, e_B^b$ survive to $\Gam_{b'-1}$ and $(ii)$ a saddle of $P$ between $\Gam_{b'-1}$ and $\Gam_{b'}$ is incident to at least one of $e_A^b, e_B^b$ (because $\Gam_{b_Q}$ is preceded by a saddle in $Q$ which results in edges connecting $A$ to $B$). Assume this occurs with $e_A^b$. Now consult  Figure~\ref{fig:b<bQ}.
This saddle must connect $e_A^b$ to an edge in the same isotopy class as $e_B^b$, and results in edges in classes $\hat{\alpha}$ or $\hat{\alpha}'$. Let $e^{b'}$ be one of these edges. Then $e^{b'}$ has the same label $x$ at both endpoints (it is parallel to an edge of $\Gam_{b'} \cap Q$). Let $e^0$ be the edge of $\Pi_0 \cap P$ that has both endpoints labelled $x$. As argued above, in $P$ we form a simple closed curve $c=e^0 \cup e^{b'} \cup \tau_1 \cup \tau_2 \subset P$, where $\tau_1, \tau_2$ are arcs running vertically along $\partial X_K$ between the endpoints of $e^0, e^{b'}$ (in the component of $\partial P$ corresponding to the fat vertex $x$). But $c$ can be projected onto $\Pi_0$ in the complement of $\A$; hence, $c$ is trivial in $\pi_1(X_K)$. By the incompressibility of $P$, we conclude that $c$ bounds a disk in $P$ on the $\alpha$ side of $e^{b'}$. But $e^{b'}$ is incident to the saddle between $\Gam_{b'-1}$ and $\Gam_{b'}$ on this side, implying that this disk contains another arc intersection with $\Gam_{b'}$. As before, this contradicts that $\tau_1, \tau_2$ are transverse to the height function on $X_K$.

%as $P$ in its Morse presentation in $T \times $ has no center singularities (e.g. Lemma~\ref{lem:nomeridians}$(2)$). 

Similarly, we argue that $b < b_Q'$, finishing the proof of the Lemma.
%\begin{clm}\label{clm:edgesamelabels}
%There is a $b' \geq b$ and an edge $e \subset \Gam_{b'} \cap P$ such that both edges of $e$ have the same label and this edge involves the saddle between $\Gam_{b'-1}$ and $\Gam_{b'}$. 
%\end{clm}
%\noindent {\em Proof of Claim:}
%Assume $\Gam_b \cap P$ has an edge connecting $A$ to $B$, that is not in isotopy class $\alpha_1$ or $\alpha_2$. By counting labels around $A$, $B$ (and using that each edge of $\Gam_b \cap Q$ has the same label on both ends) one sees that there must be a second edge of $\Gam_b$ connecting $A,B$ that is not isotopic in $\Gam_b$. That is, $\Gam_b \cap P$ has four isotopy classes of arcs connecting $A$ to $B$. The argument of Case $(I)$ above gives a contradiction. We conclude that all edges of $\Gam_b \cap P$ that connect $A$ to $B$ lie in isotopy classes $\alpha_1, \alpha_2$. 
%Let $\mathcal{A}$ be an annulus containing the edges of $\Gam_b \cap (Q \cup P)$ in classes $\alpha_1$ and $\alpha_2$. Then $e_A^b,e_B^b$ must connect one side of $\mathcal{A}$ to the other. But then there must be a $b'>b$ where a saddle is incident to one of these two edges. This saddle must connect $e_A^b$ to and edge in the same isotopy class as $e_B^b$, so in fact, by looking at labels, to $e^b_B$. The resulting edges are the desired edges.
%\qed{(Claim)}
\end{proof}

\begin{proof}
[Proof of Theorem~\ref{thm:mainNextended}]
Apply Lemma \ref{lem:annulusr>1}.
Extract the union of graphs $\Gamma$ out of the union of graphs $\Gam$ by forgetting $P$ and eliminating repetitions.
The value $b_Q$ shifts into $0$, the value $b_Q'$ shifts into $2q-2$, and the value $b$ shifts, but it remains strictly between these values.
We thereby obtain the desired value $0 < b < 2q-2$ such that all arcs of $\Gamma_b$ belong to two isotopy classes.
Apply Theorem~\ref{thm:annulus} to conclude the proof. Note that $k>1$ since $K(\lambda)$ is not prime.
\end{proof}

\section{Weighting the arcs of \pdfmath{\Gamma}{Gamma}, \pdfmath{\Gamma'}{Gamma'}
}\label{sec:alphabetawts}

Assume $K$ is a simple knot with a planar Seifert surface $Q$ of longitudinal slope $\lambda$ and $q$ boundary components. Then $K(\lambda) = (S^1 \times S^2) \# N$ for some rational homology sphere $N$. In Section~\ref{sec:r>1}, we proved Theorem~\ref{thm:mainNextended}, which classifies those knots $K$ for which $N \ne S^3$. Theorem~\ref{thm:Gr-Li-Lu} classifies those simple $K$ for which $K(\lambda)=S^1 \times S^2$.  Its proof occupies the next six sections. 

We begin by outlining the argument for Theorem~\ref{thm:Gr-Li-Lu}. Assume $K(\lambda)=S^1 \times S^2$. By Theorem~\ref{thm:lisca-rasm} (rephrased as in Theorem~\ref{thm:lisca}), $K \subset L(p,m)=L(q^2,m)$ where $q,m$ are as described. Recall the Morse presentation of $Q$ from Section~\ref{sec:setup} and its intersections with transverse levels given by $\Gamma$.
%Let $V_{\alpha} \cup_{\hat{T}} V_{\beta}$ be the standard Heegaard splitting of $L(q^2,m)$, where $\hat{T}$ is a 2-torus and $V_{\alpha},\, V_{\beta}$ are solid tori. View the Seifert surface $Q$ for $K$ in $L(q^2, m)$ in a neighborhood of $\hat{T}$, $\hat{T} \times [0,1]$.
%under a Morse function with levels given by Heegaard 2-tori. 
%As $K$ is simple, it intersects each of these tori, $\hat{T} \times \{b\}$, in two points, and we call the resulting twice-punctured torus in the knot exterior $T \times \{b\}$. 
%That is, these level punctured tori in the knot exterior are identified with $T \times [0,1]$. 
%We take $Q$ as standard outside of $T \times [0,1]$ in the knot exterior and focus on $Q \cap (T \times [0,1])$ and the resulting Morse function thereon. We may assume the only Morse singularities for $Q$ in $T \times [0,1]$ are saddle singularities, and that we have only $2q-2$ such singularities. $Q$ intersects each transverse level $T \times \{b\}$ in $q$ arcs. For example, in $T \times \{0\}$, all $q$ arcs are in two isotopy classes, and in $T \times \{1\}$ all $q$ arcs are in two isotopy classes. 
Let $\tau$ be an arc of $\Gamma_b$ for some $b=1, \dots, 2q-2$. Then $\tau$ divides $K$ into two arcs---one above $\Gamma_i$ and one below. Capping each of these arcs with $\tau$ writes $K$ as the sum of two loops---one in each of the Heegaard solid tori. This allows us to write $[K]=c \mu+f \mu' \in H_1(L(q^2,m))$ where $\mu,\, \mu'$ are the homology classes corresponding to the cores of the two Heegaard solid tori, $V_{\alpha}, V_{\beta}$, and where the integers $c,\, f$ are determined by the isotopy class of $\tau$ in $\Gamma_b$. Specifically, $c$ is called the $\beta$-weight of $\tau$, the homology class represented by $\tau \cup (K \cap V_{\beta})$ in $H_1(V_{\beta},\Z)$. Similarly, $f$ is the $\alpha$-weight of $\tau$, the homology class that $\tau \cup (K \cap V_{\alpha})$ represents in $H_1(V_{\alpha},\Z)$ (see Lemma~\ref{lem:decomposeK} and Definitions~\ref{def:alphaweights} and \ref{def:betaweights}). These weights $c,\, f$ should be small in terms of $q$ because the slopes of the arcs in $\Gamma_b$ change at most $2q-2$ times under the Morse function. This theme is exemplified in Theorem~\ref{thm:q+y,m_b=1} (Engine Theorem)  where $m_b$ ($m_b'$) is the minimum and $M_b$ ($ M_{b}'$) the maximum of the absolute values of the $\alpha$-weights ($\beta$-weights, resp.) of the arcs of $\Gamma_b$. Note that as described in section~\ref{sec:setup}, $y<q/2$
%and $q-y$ are the number of arcs of $Q \cap (T \times \{0\})$ in the two isotopy classes of $T \times \{0\}$ 
corresponds to $[K]$ in terms of $\mu$, and $y'<q/2$
%$q-y'$ the number of arcs in $T \times \{1\}$ 
corresponds to $[K]$ in terms of $\mu'$. We then use restrictions on the size of $c,\, f$ given by Theorems~\ref{thm:M:b<x+1} and \ref{thm:q+y,m_b=1}, the fact that $\mu'=m \mu$, the numerics of $L(q^2, m)$ as described by Theorem~\ref{thm:lisca-rasm}, and the fact that $[K]$ has order $q$, to show that only the homology classes claimed in Theorem~\ref{thm:Gr-Li-Lu} occur (see Theorems~\ref{thm:main:dq+1}, \ref{thm:main1-nd}, and \ref{thm:main1-nd:2q+1}). To understand this last step in more detail, we recommend the reader look through section~\ref{sec:dq+1}. Note that by Theorem~\ref{thm:q<600}, done via computer calculation, we may assume $q > 600$. 

\bigskip

Let $Q$ be a planar Seifert surface of $K \subset L(p,m)$ with $q$ boundary components as in section~\ref{sec:setup}, with its level presentation in $T \times [0,1]$ and with $\Gamma, \Gamma'$ as defined there. In this section, we define the $\alpha, \beta$-weights of level arcs described in the preceding paragraph. We assign the $\alpha$-weights to each of the arcs of $\Gamma_i$. By abuse of notation, we define $\Gamma_i$ be its level arcs along with their $\alpha$-weights. We also define $\Gamma'$, an upside down version of $\Gamma$, where we assign to the arcs of intersection their $\beta$-weights and re-index the levels (going from $V_\beta$ to $V_{\alpha}$).

We think of $L(p,m)$ as the oriented solid tori $V_{\alpha},\, V_{\beta}$, where $\alpha, \beta$ are the corresponding meridians, identified by an orientation-reversing diffeomorphism of their boundaries. Let $\widetilde{\mu'}$ be a simple closed curve in $\partial V_{\alpha}$ such that $[\alpha],\, [\widetilde{\mu'}]$ form an oriented basis for $H_1(\partial V_{\alpha})$. Note that $[\widetilde{\mu'}]=\mu' \in H_1(L(p,m))$. % After shifting to boundary?
%Then $\lambda_{\alpha}$ is isotopic to $\mu'$ in $V_{\alpha}$. 

Recall from section~\ref{sec:setup} that $A,B$ are the boundary components of $T$ gotten by removing neighborhoods of $x_0, x_\kappa$ (resp.) from $\hat{T}=\partial V_{\alpha}$ and that $K \cap V_{\alpha}$ runs from $x_\kappa$ to $x_0$.

\begin{figure}
\includegraphics[width=4in]{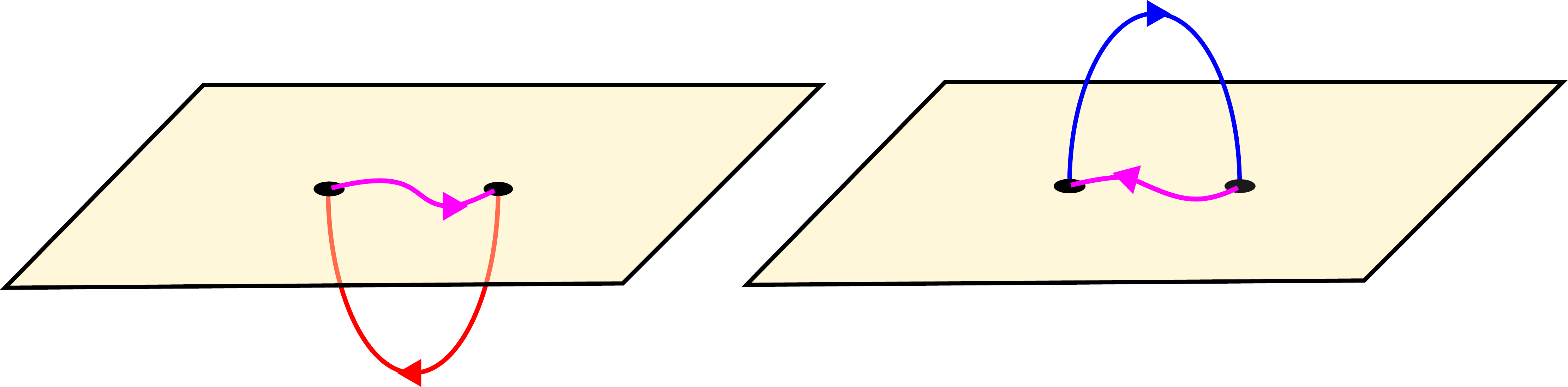}
\put(-215, 42){\textcolor{VioletRed}{$\tau_\alpha$}}
\put(-242, 30){\small$x_0$}
\put(-192, 30){\small$x_\kappa$}
\put(-85, 28){\textcolor{VioletRed}{$\tau_\beta$}}
\put(-105, 30){\small$x_0$}
\put(-57, 30){\small$x_\kappa$}
\put(-260, 0){\textcolor{red}{$K\cap V_\alpha$}}
\put(-120, 70){\textcolor{blue}{$K\cap V_\beta$}}
\put(-269, 25){$\hat{T}$}
\put(-132, 25){$\hat{T}$}
\caption{$\alpha$-weight of $\tau$ is $n'$ where $n' \mu'=[\tau_\alpha \cup (K \cap V_\alpha] \in H_1(V_\alpha)$); $\beta$-weight of $\tau$ from $H_1(V_\beta$).}
\label{fig:weights}
\end{figure}

\begin{defin}[$\alpha$-weights]\label{def:alphaweights}
Let $\tau$ be a properly embedded arc in $T$ running between boundary components $A$ and $B$.  
Define $\tau_\alpha$ to be the oriented arc in $\hat{T}=\partial V_{\alpha}$ from $x_0$ to $x_\kappa$ gotten from $\tau$ by shrinking $A$ and $B$ to the points  $x_0$ and $x_\kappa$. Define the {\em $\alpha$-weight} of $\tau$ to be $n'$, where $n' [\widetilde{\mu'}] = [\tau_\alpha \cup (K \cap V_{\alpha})] \in H_1(V_{\alpha})$. See Figure~\ref{fig:weights}. The $\alpha$-weight of the isotopy class of $\tau$ is defined to be the $\alpha$-weight of $\tau$.
\end{defin} 

%$\spadesuit$ {\em (
%\begin{defin}[$\alpha$-weights]\label{def:alphaweights'}
%We have identified each capped off level surface in $X_K$ with $\partial V_{\alpha}$ via the product structure.
%For every $T \times \{i\}$ transverse to $Q$ we orient the arcs of $Q \cap (T \times \{i\})$ from $A$ to $B$ under this identification. 
%$Note that in $\Gamma_0$ each arc is oriented to be isotopic to the oriented $K \cap V_{\alpha}$. 
%To each such oriented arc in $T \times \{i\}$, $\tau$, we let $\tau_\alpha$ be the oriented arc in $\partial V_{\alpha}$ from $x_0$ to $x_\kappa$ gotten from $\tau$ by shrinking $A$ and $B$ to points. We then define the {\em $\alpha$-weight} of $\tau$ to be $n'$, where $n' [\widetilde{\mu'}] = [\tau_\alpha \cup (K \cap V_{\alpha})] \in H_1(V_{\alpha})$. See Figure~\ref{fig:weights}.
%\end{defin} )}

%For example, in $T \times \{0\}$ we see $q$ arcs in isotopy classes $\alpha$ and $\alpha_2$ (with $y$ in one, $\hat{\alpha}$, and $x$ in the other, $\hat{\alpha'}$). 
For example, if the isotopy class of $\tau$ is either $\hat{\alpha}$ or $\hat{\alpha}'$ then the $\alpha$-weight of $\tau$ is $0$, as the corresponding $\tau_\alpha$ is isotopic rel endpoints in $V_\alpha$ to the oriented $-(K \cap V_\alpha)$---see Definition~\ref{def:alpha1} and the sentences following it. 

Let $T \times \{i\}$ be transverse to $Q$, and let $W_{\alpha}$ be the $\alpha$-side of $T \times \{i\}$ in $X_K$ (i.e. that intersecting $C_\alpha$ of section~\ref{sec:setup}).
Let $C$ be a component of $\partial Q$. 
At $C$ we locally see arcs of $Q \cap (T \times \{i\})$ incident, one at $C \cap A$ and the other at $C \cap B$. Furthermore the arc $C \cap W_{\alpha}$, which we refer to as the {\it {$\alpha$}-corner} of $C$ is isotopic to $K \cap V_{\alpha}$. Thus in $Q$ we see the $q$ arcs of its intersection with $T \times \{i\}$ forming disjoint cycles  with the components of $\partial Q$. In particular, the boundary of the surface $Q \cap W_{\alpha}$ is exactly these cycles of level arcs (along with their connecting $\alpha$-corners) and the level simple closed curves from $T \times \{i\}$.  Because the components of $\partial Q$ are coherently oriented on $\partial X_K$ we can find an orientation on $Q \cap W_{\alpha}$ that induces an orientation of its boundary from $A$ to $B$ along a level arc and from $B$ to $A$ along an $\alpha$-corner. Thus, when there are no simple closed curves in $Q \cap (T \times \{i\})$, $Q \cap W_{\alpha}$ describes a null-homology in $V_{\alpha}$ between the curves defining the $\alpha$-weights of the $q$ level arcs in $T \times \{i\}$. This proves the following.

\begin{lem}\label{lem:alphaweightTxi}
Let $T \times \{i\}$ be transverse to $Q$ with no simple closed curves of intersection, then the sum of the $\alpha$-weights of the arcs of $T \times \{i\}$ is $0$ in ${\mathbb Z}$.
\end{lem}

%\begin{lem}
%Let $T \times \{i\}$ be transverse to $Q$. 

%\end{lem}

%\noindent {\bf Remark:} Because the boundary components of $Q$ are coherently oriented along $\partial X_K$, there is an orientation on $Q$ such that at every level surface $T \times \{j\}$ that is transverse to $Q$ the orientation from $A$ to $B$ corresponds to the orientation that $Q \cap (T \times [0,j])$ induces on the arcs of $Q\cap (T \times \{j\})$. In particular, the banding occurring at critical levels of $Q$ must cancel the orientations coming from orienting the corresponding arcs from A to B. The arcs before and after the saddle will be oriented arcs from A to B (e.g., A is a source before and after).  That a banding corresponding to a saddle must locally cancel orientations implies that the banding must connect arcs in different isotopy classes in the level surface. 

%\begin{proof} 
%As noted in the Definition above, the sum of the $\alpha$-weights of the arcs in $\Gamma_0$ is $0$.
%In going from $\Gamma_j$ to $\Gamma_{j+1}$, the saddle singularity corresponds to a banding in $\partial V_{\alpha}$ that cancels orientation. Thus, the sum of closed curves, $\tau \cup \hat{\alpha}$, in $H_1(V_{\alpha})$ over the arcs $\tau$ in $\Gamma_j$ is the same as the sum over arcs $\tau$ in $\Gamma_{j+1}$.
%\end{proof}

\begin{defin}[Total $\alpha$-weight]\label{def:totalalpha} For each $i$, the {\em total $\alpha$-weight of $\Gamma_i$} is the formal linear combination over isotopy classes of arcs in $\Gamma_i$ of the $\alpha$-weight of an arc in that class with coefficient the number of arcs in that class. That is, if ${\wt{w_1}, \dots, \wt{w_n}}$ are the $\alpha$-weights of $\Gamma_i$, one for each isotopy class of arc, then the total $\alpha$-weight of $\Gamma_i$ is $c_1 \wt{w_1} + \dots + c_n \wt{w_n}$ where $c_i$ is the number of arcs in the isotopy class corresponding to $\wt{w_i}$. Note that the sum of the coefficients $c_i$ is $q$, the number of arcs in $\Gamma_i$, and that $n \leq 4$ as $T$ is a twice-punctured torus.
\end{defin}

% {\bf NOTE: We work out the total weights for $K(49,18,21)$, at the end of the paper in section~\ref{sec:exampletotalweights}.}

For example, the total weight of $\Gamma_0$ is $y [{\bf 0}] + x [{\bf 0}]$. 

We now weight arcs from the $V_{\beta}$ perspective by reversing the direction of increasing height. Recall from section~\ref{sec:setup} that $\Gamma_i'=\Gamma_{(2q-2)-i}$ and $\Gamma'=\Gamma_0' \cup \dots \cup \Gamma_{2q-2}'$. Take a simple closed curve $\widetilde{\mu}$ in $\partial V_{\beta}$ such that $[\beta],\, [\widetilde{\mu}]$ form an oriented basis for $H_1(\partial V_{\beta})$. Note that $[\widetilde{\mu}]=\mu \in H_1(L(p, m))$.

\begin{defin}[$\beta$-weights]\label{def:betaweights}
Let $\tau$ be a properly embedded arc in $T$ running between boundary components $A$ and $B$.  
Define $\tau_\beta$ to be the oriented arc in $\hat{T}=\partial V_{\beta}$ from $x_\kappa$ to $x_\beta$ gotten from $\tau$ by shrinking $A$ and $B$ to the points  $x_0$ and $x_\kappa$. Define the {\em $\beta$-weight} of $\tau$ to be $n$, where $n [\widetilde{\mu}] = [\tau_\beta \cup (K \cap V_{\beta})] \in H_1(V_{\beta})$. See Figure~\ref{fig:weights}. The $\beta$-weight of the isotopy class of $\tau$ is defined to be the $\beta$-weight of $\tau$.
\end{defin} 

%\begin{defin}[$\Gamma'$ and $\beta$-weight]\label{def:betaweights'}
%Recall that $\Gamma_i'=\Gamma_{(2q-2)-i}$ and $\Gamma'=\Gamma_0' \cup \dots \cup \Gamma_{2q-2}'$.

%We identify each capped off level surface in $X_K$ with $\partial V_{\beta}$ via the product structure.
%For every $\Gamma_i'$ we orient the arcs of intersection from $B$ to $A$. 
%Note that in $\Gamma_0'$ each arc is oriented to be isotopic to the oriented $K \cap V_{\beta}$. 
%Take a simple closed curve $\widetilde{\mu}$ in $\partial V_{\beta}$ such that $[\beta],\, [\widetilde{\mu}]$ form an oriented basis for $H_1(\partial V_{\beta})$. Note that $[\widetilde{\mu}]=\mu \in H_1(L(p, m))$.
%Then $\lambda_{\beta}$ is isotopic to $\mu$ in $V_{\beta}$.
%To each oriented arc in $\Gamma_i'$, $\tau$, we let $\tau_\beta$ be the oriented arc in $\partial V_{\beta}$ from $x_\kappa$ to $x_0$ gotten from $\tau$ by shrinking $A$ and $B$ to points. We then assign a {\em $\beta$-weight} $n$, such that $n [\widetilde{\mu}]=[\tau_\beta \cup (K \cap V_{\beta})] \in H_1(V_{\beta})$.
%\end{defin}

For example, in $\Gamma_0'$ we see $q$ arcs with $\beta$-weight $0$.

Arguing exactly as in the proof of Lemma~\ref{lem:alphaweightTxi} but now from the point of view of the $\beta$-side of $T \times \{i\}$ (that containing $C_\beta$), we get the 
following.

\begin{lem}\label{lem:betaweightTxi}
Let $T \times \{i\}$ be transverse to $Q$ with no simple closed curves of intersection, then the sum of the $\beta$-weights of the arcs of $T \times \{i\}$ is $0$ in ${\mathbb Z}$.
\end{lem}

\begin{defin}[Total $\beta$-weight]\label{def:totalbeta} For each $i$, the {\em total $\beta$-weight of $\Gamma_i'$} is the formal linear  combination over isotopy classes of arcs in $\Gamma_i'$ of the $\beta$-weight of an arc in that class with coefficient the number of arcs in that class. That is, if ${\wt{w_1}, \dots, \wt{w_n}}$ are the $\beta$-weights of $\Gamma_i'$, one for each isotopy class of arc, then the total $\beta$-weight of $\Gamma_i'$ is $c_1 \wt{w_1} + \dots + c_n \wt{w_n}$. Note that the sum of the coefficients $c_i$ is $q$, the number of arcs in $\Gamma_i$, and $n \leq 4$. 
\end{defin}

\begin{remark}\label{rmk:gammaconvention}
$\Gamma'_i=\Gamma_{2q-2-i}$ means that they correspond to the same level surfaces with the same arcs of intersection with $Q$. The difference is the weighting of these arcs: $\Gamma_i$ from the point of view of $V_{\alpha}$, $\Gamma_i'$ from the point of view of $V_{\beta}$. Thus, we take the following convention: In the context of $\Gamma_i$, (total) weight will mean (total) $\alpha$-weight. Similarly, in the context of $\Gamma_i'$, (total) weight will mean (total) $\beta$-weight.
\end{remark}

\begin{figure}
\includegraphics[width=4in]{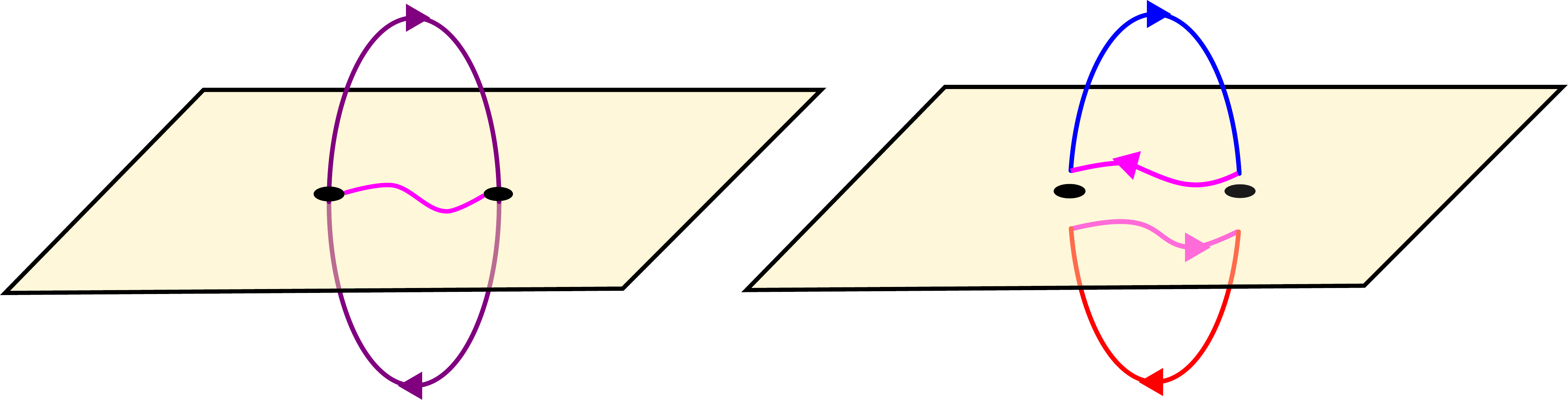}
\put(-85, 25){\textcolor{VioletRed}{$\tau_\alpha$}}
\put(-215, 42){\textcolor{VioletRed}{$\tau$}}
\put(-77, 45){\textcolor{VioletRed}{$\tau_\beta$}}
\put(-120, 0){\textcolor{red}{$K\cap V_\alpha$}}
\put(-120, 70){\textcolor{blue}{$K\cap V_\beta$}}
\put(-262, 25){$\Gamma_i$}
\put(-132, 25){$\Gamma_i$}
\put(-60, 70){$c\cdot \mu\in H_1(L(p,m))$}
\put(-60, 5){$f\cdot \mu'\in H_1(L(p,m))$} 
\caption{Banding across $\tau$ writes $[K]$ as sum.}
\label{fig:write_K_as_sum}
\end{figure}

\begin{lem}\label{lem:decomposeK}
For some $i$, let $\tau$ be an arc in $\Gamma_i$ with $\alpha$-weight $f$ and assume that the same arc in $\Gamma_{(2q-2)-i}'$ has $\beta$-weight $c$. Then $[K]=c \mu + f \mu' \in H_1(L(p,m))$ where $\mu, \mu'$ are the canonical generators of $H_1(L(p,m))$. 
\end{lem} 

\begin{proof} Let $\tau$ is an unoriented arc on $\Gamma_i$. Note that $\tau_\alpha = - \tau_\beta$.
Thus, $[K] = [(K \cap V_{\alpha}) \cup (K \cap V_{\beta})] = [(K \cap V_{\alpha}) \cup \tau_\alpha \cup \tau_\beta \cup (K \cap V_{\beta})] = [(K \cap V_{\alpha}) \cup \tau_\alpha] + [\tau_\beta \cup (K \cap V_{\beta})]= f [\widetilde{\mu'}] + c [\widetilde{\mu}] = f \mu' + c \mu \in H_1(L(p,m))$. See Figure~\ref{fig:write_K_as_sum}.
\end{proof}

%$\spadesuit$ {\em  It would make statements simpler if we could show that you can't do sequence of bandings, then a sequence of companions of the reverse bandings---as it done below when the weight is $0$.}

\begin{defin}
The smallest (largest) among the absolute values of the weights of arcs in $\Gamma_b$ is defined to be $m_b$ ($M_b$).  
Similarly, the smallest (largest) among the absolute values of the weights of arcs in $\Gamma_b'$ is defined to be $m_b'$ ($M_b'$). 
\end{defin}

\begin{cor}\label{cor:Mmdecom}
    Assume $0\le b\le 2q-2$. Let $b'=2q-2-b$. Then up to sign $[K]$ can be written as $c\mu\pm f\mu'$, where $c,f$ can take either
 of the following forms   
 \begin{enumerate}
        \item $c=m_b'$, $0\le f\le M_{b'}$; or
        \item $f=m_b$, $0\le c\le M_{b'}'$.
    \end{enumerate}
    We indicate this by $\pm[K]=c\mu\pm f\mu'$.
\end{cor}

\begin{proof}
    We give the proof of $(1)$. The proof of $(2)$ is identical to this, with $\Gamma$ and $\Gamma'$ switched.
    
    By the definition of $m_b'$, there is an arc $\tau$ in $\Gamma_b'$ whose absolute value of $\beta$-weight is $m_b'$. On the other hand, any arc in $\Gamma_{b'}$ has $\alpha$-weight with absolute value at most $M_{b'}$. Note that $\Gamma_{b}'=\Gamma_{b'}$; that is, we may consider $\tau$ as an arc in $\Gamma_{b'}$. Thus, $\tau$ has $\alpha$-weight with absolute value at most $M_{b'}$. Lemma~\ref{lem:decomposeK} says $[K]=\pm c\mu\pm f\mu'$ with appropriately chosen signs, which proves the corollary.
\end{proof}

\begin{defin}
$\Gamma_b$ is called a {\em $\mathcal{T}$-level} iff $\Gamma_b$ has exactly three isotopy classes of arcs. Similarly, $\Gamma_b'$ is a {\em $\mathcal{T}$-level} iff $\Gamma_b'$ has exactly three isotopy classes of arcs.
\end{defin}

The following is used repeatedly in the proof of the Theorem~\ref{thm:q+y,m_b=1}, the Engine Theorem. It describes how weights change moving down levels of $\Gamma, \Gamma'$, starting at a $\mathcal{T}$-level.

\begin{lem}\label{lem:b'>b}
Assume $b$ is such that $\Gamma_b$ is a $\mathcal{T}$-level. Let the total weight of $\Gamma_b$ be $c_1 \wt{w_1} + c_2 \wt{w_2} + c_3 \wt{w_3}$ where $0<c_1 \leq c_2 \leq c_3$. Then if $b'>b$ is the next level such that $\Gamma_{b'}$ has at most three isotopy classes of arcs then either

\begin{enumerate}
\item $b' \geq b + c_1$ and the total weight of $\Gamma_{b'}$ is $(c_2-c_1) {\wt{w_2}}  + (c_3 + c_1) {\wt{w_3}} + c_1 {\wt{w_1 + w_2 - w_3}}$. Furthermore, if $b<b''<b'$, then weights of arcs in $\Gamma_{b''}$ are $\wt{w_1}, \wt{w_2}, \wt{w_3}$, and $\wt{w_1 + w_2 - w_3}$. %If $b+c_1 \leq x$, then $b'=b+c_1$. 
\item $b' \geq b + c_1$ and the total weight of $\Gamma_{b'}$ is $(c_3-c_1) {\wt{w_3}}  + (c_2 + c_1) {\wt{w_2}} + c_1 {\wt{w_1 + w_3 - w_2}}$. Furthermore, if $b<b''<b'$, then weights of arcs in $\Gamma_{b''}$ are $\wt{w_1}, \wt{w_2}, \wt{w_3}$, and $\wt{w_1 + w_3 - w_2}$. %$If $b+c_1 \leq x$, then $b'=b+c_1$. 
\item $b' \geq b + c_2$ and the total weight of $\Gamma_{b'}$ is $(c_3-c_2) {\wt{w_3}}  + (c_1 + c_2) {\wt{w_1}} + c_2 {\wt{w_2 + w_3 - w_1}}$. Furthermore, if $b<b''<b'$, then weights of arcs in $\Gamma_{b''}$ are
$\wt{w_1}, \wt{w_2}, \wt{w_3}$, and $\wt{w_2 + w_3 - w_1}$. %If $b+c_2 \leq x$, then $b'=b+c_2$. 
\item The total weight of $\Gamma_{b'}$ is the same as for $\Gamma_b$. Furthermore, if $b<b''<b'$, then the weights of arcs in $\Gamma_{b''}$ are $\wt{w_1}, \wt{w_2}, \wt{w_3}$ along with one of
\[\bc{\wt{w_1+w_2-w_3}, \wt{w_2 + w_3 - w_1}, \wt{w_3+w_1-w_2}}.\]
The isotopy classes for $\Gamma_{b'}$, as well as the number of arcs in each class, are the same as for $\Gamma_{b}$.
\end{enumerate}

The above holds if $\Gamma$ is replaced by $\Gamma'$.

\end{lem}

\begin{proof}
We assume $\Gamma_b$ has three isotopy classes with corresponding weights  ${\wt{w_1}, \wt{w_2}, \wt{w_3}}$. $\Gamma_b$ has no simple closed level curves. By Lemma~\ref{lem:bandingdiffclasses}$(1),(3)$,  
the banding then connects arcs in two different classes of $\Gamma_b$, removes one arc from each of those classes and creates two new arcs in $\Gamma_{b+1}$.  In going from $\Gamma_b$ to $\Gamma_{b+1}$, assume the isotopy classes of $\Gamma_b$ corresponding to $\wt{w_i}, \,\wt{w_j}$ are banded and denote by $\wt{w_k}$ the third class ($\{i, j, k\}=\{1,2,3\}$).  The classes of $\wt{w_i}, \,\wt{w_j}$ each lose an arc, $\wt{w_k}$ gains an arc, and a new isotopy class is created with weight $\wt{w_4}$, say. By Lemma~\ref{lem:alphaweightTxi}, the sum of the $\alpha$-weights of all arcs in $\Gamma_{b+1}$ must be $0$. That is, we have that 
$c_i w_i+c_j w_j+c_k w_k=0=(c_i-1)w_i+(c_j-1)w_j+(c_k+1)w_k+w_4$.
Thus $\wt{w_4}=\wt{w_i+w_j-w_k}$. If there are four isotopy classes in $\Gamma_{b+1}$, a banding from $\Gamma_{b+1}$ to $\Gamma_{b+2}$ does one of two things (recalling that a level arc must connect different boundary components of $T$)
\begin{itemize}
\item removes an arc from each of the classes of $\wt{w_i}, \,\wt{w_j}$ and adds an arc to each of the classes of $\wt{w_k}, \,\wt{w_4}$; or
\item removes an arc from each of the classes of $\wt{w_k}, \,\wt{w_4}$ ands adds an arc to each of the classes of $\wt{w_i}, \,\wt{w_j}$.
\end{itemize}
This continues as long as level surfaces with four isotopy classes are produced. 
Let $b'>b$ be the next level such that $\Gamma_{b'}$ has at most three isotopy classes. Then for any $b \leq z \leq b'$ there is an $\delta \geq 0$ such that the total weight of $\Gamma_z$ is 
\[(c_i-\delta)\wt{w_i}+(c_j-\delta) \wt{w_j}+(c_k+\delta)\wt{w_k}+(\delta) \wt{w_4}.\]
where $\delta$ is the difference between the increases and decreases in the isotopy class with weight $\wt{w_4}$ in going from $b$ to $b'$. Note that $b'-b \geq \delta$ and $c_i,c_j \geq \delta \geq 0$. As $\Gamma_{b'}$ has at most three isotopy classes, we have that $\delta$ is either $c_i$, $c_j$, or $0$. We arrive at the conclusions of the Lemma, noting that if $\delta=0$ the isotopy classes, as well as their cardinality, in $\Gamma_b$ and in $\Gamma_{b'}$ are the same. As claimed, if $b < b'' < b'$, then the weights of $\Gamma_{b''}$ are $\wt{w_i}, \,\wt{w_j}, \,\wt{w_k}, \,\wt{w_i+w_j-w_k}$.
\end{proof}

\begin{remark}
The change in total weights in Lemma~\ref{lem:b'>b} can be stated as follows. It may be the total weight of $\Gamma_{b'}$ is the same as of $\Gamma_{b}$ (conclusion $(4)$), even though saddles have occurred. Otherwise,
one picks a pair of weights $\wt{w_i}, \wt{w_j}$ of the total weight of $\Gamma_b$. Let $c$ be the minimum of the coefficients of these two weights. The total weight of $\Gamma_{b'}$ loses the weight with smallest coefficient, or loses both weights if the coefficients are the same. The third weight, $\wt{w_k}$, remains in the total weight of $\Gamma_{b'}$ with coefficient increased by $c$. Finally, a new weight, $\wt{w_i}+\wt{w_i}-\wt{w_k}$, is produced with coefficient $c$. 

For example, assume the total weight of $\Gamma_{b}$  is $3 \, \wt{-4} + 4 \, \wt{3} + 6 \, \wt{0} $. Then the total weight of $\Gamma_{b'}$ is one of 

\begin{enumerate}
    \item $1 \, \wt{3} + 9 \, \wt{0} + 3 \,\wt{-1} $
    \item $3 \, \wt{0} + 7 \, \wt{3} + 3 \, \wt{-7} $
    \item $2 \, \wt{0} + 7 \, \wt{-4} + 4 \, \wt{7} $
    \item $3 \, \wt{-4} + 4 \, \wt{3} + 6 \, \wt{0} $ (no change)
\end{enumerate}

As another example, if the total weight of $\Gamma_b$ is $3 \, \wt{-4} + 1 \, \wt{9} + 3 \, \wt{1} $, note that $\Gamma_{b'}$ can have total weight $4 \, \wt{9} + 3 \, \wt{-12}$.
    
\end{remark}

\begin{remark} In our current context of Theorem~\ref{thm:Gr-Li-Lu}, Corollary~\ref{lem:alphaweightr} allows us to assume that each $\Gamma_i, i \neq 0, 2q-2$ has either three or four isotopy classes of arcs. In particular if $b'<2q-2$, $\Gamma_{b'}$ in Lemma~\ref{lem:b'>b} is a $\mathcal{T}$-level. See the remarks following Corollary~\ref{cor:Korder}.
\end{remark}

We take a moment to point out the following, which will also be used repeatedly in the proof of  Theorem~\ref{thm:q+y,m_b=1}, in particular to show the persistence of low weights as the levels increase.

\begin{lem}~\label{lem:ztoz+1}
In going from $\Gamma_z$ to $\Gamma_{z+1}$, an isotopy class loses at most one arc. In particular, if the total weight of $\Gamma_b$ contains the term $c \wt{w}$, then for $b \le j \le b+c-1$, the total weight of $\Gamma_{j}$ will have a term of weight $\wt{w}$.
\end{lem}

\begin{proof}

This follows immediately from Lemma~\ref{lem:bandingdiffclasses}$(4)$.
\end{proof}

As an example of Lemma~\ref{lem:ztoz+1}, assume the total weight of $\Gamma_b$ is $3 \, \wt{-4} + 1 \, \wt{9} + 3 \, \wt{1} $. Then we know that for $b \leq j \leq b+2$, the total weight of $\Gamma_j$ will include a term of weight $\wt{1}$. That is, $m_j \leq 1$ for any $b \leq j \leq b+2$.

For $b<x$ the following gives a fairly explicit description of the total weight.

\begin{lem}\label{thm:m:b<x} \hfill

\begin{enumerate}[label=(\Alph*)]
\item
Assume $0<b<x$ and $\Gamma_b$ has at most three isotopy classes of arcs. Then $\Gamma_b$ is a $\mathcal{T}$-level, and there are positive integers $w_1,w_2$ such that $b=(w_1+w_2-1)y$ and the total weight of $\Gamma_b$ is
\[(x-b) \wt{0} + (w_2 y) \wt{w_1} + (w_1 y) \wt{-w_2}.\]

Then $M_b \leq b/y$ and $m_b=0$.  If $b>b'$ and $\Gamma_{b'}$ is also a $\mathcal{T}$-level then $M_b>M_{b'}$. 
\item If $\iota<x$, then $M_\iota < \iota/y+1$ and $m_\iota=0$.
\end{enumerate}

Similarly,
\begin{enumerate}[label=(\Alph*)]
\item Assume $0<b<x'$ and $\Gamma_b'$ has at most three isotopy classes of arcs. Then $\Gamma_b'$ is a $\mathcal{T}$-level and there are positive integers $w_1,w_2$ such that $b=(w_1+w_2-1)y'$ and the total weight of $\Gamma_b'$ is
\[(x-b) \wt{0} + (w_2 y') \wt{w_1} + (w_1 y') \wt{-w_2}.\]

Then $M_b' \leq b/y'$ and $m_b'=0$.  If $b>b'$ and $\Gamma_{b'}'$ is also a $\mathcal{T}$-level then $M_b'>M_{b'}'$. 
\item If $\iota<x'$, then $M_\iota' < \iota/y'+1$ and $m_\iota'=0$.
\end{enumerate}

\end{lem}

\begin{proof}
We give the proof for $\Gamma$. The proof of the second half, for $\Gamma'$, is the same.

The total weight at $\Gamma_0$ is 
\[(y) \wt{0} + (x) \wt{0}.\]
In $\Gamma_0$, the union of arcs from each of the two isotopy classes, $\hat{\alpha}$, $\hat{\alpha}'$, lie in a neighborhood of  
%recall that $\hat{\alpha}$ is the isotopy class of the $y$ arcs and let $\hat{\alpha}'$ be the class of the $x$ arcs. Then $\hat{\alpha} \cup \hat{\alpha}$ 
$\alpha$, the meridian in $\partial V_{\alpha}$. Note that as $K \cap V_{\alpha}$ lies in a disk, $\alpha$ is trivial in $X_K$. Thus by Lemma~\ref{lem:morseposition}$(b)$ $\Gamma_0$ has no simple closed level curves. Thus any banding is between arcs of intersection. It cannot be a banding of an arc with itself as the resulting simple closed curve in $\Gamma_1$ would also be trivial in $X_K$. We then check that the total weight of $\Gamma_1$ is
\[(y-1) \wt{0} + (x-1) \wt{0} + (1) \wt{1} + (1) \wt{-1}\]
If $y=1$, then $\Gamma_1$ is a $\mathcal{T}$-level of the desired form. If $y>1$, in going from $\Gamma_1$ to $\Gamma_2$ we lose an arc in each of a pair of classes and gain an arc in each of the other.  By Corollary~\ref{cor:xdisappears}, we must lose arcs in classes with weight $\wt{0}$, hence gain in those with weight $\wt{1}, \wt{-1}$. Continuing in this way, we see that if $0<j \leq y$ then $\Gamma_j$ has total weight $(y-j) \wt{0} + (x-j) \wt{0} + (y) \wt{1} + (y) \wt{-1}$. The first $\mathcal{T}$-level ($x>y$ are relatively prime) occurs at $\Gamma_y$ with total weight
\[(x-y) \wt{0} + (y)\wt{1} + (y) \wt{-1}.\]
For the first $\mathcal{T}$-level $\Gamma_y$, we have that $m_y=0$, $M_y=1=y/y$, and $y=(1+1-1)y$.

Now we proceed inductively. Assume the Lemma holds up to and including the  $\mathcal{T}$-level $j$ with total weight
\[(x-j) \wt{0} + (w_2 y) \wt{w_1} + (w_1 y) \wt{-w_2}\]
where $w_1, w_2>0$ and $j=(w_1+w_2-1)y$. Assume that the next level with at most three isotopy classes occurs at level $j'<x$. Applying Lemma~\ref{lem:b'>b} along with Corollary~\ref{cor:xdisappears} (that $x-j'$ is the coefficient of weight $\wt{0}$), we conclude one of the following:
\begin{enumerate}
\item $j'=j+w_2y$ and that the total weight of $\Gamma_{j'}$ is
\[(x-j-w_2y) \wt{0} + ((w_1+w_2) y) \wt{-w_2} + (w_2 y) \wt{w_1+w_2}.\]
\item $j'=j+w_1y$ and that the total weight of $\Gamma_{j'}$ is
\[(x-j-w_1y) \wt{0} + ((w_2+w_1) y) \wt{w_1} + (w_1 y) \wt{-(w_1+w_2)}.\]
\end{enumerate}
In either case we see that $\Gamma_{j'}$ is a $\mathcal{T}$-level of the desired form with $j'=j+w_iy=(w_1+w_2-1)y+w_iy=(w_1+w_2+w_i-1)y$. Then $m_{j'}=0,M_{j'}=w_1+w_2 \leq j'/y$ and $M_j < w_1+w_2= M_{j'}$, as desired. By this inductive procedure we conclude that if $b<x$ and $\Gamma_b$ has at most three isotopy classes, then it is a $\mathcal{T}$-level of the desired form. We also conclude that if $b<b'<x$ are $\mathcal{T}$-levels, then $M_b<M_{b'}$ as claimed. This is part $(A)$.
 
For $(B)$, let $\iota<x$. If $\iota=0$, then $m_\iota=M_\iota=0<\iota/y+1$ as claimed. If $0<\iota \leq y$, from the above we check that $M_\iota =1 < \iota/y+1$ and $m_\iota=0$ as claimed. Otherwise, let $j \leq \iota$ be the largest value such that $\Gamma_j$ is a $\mathcal{T}$-level. By the above, there are positive integers $w_1, w_2$ with  $j=(w_1+w_2-1)y$ such that 
\[(x-j) \wt{0} + (w_2 y) \wt{w_1} + (w_1 y) \wt{-w_2}.\]
If $j=\iota$, then by the above $M_\iota \leq \iota/y < \iota/y+1$ and $m_\iota=0$. Thus we assume $\iota>j$.
By Lemma~\ref{lem:b'>b}, as $\iota$ is between $j$ and the next level with at most three isotopy classes, $M_\iota \leq w_1 + w_2 =j/y+1 <\iota/y+1$ and $m_\iota=0$, as claimed.
\end{proof}

The total weight of $\Gamma_x$, defines the useful constants $r,s,e$.

\begin{thm}\label{thm:=x}
For some $r,s,e>0$ with $ry+sy+e=q$, the total weight for $\Gamma_x$ is either 

\begin{enumerate}
\item $(ry-e) \wt{s} + (sy+e) \wt{-r} + (e) \wt{s+r}$ with $ry>e>0, s>0$; or
\item $(ry+e) \wt{s} + (sy-e) \wt{-r} + (e) \wt{-(s+r)} $ with $sy>e>0,\, r>0$.
\end{enumerate}

Furthermore, if $b \leq x$, $M_b \leq s+r$.

Similarly for some $r',s',e'>0$ with $r'y'+s'y'+e'=q$, the total weight for $\Gamma_{x'}'$ is either 

\begin{enumerate}
\item $(r'y'-e') \wt{s'} + (s'y'+e') \wt{-r'} + (e') \wt{s'+r'}$ with $r'y'>e'>0, \,s'>0$; or
\item $(r'y'+e') \wt{s'} + (s'y'-e') \wt{-r'} + (e') \wt{-(s'+r')}$ with $s'y'>e'>0, \,r'>0$.
\end{enumerate}

Furthermore, if $b \leq x'$, $M_b' \leq s'+r'$.

\end{thm}

%$\spadesuit$ {\em Shouldn't the first part, with $\Gamma$, be instead and up to sign of the weights: $\,\, (ry-\beta e) \,\, \wt{s} + (sy+\beta e) \,\, \wt{-r} + e \,\, \wt{(s+r)}}$ with positive coefficients and $|\beta|=1$? The two listed are the same up to changing the sign of the weights and $r,s$.}

%$\spadesuit$ {\em (Note that we don't really need b to be a $\mathcal{T}$-level, but that is all we need since all of our estimates will be done at such levels.)}

\begin{proof}
We prove the statement for $\Gamma$. The statement for $\Gamma'$ follows similarly.

Let $j <x$ be the largest value such that $\Gamma_j$ is a $\mathcal{T}$-level (as $\Gamma_y$ is such a level, there is one). By Lemma~\ref{thm:m:b<x} there are positive integers $r,s$ such that $j=(r+s-1)y$ and the total weight of $\Gamma_j$ is
\[(x-j) \wt{0} + (r y) \wt{s} + (s y) \wt{-r}.\]
Note that $q=(x-j)+ry+sy$. As $y$ is relatively prime to $q$, $x-j$ is relatively prime to $y$.  
Furthermore, Lemma~\ref{thm:m:b<x} implies that for any $j<i<x$, $\Gamma_i$ has four isotopy classes of arcs. 
By Corollary~\ref{cor:xdisappears}, $\Gamma_x$ is the first level for the class $\hat{\alpha}'$ to vanish. Thus $x$ is the first level past $j$ that has at most three isotopy classes of arcs. By Lemma~\ref{lem:b'>b} and the fact that $x-j$ is relatively prime to $y$, $\Gamma_x$ must be a $\mathcal{T}$-level with total weight either
\[\tag{1}(r y - (x-j)) \wt{s} + (s y + (x-j)) \wt{-r} + (x-j) \wt{s+r};\]
\[\tag{2}(s y - (x-j)) \wt{-r} + (r y + (x-j)) \wt{s} + (x-j) \wt{-(s+r)}.\]

Setting $e=x-j$ gives us the desired form of $\Gamma_x$. Then $M_x=r+s > M_j \geq M_b$ for any $b<x$ such that $\Gamma_b$ is a $\mathcal{T}$-level, by Lemma~\ref{thm:m:b<x}. Furthermore, Lemma~\ref{lem:b'>b} shows that if  $b'>b$ are consecutive $\mathcal{T}$-levels and $b'>j>b$, then $M_j \leq \max\{M_b,M_{b'}\}$. If $b' \leq x$, then $ \max\{M_b,M_{b'}\} \leq r+s$. Thus $j \leq x$ implies $M_j \leq r+s$ as claimed.
%By Theorem~\ref{thm:annulus}(2), the coefficients in the total weight of $\Gamma_x$ are non-zero.
\end{proof}

\section{Proof of Theorem~\ref{thm:Gr-Li-Lu} (Spherical Theorem)}\label{subsec:specialize}

%Results when \pdfmath{K(\lambda)=S^1 \times S^2}. The proof of Theorem~\ref{thm:Gr-Li-Lu} (Spherical Theorem)

In this section, we give the proof of Theorem~\ref{thm:Gr-Li-Lu}, based on the technical results of Sections~\ref{sec:estimates}, \ref{sec:dq+1}, and \ref{sec:d(q-ep)}.  

%establish the results needed to prove Theorem~\ref{thm:Gr-Li-Lu}, that is, when $K(\lambda)=S^1 \times S^2$. In this section, after gathering some theorems setting the context for these sections and then finish with a proof of 

We begin with a restatement of Theorem~\ref{thm:lisca-rasm} as discussed in the proof of Theorem~\ref{thm:Gr-Li-Lu} from subsection~\ref{subsec:pfandoutline}.

\begin{thm}\label{thm:lisca} Let $M$ be a lens space. $M$ is obtained by a longitudinal Dehn surgery on a knot in $S^1 \times S^2$ if and only if there are $d,q \in \Z$ such that up to diffeomorphism, possibly orientation-reversing, $M$ is $L(q^2,m)$ where either
\begin{enumerate}
\item $m=dq + 1$ with $d,q>0$ and $\gcd(d,q)=1,2$. 
\item There is an $\epsilon=\pm 1$ such that $m=d(q-\epsilon)$  with $q>0,\, d \geq 3$, $d$ odd, and $d | 2q + \epsilon$ (properly). We write this $2q + \epsilon=nd$ with $n$ odd and $n \geq 3$. 
\item There is an $\epsilon=\pm 1$ such that $m=d(q-\epsilon)$  with $q>0,\, d \geq 3$, $d$ odd, and $d | q - \epsilon$ (properly). We write this $q - \epsilon =nd$ with $n>1$.
\end{enumerate}
\end{thm}

%\begin{thm}\cite{lisca2007lensSpacesRibbonConjecture, greene2013lensSpaceRealization}
%Let $M$ be a lens space. $M$ contains a knot that admits longitudinal surgery to $S^1 \times S^2$ if and only if up to diffeomorphism, possibly orientation-reversing, there are $d,q \in \Z$ such that $M$ is $L(q^2,m)$ where
%\begin{itemize}
%\item $m=dq+1$ with $\gcd(d,q)=1,2$, or
%\item $m=d(q-1)$ when $d$ odd, and $d$ divides $q-1$ or $2q+1$.
%\end{itemize}
%\end{thm}

\begin{proof}
Let $K$ be a simple knot in a lens space $M$ that admits a longitudinal surgery to $S^1 \times S^2$. Lisca \cite{lisca2007lensSpacesRibbonConjecture} shows that such a lens space (up to a diffeomorphism which is possibly orientation-reversing) is of the form $L(q^2, m)$ where either $(a)$ $m=dq+1$ or $(b)$ $m=d(q-1)$ as described in Theorem~\ref{thm:lisca-rasm}.  Note that in case $(b)$, $d$ is odd. If $d=\pm 1$, then we are in the case $(a)$. We will hereafter assume that in case $(b)$, $|d| \geq 3$. Similarly, in case $(b)$ when $d | 2q+1$, if $d=\pm (2q+1)$, then we are in case $(a)$. We assume then that $d \neq \pm (2q+1)$. In case $(b)$ when $d | q-1$, if $d= \pm (q-1)$, then we are again in case $(a)$. We assume then that $d \neq \pm (q-1)$.

$K$ admits a longitudinal surgery if and only if $-K$ does. For this reason, all homology classes listed are in terms of $\pm [K]$. These homology classes are defined in terms of the canonical generator $\mu$ for the homology. See Definition~\ref{def:mu} for the definition of the canonical generators $\mu, \,\mu'$ of the first homology coming from the cores of the solid tori of the Heegaard splitting. We have the relationship $\mu' \equiv m \mu$ in $H_1(L(q^2,m))$. There is a diffeomorphism between $L(q^2,m)$ and $L(q^2,-m)$ identifying $\mu$ ($\mu'$) of $H_1(L(q^2,m))$ to the generator $\mu$ ($-\mu'$, resp.) of $H_1(L(q^2,-m))$.  This means that if we prove the Theorem for $L(q^2, m)$, we have done so for $L(q^2, -m)$. As a result, we may assume that $m>0$. Taking $q=\epsilon |q|$, case $(b)$ becomes conclusions $(2)$ and $(3)$. Case $(a)$ gives rise to $(1)$ where, taking $d,q>0$, we have to allow the possibility that $m=dq-1$. 

Similarly, there is a diffeomorphism taking $L(q^2,m)$ to $L(q^2,m^{-1})$ identifying the canonical generator $\mu$  ($\mu'$) of $H_1(L(q^2,m))$ to the canonical generator $\mu'$ ($\mu$, resp.) of $H_1(L(q^2,m^{-1}))$, where $m^{-1}m \equiv 1 \pmod{q^2}$. In particular, as $-(dq+1)(dq-1) \equiv 1 \pmod{q^2}$, we get the composition of diffeomorphisms $L(q^2,dq+1)  \rightarrow L(q^2,-(dq-1)) \rightarrow L(q^2,dq-1)$ identifies $\mu$ ($\mu'$) of $L(q^2,dq+1)$ with $\mu'$ ($-\mu$, resp.). On the other hand, in $H_1(L(q^2,dq \pm 1))$ we have that $q\mu' \equiv q(m\mu) \equiv \pm q\mu$. Thus this diffeomorphism between $L(q^2,dq+1)$ and $L(q^2,dq-1)$ identifies the classes $\pm n q\mu$ for any integer $n$. Thus proving the Theorem for $L(q^2,dq+1)$ will verify it for $L(q^2,dq-1)$. Thus case $(a)$ can be written as stated in conclusion $(1)$.
\end{proof}

The next three results appeared first in one form or another in Cebanu's PhD thesis~\cite{cebanu2012propertyR}. 

%\begin{thm}\label{thm:homdetection} Assume a longitudinal surgery along a knot $K$ in a lens space yields $S^1\times S^2$. The unique simple knot in the same homology class as $K$ admits a longitudinal surgery to $S^1\times S^2$.
%\end{thm}

%This section is devoted to the proof of Theorem~\ref{thm:homdetection}. 
%We start with the following fibration result by Ni and Vafaee.

\begin{thm}(\cite{cebanu2012propertyR}, Theorem 3.7.1)\label{thm:cebanufiber}  Suppose $K' \subset S^1 \times S^2$ is a knot with a lens space surgery. Then the exterior of $K'$ in $S^1 \times S^2$ admits a fibration over $S^1$.
\end{thm}

Note that Theorem 1.1 of \cite{nivafaee2019nullSurgery} generalizes this by replacing lens space with L-space.

A knot in $S^1 \times S^2$ is called a {\em spherical braid} if it is transverse to $\theta \times S^2$ for every $\theta \in S^1$. 
%Note that this is equivalent to the exterior of the knot being fibered over $S^1$ with fiber a planar surface, as Dehn filling the exterior of the knot along the unique fiber slope yields an orientable $S^2$-bundle over $S^1$, i.e. $S^1\times S^2$.

\begin{cor}\label{cor:braid} Suppose $K' \subset S^1 \times S^2$ is a knot which admits a lens space surgery. Then $K'$ is a spherical braid in $S^1 \times S^2$.
\end{cor}

\begin{proof}
Let $K' \subset S^1 \times S^2$ be a knot which admits a surgery to a lens space, $M$.
By Theorem~\ref{thm:cebanufiber}, the exterior of $K'$,  admits a fibration over $S^1$ with fiber a surface $Q$, say. As argued in Lemma 1.18 of \cite{baker2016knotsS1xS2}, $\partial Q$ must be a collection of meridians of $K'$ and Dehn filling the exterior of $K'$ along this meridian gives a fibration of $S^1 \times S^2$ over $S^1$ with fiber $\hat Q$---where $\hat Q$ is $Q$ capped off with disks. We conclude that $\hat Q$ is a sphere, writing $K'$ as a spherical braid. 
\end{proof}

From a different perspective, let $K$ be a knot in a lens space with a planar Seifert surface of longitudinal slope $\lambda$ which is the fiber of a fibration of the knot exterior over $S^1$. Then $K(\lambda)=S^1 \times S^2$. The following gives a converse.

\begin{cor}\label{cor:Korder} (See Lemma 3.2.2, Theorem 3.4.1, Section 3.5.4, and Theorem 3.7.1 of \cite{cebanu2012propertyR}) Let $K$ be a knot in a lens space with $K(\lambda)=S^1 \times S^2$ for a longitudinal slope $\lambda$.  Then, up to diffeomorphism, the lens space must be $L(q^2, m)$ for some $q,m$ as described in Theorem~\ref{thm:lisca}.  The exterior of $K$ fibers over the circle with fiber a planar surface. This fiber is a planar Seifert surface for $K$ of slope $\lambda$ and $q$ boundary components.  $[K] \in H_1(L(q^2,m))$ has torsion $q$.
\end{cor}
\begin{proof}
%Let $K$ be a knot in $S^1 \times S^2$ upon which a longitudinal surgery yields a lens space. 
Let $K' \subset S^1 \times S^2$ be the dual knot to the surgery $K(\lambda)$. 
Applying Theorem~\ref{thm:lisca} and Theorem~\ref{thm:cebanufiber}, we conclude that $K \subset L(q^2,m)$ with $q,m$ as stated and that the exterior of $K$ fibers over the circle, with fiber $Q$, say, a planar surface of slope $\lambda$. $Q$ is a Seifert surface for $K$ and, as the fiber of a fibration of the exterior of $K$, is the unique incompressible Seifert surface up to isotopy. If $|\partial Q|=1$, then $Q$ is disk and $X_K$ is a solid torus. As $\lambda$ is longitudinal, $L(q^2,m)$ is $S^3$ and $q=1$. The Corollary holds. 

So assume $|\partial Q|>1$. By Lemma 2.13 of \cite{calegarigordon2013smallRationalGenus}, $[K] \in H_1(L(q^2,m))$ has order $|\partial Q|>1$. Thus $L(q^2,m)$ is not $S^3$, and the fact that $K(\lambda)=S^1 \times S^2$ implies that $X_K$ is irreducible. Theorem~\ref{thm:r>0reduciblemeansreductiblesimple}  says $|\partial Q|=q$. 
%As $K'$ surgers to a rational homology sphere, $H_1(X_{K'})=Z \oplus H$ where $H$ is a finite. As filling $X_{K'}$ along $\mu$ has first homology $Z$, we conclude that $H$ is generated by $[\mu]$. 
%As $K'$ is a spherical braid,
%H_1(X_{K'})$ is generated by $[\mu]$ and $[\nu]$, where $\nu$ is a simple closed curve intersecting $Q$ once. 
%On the other hand, $[\mu]$ has torsion $n$ in $H_1(X_{K'})$
%($Q$ shows that the torsion is a factor of $n$, and $[\mu]=[K] \in H_1(L(q^2,m))$ has torsion $n$). Thus $H_1(X_{K'}) \cong \Z\oplus {\Z/n\Z}$.
%with free factor generated by $[\nu]$ and torsion subgroup generated by $\mu$. 
%Apply Lemma 3.2 of \cite{watson2012khovanovObstructions}, where $M=X_{K'}$, $\lambda_M=\mu$, $|l|=n$, and $M(\alpha)=X_{K'}(\alpha)=L(q^2,m)$. We conclude that  $q^2=|H_1(L(q^2,m)|=c_M \cdot 1$ where $c_M$ is a product of $|l|=n$ and the order, $n$, of the torsion subgroup of $H_1(X_{K'})$. We conclude that $n=q$.
\end{proof}

In summary, let $K$ be a knot (not necessarily simple) in a lens space and $\lambda$ be a longitudinal slope on its exterior, $X_K$. $X_K$ fibers over $S^1$ with planar fiber of slope $\lambda$ if and only $K(\lambda)=S^1 \times S^2$. Note that the dual knot, $K'$, to the $S^1 \times S^2$ surgery is a spherical braid in $S^1 \times S^2$---hence we refer to Theorem~\ref{thm:Gr-Li-Lu} as the Spherical Theorem. 

In sections~\ref{sec:estimates}, \ref{sec:dq+1}, and \ref{sec:d(q-ep)}, we assume $K \subset L(q^2,m)$ with $q,m$ as in Theorem~\ref{thm:lisca} and that $K$ has a planar Seifert surface, $Q$, of longitudinal slope and with $q$ boundary components.

\bigskip

When $K(\lambda)=S^1 \times S^2$, equation~(\ref{eq:yhomology}) of section~\ref{sec:setup} becomes
\[ \pm [K]=yq \mu, \, \pm [K]=y'q \mu' \]
This and the fact that $y<q/2$ says that {\bf $y=1$ iff $\pm [K]=q \mu$}. Similarly, {\bf $y'=1$ iff $\pm [K]=q \mu'=mq \mu$}.
In particular, $y'=1$ iff $\pm [K]=q \mu$ in case $(1)$ of Theorem~\ref{thm:lisca}, and $y'=1$ iff $\pm [K]=dq \mu$ in cases $(2),(3)$ of Theorem~\ref{thm:lisca}. Thus in the proof of Theorem~\ref{thm:Gr-Li-Lu}, we may assume 

\[y,y' > 1\]

Similarly, in the proof of Theorem~\ref{thm:Gr-Li-Lu}, Theorem~\ref{thm:annulus} and Corollary~\ref{cor:nocircletangciesr} allow us to assume

\begin{itemize}\label{lem:alphaweight}
\item For every $i \neq 0, 2q-2$, $\Gamma_i$ has at least three isotopy classes of level arcs; hence also $\Gamma_i'$.
\item For every $i=0, \dots, 2q-2$ there are no simple closed level curves in $\Gamma_i$; hence so for $\Gamma_i'$. 
\item The sum of the ($\alpha$-)weights of the arcs of $\Gamma_i$ is $0$ in ${\mathbb Z}$. The sum of the ($\beta$-)weights of the arcs of $\Gamma_i$ is $0$ in ${\mathbb Z}$.
\item There are no circle tangencies in the level presentation of $Q$.
\end{itemize} 

In sections~\ref{sec:estimates}, \ref{sec:dq+1}, and \ref{sec:d(q-ep)}, we will make these assumptions.

\bigskip

We finish this section with the formal proof of our Spherical Theorem, which we restate here for convenience. Recall that an outline of the ideas of this proof is given at the beginning of Section~\ref{sec:alphabetawts}.

\newcounter{savedthm1}
\setcounter{savedthm1}{\value{thm}}
\newcounter{savedsec6}
\setcounter{savedsec6}{\value{section}}
\setcounter{section}{\value{mainsec}}
\setcounter{thm}{\value{mainthm}}
\begin{thm}
Let $K$ be a simple knot in a lens space that admits a longitudinal surgery to $S^1 \times S^2$. Then $K$ is listed in Theorem~\ref{thm:BBL}. 
\end{thm}

\begin{proof}

By Theorem~\ref{thm:q<600}, if $q \le 600$ then $\pm [K]$ is of the desired form. We hereafter assume that $q>600$. 

The argument now considers the above three cases for $m$. 

Assume we are in case $(1)$ with $m=dq + 1$. Then Theorem~\ref{thm:main:dq+1} shows that $\pm [K] =q \mu, dq \mu$ as claimed. Note that by Corollary~\ref{cor:Korder}, $[K] \in H_1(L(q^2,m))$ as order $q$. Thus $\pm [K] \neq dq \mu$ when $\gcd(d,q)=2$.

Assume we are in case $(2)$ with $m=d(q-\epsilon)$ and $d \, | \, 2q + \epsilon$. Then Theorem~\ref{thm:main1-nd:2q+1} shows that $\pm [K] =q \mu, dq \mu$ has claimed.

Finally, assume we are in case $(3)$ with $m=d(q-\epsilon)$  and $d \, | \, q - \epsilon$. Then Theorem~\ref{thm:main1-nd} implies that $\pm [K] =q \mu, 2q \mu, dq \mu$. Furthermore, by Theorem~\ref{thm:q=1+nd,nmustbe2}, if $\pm [K]=2q \mu$, then $q=2d \pm 1$ as in Theorem~\ref{thm:BBL}.
\end{proof}
\setcounter{thm}{\value{savedthm1}}
\setcounter{section}{\value{savedsec6}}

\section{Estimates on weights}\label{sec:estimates}

Sections~\ref{sec:estimates}, \ref{sec:dq+1}, and \ref{sec:d(q-ep)} gather results needed to prove Theorem~\ref{thm:Gr-Li-Lu}.
In Sections~\ref{sec:estimates}, \ref{sec:dq+1}, and \ref{sec:d(q-ep)}, recall that we are making the assumptions of subsection~\ref{subsec:specialize}, stated after Corollary~\ref{cor:Korder}.
Note that $y,y'<q/2$ and that $x=q-y, \, x'=q-y'$. 

For $b \neq 0, 2q-2$, $\Gamma_b$ has at least three isotopy classes of arcs. Thus the change of total weight of $\Gamma_b$ as $b$ increases is given by Lemma~\ref{lem:b'>b}. 
Indeed, Lemma~\ref{thm:m:b<x} and Theorem~\ref{thm:=x}, describe explicitly the total weight of $\Gamma_b$ when $x \leq b$. The goal of section~\ref{sec:estimates} is to get good bounds on $M_b$, the largest weight (in absolute value) of $\Gamma_b$, when $b>x$. We are able to do this initially, that is approximately through the first half of the levels of $\Gamma$. This is Theorem~\ref{thm:M:b<x+1}. When we try to extend this for $b>q$ we don't always have the control on $M_b$ we would like. However, the loss of appropriate bounds on the largest weight corresponds to the smallest weight of $\Gamma_b$, $m_b$, remaining very small. This is the Engine Theorem, Theorem~\ref{thm:q+y,m_b=1}. The analogous discussion holds for $\Gamma'$. The estimates from the Engine Theorem will give enough restriction that the order of $[K]$ along with the numerics of $L(q^2,m)$ will force $\pm [K]=q\mu, dq\mu, 2q\mu$ (sections~\ref{sec:dq+1} and \ref{sec:d(q-ep)}). Note that the special cases when $q<4y, q<4y$ pop up in the Engine Theorem. The extension of the appropriate estimates in these cases are done separately, though analogously (e.g. Addendum~\ref{thm:b=<q}, Theorem~\ref{thm:q<3y,b>q+2e} ). 

\subsection{Initial Growth Theorem}\label{subsec:b<q}
\hfill

In this subsection, we bound the initial growth (approximately $b \leq q$) of the weights of $\Gamma,\, \Gamma'$.

\begin{thm}[\textbf{Initial Growth Theorem}]\label{thm:M:b<x+1} Let $e, \,e'$ be as in Theorem~\ref{thm:=x} and denote by $[q]_y$ the least positive residue of $q$ modulo $y$.   %Assume $b,b'$ are $\mathcal{T}, \mathcal{T}'$-levels.
\begin{enumerate}
\item If $b \leq x$, $M_b < q/y$. %$\spadesuit$ {\em Also, $m_b \leq \min\{r,s\}$?}
\item If $b' \leq x'$, $M'_{b'} < q/y'$. 
%\item If $b \leq 4y$, $M_b \leq 4$ or $q \leq 5y$. 
\item For every $b \leq q-[q]_y$, $M_b < (y/e+2)(q-e)/y$. %For every $x \leq b \leq q-[q]_y, m_b \in \{r,s\}$. $\spadesuit$ {\em Or, For $b \leq q-[q]_y,  m_b \in \{r,s\}$?}
\item For every $b' \leq q-[q]_{y'}$, $M'_{b'} < (y'/e'+2)(q-e')/y'$.
\end{enumerate}
\end{thm}

\begin{remark}\label{rmk:[q]y} By Theorem~\ref{thm:=x},
$(r+s)y+e=q$; hence, $[e]_y=[q]_y$ and $e=[q]_y+\lfloor e/y \rfloor y$. Here $\lfloor \cdot \rfloor$ denotes the floor function.
\end{remark}

%$\spadesuit$ {\em (Do we need to assume that $b,b'$ are $\mathcal{T}, \mathcal{T}'$-levels for parts $(3),(4)$? If so, can we say for the Engine Theorem that the $b$-levels given there are $\mathcal{T}$?)}

%$\spadesuit$ {\em (I don't think we use the statement about $m_b$ in $(3)$---or the analog in $(4)$.)}

% For an example of Theorem~\ref{thm:M:b<x+1}(1), see Table~\ref{tab:alphacomplexities} where $M_x=M_4=2<7/3=q/y$. For an example of Theorem~\ref{thm:M:b<x+1}(2), see Table~\ref{tab:betaweights} where $M_x'=M_6'=5<7/1=q/y'$.

%$\spadesuit$ {\em (Do we still need $(3),(4)$?)}
The proof of Theorem~\ref{thm:M:b<x+1}(3)(4) with a little further attention will also give the following

\begin{addendum}\label{thm:b=<q} Let $e, \,e'$ be as in Theorem~\ref{thm:=x}.

Assume $3y<q<4y$. One of the following holds:
\begin{enumerate} 
\item For every $j \leq q$, $M_j < (y/e+2)(q-e)/y$. %Furthermore, for every $j \leq q, m_j \leq \max\{r,s\}$.
\item There is a $b \geq q+y$ with $m_b \leq 1$.
\end{enumerate}

Similarly, assume $3y'<q<4y'$. Then one of the following holds:
\begin{enumerate}
\item For every $j \leq q$, $M_j' < (y'/e'+2)(q-e')/y'$. %Furthermore, for every $j \leq q, m_j' \leq \max\{r',s'\}$.
\item There is a $b \geq q+y'$ with $m_b' \leq 1$.
\end{enumerate}

\end{addendum}

\begin{addendum}\label{lem:r=1=sM} Assume $q=2y+1>11$. One of the following holds:
\begin{enumerate}
    \item For every $j \le q-1$,  $M_{j} < q-1$.
    \item There is a $b \ge q+y-2$ such that $m_b\le 1$.
\end{enumerate}

Similarly, assume $q=2y'+1>11$. Then one of the following holds:
\begin{enumerate}
    \item For every $j \le q-1$, $M_{j}' < q-1$.
    \item There is a $b \ge q+y'-2$ such that $m_b'\le 1$.
\end{enumerate}
% In fact,
% \begin{enumerate}
%     \item If $q=2y+1$, then for every $j \le q-1$,  $M_{j} \le q-1$. Furthermore, if $M_{j}=q-1$ for any $j \le q-1$, then there is a $b \ge q+y-2$ such that $m_b=1$.
%     \item If $q=2y'+1$, then for every $j \le q-1$, $M_{j}' \le q-1$. Furthermore, if $M_{j}'=q-1$ for any $j\le q-1$, then there is a $b \ge q+y'-2$ such that $m_b'=1$.
% \end{enumerate} 
\end{addendum}

\noindent {\bf Proof of Theorem~\ref{thm:M:b<x+1} and Addendum~\ref{thm:b=<q}}.

We give the proof of Theorem~\ref{thm:M:b<x+1} and include parenthetical remarks when necessary to imply Addendum~\ref{thm:b=<q}. 

Let $r,s,e>0$ be as in Theorem~\ref{thm:=x}. Parts $(1)$, $(2)$ of Theorem~\ref{thm:M:b<x+1} follow immediately from Theorem~\ref{thm:=x}, as $q=(r+s)y+e=(r'+s')y'+e'$.
%For $(3)$, note that if $q\geq 5y$ then $x \geq 4y$ so we can apply the proof of Theorem~\ref{thm:=x} and the fact there we are going up at most four jumps by $y$.}
%%%%%%%%%%%%%

%\noindent {\bf Proof of Theorem~\ref{thm:M:b<x+1} (3)(4).} 

We prove Theorem~\ref{thm:M:b<x+1} $(3)$ along with the first half of Addendum~\ref{thm:b=<q}. The proof of Theorem~\ref{thm:M:b<x+1} $(4)$ and the second half of Addendum~\ref{thm:b=<q} is identical to this, replacing $\Gamma$ with $\Gamma'$.  

\begin{defin}
$R(b,\alpha,\beta)$ is the statement: There is a $b \geq x+\alpha e$ such that 

\begin{itemize}
\item the total weight of $\Gamma_b$ is 
\[(ry - \beta e) \wt{s} + (sy+ \beta e) \wt{-r} + (e) \wt{\beta (s+r)}\]
with $ \alpha \geq 0, \,|\beta| \leq \alpha+1, \,\alpha e < y$;
\item if $j \leq b$ then $M_j < (y/e+2)(q-e)/y$. 
%\item For every $x \leq j \leq b, m_j \leq \max\{r,s\}$.  $\spadesuit$ {\em (Probably with more work can say $m_j \in \{r,s\}$.)}
\end{itemize}
%$\spadesuit$ {\em Change to $j \leq b$?}
\end{defin}

By Theorem~\ref{thm:=x},  when $j \leq x$, $M_{j} \leq s+r = (q-e)/y <  (y/e+2) (q-e)/y$. Thus from that Theorem we conclude that $R(x,0, \beta)$ holds for some $\beta \in \{\pm1\}$. 

\noindent {\bf Let $b$ be the largest value such that $R(b, \alpha, \beta)$ holds for some $\alpha, \beta$.} 

We may assume that $b <q$, otherwise Theorem~\ref{thm:M:b<x+1}(3) and Addendum~\ref{thm:b=<q}(1) hold. Each of the coefficients in the total weight of $\Gamma_{b}$ are non-zero, otherwise by Theorem~\ref{thm:annulus}, $b=2q-2>q$. The argument now breaks into cases: {\bf (A)} $e = \min \{e, ry-\beta e, sy + \beta e \}$; {\bf (B)} $ry-\beta e = \min \{e, ry-\beta e, sy + \beta e\}$; {\bf (C)} $sy+\beta e= \min \{e, ry-\beta e, sy + \beta e\}$.

We give the proof in cases (A) and (B). The argument for (C) is the same as for (B) with the roles of $s, \,r$ interchanged and with $\beta$ replaced with $-\beta$.

\noindent {\bf Proof in Case (A).}

Applying Lemma~\ref{lem:b'>b} to $\Gamma_b$ we have four possibilities:

\begin{enumerate}
\item There is $b' \geq b+e \geq x+(\alpha+1)e$ and the total weight of $\Gamma_{b'}$ is
\[(ry - (\beta + 1)e) \wt{s} + (sy+ (\beta + 1)e) \wt{-r} + (e) \wt{(\beta+1)(s+r)}.\]
\item There is $b' \geq b+e \geq x+(\alpha+1)e$ and the total weight of $\Gamma_{b'}$ is
\[(sy+ (\beta - 1)e) \wt{-r} + (ry - (\beta - 1)e) \wt{s} +   (e)\wt{(\beta-1)(s+r)}.\]
\item Define $\Delta := \text{min}\{ry-\beta e, sy + \beta e\}$.
\begin{enumerate}
\item $\Delta =ry-\beta e$. There is a 
$b' \geq b+ \Delta \geq x+\alpha e + \Delta$ and the total weight of $\Gamma_{b'}$ is
\[ (sy+ \beta e - \Delta) \wt{-r} + (e+\Delta) \wt{\beta (s+r)} + \Delta \wt{-\beta (s+r)+s-r}.\]
\item $\Delta =sy+\beta e$. There is a 
$b' \geq b+ \Delta  \geq x+\alpha e +\Delta$ and the total weight of $\Gamma_{b'}$ is 
\[ (ry - \beta e - \Delta) \wt{s}  + (e+\Delta) \wt{\beta (s+r)} + \Delta \wt{-\beta (s+r)+s-r}.\]
\end{enumerate}

\item There is a $b'>b$ such that the total weight of $\Gamma_{b'}$ is the same as for $\Gamma_b$. Furthermore if $b<b''<b'$, then weights of arcs in $\Gamma_{b''}$ are in
\[\bc{\wt{-r}, \wt{s}, \wt{\beta (r+s)}, \wt{s-r-\beta (s+r)}, \wt{(\beta+1) (s+r)}, \wt{(\beta-1)(s+r)}}.\] 

\end{enumerate}

%We may assume that each of the coefficients in the total weight of $\Gamma_{b'}$ are non-zero. Otherwise by Theorem~\ref{thm:annulus}, $b'=2q-2$, contradicting that $b <q$. For this reason, as $e\le \Delta$, we exclude the possibility from Lemma~\ref{lem:b'>b} that the coefficient $e-\Delta$ arises in the cases above. 

%In each case above, we take $b'$ to be the largest value satisfying the given weight condition on $\Gamma_{b'}$ and bound on $M_j$ for $j \leq b. 

Lemma~\ref{lem:b'>b} guarantees in all cases that for any $b \leq j \leq b'$, $M_j \leq (|\beta|+1)(s+r) \leq (\alpha+2)(s+r) < (y/e+2)(q-e)/y$ recalling that from $R(b, \alpha, \beta)$ that $\alpha e < y$ and that $ry  +  sy + e=q$. Combined with the assumption of $R(b, \alpha, \beta)$, for $j \leq b'$ we conclude that $M_j < (y/e+2)(q-e)/y$. 

%Similarly, Lemma~\ref{lem:b'>b} guarantees in all cases that for any $b \leq j \leq b'$,
%either ${\bf -r}$ or ${\bf s}$ is included among the weights of $\Gamma_j$ (need to add that to conclusion $(4)$ of Lemma~\ref{lem:b'>b}) . Hence, along with the definition of $b$, for $x \leq j \leq b'$,  $m_j \leq \max\{r,s\}$ as desired. 

We show in each case that $b' \geq q-[q]_y$ to prove Theorem~\ref{thm:M:b<x+1}(3); we show that, in fact, except as described in Addendum~\ref{thm:b=<q}(2), $b' \geq q$---thereby verifying the Addendum. 

Assume $\Gamma_{b'}$ is as in $(1)$, $(2)$ or $(4)$ above. In $(1)$ or $(2)$, if $(\alpha + 1)e\ge y$, then $b' \geq x + (\alpha +1)e \geq x + y=q$. Thus we assume in $(1),(2)$ that $(\alpha + 1)e < y$. Now in each of $(1)$, $(2)$ and $(4)$, $b'>b$ and $R(b', \alpha+1, \beta+1)$, $R(b', \alpha+1, \beta-1)$ or $R(b', \alpha, \beta)$ holds respectively, which contradicts the definition of $b$.  

Assume $\Gamma_{b'}$ is as in $(3)$ above. First assume $\Delta=ry-\beta e$. Then $b' \geq x + \alpha e+ \Delta =x+ ry +(\alpha-\beta)e \geq x + ry -e$. If $\beta \leq 0$, then $b' \geq q$ as desired. So we assume $\beta >0$.  By Remark~\ref{rmk:[q]y}, $e=[q]_y+\lfloor e/y \rfloor y$. By assumption, $ry-e \geq ry-\beta e>0$. Thus $r> \lfloor e/y \rfloor$ and $ry-e =y +(r-1)y-e \geq y+\lfloor e/y \rfloor y-e = y-[q]_y$. That is $ b' \geq x + ry -e \geq x+y-[q]_y = q - [q]_y$ as desired. %This verifies Theorem~\ref{thm:M:b<x+1}(4). 

%If $r \geq 2$, then $b' \geq q+y-e \geq q$. 
%Again the estimates on $M_j, m_j$ verify Theorem~\ref{thm:M:b<x+1}(4). Similarly, if $\beta <0$, then $b' \geq q+e > q$, and the Theorem holds. So we assume $r=1$ and $\beta \geq 0$. Then for $b''=b' + (sy+\beta e-\Delta)-1$, we are assured that $m_{b''}=1$. But $b'' \geq x + sy + (\alpha + \beta)e -1$. If $s>1$, then $b'' \geq q + y-1\geq q+e$ as desired. If $s=1$, then $ry-\beta e < sy + \beta e$ (the coefficents of the total weight are all non-zero) implies that  $\beta >0$. Consequently, $b'' \geq q + e-1$ as desired. 

\noindent { (Addendum~\ref{thm:b=<q}: Recall that $3y<q<4y$. Note that above we have that $b' \geq q$ unless $\beta>0$ and $r=1+\lfloor e/y\rfloor$. If $e \geq y$, then $r \geq 2$ and $q=(r+s)y+e\ge 4y$, a contradiction. Thus we assume $e<y$ and conclude $r=1$. Then $q=(r+s)y+e>3y$ implies that $s\ge 2$. The total weight of $\Gamma_{b'}$ is  $ (sy+ \beta e - \Delta) \wt{-r} + (e+\Delta) \wt{\beta (s+r)} + \Delta \wt{-\beta (s+r)+s-r}$ and $m_{b'} \leq r=1$.
Because an isotopy class of $\Gamma_z$ loses at most one arc as we move to $\Gamma_{z+1}$, we have that $m_{b''} \leq 1$ for $b''=b'+ (sy + \beta e - \Delta -1) \geq x + \alpha e + \Delta +  (sy + \beta e - \Delta -1) = x +sy + (\alpha + \beta)e -1$. As $s\geq 2$, $b'' \geq q+y$, satisfying conclusion $(2)$ of Addendum~\ref{thm:b=<q}.)

In case $(3)$ when $\Delta=sy+\beta e$, we apply the same argument replacing $r$  with $s$ and $\beta$ with $-\beta$ to verify the Theorem and the Addendum.

This completes the proof of  Theorem~\ref{thm:M:b<x+1}(3) and the Addendum~\ref{thm:b=<q} in Case (A).

\noindent {\bf Proof in Case (B).} 

From Case (A) we may assume $e>ry-\beta e$, that is $(\beta + 1)e>ry$. In particular, $\beta \geq 0$. 

Applying Lemma~\ref{lem:b'>b} to $\Gamma_b$ we have four possibilities:

\begin{enumerate}
\item There is $b' \geq b+ ry - \beta e \geq x+ ry + (\alpha-\beta)e$ and the total weight of $\Gamma_{b'}$ is
\[(sy+ ry) \wt{-r} + ((\beta + 1)e-ry) \wt{\beta (s+r)} + (ry - \beta e) \wt{(\beta+1)(s+r)}.\]
\item  There is a 
$b' \geq b+ ry - \beta e \geq x+ ry + (\alpha-\beta)e$ where the total weight of $\Gamma_{b'}$ is 
\[(sy -ry + 2\beta e) \wt{-r} + (ry- (\beta-1)e) \wt{\beta (s+r)} + (ry-\beta e) \wt{-\beta (s+r)+s-r}.\]
\item Define $\Delta := \text{min}\{e,sy+\beta e\}$.
\begin{enumerate}
\item $\Delta =e$. There is $b' \geq b+e \geq x+(\alpha+1)e$ and the total weight of $\Gamma_{b'}$ is
\[(sy+ (\beta - 1)e) \wt{-r} + (ry - (\beta - 1)e) \wt{s} +   (e) \wt{(\beta-1)(s+r)}.\]
\item $\Delta =sy+\beta e$. Then $\beta=0$. There is a 
$b' \geq b+ \Delta  \geq x+\alpha e + sy$ and the total weight of $\Gamma_{b'}$ is 
\[(e-sy) \wt{0}  + (ry+sy) \wt{s} + (sy) \wt{-(s+r)}.\]
\end{enumerate}
\item There is a $b'>b$ such that the total weight of $\Gamma_{b'}$ is the same as for $\Gamma_b$. Furthermore if $b<b''<b'$, then weights of arcs in $\Gamma_{b''}$ are in 
\[\bc{\wt{-r}, \wt{s}, \wt{\beta (r+s)}, \wt{s-r-\beta (s+r)}, \wt{(\beta+1) (s+r)}, \wt{(\beta-1)(s+r)}}.\] 

\end{enumerate}

%Note again that we may assume that the coefficients in the total weight of $\Gamma_{b'}$ are non-zero by Theorem~\ref{thm:annulus}.

%In each case above, we take $b'$ to be the largest value satisfying the given weight condition on $\Gamma_{b'}$ and bound on $M_j$ for $j \leq b. 

As in Case (A), Lemma~\ref{lem:b'>b} guarantees in all cases that for any $b \leq j \leq b'$, $M_j \leq (|\beta|+1)(s+r) < (y/e+2)(q-e)/y$ recalling that from $R(b, \alpha, \beta)$ that $\alpha e < y$, $|\beta| \le \alpha+1$ and $ry  +  sy + e=q$. Combined with $R(b, \alpha, \beta)$, for $j \leq b'$ we conclude that $M_j <  ((y+2e)/e)(q-e)/y$. 

%Similarly, Lemma~\ref{lem:b'>b} guarantees in all cases that for any $b \leq j \leq b'$,
%either ${\bf -r}$ or ${\bf s}$ is included among the weights of $\Gamma_j$ (need to add that to conclusion $(4)$ of Lemma~\ref{lem:b'>b}) . Hence, along with the definition of $b$, for $x \leq j \leq b'$,  $m_j \leq \max\{r,s\}$ as desired. 

First assume $\Gamma_{b'}$ is as in $(4)$ above. The restriction on $M_j$ shows that $R(b', \alpha, \beta)$ holds, contradicting the definition of $b$.

We show that in cases (1)--(3), $b' \geq q - [q]_y$ to verify Theorem~\ref{thm:M:b<x+1}(3); furthermore, we show that in fact, except as described in Addendum~\ref{thm:b=<q}(2), that $b' \geq q$.

In cases $(1)$ and $(2)$ we have that $b' \geq x +ry + (\alpha-\beta)e \geq x+ry-e$. If $\beta \le 0$, then $b' \geq x+y=q$. So we assume $\beta>0$. By Remark~\ref{rmk:[q]y}, $e=[q]_y+\lfloor e/y \rfloor y$. As $ry-e \geq ry-\beta e>0$,  $r> \lfloor e/y \rfloor$ and $ry-e \geq y-[q]_y$. That is, $ b' \geq x + ry -e \geq x+y-[q]_y = q - [q]_y$ as claimed by Theorem~\ref{thm:M:b<x+1}.

So assume $\Gamma_{b'}$ is as in $(3)$ above. In $(3)(b)$, $b' \geq x+y=q$. Thus we assume $(3)(a)$ occurs. If $(\alpha+1)e < y$ then $R(b', \alpha+1, \beta-1)$ holds, contradicting the definition of $b$.  So we assume
$(\alpha+1)e \geq y$, then $b' \geq x + (\alpha +1)e \geq x + y =q$.

This completes the proof in Case (B) of Theorem~\ref{thm:M:b<x+1}$(3)$. 

\noindent (Addendum~\ref{thm:b=<q}: Note that above we have that $b' \geq q$ as desired unless we are in cases $(1), (2)$, with $\beta>0$ and $r=1+\lfloor e/y\rfloor$, which we now assume. We proceed to show that conclusion $(2)$ of the Addendum must hold. Recall that $3y<q<4y$. If $e \geq y$, then we conclude that $r \geq 2$ and $q=(r+s)y+e\ge 4y$, a contradiction. Thus we assume $e<y$ and conclude $r=1$. Then, $q=(r+s)y+e > 3y$ implies $s \geq 2$. 

An isotopy class of $\Gamma_z$ loses at most one arc as we move to $\Gamma_{z+1}$. Thus looking at the total weight of $\Gamma_{b'}$, we have that $m_{b''} \leq 1$ when
\begin{enumerate}[label=Case (\arabic*):, left=1em]
    \item $b''=b'+ (s+r)y -1 \geq x + \alpha e + ry-\beta e +  (s+r)y -1 \geq x +(s+r)y \ge q + 2y$;
    \item $b''=b'+ sy - ry + 2\beta e - 1 \geq x + \alpha e + ry - \beta e + sy - ry + 2\beta e - 1 = x +sy + (\alpha + \beta)e -1 \ge q+y$.
\end{enumerate}

These give conclusion $(2)$ of Addendum~\ref{thm:b=<q}.)

\qed{ (Theorem~\ref{thm:M:b<x+1} and Addendum~\ref{thm:b=<q} )}

\noindent {\bf Proof of Addendum~\ref{lem:r=1=sM}}.

We prove the first half of Addendum~\ref{lem:r=1=sM}. The second half follows replacing $\Gamma$ with $\Gamma'$. Note that $q=2y+1$ implies $r=s=e=1$. We follow the proof above for Theorem~\ref{thm:M:b<x+1} and Addendum~\ref{thm:b=<q} using these values. We modify $R(b, \alpha, \beta)$ by changing the conditions $\alpha e<y$ to $\alpha \le y-2$ and $M_j < (y/e+2)(q-e)/y$ to $M_j < q-1$. As $q>11$, Theorem~\ref{thm:=x} shows that $R(x,0, \beta)$ holds for some $\beta \in \{\pm1\}$. 
As $e=1$, only case (A) occurs. Again consider the possibilities $(1)$--$(4)$ for $\Gamma_{b'}$. In all cases, for all $b\le j\le b'$, $M_j \le (|\beta|+1) |(r+s) \le 2(\alpha+2) \le 2(y-2+2)=q-1$ since $\alpha \leq y-2$, with $M_j=q-1$ only if $\alpha=y-2$. 

In cases $(1)$, $(2)$, and $(4)$, first assume $\alpha=y-2$. In the total weight of  $\Gamma_{b'}$ we find that at least one of the coefficients associated to weight $\wt{\pm 1}$ is at least $y$, so for $j=b'+y-1 \geq x+\alpha+1+y-1= q+\alpha=q+y-2$ that $m_j \le 1$ as claimed in conclusion $(2)$ of the Addendum. Thus we assume in $(1)$, $(2)$, and $(4)$ that $\alpha \leq y-3$ and consequently for $b\le j\le b'$, $M_{j} \leq 2(|\beta|+1)<q-1$. But then in $(1)$, $(2)$, and $(4)$, $b'>b$ again contradicts the definition of $b$.

Thus we assume $(3)$ holds. Then $b' \geq b+ \Delta \geq x+\alpha +y-|\beta|=q+(\alpha-|\beta|) \geq q-1$ since by $R(b,\alpha,\beta)$ we have $|\beta| \leq \alpha+1$. Then if $b<j \leq b'$, $M_{j} \leq 2|\beta| \leq 2(\alpha+1) \leq 2(y-1) = q-3$, recalling that $\alpha \leq y-2$. By $R(b,\alpha,\beta$) then, for any $j \leq q-1\leq b'$, we have $M_j<q-1$ establishing the Addendum.

\qed{ (Addendum~\ref{lem:r=1=sM})}

\subsection{Engine Theorem }\label{subsec:engine}

In this subsection, we give extensions of Theorem~\ref{thm:M:b<x+1}(3)(4) and Addendum~\ref{thm:b=<q} to levels $b >q$: Theorems~\ref{thm:q<3y,b>q+2e} and \ref{thm:q+y,m_b=1}. Their proofs follow the same line of reasoning. As a warm-up, we start with the Theorem~\ref{thm:q<3y,b>q+2e}, which gives useful estimates in the situation that $q<3y$. Most of the section is devoted to Theorem~\ref{thm:q+y,m_b=1}, which is an engine for the proof of the Spherical Theorem, Theorem~\ref{thm:Gr-Li-Lu}.

\begin{thm}\label{thm:q<3y,b>q+2e}
If $q<3y$ then either 
\begin{enumerate} 
\item there is a $b \geq q +2e$ with $m_b \leq 1$; or
\item there is a $b \geq q + \min\{e, y-e \}$ with $M_b<2q$; or
\item there is a $b \geq q +y$ with $m_b \leq 1$.
\end{enumerate}

Similarly, if $q<3y'$ then either 
\begin{enumerate} 
\item there is a $b \geq q +2e'$ with $m_b' \leq 1$; or
\item there is a $b \geq q + \min\{e', y'-e' \}$ with $M_b'<2q$; or
\item there is a $b \geq q +y'$ with $m_b' \leq 1$.
\end{enumerate}
\end{thm}

\begin{proof} 

We prove the first part, about $m_b, M_b$. The argument of the second part is the same replacing $\Gamma$ with $\Gamma'$.

Let $s, r, e>0$ as in Theorem~\ref{thm:=x}. As $y(s+r)+e=q<3y$, we conclude that $s=r=1$ and $e<y$. Then the total weight of $\Gamma_x$ is either
\begin{enumerate}
\item $(y+e) \wt{1} + (y-e) \wt{-1} + (e) \wt{-2} $; or
\item $(y-e) \wt{1} + (y+e) \wt{-1} + (e) \wt{2}$.
\end{enumerate}

We assume that the total weight of $\Gamma_x$ is given by $(1)$
as the argument for $(2)$ is the same (one gets $(1)$ by taking the negative of all individual weights in $(2)$). 

We choose $b$ to be the largest value such that the total weight of $\Gamma_b$ is
\[(y+e) \wt{1} + (y-e) \wt{-1} + (e) \wt{-2}.\]
Then $b \geq x$.

The definition of $b$ along with  Lemma~\ref{lem:b'>b} implies that either:
\begin{enumerate}[label=(\Alph*)] 
\item there is a $b' \geq b+(y-e) \geq x + (y-e)$ such that $\Gamma_{b'}$ has total weight
\[(2e) \wt{1} + (y) \wt{-2} + (y-e) \wt{2};\]
\item there is a $b' \geq b +e \geq x+e$ such that $\Gamma_{b'}$ has total weight
\[(y) \wt{1} + (y) \wt{-1} + (e) \wt{0};\]
\item there is a $b' \geq b +\Delta'  \geq x+\Delta'$, where $\Delta'=\min\{e, y-e\}$, such that $\Gamma_{b'}$  has total weight 
\begin{enumerate}[label=(\roman*)]
\item $ (y-2e) \wt{-1} + (y+2e) \wt{1} + (e) \wt{-4}$ when $\Delta'=e$,
\item $(2e-y) \wt{-2} + (2y) \wt{1} + (y-e) \wt{-4}$ when $\Delta'=y-e$.
\end{enumerate}
\end{enumerate}

We treat these in turn.

\noindent {\bf I. Assume $(A)$ holds.}

\begin{defin}
$P(j,\alpha)$ is the statement: There is a $j \geq b+\alpha (y-e)$ such that the total weight of $\Gamma_j$ is
\[((y+e) - \alpha (y- e)) \wt{1} + (e + \alpha (y-e)) \wt{-2} + (y-e) \wt{3\alpha-1},\]
with $\alpha \ge 1,  \,(\alpha -1) (y-e) < y$.
\end{defin}

Then $P(b',1)$ holds. Assume $j$ is the largest value such that $P(j, \alpha)$ holds for some $\alpha$. 

We again apply Lemma~\ref{lem:b'>b}, along with the definition of $j$, to conclude one of the following:
\begin{enumerate}[label=(\alph*)] 
\item Let $\Delta'=\min \{(y+e)-\alpha (y-e), e + \alpha (y-e)\}$. There is a $j' \geq j+\Delta' $ such that the weights of
$\Gamma_{j'}$ are in $\{ 1,-2,3\alpha-1, -3\alpha \}$
%\begin{enumerate}[label=(\roman*)]
%\item when $\Delta'=(y+e)-\alpha (y-e)$
%$$(2\alpha (y-e)-y) \wt{-2} + (2y-\alpha(y-e)) \wt{2+3(\alpha-1)} + ((y+e)-\alpha (y-e)) \wt{-(3+3(\alpha-1))}$$
%\item when $\Delta'=e+\alpha (y-e)$
%$$(y-2\alpha (y-e)) \wt{1} + (e+(\alpha+1)(y-e)) \wt{2+3(\alpha-1)} + (e+\alpha (y-e)) \wt{-(3+3(\alpha-1))}$$
%\end{enumerate}
\item Let $\Delta'=\min \{(y+e)-\alpha (y-e), y-e\}$. There is a $j' \geq j+\Delta' $ such that the weights of $\Gamma_{j'}$ lie in $\{1,-2,3\alpha-1, 3\alpha+2\}$.
%\begin{enumerate}[label=(\roman*)]
%\item when $\Delta'=(y+e)-\alpha (y-e)$
%$$ ((\alpha+1)(y-e)-(y+e)) \wt{2+3(\alpha-1)}+(y+2e) \wt{-2}+ ((y+e)-\alpha (y-e)) \wt{2+3\alpha}$$
%\item when $\Delta'=y-e$
%$$ ((y+e)-(\alpha+1)(y-e)) \wt{1}+((e+(\alpha+1) (y-e)) \wt{-2}+ (y-e) \wt{2+3\alpha}$$
Furthermore, when $\Delta'=y-e$, the definition of $j$ implies that $\alpha(y-e) \geq y$. 
%\end{enumerate}
\item there is a $j' \geq j +(y-e)  \geq x+ (\alpha+1)(y-e)$ such that $\Gamma_{j'}$ has total weight
\[(e+(\alpha-1)(y-e)) \wt{-2} + ((y+e)-(\alpha-1)(y-e)) \wt{1} + (y-e) \wt{3\alpha-4}.\]
\end{enumerate}

\bigskip

\begin{clm}\label{clm:j'q<3y}
With $j'$ as above, $M_{j'}<2q$.
\end{clm}

\noindent {\em Proof of Claim:}
$M_{j'} \leq 2+3 \alpha$. As $\alpha-1 <{y}/\pth{y-e}$ by the definition of $P(j,\alpha)$, we get that $M_{j'} < 5+{3y}/\pth{y-e}$. Recall that $q=(r+s)y+e=2y+e$. If $y-e \geq 2$, then $M_j \leq 5+{3y}/{2}<5+(q-e)<2q$. If $y-e=1$, then $\alpha\le y$ and $M_j\le 2+3y = 3+2y+e = 3+q < 2q$.

\qed {\em $\,$ (Claim)}
\medskip

Assume $(a)$ above holds. Then $\Delta'=\min \{(y+e)-\alpha (y-e), e + \alpha (y-e)\}$ and $j' \geq j+\Delta' \geq b'+\Delta' \geq x + \alpha (y-e) +\Delta' =x+\min \{y+e, y + (2\alpha-1) (y-e)\} \ge q + \min \{e,y-e\}$. Together with Claim~\ref{clm:j'q<3y}, we obtain conclusion $(2)$ of Theorem~\ref{thm:q<3y,b>q+2e}.

Assume $(b)$ holds. If $\Delta'=(y+e)-\alpha (y-e)$, then $j' \geq j+ (y+e)-\alpha (y-e) \geq x+y+e=q+e$. If $\Delta'= y-e$, then $\alpha(y-e) \geq y$ and $j' \geq j+(y-e) \geq x+ \alpha(y-e)+(y-e) \geq x+y+y-e=q+y-e$. In both cases, conclusion $(2)$ of Theorem~\ref{thm:q<3y,b>q+2e} holds.

Assume $(c)$ holds. Then $P(j', \alpha-1)$ holds, contradicting the definition of $j$ unless $\alpha=1$. In that case, we have $j' > j > b$ and $\Gamma_{j'}$ with total weight $(e) \wt{-2} + (y+e) \wt{1} +(y-e) \wt{-1}$, contradicting the definition of $b$.

\noindent {\bf II. Assume $(B)$ holds.}

Taking $b'$ as largest with the given total weight, we apply Lemma~\ref{lem:b'>b} again to conclude one of the following.
\begin{enumerate}[label=(\alph*)] 
\item There is a $b'' \geq b'+y \geq x+y+e = q+e$ such that $\Gamma_{j'}$ has total weight
\[ (y+e) \wt{0} + (y) \wt{0}.\]
%\begin{enumerate}[label=(\roman*)]
%\item when $\Delta'=(y+e)-\alpha (y-e)$
%$$(2\alpha (y-e)-y) \wt{-2} + (2y-\alpha(y-e)) \wt{2+3(\alpha-1)} + ((y+e)-\alpha (y-e)) \wt{-(3+3(\alpha-1))}$$
%\item when $\Delta'=e+\alpha (y-e)$
%$$(y-2\alpha (y-e)) \wt{1} + (e+(\alpha+1)(y-e)) \wt{2+3(\alpha-1)} + (e+\alpha (y-e)) \wt{-(3+3(\alpha-1))}$$
%\end{enumerate}

\item  There is a $b'' \geq b'+e \geq x+2e $ such that $\Gamma_{b''}$ has total weight
\[ (y-e) \wt{1}+(y+e) \wt{-1}+ (e) \wt{2}.\]
\item  There is a $b'' \geq b'+e \geq x+2e $ such that $\Gamma_{b''}$ has total weight
\[ (y-e) \wt{-1}+(y+e) \wt{1}+ (e) \wt{-2}.\]
\end{enumerate}

If $(a)$ holds, then $M_{b''}<2q$ and $b'' \geq q+e$,  and conclusion $(2)$ of Theorem~\ref{thm:q<3y,b>q+2e} holds. If $(b)$ or $(c)$ holds, then $m_{b''+y+e-1} \leq 1$. And $b''+y+e-1\geq x+y+3e-1\geq q+2e$ verifies conclusion $(1)$ of the Theorem. 

\noindent {\bf III. Assume $(C)$ holds.}

Assume $(C)(i)$ holds. Then $m_{b'+y+2e-1} \leq 1$. As $b'+y+2e-1 \geq x+y+3e-1\geq q+2e$, conclusion $(1)$ of the Theorem holds. If $(C)(ii)$ holds, then $m_{b'+2y-1} \leq 1$. As $b'+2y-1 \geq x+(y-e)+2y-1=q+2y-e-1\geq q+y$, conclusion $(3)$ of the Theorem holds.
\end{proof}

\bigskip

We now come to the main theorem of this subsection, Theorem~\ref{thm:q+y,m_b=1}. Its proof follows the theme of the above, but significantly more possibilities must be pursued.

\begin{thm}[\textbf{Engine Theorem}]\label{thm:q+y,m_b=1}
%Let $e, \,r, \,s, \,e', \,r', \,s'$ be as in Theorem~\ref{thm:=x}. %Assume $q > 8$.

\noindent One of the following holds.

\begin{enumerate}
%\item  there is a $b \geq q+y-1$ such that $m_b =1$; 
\item There is a $b \geq 2q-y-2$ and $m_b \leq 3$.

\item There is a $b \geq q+y-1$ and $M_b < 2q$.

%\item  $q>4y$ and there is a $b \geq 2q-2y-2$ such that $m_b \leq 1$;

\item There is a $b \geq q+y-1$ and $m_b \leq 1$.
%\item  $q<4y$ and there is a $b \geq q+ \max \{e,y-e \} -2$ such that $m_b \leq 2$.
%\item  $q<4y$ and there is a $b \geq 2q-y-2$ such that $m_b \leq 3$.
%Before we had: \item  there is a $b \geq q+x-2y$ s.t. $m_b \leq 3$. $\spadesuit$ {\em (Maybe all we need here is $b \geq q+y-2$ as for $(1)$.)}

\item $q < 4y$. %($q < 3y$)?
\end{enumerate}

\noindent Similarly, one of the following holds.

\begin{enumerate}
\item There is a $b \geq 2q-y'-2$ and $m_b' \leq 3$.
%\item  there is a $b \geq q+y'-1$ such that $m_b' =1$;
\item There is a $b \geq q + y' - 1$ and $M_b' < 2q$.
\item There is a $b \geq q+y'-1$ and $m_b' \leq 1$.
\item $q < 4y'$. %($q < 3y$)?
%\item  $q<4y'$ and there is a $b \geq q+ \max \{e',y'-e' \}-2$ such that $m_b' \leq 2$;
%\item  $q<4y'$ and there is a $b \geq 2q-y'-2$ such that $m_b' \leq 3$;
%\item Before we had: there is a $b \geq q+x'-2y'$ such that $m_b' \leq 3$. 
\end{enumerate}

\end{thm}

% $\spadesuit$ {\em (For (1), we might be able to sharpen to $b \geq 2q-y-1$ if that turns out to be helpful. For (2), we might be able to sharpen $M_{b}<2q$ to $M_{b}<2q - 2$, if that turns out to be helpful. (John: Doesn't look like we need any of these improvements.))}

% If we can improve $(4)$ to q < 3y, then here $q \ge 6$, we might need $7$ though.
\begin{proof}

We prove the first part, about $m_b, \,M_b$. The argument of the second part is the same replacing $\Gamma$ with $\Gamma'$. Recall that $y,y'>1$.

We assume for contradiction that Theorem~\ref{thm:q+y,m_b=1} does not hold. Since $y$ is relatively prime to $q$, we hereafter assume $q > 4y \ge 8$.

Let $e, \,r, \,s, \,e', \,r', \,s'$ be as in Theorem~\ref{thm:=x}.

\begin{defin}
For integers $b,\alpha, \beta$, let $P(b,\alpha,\beta)$ is the statement: There is a $b \geq x+\alpha e$ such that the total weight of $\Gamma_b$ is 
\[ (ry - \beta e) \wt{s} + (sy+ \beta e) \wt{-r} + (e) \wt{\beta (s+r)}, \]
with $r,s >0,\, \alpha \ge 0,\, \abs{\beta} \leq \alpha+1,\,  \alpha e < 2y - 1$.
\end{defin}

%As noted in Definition~\ref{def:totalalpha} the sum of the coefficients at any level is $q$. Thus from the total weight of $\Gamma_b$, $q=ry+sy+e$. Thus $4y>(r+s)y+e$ with $r,s>0$ gives the desired results.
%\end{proof}
 
Theorem~\ref{thm:=x}$(1)$ implies $P(x,0,1)$, and Theorem~\ref{thm:=x}$(2)$ implies $P(x,0,-1)$. 
Let $2q-2 \geq b \geq x$ be the largest $b$ for which $P(b, \alpha, \beta)$ holds for some $\alpha, \beta$.

% \begin{clm}
% $M_{b} < 2q$
% \end{clm}

% \noindent {\em Proof of Claim:}
% If $\beta=0$, then the identity $y(s+r)+e=q$ shows that $M_b<q$. So we assume
% $\beta \neq 0$ and $M_b \geq 2q$. Then $ |\beta| (s+r) =M_b \geq 2q=2(y(s+r)+e)$. Then ($|\beta| -2y) (s+r) \geq 2e>0$. By $P(b,\alpha,\beta)$, $2y +1 > \alpha +1 \geq |\beta| \geq 2e/(s+r) +2y$ implying that $2y+1>|\beta| > 2y$, a contradiction. \qed {\em $\,$ (Claim)}

% The Claim, along with our assumption that Theorem~\ref{thm:q+y,m_b=1}(2) does not hold, implies
% \[b < q+y-1.\] 
% $\spadesuit$ {\em (I don't think we need the above fact about $b$. Check. If so, delete the above Claim. Don't need it for Case (I).) }

Applying Lemma~\ref{lem:b'>b} for $\Gamma_{b'}$, the next level after $\Gamma_b$ with at most three isotopy classes of arcs, we find the weights of arcs in $\Gamma_{b'}$ are among
\[\bc{\wt{s}, \wt{-r}, \wt{\beta(s+r)}, \wt{-\beta(s+r)+s-r}, \wt{(\beta+1)(s+r)}, \wt{(\beta-1)(s+r)}}.\]

\begin{clm}\label{clm:Mb'}
$M_{b'} < 2q$.
\end{clm}

\noindent {\em Proof of Claim:}

By $P(b,\alpha, \beta)$, $\alpha e < 2y - 1$, so $2y \ge \alpha e + 2 \ge \alpha + 2$. Therefore, $M_{b'} \leq (|\beta|+1) (s+r) \leq (\alpha+2) (q-e)/y \le 2(q - e) < 2q$. Here again we are using that $ry+sy+e=q$. %Significantly shortened after changing $\alpha 2 < 2y$ to $2y - 1$.
\qed {\em $\,$ (Claim)}
\medskip

By Claim~\ref{clm:Mb'}, if $b'\ge q+y-1$ then Theorem~\ref{thm:q+y,m_b=1}(2) holds.
We hereafter assume \[b' < q+y-1.\]
Note that this implies, by Theorem~\ref{thm:annulus}, that $\Gamma_{b}$ and $\Gamma_{b'}$ each have three isotopy classes of arcs. In particular, we have that \[ry-\beta e, sy+\beta e >0\]

%The argument of Claim~\ref{clm:Mb'} shows that $M_b <2q$. Thus we may assume $ry-\beta e, sy+\beta e >0$, otherwise Theorem~\ref{thm:annulus} implies that $b=2q-2$ verifying conclusion $(2)$ of Theorem~\ref{thm:q+y,m_b=1}.  Similarly, by Claim~\ref{clm:Mb'}, we may assume that $\Gamma_{b'}$ has three isotopy classes of arcs.

For specificity, we apply Lemma~\ref{lem:b'>b} under the following three situations: 

\begin{enumerate}[label=\textbf{Case \Roman*:}, left=1em]
    \item $\textbf{min}\{ ry-\beta e, sy+\beta e, e \}=e$.
    \item $\textbf{min}\{ ry-\beta e, sy+\beta e, e \}=ry-\beta e$.
    \item $\textbf{min}\{ ry-\beta e, sy+\beta e, e \}=sy+\beta e$.
\end{enumerate}

The proof of Case III is that of Case II with $s$ and $r$ interchanged and $\beta$ and $-\beta$ interchanged. So we consider below Case I and Case II. 

\noindent {\bf Proof in Case I: $\textbf{min}\{ ry-\beta e, sy+\beta e, e \}=e$}  

Applying Lemma~\ref{lem:b'>b} to $\Gamma_b$, using the definition of $b$ to rule out Lemma~\ref{lem:b'>b}(4) (where $P(b',\alpha,\beta)$ would arise), we have three possibilities:
\begin{enumerate}
\item There is $b' \geq b+e \geq x+(\alpha+1)e$ and the total weight of $\Gamma_{b'}$ is
\[(ry - (\beta + 1)e) \wt{s} + (sy+ (\beta + 1)e) \wt{-r} + (e) \wt{(\beta+1)(s+r)}.\]
\item There is $b' \geq b+e \geq x+(\alpha+1)e$ and the total weight of $\Gamma_{b'}$ is
\[(sy+ (\beta - 1)e) \wt{-r} + (ry - (\beta - 1)e) \wt{s} + (e) \wt{(\beta-1)(s+r)}.\]
\item Define $\Delta := \text{min}\{ry-\beta e, sy + \beta e\}$.
\begin{enumerate}
\item $\Delta =ry-\beta e$. There is a 
$b' \geq b+ \Delta \geq x+\alpha e + \Delta$ and the total weight of $\Gamma_{b'}$ is
\[(sy+ \beta e - \Delta) \wt{-r} + (e+\Delta) \wt{\beta (s+r)} + \Delta \wt{-\beta (s+r)+s-r}.\]
\item $\Delta =sy+\beta e$. There is a 
$b' \geq b+ \Delta  \geq x+\alpha e +\Delta$ and the total weight of $\Gamma_{b'}$ is
\[(ry - \beta e - \Delta) \wt{s} + (e+\Delta) \wt{\beta (s+r)} + \Delta \wt{-\beta (s+r)+s-r}.\]
\end{enumerate}

\end{enumerate}

\noindent In each case above, we take $b'$ to be the largest value satisfying the given weight condition. %Note that by Theorem~\ref{thm:annulus}, {\bf we may assume that the coefficients in the total weight in each of the above cases is non-zero}.

% \begin{clm}\label{clm:mainI}
% In each of cases $(1)$--$(3)$ above, $M_{b'} < 2q$. 
% \end{clm}

% If $b' \geq q+y-1$, then Theorem~\ref{thm:q+y,m_b=1}(2) holds, so we assume that
% \[b' < q+ y-1 <2q-2.\]

\noindent {\bf Consider Case I(1).} 

If $(\alpha+1)e<2y$, then $P(b',\alpha+1, \beta +1)$ holds, contradicting the definition of $b$. On the other hand, if $(\alpha + 1) e \geq 2y$, then $b' \geq x+2y=q+y$ and Claim~\ref{clm:Mb'} implies  Theorem~\ref{thm:q+y,m_b=1}(2).

\noindent {\bf Consider Case I(2).} 

We replace $\beta$ with $-\beta$, switching the role of $r,\, s$, and apply the argument from Case I(1) to show that Theorem~\ref{thm:q+y,m_b=1}(2) must hold.

% $\spadesuit$ {\em (Here is an explicit argument. If $(\alpha+1)e<2y$ then $P(b',\alpha+1, \beta - 1)$ holds, contradicting the definition of $b$. 

% So we assume that $(\alpha + 1) e \geq 2y$ and hence $b' \geq x+2y=q+y$. On the other hand, we have that $\alpha e < 2y$ and $y>1$. But $M_{b'} \leq (|\beta|+1) (s+r) \leq (|\alpha|+2) (q-e)/y < (2y/e+2) (q-e)/y < 3q$.). Then $M_b < 2q$ unless $e=1$ which means $q=(r+s)y+1$. Then $2(y(s+r)+1)=2q=|\beta|(s+r)$ means that $2 \equiv 0 \,\, (\mod s+r)$ which means that $s=1=r$, contradicting that $q>4y$.)}

%This finishes the proof of Theorem~\ref{thm:q+y,m_b=1} in Case I(2).

\noindent {\bf Consider Case I(3).} 

By interchanging the role of $r$ with $s$ and $\beta$ with $-\beta$, the proof under assumption I$(3)(a)$ gives a proof under assumption I$(3)(b)$. So we assume I$(3)(a)$ holds.

Recall that $\alpha e \leq 2y-2$ ($P(b,\alpha,\beta)$) and $y>1$. Furthermore, Claim~\ref{clm:Mb'} shows that either $b'\ge q+y-1$ and Theorem~\ref{thm:q+y,m_b=1}(2) holds, or $x + \alpha e + \Delta \leq b' < q+ y-1 =x +2y-1$. Thus, we hereafter assume
\[ry +  (\alpha-\beta)e = \alpha e + \Delta < 2y-1.\]

First assume $r=1$. Because an isotopy class of $\Gamma_z$ loses at most one arc as we move to $\Gamma_{z+1}$, we have that $m_{b''} \leq 1$ for $b''=b'+ (sy + \beta e - \Delta -1) \geq x + \alpha e + \Delta +  (sy + \beta e - \Delta -1) = x +sy + (\alpha + \beta)e -1$. Thus, either Theorem~\ref{thm:q+y,m_b=1}(3) holds, or we have
\[sy+ (\alpha + \beta)e < 2y.\]

Recall that $\abs{\beta}\le \alpha+1$. First assume further that $\beta\ge -\alpha$. Then $sy \le sy+ (\alpha + \beta)e < 2y$ implies that $s=1$. As the coefficients of the total weight of $\Gamma_{b'}$ are positive,  $0<sy+ \beta e - \Delta=sy+\beta e - ry +\beta e = 2\beta e$. Thus $\beta>0$. Similarly, $0<\Delta=y-\beta e$ then implies that $e<y$. Then $q=ry+sy+e=2y+e<3y$, contradicting that $q>4y$.

Then we assume $\beta = -\alpha-1$. Thus, $y+e \le ry+(\alpha-\beta)e < 2y-1$, so $e < y$. Moreover, $(s-1)y < sy-e = sy+(\alpha+\beta)e < 2y$, so $s\le 2$. Then again $q=ry+sy+e \le 3y+e < 4y$, contradicting that $q>4y$.

Therefore, Theorem~\ref{thm:q+y,m_b=1} holds whenever $r=1$. In particular, if $\abs{\beta}\le \alpha$, $ry\le ry+(\alpha-\beta)e < 2y-1$ implies $r=1$. In Case I(3), we hereafter assume $r \geq 2$ along with $|\beta|=\alpha +1$.

Recall that
$ry +  (\alpha-\beta)e  < 2y-1$. Thus $\beta = \alpha +1>0$ and $(r-2) y < e-1$. Note that the latter inequality says that $e>1$. Recall also that in Case (I), we have $e \leq \Delta=ry-\beta e$. Thus $(\beta+1)e \leq ry$. Summing up, we have
\[\beta = \alpha +1, \quad (\beta+1)e \leq ry, \quad (r-2) y < e-1, \quad e>1.\]

As $\beta=\alpha+1 >0$, the above says $(r-2)y\leq e-2 \leq ry/(\beta+1)-2 \leq ry/2-2$. This implies $r \leq 3$ and $m_{b'} \leq 3$. 
Because an isotopy class of $\Gamma_z$ loses at most one arc as we move to $\Gamma_{z+1}$, we have that $m_{b''} \leq 3$ for $b''=b'+ (sy + \beta e -\Delta-1)= x+(\alpha+\beta) e + sy -1$. Assuming Theorem~\ref{thm:q+y,m_b=1}(1) does not hold, this implies that $2q-y-2 > b'' \geq x+(\alpha+\beta) e + sy -1=(q-y)+(\alpha+\beta-1)e+(q-ry)-1$. We conclude $ry>(\alpha+\beta-1)e+1$. We record then that 
\[r \leq 3, \quad ry > (\alpha+\beta-1)e+1, \quad (r-2) y \le e-2.\]

%$\spadesuit$ {\em (Tuesday) Here we could say that if $2q-y-2>b$ then $ry>1+(\alpha + \beta -1)e$. This will allow us to prove the Claim below. Thus we could conclude the stronger version of Theorem~\ref{thm:q+y,m_b=1}(1) with $b \geq 2q-y$ in Case I(3).)}

%$\spadesuit$  {\em (If we want $b \geq 2q-2y$, we need to modify the above inequality to $(r-1)y>(\alpha + \beta-1)e-1$, equivalently change the strict inequality to $(r-1)y \geq (\alpha+\beta-1)e$. Then the argument of the Claim below goes through except when $r=2,s=2,e=1$ so $q=4y+1$ and then if the Claim fails to hold we get that $\beta \geq y$ which implies that $\beta=y$. But then $(r-1)y \geq (\alpha+\beta-1)e$ becomes $y \geq y-1+y-1=2y-2$, implying that $y=2$. But then $q=9$, a contradiction.)}

We apply Lemma~\ref{lem:b'>b} to $\Gamma_{b'}$. By the definition of $b$ and $b'$, there are two possibilities.

\begin{enumerate}
\item Let $\Delta'=\text{min}\{sy+\beta e - \Delta, e+\Delta\}$. Then there is a $b'' \geq b' + \Delta' \geq x + \alpha e + \Delta + \Delta'$ such that the weights of $\Gamma_{b'}$ lie in $\{ -r, \beta (s+r), -\beta (s+r)+s-r, 2\beta (s+r)-s\}$.
\item Let $\Delta'=\text{min}\{sy+\beta e - \Delta, \Delta \}$. Then there is a $b'' \geq b' + \Delta' \geq x + \alpha e + \Delta + \Delta'$ such that the weights of $\Gamma_{b'}$ lie in $\{ -r, \beta (s+r), -\beta (s+r)+s-r, -2\beta (s+r)+s-2r\}$.
\end{enumerate}

\begin{clm} $M_{b''} < 2q$. \end{clm}
% Either \begin{enumerate}
% \item
%; or
%\item $M_{b'} < 4q, e=1,r=2$  and if $y \leq 8$ then none of the weights of $\Gamma_{b''}$ are $\pm 2q, \pm 3q$.
%\end{enumerate}

\noindent {\em Proof of Claim:}

From the list of weights of $\Gamma_{b''}$ we have
$ M_{b''} \leq 2 \beta (s+r)+2r-s =(2 \beta -1) (s+r) +3r$.

First assume $\beta >1$. Then $(r-2)y\leq e-2 \leq ry/(\beta+1)-2 \leq ry/3-2$ and we conclude $r \leq 2$. That is, $r=2$. Also $\beta >1$ implies $\alpha \geq 1$ and 
$2y=ry \geq (\alpha+\beta-1)e+2 \geq \beta e+2$. That is, $\pth{2y-2}/{e} \geq \beta \geq 2$. Note this implies $y-1\geq e$.  
Using the above,  $ M_{b''}  \leq (2 \beta -1) (s+r) +3r \leq  (2(2y-2)/e -1) (q-e)/y + 6 < (2y-3)(q/y) + 2q/y \leq (2-3/y+2/y)q <2q$. In the third inequality, we used the fact that $e \geq 2$ and that $q>4y$ implies $2q/y>8$. 

%((2y-4-e)/e)  ((q-e)/y) + 6 \leq (2/e-4/(ey)-1/y)q - ((2y-4-e)/y)+6= (2/e-4/(ey)-1/y)q +(4+e)/y)+4 \leq  (2/e-4/(ey)-1/y)q +(3+y)/y)+4 \leq $ using that $e \geq 2$.
%
%and we conclude $y+2=(r-1)y+2 \geq \beta e$. That is, $(y+2)/e \geq  \beta$. Then we have $ M_{b''}  \leq (2 \beta -1) (s+r) +3r \leq  (2(y-2)/e -1) (q-e)/y + 6= ((2y-4-e)/e)  ((q-e)/y) + 6 \leq ((2y-4-e)/e)(q/y) -((2y-4-e)/e) (e/y)+6=(2-4/(ey)-1/y)q-(2-(4+e)/y)+6$ using that $e \geq 2$. Thus  $M_{b''} \leq (2/e+2/(ey))q-2-2/y+6 < (2/e+2/(ey))q +4 < (2/e+2/(ey)+4/q)q<2q$ as $e\geq 2, y \geq 2, q>8$. 
%If $e >1$, then $M_{b''}<2q$. So we assume, $e=1$ and $M_{b''} < 2q +4-(q-1)/y$. As $q>4y$, $0 \geq 4- (q-1)/y$, and we once more conclude that 
%$M_{b''} <2q$.

It must that $\beta=1, \alpha=0$. We also know that $2(q - e) < 2q \le yq$. Then, $ M_{b''} \leq (s+r) +3r \leq (q-e)/y + 9 < q/2 + 9 \le 2q$, recalling that $r \leq 3$ and $q > 8$.
\qed{\em $\,$ (Claim)}

Assume $(1)$ or $(2)$ holds, so either $\Delta'=\text{min}\{sy+\beta e - \Delta, e+\Delta\}$ or $\Delta' = \text{min}\{sy+\beta e - \Delta, \Delta\}$. In both cases, $b'' \geq b' + \Delta' \geq x + \alpha e + \Delta + \Delta' \ge x + \text{min}\{(\alpha + \beta) e+sy, \alpha e + 2\Delta\}$.  If $e \geq y$, then $b'' \geq x+2y$ recalling that $\beta, s>0$, and that we are in Case I so $\Delta \geq e$ and Theorem~\ref{thm:q+y,m_b=1}(2) holds. So we assume $e <y$. Then $(r-2)y<e$ and $r \geq 2$ imply that $r=2$. And $4y < q=y(s+r)+e$ implies $s \geq 2$. Also, $\alpha e + 2\Delta \ge (\alpha+1)e + \Delta=ry$ as $\beta=\alpha+1$. Then $b'' \geq x + \text{min}\{sy, ry\}=x+2y$, so Theorem~\ref{thm:q+y,m_b=1}(2) holds.
% Above is the second place to use $q > 4y$ instead of $3y$. Again maybe it can be improved.

This finishes the proof of I(3), and hence of Case I.

%%%%%%%%%%%%%%%%%%%%

\noindent {\bf Proof in Case II: $\textbf{min}\{ ry-\beta e, sy+\beta e, e \}=ry-\beta e$}  

By Case I, we may take $e > ry-\beta e>0$ which means that $(\beta+1)e>ry>\beta e$. In particular, as $r>0$, we have $\beta \geq 0$ and $e>1$.
%and if $e=2$ then $2(\beta +1)>ry>2\beta$ and $ry=2\beta+1$ which implies that $y$ is odd. Recapitulating, in Case II:
%$$ e>1,  \text{and} \,\, e=2 \,\, \text{implies} \,\, ry=2\beta+1.$$

Let $\Delta= ry-\beta e$. Applying Lemma~\ref{lem:b'>b} to $P(b, \alpha, \beta)$, using the definition of $b$ to rule out Lemma~\ref{lem:b'>b}(4) where we have $P(b',\alpha,\beta)$, we have three possibilities:

\begin{enumerate}
\item There is $b' \geq b+\Delta \geq x+\alpha e+ \Delta=x+ry+(\alpha-\beta) e$ and the total weight of $\Gamma_{b'}$ is
\[ (sy+ \beta e - \Delta) \wt{-r} + (e + \Delta) [{\bf \beta (s+r)}] + \Delta \wt{(-\beta (s+r)+s-r)}.\]
\item There is $b' \geq b+\Delta \geq x+\alpha e + \Delta$ and the total weight of $\Gamma_{b'}$ is
\[ (e-\Delta) [{\bf \beta(s+r)}] + (sy+ry)  \wt{-r} + \Delta \wt{(\beta+1) (s+r)}.\]
%\item Let $\Delta' = \text{min}\{e, sy + \beta e\}$.
%\begin{enumerate}
%\item $\Delta' =sy + \beta e$. There is a 
%$b' \geq b+ \Delta' \geq x+\alpha e + \Delta'$ and the total weight of $\Gamma_{b'}$ is
%$$(e-\Delta') \,\, {\bf \beta (s+r)} +  (\Delta+\Delta') \,\, \wt{s} + \Delta' \wt{(\beta-1)(s+r)}.$$
\item Define $\Delta' := \text{min}\{e, sy + \beta e\}$.
\begin{enumerate}

\item $\Delta' =e$. There is a 
$b' \geq b+ \Delta'  \geq x+(\alpha+1) e$ and the total weight of $\Gamma_{b'}$ is
\[(sy+(\beta-1)e) \wt{-r} + (ry-(\beta-1)e) \wt{s} + (e) \wt{(\beta-1)(s+r)}.\]
\item $\Delta' =sy+\beta e$. Then $\beta=0$ and there is a 
$b' \geq b+ \Delta' \geq x+\alpha e + sy$. The total weight of $\Gamma_{b'}$ is
\[(e-sy) \wt{0} + (sy+ry) \wt{s} +  (sy) \wt{-(s+r)}.\]
\end{enumerate}

%\item $\spadesuit$ There is a 
%$b' \geq b+ e \geq x+(\alpha+1) e$ and the total weight of $\Gamma_{b'}$ is
%$$(sy + (\beta-1) e) \,\, \wt{-r} + (\Delta + e) \,\, \wt{s} + e\,\,\wt{(\beta-1)(s+r)}.$$

\end{enumerate}

\noindent In each case above, we take $b'$ to be the largest value satisfying the given weight condition (with $b, \alpha, \beta$ defined by $P(b,\alpha,\beta)$). 

% \begin{clm}\label{clm:Mb'}
% In each of cases $(1)$--$(3)$ above, $M_{b'} < 2q$. 
% \end{clm}

% \noindent {\em Proof of Claim:}

% By $P(b,\alpha, \beta)$, $\alpha e < 2y - 1$, so $2y \ge \alpha e + 2 \ge \alpha + 2$. Therefore, $M_{b'} \leq (|\beta|+1) (s+r) \leq (\alpha+2) (q-e)/y \le 2(q - e) < 2q$.  Here again we are using that $ry+sy+e=q$. %Significantly shortened after changing $\alpha 2 < 2y$ to $2y - 1$.
% \qed {\em $\,$ (Claim)}
% \medskip

% If $b' \geq q+y-1$, then Theorem~\ref{thm:q+y,m_b=1}(2) holds, so we conclude that 
% \[b' < q+ y-1 <2q-2.\]
%Note that by Theorem~\ref{thm:annulus}, {\bf we may assume that the coefficients in the total weight in each of the above cases is non-zero}.

\noindent {\bf Consider Case II(3).} 

Assume $II(3)(a)$ holds. If $b'\ge q+y-1$, then Claim~\ref{clm:Mb'} implies that Theorem~\ref{thm:q+y,m_b=1}(2) holds. Thus, we assume $x+(\alpha+1)e\le b'< q+y-1$, so $(\alpha+1) e < 2y-1$. Furthermore, $|\beta-1| \le \alpha+2$. Thus, $P(b', \alpha+1, \beta-1)$ holds. As $b'>b$ this contradicts the definition of $b$.

Assume $II(3)(b)$ holds. By Claim~\ref{clm:Mb'}, either $b'\ge q+y-1$ and Theorem~\ref{thm:q+y,m_b=1}(2) holds, or $x+\alpha e + sy\le b' < q+y-1$, so $s=1$. Because an isotopy class of $\Gamma_z$ loses at most one arc as we move to $\Gamma_{z+1}$, $m_{b'+sy+ry-1}\le 1$. Since $b'+sy+ry-1 \ge b'+2y-1 > x+2y-1=q+y-1$, Theorem~\ref{thm:q+y,m_b=1}(3) holds.

This completes the proof of Theorem~\ref{thm:q+y,m_b=1} for Case II(3). 

\bigskip

We now make some observations common to both Cases II(1) and II(2).

\medskip

\begin{clm}~\label{clm:II(1)(2)}
In both Cases II(1) and II(2), either Theorem~\ref{thm:q+y,m_b=1} holds, or we have
\[e>\Delta>0, \quad \beta=\alpha+1>0, \quad 1+(r-2)y<e, \quad \Delta+\alpha e < 2y-1,\quad r \ge 2, \quad e>2.\]
\end{clm}

\noindent {\em Proof of Claim:}
Recall that for Case II generally we assume $e>\Delta>0$. Assume we are in Cases II(1), II(2). Recall that $q+y-1>b' \ge x+\alpha e+\Delta$. We conclude that $\alpha e + \Delta = ry + (\alpha-\beta)e < 2y-1$. 

First assume $r = 1$ in Case II(1). Because an isotopy class of $\Gamma_z$ loses at most one arc as we move to $\Gamma_{z+1}$, we have that $m_{b''} \leq 1$ for $b''=b'+ (sy + \beta e - \Delta -1) \geq x + \alpha e + \Delta +  (sy + \beta e - \Delta -1) = x +sy + (\alpha + \beta)e -1$. 
%As $(\beta + \alpha) \geq (\beta+1)e>ry$, 
If $b'' \ge x+2y-1$, then Theorem~\ref{thm:q+y,m_b=1}$(3)$ holds. Thus, we assume that $sy + (\alpha+\beta)e < 2y$. Thus, $s=1$ and $(\alpha+\beta)e < y$. Then $4y<q=y(s+r)+e$ implies that $e > 2y$. %If $e<2y$ then $q=y(s+r)+e<4y$, contradicting our assumption that $q>4y$. Thus assume $e \geq 2y$. 
Then $\alpha=\beta=0$, and $m_j=0$ for $j=b'+e+\Delta-1 \geq x+ e+2\Delta -1 >x+2y-1=q+y-1$, verifying Theorem~\ref{thm:q+y,m_b=1}$(3)$.
Thus, we assume $r \ge 2$ in Case II(1). 

Assume $r=1$ in Case II(2). Similarly, we have that $m_{b''} \le 1$ for $b'' = b' + sy+ ry -1 \ge x + 2y = q +y$, verifying Theorem~\ref{thm:q+y,m_b=1}$(3)$. Thus in Case II(2) we assume that $r \ge 2$. 

In either case $r \ge 2$, and $ry + (\alpha-\beta)e < 2y-1$ implies that $\beta=\alpha+1$ and $1+(r-2)y<e$. Recall that $e>1$. If $e = 2$ then $1+(r-2)y<e = 2$ implies $r=2$ ($y>1, r \geq 2$). As $2 = e>\Delta>0$, we conclude that $1=\Delta=2y-2 \beta $, a contradiction. We conclude $e>2$.

\qed {\em $\,$ (Claim)}

\medskip

The following estimate will be useful in both cases.

\begin{clm}\label{clm:2beta+1}
$(2 \beta+1) (s+r)+r<2q$
\end{clm}

\noindent {\em Proof of Claim:}
% First assume $\beta=1$. As $2e=(\beta+1)e>ry$, $2q=2(y(r+s)+e) > 4(r+s)+ry >3(r+s)+r$, verifying the Claim.

First assume $\beta\le 2$. Then $\alpha >0$ and $1+ \alpha e + \Delta < 2y$ along with $e \geq 3$ implies $y \geq 3$. Thus, $(2 \beta +1) (s+r)+r \le 5s + 6r < 6(s+r) < 2y(s+r)+2e=2q$ as claimed. 

We hereafter assume $\beta >2$ and hence that $\alpha \geq 2$. Then $1+ \alpha e + \Delta < 2y$ along with $e \geq 3$ implies that $y \geq 5$. As $1+ \alpha e + \Delta < 2y$, we have $2y/e>\alpha$. Thus,
% \[
% \begin{aligned}
%     (2 \beta +1) (s+r)+r \leq (2 \alpha+3) (s+r) +r 
%     &< \pth*{\frac{4y}{e} + 3} \frac{q-e}{y} + r\\ 
%     &= \pth*{\frac{4}{e}+\frac{3}{y}} \pth*{q-e} +r\\
%     &< \pth*{\frac{4}{e}+\frac{3}{y}} q - \pth*{\frac{4}{e}+\frac{3}{y}} e + \pth*{2+\frac{e}{y}-\frac{1}{y}}\\ 
%     &= \pth*{\frac{4}{e}+\frac{3}{y}} q - \pth*{2+\frac{2e}{y} + \frac{1}{y}},
% \end{aligned}
% \]
\[
\begin{aligned}
    2q
    &> (3/y+4/e) q
    = (3/y+4/e)((r+s)y+e)
    = 3(r+s)+2(r+s)2y/e + 3e/y +4\\
    &> 3(r+s)+2\alpha (r+s) + r
    = (2\alpha+3)(r+s)+r\\
    &\ge (2 \beta +1) (s+r)+r,
\end{aligned}
\]
where above we use that $1+(r-2)y<e$. 
\qed {\em $\,$ (Claim)}

\bigskip

\noindent {\bf Consider Case II(1).}

Applying Lemma~\ref{lem:b'>b} to $\Gamma_{b'}$ and using the definitions of $b,\, b'$ as the largest values satisfying their defining conditions, we are left with two possibilities:

\begin{enumerate}
\item Let $\Delta''=\min\{sy+ \beta e - \Delta, e+\Delta\}$. There is a $b'' \geq b' + \Delta'' \geq x + \alpha e + \Delta + \Delta''$ such that the total weights of $\Gamma_{b''}$ lie in $\{-r, \beta (s+r), -\beta (s+r) + s -r, 2\beta (s+r) -s\}$. 

\item Let $\Delta''=\min\{sy+ \beta e - \Delta, \Delta\}$. There is a $b'' \geq b' + \Delta'' \geq x + \alpha e + \Delta + \Delta''$ such that the total weights of $\Gamma_{b''}$ lie in $\{-r, \beta (s+r), -\beta (s+r)+s-r, -(2\beta+1) (s+r)-r+2s\}$. 
\end{enumerate}

\begin{clm}$M_{b''} < 2q$.
\end{clm}

\noindent {\em Proof of Claim:}
$M_{b''} \leq (2\beta+1)(s+r)+r-2s<2q$ by Claim~\ref{clm:2beta+1}.
\qed {\em $\,$ (Claim)}

Assume $(1)$ holds: $\Delta''=\min\{sy+\beta e - \Delta, e+\Delta\}$ and $b'' \geq b' + \Delta'' \geq x + \alpha e + \Delta + \Delta'' = x + \min\{sy+(\alpha + \beta) e, (\alpha +1) e + 2 \Delta\} = x + \min\{sy+(\alpha + \beta) e, \beta e + \Delta +\Delta\} = x+\min\{sy + (\alpha + \beta) e, ry + \Delta\}$.  Assume $e \geq y$. Then $b'' \geq x+2y$ recalling Claim~\ref{clm:II(1)(2)}. This implies that Theorem~\ref{thm:q+y,m_b=1}(2) holds. So we assume $e <y$. Then $(r-2)y<e$ and $r \geq 2$ imply that $r=2$. And $4y < q=y(s+r)+e$ implies $s \geq 2$. Then $b'' \geq x + \text{min}\{sy, ry\}=x+2y$, again implying that Theorem~\ref{thm:q+y,m_b=1}(2) holds.

%$\spadesuit \spadesuit$ {\em (In what follows is where we need to argue for  Theorem~\ref{thm:q+y,m_b=1}$(2)$ that there is a $b \geq q+y-1$ rather than $b \geq q+y$.)}

Finally, consider $(2)$: $\Delta''=\min\{sy+\beta e - \Delta, \Delta\}$. If $\Delta''=sy+\beta e - \Delta$, the arguments in the paragraph above show that the Theorem holds. We therefore assume that $\Delta''=\Delta$ and (by the preceding Claim) that $\alpha e + 2\Delta <2y-1$. The latter can be written $2y-1>
(\beta e + \Delta) +\Delta -e = ry + \Delta - e$ from which we conclude $e-\Delta-1>(r-2)y$. Note that this implies that $e>2$ as $r \geq 2$. In sum, for what follows, we have that $b'' \geq x+ \alpha e + 2 \Delta$ and the weighting of $\Gamma_{b''}$ is
\[\tag{$\ast$}(sy+ \beta e - 2\Delta) \wt{-r} + (e + 2\Delta) \wt{\beta (s+r)} + \Delta \wt{-2\beta (s+r)+s-2r}.\]
And furthermore, 
\[\Delta < e, \quad \beta=\alpha+1, \quad (r-2)y<e-\Delta-1, \quad \alpha e + 2\Delta <2y-1, \quad r \ge 2, \quad e>2.\]

Assume $\alpha>0$. Then $e+2\Delta < 2y-1$ implies that $(r-2)y < e-\Delta-1 < 2y-3\Delta-2$. This means $r \leq 3$.
Because an isotopy class of $\Gamma_z$ loses at most one arc as we move to $\Gamma_{z+1}$, we have, from the total weight of $\Gamma_{b'}$, that $m_{j} \leq 3$ for $j=b''+ sy + \beta e - 2\Delta -1 \geq x+ \alpha e +sy+\beta e-1=x +(\alpha-1)e + (e+sy+ry)-\Delta-1=2q-y+(\alpha-1)e-\Delta-1$ using that $\beta e=ry-\Delta$. If $\Delta \le (\alpha-1)e+1$, then Theorem~\ref{thm:q+y,m_b=1}$(1)$ holds. Thus, we assume $\Delta >(\alpha-1)e+1$. As $e>\Delta$, we conclude that $\alpha=1, \Delta \geq 2$. 

We assume hereafter in Case II(1), $(\ast)$ holds and either
\begin{itemize}
\item $\alpha=1, \beta=2, r \leq 3, \Delta \geq 2$, or
\item $\alpha=0, \beta=1$.
\end{itemize}

%As $q=y(s+r)+e, x=q-y, \Delta<2y/3$, we continue,  $b''' \geq x+ q-\Delta-1>q-y +q-(2/3)y-1 > 2q-2y-1$, contradicting that Theorem~\ref{thm:q+y,m_b=1}$(1)$ does not hold.

%$\spadesuit$ {\em (Note that here we get $b'' \geq 2q-2y$ automatically. Furthermore we get $b'' \geq 2q-2y+2$ as long as $-2y/3-1 \geq -y+1$ which means $y/3 \geq 2$. Need to try to elevate this to Superversion)}

\begin{defin}\label{defin: Q1}
Let $s,r,e, \Delta, \alpha, \beta, b$ be as above.
$Q_1(i, \omega)$ is the statement: $i \geq b' + (\omega-1) \Delta \geq b + \omega \Delta  \geq x + \alpha e+ \omega \Delta$ and the total weight of $\Gamma_i$ is 
\[(sy+ \beta e - \omega \Delta)  \wt{-r}  + (e+\omega \Delta) \wt{\beta(s+r)} +\Delta \wt{-\beta \omega (s+r) + s -\omega r}\]
with $\omega\geq 2, \alpha e+\omega \Delta< 2y-1$.
\end{defin}
Note that if $i, \omega$ are such that $Q_1(i, \omega)$ holds, then $2y - 1 > \alpha e + \omega\Delta \ge \alpha e + 2(ry - \beta e) = 2ry - (\beta + 1)e \ge 4y - (\beta + 1)e$, so either $\beta = 1$ and $y < e$ or $\beta = 2$, $\alpha  = 1$ and $y < 3e/2$, and hence $sy + \beta e - \omega\Delta > sy + \beta e + (\alpha e- 2y + 1) = (s - 1)y + ((2\alpha + 1)e - y ) + 1 > 0$. Namely, the coefficient of $\wt{-r}$ is always positive (as are clearly the coefficients $\Delta, e+\omega \Delta$). 

Also note that $Q_1(i, \omega)$ implies that $2y-1>\alpha e + \omega \Delta=(\alpha +1) e + \Delta +(\omega-1) \Delta -e=ry + (\omega-1) \Delta -e$, implying that 
\[e>1+(r-2)y+(\omega-1) \Delta.\]

Now, $Q_1(b'', 2)$ holds.
Let $i$ be the largest value such that $Q_1(i, \omega)$ holds for some $\omega$.
Thus $i \geq b''$. Applying Lemma~\ref{lem:b'>b} to $\Gamma_{i}$ and using the definitions of $b',\, i$ as the largest values satisfying their defining conditions, we are left with two possibilities:

\begin{enumerate}[label=(\alph*)] 
\item Let $\Delta'''=\text{min}\{sy+ \beta e - \omega \Delta, e+ \omega \Delta\}$. There is an $i' \geq i + \Delta''' \geq x + \alpha e+ \omega \Delta + \Delta'''$ such that the total weights of $\Gamma_{i'}$ lie in $\{-r, \beta (s+r), - \beta \omega (s+r) + s - \omega r, \beta (\omega +1)  (s+r) -s + (\omega -1) r\}$. 

\item Let $\Delta'''=\text{min}\{sy+ \beta e - \omega \Delta, \Delta\}$. There is a $i' \geq i + \Delta''' \geq x + \alpha e +\omega \Delta + \Delta'''$ such that the total weights of $\Gamma_{i'}$ lie in $\{-r, \beta (s+r), - \beta \omega (s+r) + s - \omega r,  - \beta (\omega +1) (s+r) +s - (\omega +1) r\}$. Furthermore, if $\Delta'''=\Delta$ then $i' \geq x+ \alpha e+(\omega+1) \Delta \geq x+2y-1=q+y-1$ (else $Q_1(i',\omega+1)$ holds, contradicting the definition of $i$). 

\end{enumerate}

\begin{clm} \label{clm: II M < 2q second}
$M_{i'} <  2q$.
\end{clm}

\noindent {\em Proof of Claim:}

{\bf We first assume $\alpha=0, \beta=1$.} 
Then $M_{i'} \leq \beta (\omega+1) (r+s)-s+(\omega+1)r=\omega (2r+s)+ 2r$. Now $\Delta=ry-e$, thus $q=y(r+s)+e=y(s+r)+ry-\Delta$, giving $2r+s=(q+\Delta)/y$. 
% Furthermore, $2+2r \leq y(s+r)<q$.
Then, since $\omega\Delta < 2y - 1$ and $y\ge 2$,
\[M_{i'} \leq \omega (s+2r)+ 2r \leq \pth*{\frac{2y-2}{\Delta}}\pth*{\frac{q+\Delta}{y}}+2r < \frac{2q}{\Delta}\pth*{1-\frac{1}{y}} + 2 + 2r\]
Thus, as $\Delta\ge 1$,
\[M_{i'} < 2q - \frac{2q}{y} + 2 + 2r = 2q - \frac{2}{y}\pth*{q - (r + 1)y} < 2q\]
noting that $y(r+1) \le y(r+s)=q-e$.

Thus we may assume $\alpha=1,
\,\beta=2$. Recall that then $2\leq r \leq 3$ and $\Delta \geq 2$.

From $Q_1(i, \omega)$ we have that $e+\omega \Delta \leq 2y-2$. We rewrite this as $\omega \leq \frac{2y-2-e}{\Delta}$. 
Thus $M_{i'} \leq 2 (\omega+1) (s+r)-s+(\omega+1)r=(2 \omega+1) (s+r)+(\omega+2)r \leq (\frac{4y-4-2e}{\Delta}+1) \frac{q-e}{y}+(\frac{2y-2-e}{\Delta}+2)r$.

%= \frac{4}{\Delta}-\frac{6+e}{y \Delta} (q-e)+(\frac{2y-4+\Delta}{\Delta})r$. 

First assume $r=3$. As $\omega \geq 2$, the inequality $e>1+(r-2)y+(\omega-1) \Delta$ shows that $e>y+\Delta+1$. Then $M_{i'} <( \frac{2y-6}{\Delta}) (\frac{q}{y})+ 3(\frac{y-3+\Delta}{\Delta}) < (1-3/y) q+(\frac{3}{2})y+3 $ remembering that $\Delta \geq 2$. As $q>4y$, note that $q/2>2y$ and $-3q/y<-12$. Thus we get $M_{i'} <q+q/2<2q$.

Now assume $r=2$. Then $2 \leq \Delta=2y-2e$. Again we have $M_{i'} \leq (\frac{4y-4-2e}{\Delta}+1) (\frac{q-e}{y})+(\frac{2y-2-e}{\Delta}+2)r \leq (\frac{2y-4}{\Delta}+2) (\frac{q-e}{y}) +2 (\frac{e-2}{\Delta}+3)$. Assume $\Delta>2$. Then $2y-2e=\Delta \geq 4$ and $M_{i'}  \leq (\frac{2y-4}{4}+2) (\frac{q-e}{y}) + 2(\frac{e-2}{4}+3) =(\frac{1}{2} + \frac{1}{y}) (q-e) + \frac{e}{2}+5 \leq q + 5-\frac{e}{2} <2q$, remembering that $e>2, q>4$. Thus we may take $2 = \Delta=2y-2e$ and we get
$M_{i'} \le (\frac{2y-4}{\Delta}+2) (\frac{q-e}{y}) + 2(\frac{e-2}{\Delta}+3) =(\frac{2y-4}{2}+2) (\frac{q-e}{y}) + 2(\frac{e - 2}{2}+3) = (q-e) + (e + 4)<2q$.
%We also have that $ (r-2)y<e-\Delta-1$. As $r \geq 2$, this implies that $e \geq \Delta +2$. Thus we have $2+(\omega +1) \Delta=\Delta +2 + \omega \Delta \leq e + \omega \Delta \leq 2y-2$. Thus $\omega+1 \leq \frac{2y-4}{\Delta}$.
%
%
%If $\Delta \geq 3$ then 
%$M_{i'} < \frac{4}{3} q+ 2y-4+\Delta= \frac{4}{3} q+ 2y-4+ry-2e<  \frac{4}{3} q+ 5y< \frac{4}{3} q+ $
\qed {\em $\,$ (Claim)}
\medskip

In $(b)$, if $\Delta'''=\Delta$, then $i'\ge q+y-1$ as stated there, verifying Theorem~\ref{thm:q+y,m_b=1}$(2)$. Therefore, whether in case $(a)$ or $(b)$ we assume $i' \geq i + \Delta''' \geq x + \alpha e + \omega \Delta + \Delta'''$ where $\Delta'''=\min\{sy+ \beta e - \omega \Delta, e+ \omega \Delta\}$. If $\Delta'''= sy+ \beta e - \omega \Delta$, then $i' \geq x+ sy + (\alpha + \beta) e$. Thus $s=1$, otherwise Theorem~\ref{thm:q+y,m_b=1}$(2)$ holds. If $e \geq y$ then $i' \geq x+2y=q+y$, which also verifies Theorem~\ref{thm:q+y,m_b=1}$(2)$. Thus we assume $e < y$. Then $1+(r-2)y+(\omega-1) \Delta <e$ implies that $r=2$. But $4y<q=y(r+s)+e < 3y+y$, a contradiction. Thus, it must be that $\Delta'''=e+ \omega \Delta$. Then $i' \geq x + (\alpha +1)e + 2 \omega \Delta =x+ \beta e + \Delta + (2 \omega-1) \Delta >x+ry \geq x+2y=q+y$, again implying that Theorem~\ref{thm:q+y,m_b=1}$(2)$ holds. 

% In $(a)$, $i' \geq i + \Delta''' \geq x + \alpha e + \omega \Delta + \Delta'''$ where $\Delta'''=\text{min}\{sy+ \beta e - \omega \Delta, e+ \omega \Delta\}$. If $\Delta'''= sy+ \beta e - \omega \Delta$, then $i' \geq x+ sy + (\alpha + \beta) e $. If $e \geq y$ then $i' \geq x+2y=q+y$, contradicting that Theorem~\ref{thm:q+y,m_b=1}$(2)$ does not hold. Thus $e < y$ and, again as Theorem~\ref{thm:q+y,m_b=1}$(2)$ does not hold, $s=1$. Then $1+(r-2)y+(\omega-1) \Delta <e$ implies that $r=2$. But $4y<q=y(r+s)+e < 3y+y$, a contradiction. Thus it must be that $\Delta'''=e+ \omega \Delta$. Then $i' \geq x + (\alpha +1)e + 2 \omega \Delta =x+ \beta e + \Delta + (2 \omega-1) \Delta >x+ry \geq x+2y=q+y$, again contradicting that Theorem~\ref{thm:q+y,m_b=1}$(2)$ does not hold. Thus $(a)$ does not occur. 

% Assume $(b)$ occurs. Then $i' \geq x + \alpha e + \omega \Delta + \Delta'''$ where $\Delta'''=\text{min}\{sy+ \beta e - \omega \Delta, \Delta\}$. As stated in $(b)$, if $\Delta'''=\Delta$, then $i' \geq q+y-1$, contradicting that Theorem~\ref{thm:q+y,m_b=1}$(2)$ does not hold. On the other hand, as argued in case $(a)$,  $x+ sy + (\alpha + \beta) e \geq q+y$, also contradicting that Theorem~\ref{thm:q+y,m_b=1}$(2)$ does not hold.

% These contradictions finish the proof of Theorem~\ref{thm:q+y,m_b=1} in Case II(1).

\bigskip

\noindent {\bf Consider Case II(2).} For the convenience of the reader, we first recall the context.
%%%%%%%%%%%%%%
\begin{defin}
$P(b,\alpha,\beta)$ is the statement: There is a $b \geq x+\alpha e$ such that the total weight of $\Gamma_b$ is
\[(ry - \beta e) \wt{s} + (sy+ \beta e) \wt{-r} + (e)\wt{\beta (s+r)}\]
with $\alpha \ge 0, \,  |\beta| \leq \alpha+1,\,   \alpha e < 2y - 1$.
\end{defin}

For Case II, we take $\Delta= ry-\beta e < e$. That is, $(\beta +1)e>ry>\beta e$ and $\beta \geq 0, e>1$. %In fact, $\beta>0$ by Lemma~\ref{lem:samecomplexity} and the fact that $b >x$. 
For Case II(2), 
there is $b' \geq b+\Delta \geq x+\alpha e + \Delta$ and the total weight of $\Gamma_{b'}$ is
\[(e-\Delta) \wt{\beta(s+r)} + (sy+ \beta e + \Delta) \wt{-r} + \Delta \wt{(\beta+1) (s+r)}.\]
%%%%%%%%%%%%%%%%

We assume that $b'$ is taken to be the largest integer so that it still satisfies the above statement (with the given $\alpha, \beta$). 

We assume the conditions of Claim~\ref{clm:II(1)(2)} hold.

%In fact, we hereafter assume $\beta >0$. Otherwise, $\beta=0$, and we see that $m_{b'} =0$. If $r >1$, then Theorem~\ref{thm:q+y,m_b=1}$(3)$ holds. If $r=1$ then $m_{x + \alpha e+ ry +(r+s)y-1 \le 1}$, and again Theorem~\ref{thm:q+y,m_b=1}$(3)$ holds.

% Then, $M_{b'}= (\beta + 1) (s+r) \leq (\alpha + 2) (q-e)/y < (2y/e + 2) (q-e)/y < 2q$ since $\alpha e < 2y, y>1, e>\Delta \geq 1 $. As Theorem~\ref{thm:q+y,m_b=1}(2) does not hold, $ry + (\alpha-\beta) e = \alpha e + \Delta < 2y-1$. 
% Then $ry = \Delta + \beta e \leq \Delta + \alpha e + e < 2y-1+e$, and we have $1+(r-2)y<e$.  
% That is
% \[1+ \alpha e + \Delta < 2y, \,\, 1+(r-2)y<e.\]

% First assume $r =1$. Then as we lose at most one arc in each isotopy class in moving from one level of $\Gamma$ to the next, we have that $m_{b''} =1$ for $b'' = b' + sy+ \beta e + \Delta -1 \geq x + \alpha e + \beta e + 2 \Delta + sy -1=x + \alpha e + \Delta +ry +sy-1 \geq x +2y = q +y$. This contradicts that Theorem~\ref{thm:q+y,m_b=1}$(3)$ does not hold. 

% {\bf We hereafter assume $r \geq 2$.} 

% Recall that in Case II we have $e>1$. If $e=2$ then $1+(r-2)y<e=2$ implies $r=2$ as we assume $y>1, r \geq 2$. As $2=e>\Delta>0$ we conclude that $\Delta=1$. But then $1=\Delta=2y-2 \beta $, a contradiction. Thus in considering Case II(2), we also have
% \[e>2.\]

% Assume $\beta=0$. Then $(\beta+1)e>ry$ becomes $e>ry$, and $\Delta=ry$. 
% But $2y>1+\alpha e + \Delta>ry$, contradicting that $r \geq 2$. Thus we have
% \[\beta >0.\]

\bigskip

Applying Lemma~\ref{lem:b'>b} to $\Gamma_{b'}$ and using the definitions of $b, b'$ as the largest values satisfying their defining conditions, we are left with two possibilities:

\begin{enumerate}[label=(\alph*)] 
\item There is a $b'' \geq b' + (e-\Delta)  \geq x + (\alpha+1) e$ such that the total weight of $\Gamma_{b''}$ is 
\[(sy + (\beta-1)e + 2 \Delta) \wt{-r} + (e) \wt{(\beta+1) (s+r)} + (e-\Delta) \wt{-(s+2r)}.\] 

\item Let $\Delta''=\text{min}\{e - \Delta, \Delta\}$. There is a $b'' \geq b' + \Delta'' \geq x + \alpha e + \Delta + \Delta''$ such that the total weight of $\Gamma_{b''}$
is % lie in $\{-r, \beta (s+r), (\beta+1) (s+r), (2 \beta +1) (s+r)+r \}$. 
\begin{enumerate}[label=(\roman*)]
\item $(2\Delta-e) \wt{(\beta+1) (s+r)} + (sy+ (\beta+1) e) \wt{-r} + (e-\Delta) \wt{(2 \beta +1) (s+r)+r}$ when $\Delta''=e - \Delta$; %(hence $e < 2 \Delta$ by Theorem~\ref{thm:annulus}); 
or
\item   $(e-2 \Delta) \wt{\beta(s+r)} + (sy+ \beta e + 2\Delta)  \wt{-r}  + \Delta \wt{(2 \beta +1) (s+r)+r}$ when $\Delta''= \Delta$. %(hence $e > 2 \Delta$ by Theorem~\ref{thm:annulus}).
\end{enumerate} 
\end{enumerate}

\medskip

\begin{clm}\label{clm:II(2)}
$M_{b''}<2q$.
\end{clm}

\noindent {\em Proof of Claim:}
$M_{b''} \le (2 \beta +1) (s+r)+r < 2q$ by Claim~\ref{clm:2beta+1}.
\qed {\em $\,$ (Claim)}
\smallskip

If $b''\ge x+2y-1$, then Claim~\ref{clm:II(2)} implies Theorem~\ref{thm:q+y,m_b=1}(2) holds. Therefore, we assume 
\[b''< x + 2y-1.\]
Theorem~\ref{thm:annulus} implies that the three coefficients of the total weight of $\Gamma_{b''}$ are positive. 

\bigskip

\noindent{\bf Consider Case II$(2)(a)$ above.}  

As $b''< x + 2y-1$, in Case II$(2)(a)$ we conclude that $(\alpha+1)e <2y-1$.

Along with Claim~\ref{clm:II(1)(2)}, in Case II(a) we have the following:
\[r \geq 2, \,\, 1+ \alpha e + \Delta < 2y, \,\, 1+(r-2)y<e, \,\, e>2, \,\, (\alpha+1)e<2y-1, \,\, \beta=\alpha+1 >0.\]

Now $1+(r-2)y<e<2y-1$ implies that $r=2,3$. Then we have that $m_i \leq 3$ for $i = b'' + sy + (\beta-1)e + 2 \Delta-1 \geq  x + \alpha e + sy+ (\beta e + \Delta) +\Delta -1=x + (\alpha-1)e + (e + sy+ry)+ \Delta-1=2q-y+(\alpha-1)e+\Delta-1$. If $(\alpha - 1)e+\Delta \ge -1$, then $m_i\le 3$ for $i\ge 2q-y-2$, and Theorem~\ref{thm:q+y,m_b=1}$(1)$ holds. Thus we assume $(\alpha - 1)e + \Delta < -1$. This implies 
\[\alpha=0, \quad \beta=1, \quad e-\Delta>1, \quad r=2,3.\]

%$(\alpha-1)e+\Delta-1 < -y+2$. 
%As $y \geq 2$, we conclude that 
%$$\alpha=0, \,\, \beta=1, \,\, y+\Delta \leq e+2$$
%
%$\spadesuit$ {\em (Going for the stronger version of Theorem~\ref{thm:q+y,m_b=1}$(1)$, we have $2q-y-2 > i \geq 2q-y+(\alpha-1)e +\Delta-1$. This implies $\alpha=0, \beta=1, e-\Delta>1$. Try to use this below in place of $y+\Delta \leq e+2$ to prove the Claim.)}

Then,  $b'' \geq x + e$ such that the total weight of $\Gamma_{b''}$ is 
\[(sy + 2 \Delta) \wt{-r} + (e) \wt{2(s+r)} + (e-\Delta) \wt{-(s+2r)}.\]

\begin{defin}\label{defin: Q2}
Let $s,r,e, \Delta$ be as above.
$Q_2(i, \omega)$ is the statement: $i \geq b'+ \omega (e-\Delta) =x + \Delta +\omega (e-\Delta)$ and the total weight of $\Gamma_i$ is 
\[(sy+ e + \Delta - \omega (e-\Delta)) \wt{-r}  + (\Delta +\omega (e-\Delta))\wt{2 (s+r)}+(e-\Delta) \wt{-\omega (2(s+r)+r) +(s+r)}\]
with $e-\Delta > 0, \omega\ge 1, \Delta + \omega (e-\Delta)<2y-1$.
\end{defin}

Then $Q_2(b'',1)$ holds. Let $i$ be largest so that $Q_2(i, \omega)$ holds for some $\omega$. Then $i \geq b''>b'$. By $Q_2(i, \omega)$,  $2y>1+ \Delta + \omega (e-\Delta)=1+ e + \Delta + \omega (e-\Delta)-e =1+ry  + \omega (e-\Delta)-e$ giving that 
\[e>1+(r-2)y+ \omega (e-\Delta).\]

In particular, if follows that $sy+e+\Delta -\omega (e-\Delta) > sy+\Delta +1>0$. That is, all the coefficients of $\Gamma_{b''}$ are positive. Applying Lemma~\ref{lem:b'>b} to $\Gamma_{i}$ and using the definitions of $i,b'$ as the largest values satisfying their defining conditions, we are left with two possibilities:

\begin{enumerate}[label=(\Alph*)] 
\item Let $\Delta'''=\min\bc{sy+ e + \Delta - \omega (e-\Delta), \Delta +\omega (e-\Delta)}$. There is an $i' \geq i + \Delta'''  \geq x + \Delta +\omega (e-\Delta)+ \Delta'''$ such that the weights of $\Gamma_{i'}$ are in 
$\bc{-r, 2 (s+r), -\omega (2(s+r)+r) +(s+r), \omega[2(s+r)+r]+s}$.

\item Let $\Delta'''=\min\bc{sy+ e + \Delta - \omega (e-\Delta), e - \Delta}$. There is an $i' \geq i + \Delta''' \geq x + \Delta + \omega (e-\Delta) + \Delta'''$ such that the weights of $\Gamma_{i'}$ are in
$\bc{-r, 2 (s+r), -\omega (2(s+r)+r) +(s+r), -\omega[2(s+r)+r]-(s+2r)}$. Furthermore, if $\Delta'''=e-\Delta$, then 
$\Delta + (\omega+1) (e-\Delta)\ge 2y-1$ (else $Q_2(i',\omega+1)$ holds, contradicting the definition of $i$).
\end{enumerate}

%Recall that $y+\Delta \leq e+2$, which we rewrite $y-2 \leq e-\Delta$. 
% In sum, in the current context we have
% $$2y>1+ \Delta + \omega (e-\Delta),  e>\Delta + 1, e \geq 3, \omega \geq 1, 1+e<2y, 1+(r-2)y+ \omega (e-\Delta) < e,  r=2,3$$ 

%Then $2(e-\Delta) \geq 2y-4> 1+ \Delta + \omega (e-\Delta)-4$ and we have $3>\Delta + (\omega-2)(e-\Delta)$. Thus $\omega \leq 2$. 

\begin{clm}\label{clm:II(2)(a)}
$M_{i'} <  2q$
\end{clm}

\noindent {\em Proof of Claim:}

%From the list of weights of $\Gamma_{i'}$ we have
%$ M_{i'} \leq \omega[2(s+r)+r]+(s+2r) \leq 5(s+r)+3r \leq 5(q-e)/y+9<(5/y+1/3)q
%$ since $q \geq 27$. On the other hand, $2y>e+1 \geq 4$, hence $y \geq 3$. Thus
%$ M_{i'} <2q$ as claimed.

By $Q_2(i,\omega)$, $2y>1+ \Delta + \omega (e-\Delta) \geq 1+\Delta+2 \omega$ (recall $e-\Delta>1$). Thus $2y-\Delta \ge 2 \omega+2$. That is, $\omega +1 \leq y-\Delta/2 $. Then, $M_{i'} \leq \omega [2(s+r)+r] +s+2r < (\omega+1) [2(s+r)+r] \le (y-\Delta/2) [2(s+r)+r] = 2y(s+r)+ry-(\Delta/2) [2(s+r)+r] < 2y(s+r)+ry-\Delta=2q-2e+ry-\Delta=2q - e < 2q$, where we have used that $\Delta = ry-e$ and $q = y(r+s)+e$.

\qed {\em $\,$ (Claim)}

\medskip

In either $(A)$ or $(B)$, if $\Delta'''=sy+ e + \Delta - \omega (e-\Delta)$, then $i'\ge  x +  \Delta + \omega (e-\Delta) + \Delta'''\ge x + sy+  (e + \Delta) + \Delta = x +(s+r)y +\Delta > x+2y$ and Theorem~\ref{thm:q+y,m_b=1}$(2)$ holds.

Thus, we assume $\Delta''' = \Delta + \omega(e-\Delta)\ge e$ in $(A)$ and $\Delta''' = e-\Delta$ in $(B)$. Then $i'\ge  x +  \Delta + \omega (e-\Delta) + \Delta'''> x +\Delta + e = x+ry\ge x+2y$ in $(A)$, and $i' \ge x+ \Delta + (\omega+1) (e-\Delta) \ge x+2y-1$ in $(B)$, recalling $\Delta + (\omega+1) (e-\Delta)\ge 2y-1$. Thus, in both cases, Theorem~\ref{thm:q+y,m_b=1}$(2)$ holds.

This completes the proof in Case II$(2)(a)$. 

\noindent{\bf Consider Case II$(2)(b)$ above.}  

Recall that in Case II$(2)(b)$, we define $\Delta''=\min\{e - \Delta, \Delta\}$ and $b'' \geq b' + \Delta'' \geq x + \alpha e + \Delta + \Delta''$. We now divide the argument according to whether $(i)$ or $(ii)$ of Case II$(2)(b)$ occurs.

% \begin{clm}\label{clm:II(2)(b)b''}
% $M_{b''} <  2q$
% \end{clm}

% \noindent {\em Proof of Claim:}
% First assume $\beta=1$. As $2e=(\beta+1)e>ry$, $2q=2(y(r+s)+e) > 4(r+s)+ry >3(r+s)+r=M_{b''}$, verifying the Claim.

% Now assume $\beta=2$. Then $\alpha >0$ and $1+ \alpha e + \Delta < 2y$ along with $e \geq 3$ implies $y \geq 3$. Thus, $M_{b''}=(2 \beta +1) (s+r)+r = 5s + 6r < 6(s+r) < 2y(s+r)+2e=2q$ as claimed. 

% We hereafter assume $\beta >2$ and hence that $\alpha \geq 2$. Then $1+ \alpha e + \Delta < 2y$ along with $e \geq 3$ implies that $y \geq 5$.  

% As $1+ \alpha e + \Delta < 2y$, we have $\alpha < \frac{2y}{e}$. Thus,
% \[
% \begin{aligned}
%     M_{b''}&=(2 \beta +1) (s+r)+r \leq (2 \alpha+3) (s+r) +r < (\frac{4y}{e} + 3) \frac{q-e}{y} + r = (\frac{4}{e}+\frac{3}{y}) (q-e) +r\\
%     &= (\frac{4}{e}+\frac{3}{y}) q - (\frac{4}{e}+\frac{3}{y}) e + (2+\frac{e}{y}-\frac{1}{y}) = (\frac{4}{e}+\frac{3}{y}) q - (2+\frac{2e}{y} + \frac{1}{y}).
% \end{aligned}
% \]
% where above we use that $1+(r-2)y<e$.
% Thus, we conclude
% \[M_{b''}<(\frac{4}{e}+\frac{3}{y}) q.\]

% Now $e \geq 3$ and $y \geq 5$, so $M_{b''} < 2q$ as claimed. 
% \qed {\em $\,$ (Claim)}

% As Theorem~\ref{thm:q+y,m_b=1}$(2)$  does not hold, the Claim implies that $b'' <q+y-1$. That is 
% $$\alpha e + \Delta + \Delta''<2y-1$$. 

\bigskip

{\bf Assume Case II$(2)(b)$(i) occurs.} That is $\Delta''=e-\Delta$. By Claim~\ref{clm:II(1)(2)}, we have $e>(r-2)y$. If $r> 3$, then $b''\ge x+\alpha e +\Delta +(e-\Delta) = x+(\alpha +1)e \ge x+e > x+(r-2)y \ge x+2y$, and Theorem~\ref{thm:q+y,m_b=1}$(2)$ holds. We hereafter assume $r\ge 3$.

The coefficient of $\wt{-r}$ in the total weight of $\Gamma_{b''}$ is $sy+ (\beta+1) e$. Thus $m_j \leq 3$ for $j=b''+(sy+ (\beta+1) e)-1\geq x+(\alpha+1)e+sy+(\beta+1)e-1=x+\alpha e+(e-\Delta)+(\Delta+\beta e)+sy+e-1=x+\alpha e+(e-\Delta)+(ry+sy+e)-1=2q-y+\alpha e+(e-\Delta)-1 \geq 2q-y$, implying that Theorem~\ref{thm:q+y,m_b=1}$(1)$ holds.

\bigskip

{\bf Assume Case II$(2)(b)$(ii) occurs.} That is $\Delta''=\Delta$ and for $b'' \geq b'+\Delta'' \geq x+\alpha e + 2 \Delta$ the weighting of $\Gamma_{b''}$ is 
\[(e-2\Delta) \wt{\beta (s+r)}  +  (sy+ \beta e + 2\Delta)  \wt{-r} + \Delta \wt{(2 \beta +1) (s+r)+r}.\]

Assume $2y-1\le \alpha e + 2\Delta$. Then $b''\ge x+\alpha e + 2\Delta \ge x+2y-1$, and Theorem~\ref{thm:q+y,m_b=1}$(2)$ holds. Thus, we assume that
\[2y-1>\alpha e +2\Delta.\]

% From the general discussion of II$(2)(b)$ above, we have $2y-1 > \alpha e + \Delta + \Delta''=\alpha e + 2\Delta$. Thus we have 
% \[e > 2\Delta, \,\, 1+ \alpha e + 2\Delta < 2y, \,\, 1+(r-2)y<e.\]
% From the general discussion of II(2), we also have
% \[r \geq 2, \,\, e>2, \,\, \beta>0.\]

\begin{clm}\label{clm:II(2)(b)alpha=0}
Either Theorem~\ref{thm:q+y,m_b=1} holds or $\alpha=0,\, \beta=1$.
\end{clm}

\noindent {\em Proof of Claim:}

If $\alpha=0$, then Claim~\ref{clm:II(1)(2)} says $\beta = \alpha + 1=1$. So assume $\alpha >0$. Claim~\ref{clm:II(1)(2)} says $2y-1>e>1+(r-2)y$, which implies $r \leq 3$. Note that $(1-\alpha)e \le 3 \Delta+1$. Then $m_j \leq 3$ if $j=b''+sy+ \beta e + 2\Delta-1=  x+\alpha e + 2 \Delta + sy+ \beta e + 2\Delta-1 \geq x+ e+sy+(\beta e+\Delta)+ (\alpha-1)e+3 \Delta-1=2q-y+ (\alpha-1)e +3 \Delta-1 \ge 2q-y-2$. Thus, Theorem~\ref{thm:q+y,m_b=1}$(1)$ holds. 

\qed {\em $\,$ (Claim)}

\begin{defin}
Let $s,r,e, \Delta, \beta$ be as above.
$Q_3(i, \omega)$ is the statement: There is an $i \geq b'+ (\omega-1) \Delta =x +\omega \Delta$ and the total weight of $\Gamma_i$ is 
\[(e- \omega \Delta) \wt{s+r}  +  (sy+ e + \omega \Delta) \wt{-r} + \Delta \wt{(\omega +1) (s+2r)-2r}\]
with $e- \omega \Delta > 0, \, \omega \ge 2, \, \omega \Delta <2y-1$.
\end{defin}

Then $Q_3(b'',2)$ holds. {\bf Let $i$ be largest so that $Q_3(i, \omega)$ holds for some $\omega$.} Then $i \geq b''>b'$. Applying Lemma~\ref{lem:b'>b} to $\Gamma_{i}$ and using the definitions of $i, b'$ as the largest values satisfying their defining conditions, we are left with two possibilities:

\begin{enumerate}[label=(\Alph*)] 
\item There is an $i'  \geq i+(e- \omega \Delta) \geq x + \omega \Delta +(e- \omega \Delta)=x + e$ such that the total weight of $\Gamma_{i'}$ is 
\[(sy + 2\omega \Delta) \wt{-r} + (e-(\omega-1)\Delta) \wt{(\omega +1) (s+2r)- 2r} + (e - \omega \Delta) \wt{-\omega (s+2r)}.\]

\item Let $\Delta'''=\text{min}\{e - \omega \Delta, \Delta\}$. There is a $i' \geq i + \Delta''' \geq x + \omega \Delta + \Delta'''$ such that the total weight of $\Gamma_{i'}$ is
\[((\omega+1)\Delta-e) \wt{(\omega+1) (s+2r)-2r} +(sy+ 2 e) \wt{-r} + (e - \omega \Delta)  \wt{(\omega+2) (s+2r)-2r}\]
when $\Delta'''= e - \omega \Delta$,
\[(e -(\omega+1)\Delta) \wt{s+r} + (sy+ e +(\omega+1) \Delta)  \wt{-r} + \Delta \wt{(\omega+2) (s+2r)-2r}\]
when $\Delta'''= \Delta<e - \omega \Delta$. Furthermore, $(\omega+1) \Delta \geq 2y-1$.
\end{enumerate}

\bigskip

% By the definition of $Q_3(i', \omega)$, $\omega \geq 2$ and $2y>1+ \omega \Delta$. Rewriting this,
% $2y + e > 1 +  \omega\Delta + e = 1 + (e + \Delta) + (\omega -1) \Delta=1+ ry + (\omega-1) \Delta$. That is $e> 1 + (r-2) y + (\omega-1) \Delta$.  %Similarly, starting from $2y>1+ \omega \Delta$ and recalling that $e > 2 \Delta$, we get that $2y > 1+ \omega \Delta$. 
% We post these for reference in the remainder of the argument:
% \[2y>1+ \omega \Delta, \, \, e> 1 + (r-2) y + (\omega-1) \Delta.\]

\begin{clm}\label{clm:II(2)(b)i'}
$M_{i'} <  2q$.
\end{clm}

\noindent {\em Proof of Claim:}
Then $M_{i'} \leq (\omega +2) (s+2r) -2r  < ({2y}/{\Delta} +2) (\pth{q+\Delta}/{y})-2r$, where we are using from $Q_3(i', \omega)$ that $2y-1>\omega \Delta$, together with the fact that $\Delta = ry-e$ so $q+\Delta=y(s+r)+e+\Delta=y(s+2r)$. Thus $M_{i'} < ({2}/{\Delta}+{2}/{y})q + (2+{2 \Delta}/{y})-2r <({2}/{\Delta}+{2}/{y})q + 4-2r \leq ({2}/{\Delta}+{2}/{y})q$ since $r \geq 2$ and since $\omega \geq 2$ means $2y>1+ 2 \Delta$.
%%using that $1+\alpha e + 2\Delta < 2y$ and that $q > 32$.  That is, 
% $$M_{i'} < (\frac{2}{\Delta} + \frac{2}{y})q$$
If $\Delta \geq 2$, then $M_{i'} < 2q$.  So we assume $\Delta=1$. As before, $M_{i'} \leq  (\omega +2) (s+2r) -2r$. Recall that $q=y(s+2r)-\Delta=y(s+2r)-1$ and that $\omega \Delta <2y-1$. The latter we rewrite $\omega+2 \leq 2y$. Thus $M_{i'} \leq  (\omega +2) (s+2r) -2r <2y(s+2r)-2=2q$, as claimed.\qed {\em $\,$ (Claim)}

\medskip
% By the same argument in the proof of Claim~\ref{clm: II M < 2q second} when $\alpha = 0$ and $\beta = 1$, since $\Delta = ry - e$ and $\omega\Delta \le 2y - 1$, we obtain
% \[M_{i'} \le \frac{2q}{\Delta} + 2 + 2r + 2s - \frac{2q}{\Delta y} - \frac{2}{y}.\]
% Thus, if $\Delta \ge 2$, since $y \ge 2$, $M_{i'} < 2q/\Delta + 1 + 2(r + s) \le 2q$ as claimed. Now assume $\Delta = 1$. From $e> 2\Delta = 2(ry - e) \ge 2(2y - e) > 3y - 2e$, we know $e > y$. Then,
% \[M_{i'} \le 2q - \frac{2}{y}(q - (r + s)y - (y - 1)) < 2q - \frac{2}{y}(q - (r + s)y - e) = 2q.\]
% \qed {\em $\,$ (Claim)}

We now argue according to whether $i'$ is as in $(A)$ or in $(B)$ above.

First we note that if case $(B)$ occurs and $\Delta'''=\Delta$, then $i' \geq x + (\omega+1) \Delta \geq x + 2y-1 = q+y-1$, and Claim~\ref{clm:II(2)(b)i'} implies that Theorem~\ref{thm:q+y,m_b=1}$(2)$ holds. Thus in Case $(B)$ we may assume that $\Delta'''=e-\omega \Delta$. 

In both Cases $(A)$ and $(B)$ we have then that $i'  \geq x + \omega \Delta +(e- \omega \Delta)=x + e$. If $e\ge 2y-1$, then Claim~\ref{clm:II(2)(b)i'} implies that Theorem~\ref{thm:q+y,m_b=1}$(2)$ holds. Thus, we may assume that
\[2y > e + 1.\]

By $Q_3(i,\omega)$, we have $e > 2\Delta = 2ry - 2e > r(e + 1) - 2e = (r - 2)e + r$. Thus it must be that $r=2$. Then, $2\Delta > (r - 2)e + r = 2$, and we have
\[r=2, \, \, \Delta \geq 2.\]

% $\spadesuit$ {\em (Or $i' \geq x + e > x + 2y=q+y$, and Claim~\ref{clm:II(2)(b)i'} contradicts that Theorem~\ref{thm:q+y,m_b=1}$(2)$ does not hold.)}

%Thus it must be that $r=2$. Then, $2\Delta > (r - 2)e + r = 2$, and we have
%\[\Delta \geq 2.\]

%$\spadesuit$ {\em (Actually, don't we have a contradiction here: $\Delta >e+1>\Delta+1$; hence, we need go no further. In any case, $e \geq 3$, so $\Delta>4$? $\spadesuit$ Don't know where $\Delta > e+1$ came from.)}

%$\spadesuit$ {\em (Note that $i' \geq x + (\alpha +1) e$ is not going to get $i' \geq x +2y$ when $r=2$ as $\alpha +1$ is likely $\beta$ and $0<\Delta=2y-\beta e$. Further note that $\Delta$ is the difference between $e$ and $2y$.}

Consider Case $(B)$, where we may assume $\Delta'''=e-\omega \Delta$. As $r =2$, $e > 2 \Delta$, and because an isotopy class of $\Gamma_z$ loses at most one arc as we move to $\Gamma_{z+1}$, we have that $m_{i''} \leq 2$ for $i''=i'+sy+ 2 e -1 \geq x + e + sy+ 2 e -1 = x + e + sy+ (e +\Delta) + (e-\Delta) -1 \ge x + (e+sy+ry) = 2q-y$, implying that Theorem~\ref{thm:q+y,m_b=1}$(1)$ holds. 

Thus we are left with Case $(A)$ with $r=2$. Recapping, there is an $i' \geq x + e$ such that the total weight of $\Gamma_{i'}$ is 
\[(sy + 2\omega \Delta) \wt{-r} + (e-(\omega-1)\Delta) \wt{(\omega +1) (s+2r)-2r} + (e - \omega \Delta) \wt{-\omega (s+2r)}.\]
We assume $i'$ to be largest satisfying the above.

\begin{clm}\label{clm:II(2)(b)2omega}
Either Theorem~\ref{thm:q+y,m_b=1} holds, or both $e>(2\omega-1)\Delta +1$ and $y\geq \omega \Delta +1$.

\end{clm}

\noindent {\em Proof of Claim:}
If $e\le (2\omega-1)\Delta + 1$, as $r=2$, looking at the the total weight of $\Gamma_{i'}$, we have that $m_j \leq 2$ if $j=i'+sy+2\omega \Delta-1 \ge x+e+ sy+2\omega \Delta-1=x+e+sy+(\Delta+e)+(2\omega-1)\Delta-e-1=x+(e+sy+ry)+(2\omega-1)\Delta-e-1\ge 2q-y-2$. Thus
Theorem~\ref{thm:q+y,m_b=1}$(1)$ holds. We assume $e>(2\omega-1)\Delta + 1$. We have $2y-\Delta=e \geq (2\omega-1) \Delta +2$, implying $2y \geq 2 \omega \Delta+2$.
That is $y\geq \omega \Delta +1$. 

\qed {\em $\,$ (Claim)}

%$\spadesuit$ {\em (If we want the stronger version of Theorem~\ref{thm:q+y,m_b=1}$(1)$, we need to change the Claim to: {\em $1+e > (2\omega-1) \Delta$ and $\omega \Delta \leq y+1$.} Because we argue $2q-y > j=2q-y+(2\omega-1)\Delta-e-1$. We continue $2+2y-\Delta=2+e \geq (2\omega-1) \Delta$, implying $2(y+1) \geq 2 \omega \Delta$.)}

Note that the inequality $e > \Delta = ry-\beta e$ becomes $2e > 2y$. Collecting our inequalities:
\[2y-1>e>y, \,  r=2, \, 2\leq \Delta=2y-e, \, \omega \geq 2, \, e>(2\omega-1)\Delta +1, \, y\geq \omega \Delta +1.\]

Once more applying Lemma~\ref{lem:b'>b} to $\Gamma_{i'}$ (by Claim~\ref{clm:II(2)(b)2omega}, we assume the coefficients above are positive) and using the definitions of $i,i'$ as the largest value satisfying its defining condition, we are left with two possibilities:

\begin{enumerate}[label=(\alph*)] 
\item Let $\Delta''''=\min\{sy + 2\omega \Delta, e-(\omega-1)\Delta\}$. There is an $i''  \geq i' + \Delta'''' \geq x + e + \Delta''''$ such that the weights of $\Gamma_{i''}$ are in $\{-r, (\omega +1) (s+2r) - 2r, -\omega (s+2r), (2 \omega+1) (s+2r) - 3r\}$.

%$$(sy+(3 \omega -1) \Delta -e) \,\, [\wt{-r}] +(2e - (2 \omega-1)\Delta) \wt{-\omega (s+r) - \omega r}$$
%$$ + (e-(\omega-1) \Delta) \wt{(2 \omega+1) (s+r) + (2 \omega-2) r}$$
%Note that here we are using that $e>(2\omega-1)\Delta +1 \geq 3 \Delta+1$. 

\item Let $\Delta''''=\min\{sy + 2\omega \Delta, e-\omega \Delta\}$. There is an $i''  \geq i' + \Delta'''' \geq x + e + \Delta''''$ such that the weights of $\Gamma_{i''}$ are in $\{-r, (\omega +1) (s+2r) - 2r, -\omega (s+2r), -(2 \omega+1) (s+2r) + r\}$.

%\item There is an $i''  \geq i' + (e- \omega \Delta) \geq x + 2e - \omega \Delta$ such that the total weight of $\Gamma_{i''}$ is 
%$$(sy+3 \omega \Delta -e) \,\, [\wt{-r}] +(2e - (2 \omega-1)\Delta) \wt{(\omega+1) (s+r) + (\omega-1) r} $$
%$$ + (e- \omega \Delta) \wt{-(2 \omega+1) (s+r) - 2 \omega r}$$

\end{enumerate}

\begin{clm}\label{clm:II(2)(b)i''}
$M_{i''} < 2q$
\end{clm}

\noindent {\em Proof of Claim:}
From the list of weights of $\Gamma_{i''}$ we have $M_{i''} \leq (2 \omega +1)(s+2r) -r$.
Noting that $q+\Delta=y(s+r)+e+\Delta=y(s+2r)$ and recalling that
$\omega \leq \pth{y-1}/{\Delta}$. Then
\[\begin{aligned}
    M_{i''} &\le \pth*{\frac{2y-2}{\Delta}+1} \pth*{\frac{q+\Delta}{y}} -2 =\pth*{\frac{2y-2}{\Delta}+1} \pth*{\frac{q}{y}} + \pth*{\frac{2y-2}{\Delta}+1} \pth*{\frac{\Delta}{y}}-2\\
    &= \pth*{\frac{2}{\Delta}-\frac{2}{\Delta y}+\frac{1}{y}}q+ \pth*{2+\frac{-2+\Delta}{y}} -2 < \pth*{\frac{2}{\Delta}+\frac{1}{y}}q+ \frac{\Delta}{y} < \pth*{\frac{2}{\Delta}+\frac{1}{y}}q+ 2.
\end{aligned}\]
In the last step we use $\Delta=2y-e<2y$.
Recall that $\Delta \geq 2$, $\omega \Delta \leq y-1$, and $\omega \geq 2$, so $y \geq 5$. Thus, $M_{i''} < 6q/5 + 2$. Since $q > 3$, we conclude that $M_{i''} < 2q$.

\qed {\em $\,$ (Claim~\ref{clm:II(2)(b)i''})}

However, in both Cases $(a)$ and $(b)$ above, $i'' \geq x+e+ \min\{sy + 2\omega \Delta, e-\omega \Delta\}$.

Observe that $e> (2\omega-1) \Delta +1$ implies $e - \omega \Delta > (\omega-1) \Delta +1 > \Delta$ as $\omega \geq 2$. We have that $i'' \geq x+e+ \min\{sy + 2\omega \Delta, e-\omega \Delta\} > x+e+\Delta=x+ry = x+2y > q + y - 1$.

Therefore, in both cases, Claim~\ref{clm:II(2)(b)i''} implies that Theorem~\ref{thm:q+y,m_b=1}$(2)$ holds. 

This completes the proof in Case II$(2)(b)$(ii). We thus complete the proof in Case II(2). We have now completed the proof in all of Case II(1),(2),(3), which completes the proof in Case II, finishing the proof of Theorem~\ref{thm:q+y,m_b=1}.
\end{proof}

\subsection{When \pdfmath{q<4y}{q<4y} or \pdfmath{q<4y'}{q<4y'}}\label{subsec:qsmall}

One of the conclusions of Theorem~\ref{thm:q+y,m_b=1} is that $q<4y$ or $4y'$. In this subsection, we explore that situation. 

\begin{lem}\label{lem:q<4y}
Let $r,s,e,r',s',e'>0$ be as in Theorem~\ref{thm:=x}.
\begin{enumerate}
\item Assume $q < 4y$. Then $r+s \leq 3$, and for any $j$ with $j \leq q+ e -1$,  $m_j \leq 2$. Furthermore,
\begin{itemize}
    \item if $r+s=2$, then $m_j \leq 1$;
    \item if $r+s=3$, then $e=[q]_y<y$ and for any $z \leq q+(y-e)-1$, $m_{z} \leq 1$.
\end{itemize}
Note then that $q<4y$ implies that if $z \leq q$, then $m_z \leq 1$.

\item Assume $q < 4y'$. Then $r'+s' \leq 3$, and for any $j$ with $j \leq q+ e' -1$, $m_j' \leq 2$. 
Furthermore,
\begin{itemize}
    \item if $r'+s'=2$, then $m_j' \leq 1$;
    \item if $r'+s'=3$, then $e'=[q]_{y'}<y'$ and for any $z \leq q+(y'-e')-1$, $m_z' \leq 1$.
\end{itemize}
Note then that $q<4y'$ implies that if $z \leq q$, then $m_z' \leq 1$.
\end{enumerate}
\end{lem}

%$\spadesuit$ {\em (We can take any $j$ up to the first $\mathcal{T}$-level at least as big as $q+e-1$ and, when $r+s=3$, for any $\iota$ up to the first $\mathcal{T}$-level at least as big as $q+e-1$.)}

\begin{proof}
We prove $(1)$. The proof for $(2)$ follows by replacing $\Gamma$ with $\Gamma'$. 

By Theorem~\ref{thm:=x}, the total weight of $\Gamma_x$ is 
\[(ry - \beta e) \wt{s} + (sy + \beta e) \wt{-r} + (e) \wt{\beta {(s+r)}}\]
for some $|\beta|=1$ and where the coefficients
$ry - \beta e, sy + \beta e, e$ are positive.
Then $4y>q=(r+s)y+e$, and, as $e>0$, we conclude that $0<r+s \leq 3$ with $r,s>0$. Note also that $[e]_y=[q]_y$. Theorem~\ref{thm:m:b<x} says that if $z < x$, then $m_z=0$. By the above we see that $m_x=1$.
Because an isotopy class of $\Gamma_z$ loses at most one arc as we move to $\Gamma_{z+1}$, we have that $m_{j} \leq 2$ for any $j$ with $j \leq q+e-1$, as $\max \{ry - \beta e, sy + \beta e\} -1 \geq y+e-1$. Furthermore, knowing $r=s=1$ sharpens the statement to $m_j \leq 1$. Now assume that $r+s=3$. WLOG assume that $r=1,s=2$. Again because an isotopy class of $\Gamma_z$ loses at most one arc as we move to $\Gamma_{z+1}$, we have that $m_{b} \leq 1$ for any $j$ with $j \leq q+(y-e)-1 = x+2y-e-1 \leq x + sy + \beta e - 1$. As $4y>q=(r+s)y+e$, note that when $r+s=3$, $y>e=[q]_y$ ($y$ is relatively prime to $q$). 

Finally, assume $z \leq q$. Assume $r+s=2$. Then as $z \leq q \leq q+ e -1$, the above shows that $m_z \leq 1$. So assume $r+s=3$. As $z \leq q \leq q+(y-e)-1$, the above shows that $m_z \leq 1$.
\end{proof}

\section{\pdfmath{L(q^2, m)}{L(q2, m)} with \pdfmath{m=dq+1}{m=dq+1}.}\label{sec:dq+1}

In this section, we assume we are in the context of conclusion $(1)$ of Theorem~\ref{thm:lisca}: $K \subset L(q^2, m)$ with $m=dq+1$ where $d, q>0$. Recall that $K=K(q^2, dq+1, \kappa)$ is a simple knot that admits longitudinal surgery to $S^1 \times S^2$. The goal of the section is to prove Theorem~\ref{thm:main:dq+1}, which is basically the Spherical Theorem, Theorem~\ref{thm:Gr-Li-Lu}, in this setting. 

Recall that $\mu'=m \mu$ where now $m = dq+1 \equiv 1\,\,\,(q)$. Note that $m^{-1}\equiv 1 -dq\,\,\, (q^2)$. Thus, $\mu=(1-dq) \mu'$. In particular $q \mu= q \mu'$. Recall that $[K]=\pm y'q \mu'=\pm y'q \mu$ and $[K]=\pm yq \mu$. As both $0<y,y'<q/2$, \textbf{we conclude that $y=y'$}. In addition, \textbf{recall that $[K] \in H_1(L(q^2,m))$ has order $q$}, which we assume throughout this section.

%In this setting, we prove Theorem~\ref{thm:main:dq+1}, as a version of the Model Theorem~\ref{thm:Gr-Li-Lu}. 

Corollary~\ref{cor:Mmdecom}, Theorem~\ref{thm:q+y,m_b=1} (Engine Theorem), and 
Theorem~\ref{thm:M:b<x+1} (Initial Growth Theorem), 
 allow us to write $[K]=c \mu \pm f \mu'$ where $c$ and $f$ are relatively small.
 Theorem~\ref{thm:c+d:dq+1} uses the order of $[K]$ to list the possible classes of $[K]$ that can arise with the kind of bounds on $c,f$ that arise. Amongst these, Theorem~\ref{thm:main:dq+1} shows that conclusions of the Engine Theorem only allow the possibilities $\pm q\mu, dq \mu$. Lemma~\ref{lem:dq+1,4y} and Theorem~\ref{thm:m_b'leq3:dq+1} address two of the four conclusions of the Engine Theorem in anticipation of Theorem~\ref{thm:main:dq+1}. 
 %completes the proof of Theorem~\ref{thm:main1-nd}.

\begin{thm}\label{thm:c+d:dq+1}
Assume $q>3$. Assume $[K]=c \mu \pm f \mu'$ and one of the following holds on $c,f$:

\begin{enumerate}
\item $c=0$ and $0 \leq f < q^2/2$;
\item $0 < c \leq 3$ and $0 \leq f < 2q$;  
\item $0 < f \leq 3$ and $0 \leq c < 2q$;
\item $f=0$ and $0 \leq c < q^2/2$.
\end{enumerate}

Then, we have the corresponding conclusions:

\begin{enumerate}
\item $[K]= \pm f \mu' = \pm f \mu$ where $f= y' q$.
\item Either
\begin{enumerate}
\item$[K]=-cdq \mu=-cdq \mu'$ with $f=c$;
\item $[K]=(1-cd)q \mu$ with $f=q-c$;
\item $[K]=(-1-cd)q \mu$ with $f=q+c$;
\item $[K]=(2-cd)q \mu$ with $f=2q-c$.
\end{enumerate}

\item Either
\begin{enumerate}
\item $[K]=-fdq \mu=-fdq \mu'$ with $c=f$;
\item $[K]=(1+fd)q \mu$ with $c=q-f$;
\item $[K]=(1-fd)q \mu$ with $c=q+f$;
\item $[K]=(2+fd)q \mu$ with $c=2q-f$.
\end{enumerate}

\item $[K]= c \mu$ and $c= yq$. 

\end{enumerate}
\end{thm}

\begin{proof}
$[K]=c \mu \pm f \mu' = (c \pm fm) \mu$. As $[K]$ has $q$ torsion, $q | c \pm fm$.
That is
\[0 \equiv c \pm fm \equiv c \pm f \,\,\, (\mod q);\ c\equiv \mp f\,\,\, (\mod q).\]

{\bf Assume $(1)$ holds.} Then $[K]=\pm f \mu'$. Recall $[K]=\pm y'q \mu'$. Therefore, $f \equiv \pm y'q \pmod{q^2}$.
As $0 \leq f, y'q < q^2/2$, we conclude that $f = y'q$, and $[K] = \pm y'q\mu' = \pm y'q\mu = \pm f\mu$.

{\bf Assume $(2)$ holds.} As $c\equiv \mp f\pmod q$, we may write $f = \eta q\mp c$, where the sign respects that of $[K] = c\mu \pm f\mu'$. Then $[K] = (c\pm fm)\mu = (c\pm(\eta q\mp c)(dq + 1))\mu = (\pm\eta - cd)q\mu$. 

When $f=c$ we have $\eta=0$ and the sign is $+$ (that is, we get $[K]=(\eta - cd)q\mu$ above), so we have $[K]=-cdq\mu$ giving conclusion $(a)$. When $f=2q-c$, $\eta=2$ with sign $-$, and $[K]=(2-cd)q\mu$---giving conclusion $(d)$.

From the assumption $0 < c \leq 3$ and $0 \leq f < 2q$, we conclude $f = c$, $q-c$, $q+c$, or $2q-c$. We have two additional cases to the above. For $f=q \mp c$, with $\eta=1$ and sign $\mp$ (resp.), we conclude $[K]=(\pm1 - cd)q\mu$ (resp.)---giving conclusions $(b)$, $(c)$.

{\bf Assume $(3)$ holds.} Similar to $(2)$, from the assumption $0 < f \leq 3$ and $0 \leq c < 2q$, we conclude that $c=f,q \pm f, 2q-f$, denoted as $c = \eta q\mp f$, where $\eta$ and the sign are in $\bc{(0,+), (1,-), (1, +), (2, -)}$. Then $[K] = (c\pm fm)\mu = ((\eta q\mp f)\pm f(dq + 1))\mu = (\pm\eta fd)q\mu$.

% The following is the orginal proof
% As $0 \leq c<2q$ we conclude that $c=f,q \pm f, 2q-f$. 

% If $c=f$, then $[K]=(c -fm) \mu=(f-f(dq+1)) \mu=-fdq \mu=-fdq \mu'$.

% If $c=q-f$, then $[K]=(c+fm) \mu=((q-f)+f(dq+1)) \mu=(1+fd)q \mu$.

% If $c=q+f$, then $[K]=(c -fm) \mu=((q+f)-f(dq+1)) \mu=(1-fd)q \mu$.

% If $c=2q-f$, then $[K]=(c +fm) \mu=((2q-f)+f(dq+1)) \mu=(2+fd)q \mu$.

{\bf Assume $(4)$ holds}. Then $[K]=c \mu=\pm yq\mu$. Therefore, $c \equiv \pm yq\,\,\, (\mod q^2)$. As $0 \leq c, yq < q^2/2$, we conclude that $c = yq$ as claimed.
\end{proof}

\begin{lem}\label{lem:dq+1,4y}
If $12\le q < 4y=4y'$, then $[K]= \pm dq \mu$.
\end{lem}

\begin{proof}
As $12\le q<4y$, we may assume $ y=y'>3$. Recall $y(r+s)+e=q$, so $2\le r+s \le 3$.

\begin{clm}\label{clm:q+y-2}
The following statements hold.
\begin{enumerate}
    \item If there is a $b \geq q+y-2$ such that $m_b \leq 1$, then $[K]= \pm dq \mu$. 
    \item If there is a $b \geq q+y'-2$ such that $m_b' \leq 1$, then $[K]= \pm dq \mu$.
\end{enumerate}
\end{clm}

\noindent {\em Proof of Claim:}
(1) Assume there is a $b \geq q+y-2$ such that $m_b \leq 1$.
Let $b'=2q-2-b \leq q-y=q-y'=x'$. By Theorem~\ref{thm:M:b<x+1}, $M'_{b'} < q/y' \le q/4$. By Corollary~\ref{cor:Mmdecom}, up to sign $[K]=c \mu \pm f \mu'$ where $0\leq c< q/4$ and $0 \leq f \leq 1$. By Theorem~\ref{thm:c+d:dq+1}(3),(4), we have $c=f=1$ and $[K]=-dq \mu$.

(2) The proof is identical to $(1)$ with $\Gamma, \,m, \,M', \,x', \,y, \,y', \,c, \,f$, Theorem~\ref{thm:c+d:dq+1}(3),(4) replaced by $\Gamma, \,m', \,M, \,x, \,y', \,y, \,f, \,c$, Theorem~\ref{thm:c+d:dq+1}(1),(2), respectively.
% The following is the original proof of (2)
% Assume there is a $b \geq q+y'-2$ such that $m_b' \leq 1$. Let $\tau$ be an arc in $\Gamma_b'$ with $\beta$-weight at most one. Let $b'=2q-2-b \leq q-y'=q-y=x$. By Theorem~\ref{thm:M:b<x+1}, $M_{b'} < q/y \le q/2$. By Lemma~\ref{lem:decomposeK}, up to sign $[K]=c \mu + f \mu'$ where $0 \leq c \leq 1$ and $0\leq f< q/4$. By Theorem~\ref{thm:c+d:dq+1}(1),(2), we again have $c=f=1$ and $[K]=-dq \mu$.
\qed {\em $\,$ (Claim)}
\medskip
%\begin{thm}\label{thm:c+d:dq+1}
%Assume  $[K]=c \mu \pm f \mu'$ and $K$ is spherical. Furthermore assume $q>18$ and one of the following holds on $c,f$:
%
%\begin{enumerate}
%\item $c=0$ and $0 \leq f < q^2/2$;
%\item $0 < c \leq 3$ and $0 \leq f < 2q$;  
%\item $0 < f \leq 3$ and $0 \leq c < 2q$; 
%\item $f=0$ and $0 \leq c < q^2/2$;
%\end{enumerate}
%
%then we have the corresponding conclusions:
%
%\begin{enumerate}
%\item $[K]= \pm f \mu' = \pm f \mu$ where $f= y' q$;
%\item Either
%\begin{enumerate}
%\item$[K]=-cdq \mu=-cdq \mu'$ with $f=c$;
%\item $[K]=(1-cd)q \mu$ with $f=q-c$;
%\item $[K]=-(1+cd)q \mu$ with $f=q+c$;
%\item $[K]=(2-cd)q \mu$ with $f=2q-c$.
%\end{enumerate}
%
%\item Either
%\begin{enumerate}
%\item $[K]=-fdq \mu=-fdq \mu'$ with $c=f$;
%\item $[K]=(1+fd)q \mu$ with $c=q-f$;
%\item $[K]=(1-fd)q \mu$ with $c=q+f$;
%\item $[K]=(2+fd)q \mu$ with $c=2q-f$;

First assume $r + s=3$. By Lemma~\ref{lem:q<4y}, $e < y$ and for $b = q+y-e-1 \ge q$, $m_{b} \leq 1$. By Claim~\ref{clm:q+y-2}(1), we may assume $e>1$, so $e'\ge [q]_{y'}=[q]_y=e>1$. Let $b'=2q-2-b \leq q-2$. By Addendum~\ref{thm:b=<q}, Claim~\ref{clm:q+y-2}(2), and the fact that $q>3y'>9$, we may assume that $M_{b'}' < (y'/e'+2)(q-e')/y' < q-e' \le q-1$. By Corollary~\ref{cor:Mmdecom}, $\pm[K]=c \mu \pm f \mu'$ where $0\leq c< q-1$ and $0 \leq f\leq 1$. By Theorem~\ref{thm:c+d:dq+1}(3),(4), we have $c=f=1$ and $\pm[K]= -dq \mu$.

Thus we assume $r=s=1$. By Lemma~\ref{lem:q<4y}, for $b = q+e-1$, $m_b \leq 1$. By Theorem~\ref{thm:=x} where the coefficients are positive, $e<y$.
% If $e \geq y$, then $b \geq q+y-1$ and Claim~\ref{clm:q+y-2} proves the Theorem. $\spadesuit$ {(\em (I think $e<y$ by Theorem~\ref{thm:=x} and the fact that the coefficients are non-zero)}. We then further assume $e<y$.
Then $q=2y+e<3y=3y'$.
 
Let $b'=2q-2-b \leq q-e-1<q-[q]_{y}=q-[q]_{y'}$. If $e'>1$, then by Theorem~\ref{thm:M:b<x+1}, $M_{b'}'< (y'/e'+2)(q-e')/y'< q-e' <q-1$; if $e'=1$, by Addendum~\ref{lem:r=1=sM} and Claim~\ref{clm:q+y-2}, we may assume for any $j \leq q-1$, $M_{j}' < q-1$, in particular, $M_{b'}' < q-1$. In either case, by Corollary~\ref{cor:Mmdecom}, $\pm[K]=c \mu \pm f \mu'$ where $0\leq c< q-1$ and $0 \leq f\leq 1$. By Theorem~\ref{thm:c+d:dq+1}(3),(4), we have $c=f=1$ and $\pm[K]=-dq \mu$ as desired.
\end{proof}

% $\blacklozenge$ We have removed Lemma~\ref{lem:r=1=sM}, which was a version of Theorem~\ref{thm:M:b<x+1} (with $b<q-1$ in $(3)$ and $(4)$ but $M_j\le q-1$ whenever $j\le b$ in the case when $r=s=e=1$). The proof was very similar but simpler, so we may mark the difference there.

%\begin{thm}\label{thm:M:b<x+1} Let $e,e'$ be as in Theorem~\ref{thm:=x}.% Assume $b,b'$ are $\mathcal{T}, \mathcal{T}'$-levels.
%
%\begin{enumerate}
%\item If $b \leq x$, $M_b <q/y$. %$\spadesuit$ {\em Also, $m_b \leq \min\{r,s\}$?}
%\item If $b' \leq x'$, $M_b' < q/y'$. 
%%\item If $b \leq 4y$, $M_b \leq 4$ or $q \leq 5y$. 
%\item For every $b \leq q-[q]_y$, $M_b < (y/e+2)(q-e)/y$. For every $x \leq b \leq q-[q]_y, m_b \in \{r,s\}$. %$\spadesuit$ {\em Or, For $b \leq q-[q]_y,  m_b \in \{r,s\}$?}
%\item For every $b \leq q-[q]_{y'}$, $M_b' < (y'/e'+2)(q-e')/y'$
%\end{enumerate}

% \begin{lem}\label{lem:m_x+sy-e}
% Let $r,s,e,r',s',e'$ be as in Theorem~\ref{thm:=x}. Then  $m_{x+sy-e-1} \leq r$, \,\,$m_{x+ry-e-1} \leq s$,\,\, $m_{x+s'y-e'-1}' \leq r'$,\,\, and $m_{x+r'y-e-1}' \leq s'$. 

\begin{thm}\label{thm:m_b'leq3:dq+1}
Assume there is a $b \geq 2q-y'-2$ with $m_b' \leq 3$. Then 
\[\pm [K] = dq \mu.\]
%Furthermore, if $\pm [K]=2dq \mu$ then $K$ is the union of two arcs on either side of a level surface that writes $\pm [K]=2\mu -2 \mu'$.

%Furthermore, 
%
%\noindent (A). If $\pm [K] = 2q \mu$ then $m_b'=1, y=2,y'=2n$ and $|d| \geq 8$. 
%Furthermore $\pm [K]=\mu + n \mu'$ as in $(d)$ of Lemma~\ref{thm:c+d}(2) holds.
%
%\noindent (B). If $\pm [K] = 4q \mu$ then $m_b'=2, y=4,y'=4n$ and $|d| \geq 16$. 
%Furthermore $\pm [K]=2 \mu + 2n \mu'$ as in $(d)$ of Lemma~\ref{thm:c+d}(2).
%Finally, we can rule out this case if we can take $b \geq 2q-2y'+2$.
\end{thm}

\begin{proof}
Let $b'=2q-2-b \leq y'$. As $y'<q/2$, we have $x'=q-y'>y'\ge b'$. Lemma~\ref{thm:m:b<x}, along with the fact that $y=y'$, says that $M_{b'} < b'/y +1 \leq y'/y +1 = 2$. That is, $M_{b'} \leq 1$.
By Corollary~\ref{cor:Mmdecom}, we can write $\pm[K]=c \mu \pm f \mu'$ with $0\le c\le 3$ and $0 \leq f =1$. By Theorem~\ref{thm:c+d:dq+1}(2), we conclude $\pm[K]=-cdq \mu$ with $f=c=1$, as claimed.
\end{proof}

\begin{thm}\label{thm:main:dq+1}
Assume $q\ge 12$. Let $K=K(q^2, dq+1, \kappa)$ be a simple knot that admits longitudinal surgery to $S^1\times S^2$. Let $[K]$ denote its class in $H_1(L(q^2, dq+1))$. Let $\mu, \mu'$ be the standard generators of $H_1(L(q^2, dq+1))$.
Then
\[\pm [K]=q \mu, dq \mu.\]
%Furthermore, if $\pm [K]=2dq \mu$ then $K$ is the union of two arcs on either side of a level surface that writes $\pm [K]=2\mu -2 \mu'$.
%$\spadesuit$ {\em (The latter is shorthand for saying there is an arc in a level surface whose $\alpha$-weight is $\pm 2$ and whose $\beta$-weight is $\mp 2$.)}
\end{thm}

%$\spadesuit$ {\em (In Theorem~\ref{thm:c+d:dq+1}, can we get either $m_b \not\equiv 0,1 \,\, \mod{q}$ or $m_b \leq q$ for some $b \geq 2q-y-2$? This would force $[K]= \pm dq \mu$ as the $\beta$-weight would have to be $-1,0,1$. Or just eliminate the $m_b \leq 3$ conclusion.)}

\begin{proof}
%$K=K(q^2,d(q-1),\kappa)$ and $[K]=\kappa \mu = l \mu'$ where $\kappa=t'q$ and $l=tq$. Recall that $y = \min \{t',q-t'\}, y'= \min \{t, q-t\}$. %By Lemma ~\ref{lem:y=1}, we may assume $y >2, y'>1$. 
Note that if $K$ admits longitudinal $S^1\times S^2$ surgery, so does $K$ with its orientation reversed. Thus, we may work with $[K]$ up to sign. Let $e,\, r,\, s,\, e',\, r',\, s'$ be as in Theorem~\ref{thm:=x}. 

If $y=y'=1$, then $[K]=\pm yq \mu=\pm q \mu$. We hereafter assume $y=y'>1$. 
Therefore, we may apply
Theorem~\ref{thm:q+y,m_b=1} to $\Gamma'$ to conclude that there is $b$ such that one of the following holds:

\begin{enumerate}
\item  $b \geq 2q-y'-2$ and $m_b' \leq 3$;
\item  $b \geq q+y'-1$ and $M_b' < 2q$;
\item  $b \geq q+y'-1$ such that $m_b' \leq 1$; or
\item  $q<4y'$.
\end{enumerate}

If $(1)$ holds, Theorem~\ref{thm:m_b'leq3:dq+1} proves the Theorem. If $(4)$ holds, then Theorem~\ref{lem:dq+1,4y} proves Theorem~\ref{thm:main:dq+1}. We hereafter assume $q >4y'$ ($y',\, q$ are relatively prime).

Now we assume either $(2)$ or $(3)$ holds. Then, $b'=(2q-2)-b \leq q-y'-1 = q-y-1 =x-1$.

Assume that $(2)$ holds where $M_b' < 2q$. Theorem~\ref{thm:m:b<x} shows that $m_{b'}=0$.
By Corollary~\ref{cor:Mmdecom}, $\pm[K]=c \mu \pm f \mu'$ where $0<c < 2q$ and $f = 0$. Then theorem~\ref{thm:c+d:dq+1}$(4)$  shows that $c=q$ and $\pm [K]=q \mu$ as claimed.

Assume that $(3)$ holds with $m_b' \leq 1$.
Thus, by Theorem~\ref{thm:M:b<x+1}$(1)$, $M_{b'} < q/y \le q/2$. By Corollary~\ref{cor:Mmdecom}, $\pm[K]=c \mu \pm f \mu'$ where $c \leq 1$ and $0 \leq f < q/2 < q - 1$. By Theorem~\ref{thm:c+d:dq+1}(1), (2), $\pm[K]= -dq \mu$ as claimed.
\end{proof}

\section{\pdfmath{L(q^2, m)}{L(q2, m)} with \pdfmath{m=d(q-\epsilon)}{m=d(q-epsilon)}.}\label{sec:d(q-ep)}

In this section, we assume we are in the context of conclusions $(2)$ and $(3)$ of Theorem~\ref{thm:lisca}: $K \subset L(q^2, m)$ with $m=d(q-\epsilon)$ and $\epsilon=\pm 1$. Recall that $K=K(q^2,d(q-\epsilon) , \kappa)$ is a simple knot that admits longitudinal surgery to $S^1 \times S^2$.
%In this section, we almost prove the main theorem of the paper, Theorem~\ref{thm:Gr-Li-Lu}, for $K \subset L(q^2,m)$ when $m=d(q-\epsilon)$ and $\epsilon=\pm 1$. 
%As pointed out in the proof of Theorem~\ref{thm:Gr-Li-Lu} have two subcases to consider.
Conclusion $(2)$ and $(3)$ are then, respectively, 
\begin{align}
    m=d(q-\epsilon), &\qquad \epsilon=\pm1, &\qquad q= nd+\epsilon, &\qquad n \ge 2,\,\, d\geq 3 \,\, \hbox{and}\,\, d \,\, \hbox{odd}. \tag{$\ast{q}$} \\
    m=d(q-\epsilon), &\qquad \epsilon=\pm1, &\qquad 2q= nd-\epsilon, &\qquad n,d\ge 3 \,\,\hbox{both odd}. \tag{$\ast{2q}$}
\end{align}

The goal of this section is to prove the main theorem of the paper, Theorem~\ref{thm:Gr-Li-Lu}, in these settings. Theorem~\ref{thm:main1-nd}, in conjunction with Theorem~\ref{thm:q=1+nd,nmustbe2} of section~\ref{sec:alexwidth} (for $\pm [K]=2q\mu$), proves Theorem~\ref{thm:Gr-Li-Lu} when $(*q)$ holds. Theorem~\ref{thm:main1-nd:2q+1}  proves Theorem~\ref{thm:Gr-Li-Lu} when $(*2q)$ holds.

Though more elaborate, the arguments for Theorems~\ref{thm:main1-nd} and \ref{thm:main1-nd:2q+1} follow that for Theorem~\ref{thm:main:dq+1} in section~\ref{sec:dq+1}. 
Again, Corollary~\ref{cor:Mmdecom},  Theorem~\ref{thm:q+y,m_b=1}, and  Theorem~\ref{thm:M:b<x+1} 
 allow us to write $[K]=c \mu \pm f \mu'$ where $c$ and $f$ are relatively small.
 Theorem~\ref{thm:c+d} for case $(*q)$ and Theorem~\ref{thm:c+d:2q+1} for case $(*2q)$ use the fact that the order of $[K]$ is $q$ to list the possible classes of $[K]$ that can arise with the kind of bounds on $c,f$ that arise. Amongst these classes, the proofs of Theorems~\ref{thm:main1-nd} for $(*q)$ and \ref{thm:main1-nd:2q+1} for $(*2q)$ show how the conclusions of Theorem~\ref{thm:q+y,m_b=1} only allow the possibilities $\pm q\mu, dq \mu$ (and $2q\mu$ in case $(*q)$). Note that subsection~\ref{subsec:b>=2q-y-2} (in particular Theorems~\ref{thm:m_bleq3}, \ref{thm:m_b'leq3} and Theorems~\ref{thm:m_bleq3:2q+1},  \ref{thm:m_b'leq3:2q+1}) and subsection~\ref{subsec:q<4yor4y'} (Theorems~\ref{thm:q<4y:1-ndprime}, \ref{thm:q<4y':1-ndprime} and Theorems~\ref{thm:q<4y:2q+1}, \ref{thm:q<4y':2q+1}) address two of the four conclusions of Theorem~\ref{thm:q+y,m_b=1} (the $\Gamma$ and $\Gamma'$ variations). 

%We assume in this section, as in Theorem~\ref{thm:Gr-Li-Lu} and the earlier sections, that $K$ is a simple knot that admits longitudinal surgery to $S^1 \times S^2$. 
We begin by recalling the following facts used freely in this section:
\begin{itemize}
\item $[K] \in H_1(L(q^2,m))$ has order $q$. 
\item The positive integers $y,y'$ associated to $K$ in section~\ref{sec:intro} are determined by the properties
\[0<y,y'<q/2 \,\,\hbox{and} \,\, [K]=\pm yq\mu,\,\, [K]=\pm y'q\mu.\]
By definition, $x=q-y,\, x'=q-y'$.
\item Theorem~\ref{thm:=x} associates to $K$ the positive integers $r,s,e,r',s',e'>0$ satisfying
\[(r+s)y+e=q=(r'+s')y+e'.\]
\item $m \equiv -\epsilon d \, \mod \, (q)$. As $\mu'=m \mu \in H_1(L(q^2,m))$, we have that $q\mu'= -\epsilon d q \mu$.
\item In case $(*q)$, $m^{-1} \equiv n \,\, \mod \, (q)$. Thus $q \mu = m^{-1} q \mu'=nq \mu$.
\item In case $(*2q)$, $m^{-1} \equiv -n \,\, \mod \, (q)$. Thus $q \mu = m^{-1} q \mu'=-nq \mu$.

\end{itemize}

\subsection{General computations}\label{subsec:gencom}
\hfill

Theorems~\ref{thm:q+y,m_b=1} and \ref{thm:M:b<x+1} allow us to write $[K]=c \mu \pm f \mu'$ where $c,\, f$ are small on the scale of $q$. In this subsection, we use these bounds on $c,\, f$ along with order $[K] = q$ to restrict $[K]$ to a few possible classes. This is done in Theorem~\ref{thm:c+d} for class $(*q)$, and in Theorem~\ref{thm:c+d:2q+1} for class $(*2q)$. We also gather some supporting lemmas for subsequent subsections. 

\subsubsection{Case $(*q)$}\label{subsubsec:gencomq}

For this subsubsection, we assume that we are in case $(*q)$ where
\[m=d(q-\epsilon),  \qquad \epsilon=\pm1, \qquad q= nd+\epsilon, \qquad n \ge 2,\,\, d\geq 3 \,\, \hbox{and}  \,\, d \, \, \hbox{odd}.\]

\bigskip

\begin{thm}\label{thm:c+d}
Assume $K \subset L(q^2,m)$ for $m$ as in $(*q)$.
Assume $[K]=c \mu \pm f \mu' \in H_1(L(q^2,m))$. Furthermore, assume $q>18$ and one of the following holds on $c,f$:

\begin{enumerate}
\item $c=0$ and $0 \leq f < q^2/2$;
\item $0 < c \leq 3$ and $0 \leq f < 2q$;  
\item $0 < f \leq 3$ and $0 \leq c < 2q$; 
\item $f=0$ and $0 \leq c < q^2/2$.
\end{enumerate}

Then, we have the corresponding conclusions:

\begin{enumerate}
\item $[K]= \pm f \mu' = \mp \epsilon fd \mu$ where $f= y' q$. 
\item Either
\begin{enumerate}
\item $\epsilon[K]= 2(c-d)q \mu$ and $f=2q-cn$;
\item $\epsilon[K]= (2c+d)q \mu$ and $f=q+cn$; 
\item $\epsilon[K]= (2c-d)q \mu$ and $f=q-cn$;
\item $\epsilon[K]= 2cq \mu$ and $f=cn$;
\item $\epsilon=-1$, $c=d=3$, and up to sign $[K]=3q\mu=dq\mu$.
\end{enumerate}
%Cases $(b),(c)$ degenerate to $[K]= \pm dq$ when $f \leq q/3$. 

\item Either
\begin{enumerate}
\item $[K]=(\eta-\epsilon fd)q \mu$ where $\eta \in  \{1,2\}$ and $c=\eta q-fd$; 
\item $[K]= (\eta+\epsilon fd)q \mu,$ where $\eta\in\bc{0, 1}$ and $c=\eta q+fd$;
\item $n=2, \,f=3$, either $[K]= (\epsilon d - 2)q \mu$ with $c=d-\epsilon$, or $[K]= (4-\epsilon d)q \mu$ with $c=2q-(d-\epsilon)$;
\item $\epsilon = -1, \,n=f=2, 3$, either $[K]=-2q\mu$ with $c=1$, or $[K]=4q\mu$ with $c=2q-1$.
\end{enumerate}

\item $[K]= c \mu$ where $c= yq$ and $0<y<q/2$. 
\end{enumerate}
\end{thm}

\begin{proof}
Recall that $\mu'=m \mu$. $[K]=c \mu \pm f \mu' = (c \pm fm) \mu$. As $[K]$ has order $q$, $q | c \pm fm$, i.e.,
\[0 \equiv c \pm fm \equiv c \mp \epsilon fd \pmod{q}.\]
Note that $c\equiv \pm \epsilon fd\ (\mod q)$ where the sign respects that of $[K] = c\mu \pm f\mu'$.

{\bf Assume $(1)$ or $(4)$ holds.} The proof is then identical to that of Theorem~\ref{thm:c+d:dq+1}.

{\bf Assume $(2)$ holds.} 
As $c\equiv \pm \epsilon fd\ (\mod q)$, $cn\equiv \pm \epsilon fnd \equiv \mp f\ (q)$. 
Thus we may write $f = \eta q\mp cn$, where the sign respects that of $[K] = c\mu \pm f\mu'$.
Then, $[K] = (c\pm fm)\mu = (c\pm(\eta q\mp cn)d(q - \epsilon))\mu = (c+(\pm \eta dq-c(q-\epsilon))(q-\epsilon))\mu = \epsilon(2c\mp \eta d)q\mu$, using that $nd=q-\epsilon$. When $f=cn$ we have $\eta=0$ and the sign is $+$ (that is, we get $[K]=\epsilon(2c + \eta d)q\mu$ above), so we have $[K]=\epsilon (2c)q\mu$ giving conclusion $(d)$. When $f=2q-cn$, $\eta=2$ with sign $-$, and $[K]=\epsilon (2c-2d)q\mu$---giving conclusion $(a)$.

First assume $q>cn$. As $0\le f < 2q$, we conclude $f = cn$, $q-cn$, $q+cn$, or $2q-cn$. We have two additional cases to the above. For $f=q \pm cn$, with $\eta=1$ and sign $\pm$ (resp.), we conclude $[K]=\epsilon (2c \pm d)q\mu$ (resp.)---giving conclusions $(b),(c)$.

%As $0\le f < 2q$, we conclude $f = cn$, $q-cn$, $q+cn$, or $2q-cn$, denoted as $f = \eta q + \eta' cn$ for the corresponding $(\eta, \eta')$. Thus, $[K] = (c\pm fm)\mu = (c\pm(\eta q +\eta' cn)d(q - \epsilon))\mu = (c \pm (\eta dq(q-\epsilon) + \eta'c(q-\epsilon)^2) \mu = (c \pm  (-\epsilon q(\eta d + 2 \eta' c) + \eta' c)) \mu$, using that $nd=q-\epsilon$. As $c \pm fm 
% \equiv 0 \, \mod (q)$, we conclude that either 
 
% \item $\eta'=-1$ and $[K] = \epsilon q(2c- \eta d) \mu$ giving conclusions $(a),(c)$; or
 
% \item $\eta'=1$ and $[K] = \epsilon q(2c+ \eta d) \mu$ giving conclusions $(b),(d)$.

Now assume $nd+\epsilon = q<cn\le 3n$. Thus, $c=d=3$ and $\epsilon=-1$ as $d\ge 3$. From $0\le f<2q$ and $q<cn<3q/2$, we conclude that $f=-q+cn$, $2q-cn$, $cn$, or $3q-cn$, giving two additional cases. When $f=-q+cn = 1$, then $\eta=-1$ and the sign is $+$. Then $-[K]= 3q\mu$ as described in conclusion $(e)$. When $f=3q-cn$, we similarly conclude $[K] = 3q\mu$ again giving conclusion $(e)$.

{\bf Assume $(3)$ holds}. As $c \equiv \pm \epsilon fd\ (q)$, we write $c = \eta q\pm \epsilon fd$, where the sign respects that of $[K] = c\mu \pm f\mu'$. Thus, $[K] = (c\pm fm)\mu = (\eta q\pm \epsilon fd\pm fd(q - \epsilon))\mu = (\eta\pm fd)q\mu=(\eta\pm \epsilon^2 fd)q\mu$. If $c=fd$,
we have $\eta=0$ and $\pm\epsilon=1$, and we conclude that $[K]=\pm \epsilon^2 fdq \mu=\epsilon fdq$ which is part of conclusion $(b)$. When $c=2q-fd$, $\eta=2,\, \pm \epsilon=-1$, and we have $[K]=(2 \pm \epsilon^2 fd)q\mu=(2-\epsilon fd)q\mu$ as in conclusion $(a)$.

First assume $q > fd$. Then, as $c \equiv \pm \epsilon fd\ (q)$ and $0\le c < 2q$, we conclude that $c = fd$, $q - fd$, $q + fd$, or $2q - fd$. When $c=q \pm fd$, we have $\eta=1, \pm \epsilon=\pm1$, and $[K]=(1 \pm \epsilon^2 fd)q\mu=(1 \pm \epsilon fd)q\mu$ as in conclusions $(a),(b)$.

%First assume $q > fd$. Then, as $c \equiv \pm \epsilon fd\ (q)$ and $0\le c < 2q$, we conclude that $c = fd$, $q - fd$, $q + fd$, or $2q - fd$, denoted as $c = \eta q\pm \epsilon fd$, where the sign respects that of  $[K] = c\mu \pm f\mu'$. Thus, $[K] = (c\pm fm)\mu = (\eta q\pm \epsilon fd\pm fd(q - \epsilon))\mu = (\eta\pm fd)q\mu$. If $c = q-fd$ or $2q-fd$, we obtain (a); otherwise, $c = fd$ or $q+fd$, and we obtain (b).

Now assume $nd+\epsilon = q < fd\le 3d$. Then, either both $n=2$ and $f=3$, or $\epsilon = -1$ and $n=f=2,3$. From $q<fd<2q$, we conclude that $c=-q+fd$, $2q-fd$, $fd$, or $3q-fd$. We have two additional cases to consider.

Consider the case where $n = 2$ and $f = 3$, $q = 2d+\epsilon$, and $fd = 3d = q+(d-\epsilon) \equiv d-\epsilon\ (\mod q)$. Assume $c=-q+fd=d-\epsilon$. Taking $\eta=-1, \pm \epsilon=1$, we get $[K]=(-1 \pm \epsilon^2 fd)q\mu=(-1+\epsilon fd)q\mu=(-1+\epsilon (d-\epsilon))q\mu = (\epsilon d - 2)q\mu$. If $c=3q-fd=2q-(d-\epsilon)$, we take $\eta=3, \pm \epsilon=-1$ and get $[K]=(3 \pm \epsilon^2 fd)q\mu=(3-\epsilon fd)q\mu=(3-\epsilon (d-\epsilon))q\mu = (4 - \epsilon d)q\mu$. These are the conclusions of $(c)$.

Finally, when $\epsilon=-1$ and $n=f=2, 3$, $fd=nd=q+1\equiv 1\ (\mod q)$. If $c=-q+fd=1$, we again conclude $[K]=(-1 + \epsilon fd)q\mu = -2q\mu$.  If $c=3q-fd = 2q-1$, we again have $[K]=(3-\epsilon fd)q\mu=4q\mu$. These are the conclusions of $(d)$.
\end{proof}

The next two lemmas determine the positive integers $y,\, y'$ associated to special classes of $[K]$.

\begin{lem}\label{lem:dne2an>2}
Let $27 < q= nd+\epsilon$. Let $r,s,e$ be as in Lemma~\ref{thm:=x}.
\begin{enumerate}[label=({\Roman*})]
    \item Let $\pm [K]=(d - \zeta a)q \mu$ for $a \in \{1,2,3,4\}$, $\zeta=\pm 1$. 
    One of the following holds: 
    \begin{enumerate}
        \item $a=4, \,d=3, \,\zeta=1$. Then $\pm [K]=q \mu, \,y=1, \,y'=n$. 
        \item  $n=2$. If $\zeta=1$, then $y=d-a, \,y'=2a+\epsilon$. If $\zeta=-1$, then $y=d-(a-\epsilon), \,y'=2a-\epsilon$.
        \item  $n>2$. Then $y=d - \zeta a$. If $2a < d$, then $y'=an + \zeta\epsilon$. When $2a \geq d$ we have the following:
        \begin{enumerate}
            \item $d=7, \,a=4, \,y'=3n+\epsilon(1-\zeta)$.
            \item $d=5, \,a=3,4, \,y'=(5-a)n+\epsilon(1-\zeta)$.
            \item $d=3, \,a=2,3,4$. When $a=2$, $y'= n+\epsilon(1-\zeta)$. 
            When $a=3$, $\zeta=-1$, $y'=2$. When $a=4$, $\zeta=-1, \,y'=n-2\epsilon$.
        \end{enumerate}
    \end{enumerate}
    \item Let $\pm [K]=aq \mu$ where $a \in \{2,3,4\}$. Then $y=a$.
    \begin{enumerate}
        \item If $d>7$, then $y'=an$.
        \item If $d=7$, then $y'=an$ for $a<4$.  When $a=4$, $y'=3n+\epsilon$.
        \item If $d=5$, $y'=2n$ if $a=2$. For $a=3,4$, $y'=(5-a)n+\epsilon$.
        \item If $d=3$, $y'=n+\epsilon$ for $a=2$. When $a=3$, $y'=1$. When $a=4$, $y'=n-\epsilon$.  
    \end{enumerate}
\end{enumerate}

\end{lem}

\begin{proof}
Recall that $[K]= \pm yq \mu$, $[K]=\pm y'q \mu'$ and $m=d(q-\epsilon)$; furthermore, $\mu=m^{-1}\mu' \,\, (q^2)$. As $m \equiv -\epsilon d$, $m^{-1} \equiv n  \,\,\, (q)$. Then,
\[\pm y \equiv y'm \equiv -\epsilon y'd \ (q), \quad \pm y' \equiv ym^{-1} \equiv yn \ (q).\]
Finally, in determining the values of $y,y'$ associated to $K$, we use the fact that they are determined by the properties
\[0<y,y'<q/2 \,\,\hbox{and} \,\, [K]=\pm yq\mu,\,\, [K]=\pm y'q\mu.\]
In particular, $y,y'$ are determined by their mod $q$ class up to sign and the fact that $0<y,y'<q/2$.

\begin{enumerate}[label=({\Roman*})]
\item {\bf Consider $\pm [K]=(d - \zeta a)q \mu$  for $a \in \{1,2,3,4\}, \zeta=\pm 1$. }

\smallskip

Then $\pm y \equiv d - \zeta a\,\,\, (q)$. First assume $\zeta=1$ and $d\le a$. If $d=3$ and $a=4$, then $\pm [K]=q \mu$, $y=1$ and $y'=n$, which is (a). On the other hand, $d=3$ and $a=3$ contradicts that $[K]$ has order $q$. We hereafter assume that if $\zeta=1$, then $d>a$. 

Assume $n=2$. Since $q = 2d+\epsilon = (d-\zeta a) + (d+\zeta a+\epsilon)$, $\pm [K] = (d - \zeta a)q \mu = \pm y q \mu$, we have $y \equiv \pm (d-\zeta a) \ (q)$ which along with  $0<y<q/2$ implies $y = \min\bc{d-\zeta a, d + \zeta a + \epsilon}$. Conclusion (b) then follows from $\pm y' \equiv yn\ (q)$ and $0< y' < q/2$.

Assume $n>2$. Since $q > 27$, $q/2-q/n\ge q/6 > 13/3$. Thus, $0 < d - \zeta a \leq (q-\epsilon)/n+4 < q/2$ and we conclude $y = d - \zeta a$ as claimed. 
Then $\pm y' \equiv -\zeta ny = -\zeta nd + an \equiv an + \zeta\epsilon\ (q)$. Thus, if $2a < d$, $0 < 2(an+\zeta\epsilon) \le n(d-1)+2\zeta\epsilon < q$, so $y' = an + \zeta\epsilon$ as desired. We hereafter assume $2a \geq d$. Then $d \leq 7$. 

Again, note that $\pm y' \equiv an + \zeta\epsilon\ (q)$ and $-(an + \zeta\epsilon) \equiv q - (an+\zeta\epsilon) = (d-a)n + \epsilon(1 - \zeta) > 0$, unless $d = 3$ and $a = 3, 4$. Moreover, $q-(an+\zeta\epsilon) \le q - ((d+1)n/2+\zeta\epsilon) = q/2 - ((n-\epsilon)/2+\zeta\epsilon) < q/2$. Thus, the conclusions follow.

Finally, assume $d = 3$ and $a = 3, 4$. Recall that in this case $\zeta = -1$. If $a = 3$, then $\mp y'\equiv 2\epsilon\ (\mod q)$ which implies $y' = 2$. If $a = 4$, $\pm y' \equiv 4n-\epsilon \equiv n-2\epsilon < q/2$, since $q = 3n+\epsilon$. Thus $y'=n - 2\epsilon$. This concludes the proof of (c).

% If $d=-7$, then $a=4$. $1+ \zeta 4n=1+\zeta [7n+1-1-3n] \equiv 1- \zeta (1+3n)  \,\, (q)$. This means that $y'=3n+1-\zeta$, as desired.  

% If $d=-5$, then $a=3,4$. $1+ \zeta an=1+\zeta [5n+1-1-(5-a)n] \equiv 1- \zeta [1+(5-a)n]$. Hence $y'=(5-a)n+1-\zeta$. 

% $d=-3$, then $a=2,3,4$. $1+ \zeta an=1+\zeta [3n-1+1-(3-a)n] \equiv 1+ \zeta [-1-(3-a)n]=1-\zeta[1+(3-a)n]$. That is, $\pm y' \equiv \zeta -[1+(3-a)n]$.  When $a=2, y'= n+1-\zeta$. Recall that when $d=-3,a=3,4$, we are assuming $\zeta=-1$. Thus, when $a=4, y'=\zeta-1+n=n-2$. When $a=3, y' \equiv \zeta-1$, and we conclude that $y'=2$.

\smallskip

\item {\bf Consider $\pm [K]=aq \mu$ where $a \in \{2,3,4\}$.} 

\medskip

Then $\pm y \equiv a\ (q)$. As $0<a<q/2$, $y=a$. Then $\pm y' \equiv yn = an$. If $d>7$, then $q = nd + \epsilon\ge 9n+\epsilon > 2an$, so $y' = an < q/2$.

In the case where $d\le 7$, $y' = an$, $q-an$, or $an - q$, defined by $0 < y' < q/2$. Then, analyzing by cases gives all the conclusions. 
% If $d=-7$, the $q=1+7n$ and $y' = an $ for $a<4$. When $a=4$, $\pm y' \equiv -4n+1+7n=1+3n$ and we conclude that $y'=1+3n$. If $d=-5$ then $q=1+5n$. Then $y'=2n$ if $a=2$. For $a=3,4$, $\pm y' \equiv -an +(1+5n) = (5-a)n+1$ implies that $y'=(5-a)n+1$. If $d=-3, q=1+3n$, then $\pm y' \equiv -an \equiv (3-a)n+1$. If $a<4$ then $y'=(3-a)n+1$. If $a=4$, then $y'=n-1$. 
\end{enumerate}
\end{proof}

\begin{lem}\label{lem:2dpmq:1-nd}
Assume $28<q=nd+\epsilon$ and $\pm [K]\!=\!(-2d+a)q \mu\!=\!(2\epsilon+an)q \mu'$ where $a\!\in\!\{-1,0,1,2,3,4\}$. Then the following holds for $y,\, y'$:
\begin{itemize}
    \item
        \begin{enumerate}
            \item If $n >4$ or both $n=4$ and $2a+\epsilon \ge 0$, then $y=2d-a$.
            \item If $n=2,\, a+\epsilon<0$, then $y=-a-\epsilon$.
            \item Otherwise, $y=(n-2)d+a+\epsilon$.
        \end{enumerate}
    % $\spadesuit$ (When $n=2$, need $a+\epsilon >0$. In particular when $a=-1,\epsilon=-1$, then $y=2$.)
    \item If $a=-1$, then $y'=n-2\epsilon$.
    \item If $a=0$, then $y'=2$.

    \item For $a \in \{1,2,3,4\}$, we have the following:
    \begin{enumerate}
        \item If $d \geq 2a$, then $y'=an+2\epsilon$.
        \item If $2a > d >a$, then $y'=(d-a)n-\epsilon$.
        \item If $d \leq a$, then either $d=3, \,a=3, \,y'=1$, or $d=3, \,a=4, \,y'=n+\epsilon$.
    \end{enumerate}
\end{itemize}
\end{lem}

\begin{proof}
If $n >4$ or both $n=4$ and $2a+\epsilon \ge 0$, then $0<2d-a < \pth{nd+\epsilon}/{2}={q}/{2}$. Then we conclude that $y=2d-a$. Otherwise, write $\pm [K]=(-2d+a+(nd+\epsilon))q \mu=((n-2)d+a+\epsilon)q \mu$. When $n=4, 2a+\epsilon < 0, q>28$, we have $0<(n-2)d+a+\epsilon<q/2$ and conclude $y=(n-2)d+a+\epsilon$. So assume $n \leq 3$. Then $3d +\epsilon \ge nd + \epsilon=q>28$ implies that $(n-2)d+a+\epsilon<\pth{nd+\epsilon}/{2}={q}/{2}$. If $n=3$ or if $n=2$ and $a+\epsilon>0$, we conclude $y= (n-2)d+a+\epsilon$ as claimed. Finally, if $n=2$ and $a+\epsilon \leq 0$, then $y=-a-\epsilon$. Note that when $n=2$, $a+\epsilon \neq 0$ since $[K]$ has order $q$.

Rewrite $\pm[K] = (2\epsilon+an)q \mu'$.  If $a=-1$, then $0 \leq y' =n-2\epsilon<{q}/{2}$. Note that $n \neq 2$ when $\epsilon = 1$, as $[K]$ has order $q$. When $a=0$, we conclude $y'=2$. We hereafter assume that $a \neq -1,0$ and have the following cases:

\begin{enumerate}
\item Assume $d \geq 2a+1$. As $n\neq 2$ when $\epsilon=1$, we have that $2(an+2\epsilon)\le (d-1)n+4\epsilon = q-(n-3\epsilon) \leq q$. We conclude that $y'=an+2$, since $an+2\epsilon >0$. ($an+2\epsilon \neq 0$ as [K] has order $q$).

\item Assume $2a > d >a$. Since $q > 28$, we have $n\ge 4$. Rewrite $\pm[K]=(-(2\epsilon+an)+(nd+\epsilon))q \mu'=((d-a)n-\epsilon)q \mu'$. Note that $0<(d-a)n-\epsilon\le \pth{dn+\epsilon+(-n-3\epsilon)}/{2}<{q}/{2}$, so we conclude that $y'=(d-a)n-\epsilon$.
%$(2d-2a)n-2\epsilon-dn-\epsilon \le -n-3\epsilon < 0$ unless $n = 2$, $\epsilon = -1, \,d = 2a-1$, contradicting $q > 28$.
\item If $d \leq a$, then $d=3$ and $a=3,4$. Recall that $\pm y'\equiv (d-a)n-\epsilon$. If $d=3$ and $a=3$, we conclude that $y'=1$. If $d=3$ and $a=4$, we conclude that $y'=n+\epsilon$. 
\end{enumerate}
\end{proof}

Here we have a variation of Theorem~\ref{thm:c+d} when sharper restrictions can be placed on the decomposition of $[K]$.

\begin{cor}\label{cor:c+d:q| d-1<4/3 q}
Assume $L(q^2,m)$ as in case $(*q)$. Assume further that $[K]=c \mu \pm f \mu'$ and $nd+\epsilon=q>18$. 

\begin{enumerate}
\item If $0 \leq c \leq 1$ and $0 \leq  f <q/2$, then $\pm [K]=2q\mu$ and $c=1, f=n$.
\item If $0 \leq f \leq 1$ and $0 \leq  c <q/2$, then $[K]=\epsilon dq \mu$.
\end{enumerate}

\end{cor}

%\begin{cor}\label{cor:c+d:q| d-1<4/3 q}
%Assume $[K]=c \mu \pm f \mu'$, $K$ is spherical, and $nd+\epsilon=q>18$. 

%\item $0 \leq c \leq 1$ and $0 \leq  f <q/2$; or 
%\item $0 \leq f \leq 1$ and $0 \leq  c <q/2$ 
%\end{enumerate}
%then $$\pm [K]=2q\mu , dq \mu$$
%\end{cor}

\begin{proof} \hfill

\begin{enumerate}
\item $0 \leq c \leq 1$ and $0 \leq  f <q/2$. Note that $c=1$, otherwise the fact that $f<q$ contradicts Theorem~\ref{thm:c+d}$(1)$. We apply 
Theorem~\ref{thm:c+d}$(2)$ and consider conclusions $(a)$--$(d)$ ($(e)$ is not relevant). As $f<q/2$ and $q=nd+\epsilon > 2n = 2cn$, conclusions $(a)$, $(b)$, $(c)$ cannot hold. We are left with conclusion $(d)$, which verifies the Corollary. 
\item $0 \leq f \leq 1$ and $0 \leq  c <q/2$. We assume $f=1$, otherwise $c<q$ contradicts Theorem~\ref{thm:c+d}$(4)$. We apply Theorem~\ref{thm:c+d}$(3)$ and consider conclusion $(a)$--$(d)$. If conclusion $(b)$ holds, then $c<q/2$ implies that $\eta=0$, $c=d$ and $[K]=\epsilon dq\mu$. Conclusions $(c)$, $(d)$ do not apply as $f=1$. 

Thus we may assume conclusion $(a)$ of Theorem~\ref{thm:c+d}$(3)$ holds. Then $c<q/2$ implies that $\eta=1, \,n=2, \,\epsilon=-1, \,c=d-1$, and $[K]=(1+d)q\mu=-(d-2)q\mu$. By Lemma~\ref{lem:dne2an>2}(I)(b) with $\zeta=1,a=2$, we conclude that $y=d-2, \,y'=3$. As $q>18$, $q=2d-1<3d-6=3y$. Then $3y>q=(r+s)y+e$, implies $r+s=2$. Then $2d-1=q=2y+e=2(d-2)+e$ implies $e=3$. But now Theorem~\ref{lem:q<4y} gives a sharper decomposition of $[K]$. Theorem~\ref{lem:q<4y}(1) says that for $j=q+e-1$, $m_j\leq 1$. Then $j'=2q-2-j=q-1-e=q-4<q-y'=x'$. By Theorem~\ref{thm:M:b<x+1},  $M_{j'}' < {q}/{y'} = {q}/{3}$. Thus we may write $\pm [K]=c' \mu \pm f' \mu'$ where $0 \leq c' <q/3, \,0 \leq f'=m_b \leq 1$. Applying Theorem~\ref{thm:c+d}$(3)$ to $c', \,f'$ we conclude as above that either $[K]=\epsilon dq\mu$ or $c'=d-1, \epsilon=-1$. Now $d-1=c'<q/3=(2d-1)/3$ would contradict that $d \geq 3$. We conclude $\pm [K]=\epsilon dq\mu$ as claimed.
\end{enumerate} 
\end{proof}

One of the conclusions of the Engine Theorem~\ref{thm:q+y,m_b=1} is that $q<4y$ or $q<4y'$. The next two theorems give further estimates on $y,y'$ in these contexts. 

\begin{thm}\label{thm:q<4y,e+1:1-nd} Let $e$ be as in Theorem~\ref{thm:=x}. Assume $26<q<4y$. Then either $\pm [K] = dq \mu$
or $y'>e+1$
\end{thm}
\begin{proof}
Let $r,s,e$ be as in Theorem~\ref{thm:=x}. Assume for contradiction that $y' \leq e+1$. Note that $y > q/4 > 3$. Furthermore, if $y' = 1$ then $\pm [K]=dq\mu$, so we may also assume $y'>1$. 

By Lemma~\ref{lem:q<4y}(1), $m_b \leq 2$ when $b=q+e-1$. Let $b'=2q-2-b = q-e-1\leq q-y'=x'$. By Theorem~\ref{thm:M:b<x+1},  $M_{b'}' < {q}/{y'} \leq {q}/{2}$. Thus we may write $\pm [K]=c \mu \pm f \mu'$ where $0 \leq c <q/2,\, 0 \leq f=m_b \leq 2$. By Corollary~\ref{cor:c+d:q| d-1<4/3 q}, we may assume $f=m_b=2$. Note then that Lemma~\ref{lem:q<4y} shows that $r+s = 3$; hence $q=3y+e$.  As $[K]$ has order $q$, $c \neq 0$. Theorem~\ref{thm:c+d}(3) applies and we consider, in turn, conclusions $(a)$--$(d)$. As $f=2$ and $y>3$ conclusions $(c)$, $(d)$ do not hold.

{\bf Assume $Theorem~\ref{thm:c+d}(3)(a)$ holds: $[K]=(\eta-2\epsilon d)q \mu$ where $\eta \in  \{1,2\}$ and $c=\eta q-2d$. } 

As $n \geq 2$, $c <q/2$ implies $\eta = 1$, and either $n = 2, 3$, or $n=4 \, \hbox{and} \, \epsilon=-1$. Then $[K]=(1-2\epsilon d)q \mu$. 

If $n=2$, then $[K]=(1-2\epsilon d +\epsilon(2d + \epsilon))q \mu=2q \mu$, contradicting $y>3$.

Assume $n=3$. Then, $[K]=(1-2\epsilon d+\epsilon(3d+\epsilon))q \mu=\epsilon(d+2\epsilon)q \mu$. As $0<2(d+2\epsilon)<3d+\epsilon = q$ ($q>10$), we conclude $y=d+2\epsilon$. Then, $0<e=q-3y=(3d+\epsilon)-3(d+2\epsilon) = -5\epsilon $, so $\epsilon = -1, \,e = 5, \,y=d-2$.
Thus, $\pm y' = yn = nd-2n\equiv 1-2n = -5\ (\mod q)$, so $y' = 5$.
Note that $y - e = (d-2)-5 = d-7 \ge 4$, as $q=nd-1>26$, $n=3$, and $d$ is odd. By Lemma~\ref{lem:q<4y}(1), for $z = q+ (y-e)-1$, $m_z\leq 1$. From $y-e\ge 4$, $z \geq q+3$, so $z'=2q-2-z \leq q-5=q-y'=x'$. By Theorem~\ref{thm:M:b<x+1}, $M'_{z'}< q/5 < q/2$. Thus we can represent $\pm [K]=c' \mu \pm f' \mu'$, where $0 \leq c' < q/2, \,0 \leq f' \leq 1$. Corollary~\ref{cor:c+d:q| d-1<4/3 q} now proves the Theorem.

Finally, assume $n=4, \epsilon=-1$. $[K]=(1+2d -(4d-1))q \mu=-(2d-2)q \mu$. As $0<2d-2<q/2$, we conclude that $y=2d-2$. But $4d-1=q=3y+e=6d-6+e$, contradicting that $e>0$.

{\bf Assume $Theorem~\ref{thm:c+d}(3)(b)$ holds.} As $c<q/2$ and $f=2$, we conclude that $[K]=2\epsilon dq \mu=2\epsilon ndq \mu'=-2q \mu'$. As $q>4$, we conclude that $y'=2$. By Lemma~\ref{lem:q<4y}(1), for $j =q$, $m_j \leq 1$. Then $j'=2q-2-j = q-2=q-y'=x'$. By Theorem~\ref{thm:M:b<x+1}, $M_{j'}'<q/2$. Thus we can represent $\pm [K]=c' \mu \pm f' \mu'$ where $0 \leq c' < q/2, \, 0 \leq f' \leq 1$. Corollary~\ref{cor:c+d:q| d-1<4/3 q} now proves the Theorem. 
\end{proof}

\begin{thm}\label{thm:q<4y',e'+1:1-nd} Let $e$ be as in Theorem~\ref{thm:=x}. Assume $31<q<4y'$. Then either $\pm [K] = q \mu, 2q \mu$ or $y>e'+1$. 
\end{thm}

\begin{proof}
Let $r',s',e'$ be as in Theorem~\ref{thm:=x}. Assume for contradiction that $y \leq e'+1$. %In particular, $\pm [K] \neq 2q \mu'=-2dq \mu$ as otherwise $q>4$ would imply that $y'=2$.  
Furthermore, if $y = 1$, then $\pm [K]=q \mu$, so we may also assume $y>1$. 

By Lemma~\ref{lem:q<4y}(2), $m_b' \leq 2$ when $b=q+e'-1$. Let $b'=2q-2-b = q-e'-1\leq q-y=x$. By Theorem~\ref{thm:M:b<x+1},  $M_{b'} < {q}/{y} \leq {q}/{2}$. Thus we may write $\pm [K]=c \mu \pm f \mu'$ where $0 \leq f <q/2,\, 0 \leq c=m_b' \leq 2$. By Corollary~\ref{cor:c+d:q| d-1<4/3 q}, we may assume $c=m_b=2$. Note then that Lemma~\ref{lem:q<4y} shows that $r'+s' = 3$; hence, $3y'+e'=q>3y'$. As $[K]$ has order $q$, $f \neq 0$. Theorem~\ref{thm:c+d}(2) applies and we consider, in turn, conclusions $(a)$--$(e)$. Note that conclusions $(a)$, $(b)$ violate the assumption that $f < q/2$, and $(e)$ that $c=2$. 

{\bf Assume $Theorem~\ref{thm:c+d}(2)(c)$ holds: $\epsilon [K]=(4-d)q \mu$ and $f = q-2n$. } 

As $d \ge 3$, $f <q/2$ implies that $d=3$. Then $\pm[K]=q \mu$ as claimed.

{\bf Assume $Theorem~\ref{thm:c+d}(2)(d)$ holds.} 

Then, $\epsilon[K]=4q \mu$. As $q > 8$, we conclude $y=4$. Note that since $[K]$ has order $q$, we obtain that $q$ is odd. As $d$ is odd and $q=nd+\epsilon$, $n$ must be even.

\begin{clm} $y'-e' \geq 3 = y-1$. \label{clm:min>2}
\end{clm}

\noindent {\em Proof of Claim:}

Assume that $d>7$. Then by Lemma~\ref{lem:dne2an>2}(II)$(a)$ with $a=4$, we conclude that $y'=4n$. As $q<4y'$, we have that $nd+\epsilon<16n$ and conclude that $d < 16$. The equation $q=3y'+e'$ becomes $nd+\epsilon=12n+e'$, or $e'=(d-12)n+\epsilon$. As $e'>0$ and $d$ odd, we conclude that $d=13, 15$. Thus $y'-e'=4n-((d-12)n+\epsilon)=(16-d)n-\epsilon \geq n-\epsilon$. As $q=nd+\epsilon>31$ and $n$ is even, we conclude that $y'-e' \geq 3$ as claimed. 

Assume $d=7$. Lemma~\ref{lem:dne2an>2}(II)(b) says that $y'=3n+\epsilon$, but this contradicts that $7n+\epsilon=q>3y'$. 

Assume $d=5$. Lemma~\ref{lem:dne2an>2}(II)(c) says that $y'=n+\epsilon$, but this contradicts that $4y'>q=5n+\epsilon$ and $q > 11$.

Finally, assume $d=3$. Lemma~\ref{lem:dne2an>2}(II)(d) says that $y'=n-\epsilon$. Then $3n+\epsilon=q=3(n-\epsilon)+e'$, and we conclude that $\epsilon=1$ and $e'=4$. Thus, $y'-e'=n-5 \geq 3$, as $q > 19$ and $n$ is even. 
\qed {\em\,(Claim)}

\medskip
By Lemma~\ref{lem:q<4y}(2), for $z = q+ (y'-e')-1$, $m_z' \leq 1$. Claim~\ref{clm:min>2} says that $z \geq q+y-2$. Then $z'=2q-2-z \leq q-y=x$. By Theorem~\ref{thm:M:b<x+1}, $M_{z'}< q/4 < q/2$. Thus we can represent $\pm [K]=c' \mu \pm f' \mu'$, where $0 \leq f' < q/2, \,0 \leq c' \leq 1$. Corollary~\ref{cor:c+d:q| d-1<4/3 q} now proves the Theorem.
\end{proof}

\subsubsection{Case $(*2q)$}\label{subsubsec:gencom2q}

For this subsubsection, we assume that we are in case $(*2q)$ where
\[m=d(q-\epsilon),  \qquad \epsilon=\pm1, \qquad 2q= nd-\epsilon, \qquad n \ge 3,\,\, d\geq 3 \,\, \hbox{both odd}.\]

\bigskip

\begin{thm}\label{thm:c+d:2q+1}
Assume $K \subset L(q^2,m)$ for $m$ as in $(*2q)$.
Assume $[K]=c \mu \pm f \mu' \in H_1(L(q^2,m))$. Furthermore, assume $q>18$ and one of the following holds on $c,f$:

\begin{enumerate}
\item $c=0$ and $0 \leq f < q^2/2$;
\item $0 < c \leq 3$ and $0 \leq f < 2q$;  
\item $0 < f \leq 3$ and $0 \leq c < 2q$; 
\item $f=0$ and $0 \leq c < q^2/2$.
\end{enumerate}
Then, we have the corresponding conclusions:

\begin{enumerate}
\item $[K]= \pm f \mu' = \pm\epsilon fd \mu$ where $f= y' q$;
\item Either
\begin{enumerate}
\item $\epsilon[K]=-(\eta d+c)q \mu=\epsilon(\eta+\epsilon cn)q \mu'$ and $f=\eta q + cn$ with $\eta=0, 1$; % old (b), (d)
\item $\epsilon[K]=(\eta d-c)q\mu=-\epsilon(\eta-\epsilon cn)q \mu'$ and $f=\eta q-cn$ with $\eta=1,2,3,4$; % old (e)
\item $\epsilon[K]=2q \mu$ and $c=3, \,d=5, \,f=\pth{n+\epsilon}/{2}$; % old (c)
% \item $\epsilon[K]=dq \mu$, with $c=d=3$ and $\epsilon=1$.
\item $\epsilon[K]=q \mu$ and $c=2, \,d=3, \,f=\pth{n+\epsilon}/{2}$; % This was missed
\item $\epsilon[K]=dq \mu$ and $c=3, \,d=3, \,\epsilon=f=1$. % old (a)
\end{enumerate}

\item Either
\begin{enumerate}
\item $[K]=(\eta +\epsilon fd)q \mu=-(\eta n+f)q \mu'$ and $c=\eta q+fd$ where $\eta=0, 1$;
\item $[K]=(\eta-\epsilon fd)q \mu=-(\eta n-f)q \mu'$ and $c=\eta q-fd$ where $\eta = 1,2,3,4$;
\item $[K]=-\epsilon 2dq \mu=2q \mu'$ and $f=3, \,n=5, \,c=\pth{d+\epsilon}/{2}$;
\item $[K]=-\epsilon dq \mu=q \mu'$ and $f=2, \,n=3, \,c=\pth{d+\epsilon}/{2}$;
\item $[K]=-\epsilon q\mu=nq \mu'$ and $f=3, \,n=3, \,\epsilon=c=1$.
\end{enumerate}

\item $[K]= c \mu$ and $c= yq$. 

\end{enumerate}

%In either case $(2),(3)$, we conclude that either $\pm [K]=q \mu, dq \mu$ or $\pm [K]=\nu d \pm \lambda \mu$ for some $\nu \in \{1,2\}$ and for some $\lambda \in \{0, \pm 1, \pm 2\}$. 
\end{thm}

\begin{proof}
Note that $m = d(q-\epsilon) \equiv -\epsilon d\ (q)$. Since $nd=2q+\epsilon$, $nd \equiv \epsilon \ (q)$ hence $m^{-1} \equiv -n \ (q)$. Recall that $\mu'=m \mu$. $[K]=c \mu \pm f \mu' = (c \pm fm) \mu$. As $[K]$ has order $q$, $q | c \pm fm$, i.e., 
\[0 \equiv c \pm fm \equiv c \mp \epsilon fd \ \mod (q).\]

{\bf Assume $(1)$ or $(4)$ holds.} The proof is then identical to that of Theorem~\ref{thm:c+d}.

{\bf Assume $(2)$ holds.}
Since $c\equiv \pm \epsilon fd\ (\mod q)$, $cn\equiv \pm \epsilon fnd \equiv \pm f\ (q)$. Thus, we write $f = \eta q\pm cn$ (where the sign respects that of $[K] =(c\pm fm)\mu$). Then $\epsilon [K] = \epsilon(c\pm fm)\mu = \epsilon(c\pm (\eta q\pm cn)d(q-\epsilon))\mu = \epsilon(c+ (\pm\eta dq+ c(2q+\epsilon))(q-\epsilon))\mu = (\mp\eta d- c)q\mu$.

First, assume $q > cn$. As $0\le f = \eta q\pm cn < 2q$, we conclude $f = cn, q-cn, q+cn, 2q-cn$. From $\epsilon [K] = (\mp\eta d- c)q\mu$, we obtain (a) when $f=cn, q+cn$ and (b) when $f=q-cn, 2q-cn$.

Now assume $cn>q$. Together with $2q=nd-\epsilon$ and $n,d\ge 3$, $c\le 3$, we conclude that either (i) $d=5, \,c=3$; (ii) $d=3, \,c=2$; (iii) $d=c=3, \,\epsilon = -1$; or (iv) $d=c=3, \,\epsilon = 1$. Note that in (i)--(iii), $q<cn<2q$; in (iv), $2q<cn<3q$. We conclude that $f=2q-cn$, $cn$ (in (i)--(iii)), $-q+cn$, $3q-cn$ (in (i)--(iv)), or $4q-cn$, $-2q+cn$ (in (iv)). Writing these options as $f=\eta q \pm cn$ for the appropriate values of $\eta$ and sign, we get as above, noting that the sign of $c\pm fm$ and of $f=\eta q \pm cn$ agree, that $\epsilon [K] = \epsilon(c\pm fm)\mu =  (\mp\eta d- c)q\mu$. Then $f=cn$ corresponds to $\eta=0$ and lands in conclusion (a). Similarly, $f=2q-cn, 3q-cn, 4q-cn$ are in conclusion (b) with $\eta=2,3,4$. We need only consider $f=-q+cn$ in (i)--(iv) and $-2q+cn$ in (iv).

Assume (i) holds. Then, $2q=5n-\epsilon$. As $f = -q+cn = (n+\epsilon)/2$, $\epsilon[K] = (d-c)q\mu = 2q\mu$, which is (c).
% if $f=3q-cn= 4n-3\epsilon/2$, $[K]=\epsilon(3d-c)q\mu = 12q\mu$.

Assume (ii) holds. Then, $2q=3n-\epsilon$. As $f = -q+cn = (n+\epsilon)/2$, $\epsilon[K] = (d-c)q\mu = q\mu$, which is (d).

Assume (iii) or (iv) holds and $f = -q+cn$. Then $\epsilon[K] = (d-c)q\mu = 0$, contradicting that $[K]$ has $q$ torsion.

Finally, assume (iv) holds and $f = -2q+cn$. Then $\epsilon[K] = (2d-c)q\mu = dq\mu$, which is (e).

{\bf Assume $(3)$ holds}.

As $c \equiv \pm \epsilon fd\ (q)$, we write $c = \eta q\pm \epsilon fd$, where the sign respects that of $[K]=c\mu\pm f\mu'$ and $\eta$ an integer. Thus, $[K] = (c\pm fm)\mu = (\eta q\pm \epsilon fd\pm fd(q - \epsilon))\mu = (\eta\pm fd)q\mu$.

First assume $fd < q$. Then $0\le c = \eta q\pm \epsilon fd < 2q$ implies $c = fd$, $q - fd$, $q + fd$, or $2q - fd$. If $c = fd$ or $q+fd$ (hence $\pm 1=\epsilon$), and we obtain (a); otherwise, $c = q-fd$ or $2q-fd$ (hence $\pm 1=-\epsilon$), and we obtain (b).

Next assume $fd > q$. Note that $\pm\epsilon[K]=\pm\epsilon(\eta\pm fd)q\mu =(\pm\epsilon\eta+\epsilon fd)(-nq\mu') = (\mp \epsilon\eta n- f)q\mu'$.
Using the analog of the arguments for $(2)$, with the roles of $c, \,n, \,\mu$ exchanged for $f, \,d, \,\mu'$, we either obtain conclusions (a), (b) or the additional cases $[K]=(n-f)q\mu',\, (2n-f)q\mu'$ corresponding to conclusions (c)--(e).
\end{proof}

The next lemma determines the positive integers $y,\, y'$ associated to special classes of $\pm [K]$.

\begin{lem}\label{lem:nupmlambda} 
Assume $q > 40$. Let $\pm [K]=(\nu d + h)q \mu=-(\epsilon\nu + h n)q \mu'$ for some $\nu \in \bc{0,1,2}$ and for some $h \in \bc{-2,-1, 0, 1, 2}$. Then:
\begin{enumerate}
    \item
        \begin{enumerate}
        \item If $n \geq 4 \nu+1$, then either $y=\nu d + h$, or $\nu=0$ and $y=|h|$.
        \item If $4\nu-1 \geq n \geq 2\nu+1$, then $y=\pth{(n-2\nu)d-(\epsilon + 2h)}/{2}$.
        \item If $2 \nu -1 \geq n$, then $n=3,\,\nu=2$ and $y=\pth{d+\epsilon + 2h}/{2}$.
        \end{enumerate}
    \item
    Moreover, letting $\rho= {h}/{|h|}$ when $h \neq 0$ and $\rho=1$ when $h=0$, we have:
        \begin{enumerate}
        \item If $d \geq 4 |h|+1$, then either $y'=|h| n + \epsilon\rho\nu$ or $h=0$ and $y'=\nu$. 
        %or $\nu=0$ and $y=|h|$.
        \item If $4|h|-1 \geq d \geq 2|h|+1$, then $y'=\pth{(d-2|h|)n-(\epsilon +2\epsilon\rho\nu)}/{2}$.
        \item If $2 |h|-1 \geq d$, then $d=3,\,|h|=2$ and $y'=\pth{n+\epsilon +2\epsilon\rho\nu}/{2}$.
    \end{enumerate}
\end{enumerate}
%In all cases above the $\pm$ are corresponding.
\end{lem}

\begin{proof}
We have $2q=nd-\epsilon$ with $n,d\ge 3$ both odd. Recall $y,y'$ are determined by their mod $q$ class up to sign and the fact that $0<y,y'<q/2$.

We first compute $y$, where up to sign $y \equiv \nu d+h \, \mod (q)$. Assume $n\ge 4\nu + 1$. If $\nu=0$, then $y=\abs{h}$. Assume $\nu \neq 0$. Then $\nu d+h > 0$. Furthermore, $q/2-(\nu d+h) = (nd-\epsilon)/4-(\nu d+h) = ((n-4\nu)d-(\epsilon+4h))/4 > 0$ as $q>40$ and $n$ is odd. Thus $y= \nu d + h$ as claimed. 

Now assume $4 \nu-1 \geq n \geq 2\nu+1$. Then, $\nu d + h > q/2$. Using $q > 40$, we have 
$q/2 > q - (\nu d + h) = \pth{(n-2\nu)d-(\epsilon +2h)}/{2} >0$.
We conclude that $y=\pth{(n-2\nu)d-(\epsilon +2h)}/{2}$. 

Finally, assume $2 \nu -1 \geq n$. Then $\nu=2$ and $n=3$. As $q>40$, 
$0 < (\nu d + h) - q = \pth{d+\epsilon+2h}/{2} < \pth{3d-1}/{4} = {q}/{2}$.
As claimed, $y=\pth{d+\epsilon+2 h}/{2}$.

For the computation of $y'$ we use that $\pm [K]=-(\epsilon\nu + h n)q \mu'$ and hence that up to sign $y' \equiv \abs{h}n+ \epsilon \rho\nu \, \,\mod (q)$.  The values of $y'$ come from the analogous argument to that for $y$ above, with the roles of $d, \,n, \,\abs{h}, \, \epsilon \rho \nu, \,y'$ playing those of $n, \,d, \,\nu, \,h, \,y$.
\end{proof}

Again we have a variation of Theorem~\ref{thm:c+d:2q+1} with sharper restriction placed on the decomposition of $[K]$.

\begin{cor}\label{cor:c+d:d|2q+1<q/2}
Assume $L(q^2,m)$ as in case $(*2q)$. Assume further that, up to sign,  $[K]=c \mu \pm f \mu'$ and $\pth{nd-\epsilon}/{2}=q > 40$. 
\begin{enumerate}
\item If $0 \leq c \leq 1$ and $0 \leq  f <q/2$, then $\pm [K]=q \mu$. 
\item If $0 \leq f \leq 1$ and $0 \leq  c <q/2$, then $\pm [K]=dq \mu$. 
\end{enumerate}
\end{cor}

\begin{proof} \hfill

\noindent {\bf Assume that $(1)$ holds:} $0 \leq c \leq 1$ and $0 \leq  f <q/2$.

If $c=0$, then Theorem~\ref{thm:c+d:2q+1}(1) contradicts $[K]$ having order $q$. Thus, we assume $c=1$, apply Theorem~\ref{thm:c+d:2q+1}(2), and consider the alternatives. If conclusion $(a)$ applies, by $c=1$ and $f<q/2$, $\eta = 0$ and $\epsilon[K]=-q\mu$, which verifies the theorem. $(c)$--$(e)$ don't apply as $c=1$. Thus, assume conclusion $(b)$ applies. Furthermore, $f<q/2$ and $n,d \geq 3$ imply that $\eta=1$, $d=3$, and $2q=3n-\epsilon$. That is
\[\epsilon[K]=2q \mu, \ d=3.\]
Applying Lemma~\ref{lem:nupmlambda} with $\nu=0,\,h=2$, we conclude that $y=2, \,y'=\pth{n+\epsilon}/{2}$. Then $q=\pth{3n-\epsilon}/{2} = 3y'-2\epsilon$. Let $e',r',s'$ be as in Theorem~\ref{thm:=x}. Then $(r'+s')y'+e'=q$.

\textbf{First assume $\epsilon = 1$.} Then $40<q=3y'-2<3y'$, so $r'=s'=1$, $y'> 4$, and
\[e'=q-2y'=y'-2.\]

Apply Theorem~\ref{thm:q<3y,b>q+2e}. Assume that conclusion $(2)$ holds: there is a $b \geq q + \min\{e',y'-e'\} > q$ with $M_b' <2q$. Letting $b'=2q-2-b < q-2 = q-y = x$, by Lemma~\ref{thm:m:b<x}, $m_{b'}=0$. Thus, we can write $\pm [K]=c' \mu$ where $c' < 2q$. As $[K]$ has order $q$, we conclude $\pm [K]=q \mu$, contradicting $\pm [K]=2q \mu$. 

Thus conclusions $(1)$ or $(3)$ of Theorem~\ref{thm:q<3y,b>q+2e} must hold. As $2e'=2y'-4 \geq y'$, there is a $b \geq q+y'$ with $m_b' \leq 1$. Then $b'=2q-2-b \leq q-y'-2 < q-2 = x$. By Lemma~\ref{thm:m:b<x}, $M_{b'}<{b'}/{2}+1 \leq \pth{q-y'}/{2} = (q-1)/3$ as $q-y'=q-(q+2)/3=2(q-1)/3$. This allows us to write $\pm [K]=c' \mu \pm f' \mu'$ where $0 \leq c' \leq 1$ and $0 \leq f' < \pth{q-1}/{3} < q$. Since $[K]$ has order $q$, $c'=1$. We rewrite $\pm[K]=(1 \pm f'm) \mu$. As $[K]$ has order $q$, we conclude that $1 \pm f'm \equiv 1\mp f'd \equiv 1 \mp 3f' \equiv 0 \,\, (\mod (q))$. This is not possible, as $0 \leq 3f'<q-1$.

\textbf{Now assume $\epsilon = -1$.} Then $q=3y'+2<4y'$, $e' \geq 2$, and $y' > 10$ as $q>40$.

Applying Addendum~\ref{thm:b=<q}, one of the following holds: $(1)$ For every $b \le q$, $M'_b < (1/e' + 2/y')(q - e')$. $(2)$ There is a $b \ge q + y'$ with $m'_b \le 1$.

Assume $(1)$ holds. Then $M'_q < (1/e' + 2/y')(q - e') < q$ and $b'=2q-2-q = q-2=x$. As $r+s \leq 3$, Theorem~\ref{thm:=x} shows that $m_{b'} = 1$. This allows us to write $\pm [K]=c' \mu \pm f' \mu'$ where $0 \leq c' < q$ and $f'= 1$.
%We rewrite $\pm[K]=(c' \pm m) \mu$. We conclude $f'=1$, as $[K]$ has order $q$.
Then $\pm[K]=(c'\pm m)\mu  = (c'\pm 3(q+1))\mu$, contradicting $\pm[K]=2q\mu$ as $q>40$.

Assume $(2)$ holds and $b$ as given. Let $b'=2q-2-b = q-y'-2 < q-2 = x$. By Lemma~\ref{thm:m:b<x}, $M_{b'}< b'/y+1 = (q-y'-2)/2+1 = y'+1$, so $M_{b'}\le y' = (q-2)/3$. With the same arguments as in the case $\epsilon = 1$ when $(1)$ or $(3)$ holds, this is impossible.

\noindent {\bf Assume that $(2)$ holds:} $0 \leq f \leq 1$ and $0 \leq  c <q/2$. 

Note that as $m=d(q-\epsilon)$, $m^{-1}\equiv\epsilon n(q-\epsilon)\pmod{q^2}\equiv -n\pmod{q}$. In particular, when $n=3$ and $\epsilon = -1$, $m^{-1}\equiv -3(q+1)\pmod{q^2}$. The proof is then identical to that we used above when $(1)$ holds up to sign of $[K]$, with the roles of $c, \,c', \,n, \,m, \,x, \,y, \,e', \,\mu, \,\Gamma_b$ exchanged with $f, \,f', \,d, \,m^{-1}, \,x', \,y', \,e, \,\mu', \,\Gamma'_b$, replacing Theorem~\ref{thm:c+d:2q+1}$(1)$, $(2)$ with Theorem~\ref{thm:c+d:2q+1}$(4)$, $(3)$.
% $\spadesuit$ {(John: It seems like you also want to switch roles of $m$ and $m^{-1}$. But in the case Thm 6.8(3)(b) with $\epsilon=-1$ and conclusion $(1)$ of Addendum 4.2, you get $\pm [K]=\mu \pm f' \mu'$. The analog is to look at $(m^{-1} \pm f')\mu'$ and show this doesn't equal $2q\mu'$. But we don't know $m^{-1} \, \mod (q^2)$. Is it the same argument as above in this case? Do we want to modify the argument above so that the analog is clear here? Or just give an explicit expression for $m^{-1}$?)}
\end{proof}

\begin{thm}\label{thm:q<4y,e+1:2q+1} Assume $L(q^2,m)$ as in case $(*2q)$. Let $e$ be as in Theorem~\ref{thm:=x}. Assume $40<q<4y$. Then either $\pm [K] = q \mu, dq \mu$ or $y'>e+1$. 
\end{thm}
\begin{proof}
Let $r,s,e$ be as in Theorem~\ref{thm:=x}. Assume for contradiction that $y' \leq e+1$. By Lemma~\ref{lem:q<4y}(1), $m_j \leq 2$ when $j=q+e-1$. Let $j'=2q-2-j = q-e-1\leq q-y'=x'$. By Theorem~\ref{thm:M:b<x+1},  $M_{j'}' < {q}/{y'} \leq {q}/{2}$. Thus we may write $\pm [K]=c \mu \pm f \mu'$ where $0 \leq c <q/2,\, 0 \leq f \leq 2$. By Corollary~\ref{cor:c+d:d|2q+1<q/2}, we may assume $f=2$. Note then that Lemma~\ref{lem:q<4y} shows that $r+s=3$. As $[K]$ has order $q$, $c \neq 0$. Theorem~\ref{thm:c+d:2q+1}(3) applies and we consider, in turn conclusions $(a)$ and $(b)$, as $f=2$ so $(c)$ or $(e)$ does not occur and $(d)$ is desired.

{\bf Assume Theorem~\ref{thm:c+d:2q+1}$(3)(a)$ holds: $\pm [K]=(\eta +\epsilon 2d)q \mu=-(\eta n+2)q \mu'$ where $c =\eta q+2d$ with $\eta = 0, 1$.} 

As $d \geq 3$, $c <q/2$ implies that $\eta =0$. Then up to sign $[K]=-2q \mu'$. As $2<q/2$, we conclude that $y'=2$. As $r+s=3$, Lemma~\ref{lem:q<4y}(1) shows that for $z = q$, $m_z \leq 1$. Then $z'=2q-2-z = q-2=q-y'=x'$. By Theorem~\ref{thm:M:b<x+1}, $M_{z'}'<q/2$. Thus by Corollary~\ref{cor:Mmdecom}, $\pm[K]=c' \mu \pm f' \mu'$ where $0 \leq c' < q/2, \,0 \leq f' \leq 1$. Corollary~\ref{cor:c+d:d|2q+1<q/2} now proves the Theorem.

{\bf Assume Theorem~\ref{thm:c+d:2q+1}$(3)(b)$ holds: $\pm [K]=(\eta -\epsilon 2d)q \mu=-(\eta n-2)q \mu'$ where $c =\eta q-2d$ with $\eta = 1,2,3,4$.}

As $n,d \geq 3$, $c <q/2$ implies that $\eta  = 1$. That is, $q-2d<{q}/{2}$. Then $\pth{nd-\epsilon}/{4}<2d$ and we conclude that $n \leq 8$. That is, $n=3,5,7$. We consider these in turn. Note that $\mp\epsilon[K]=(2d-\epsilon)q \mu$; consequently Lemma~\ref{lem:nupmlambda} applies with $\nu=2,\,h=-\epsilon$.

If $n=7$, then Lemma~\ref{lem:nupmlambda} says that $y=\pth{3d+\epsilon}/{2}$. 
%Then $d>6$ and $Lemma~\ref{lem:nupmlambda}$ says that $y=n-2=5$. 
Then $3y=\pth{9d+3\epsilon}/{2}>\pth{7d-\epsilon}/{2}=q$. Recall that $q=(r+s)y+e$ where $e>0$. This contradicts our assumption that $r+s=3$. 

If $n=5$, then Lemma~\ref{lem:nupmlambda} says that $y=\pth{d+\epsilon}/{2}$. Then $5d-\epsilon = 2q < 8y = 4(d+\epsilon)$. Thus $d=3$ and $\epsilon = 1$, so $q=7$, contradicting $q>40$.
% $10 < {q}/{4}<y=\pth{d+\epsilon}/{2}$. That is, $d>19$. 
%Then $d>19$ and $Lemma~\ref{lem:nupmlambda}$ says that $y=n-2=3$, again a contradiction.
% But $4y=\frac{4d+4}{2}<\frac{5d-1}{2}=q$, contradicting an assumption on $q$.

If $n=3$, then $\pm [K]=-(n-2)q \mu'=-q \mu'=\epsilon dq \mu$, as desired.
\end{proof}

\begin{thm}\label{thm:q<4y',e'+1:2q+1} Assume $L(q^2,m)$ as in case $(*2q)$. Let $e$ be as in Theorem~\ref{thm:=x}. Assume $40<q<4y'$. Then either $\pm [K] = q \mu, dq \mu$ or $y>e'+1$. 
\end{thm}

\begin{proof}
The proof is identical to that of Theorem~\ref{thm:q<4y,e+1:2q+1} up to the sign of $[K]$, with the roles of $c, \,c', \,n, \,y, \,e, \,\mu, \,\Gamma$ switched with $f, \,f', \,d, \,y', \,e', \,\mu', \,\Gamma'$, replacing Theorem~\ref{thm:c+d:2q+1}(3), Lemma~\ref{lem:q<4y}(1) with Theorem~\ref{thm:c+d:2q+1}(2), Lemma~\ref{lem:q<4y}(2).
\end{proof}

\subsection{Spherical Theorem}\label{subsec:mainthm}

\subsubsection{Spherical Theorem when $(*q)$ holds.}\label{subsubsec:mainthmq}

The following theorem almost proves the Spherical Theorem, Theorem~\ref{thm:Gr-Li-Lu}, in the case $(*q)$ where
\[m=d(q-\epsilon),  \qquad \epsilon=\pm1, \qquad q= nd+\epsilon, \qquad n \ge 2,\,\, d\geq 3 \,\, \hbox{and}  \,\, d \, \, \hbox{odd}.\]
The conclusion that $\pm [K]=2q\mu$ is dealt with in section~\ref{sec:alexwidth}.

\begin{thm}\label{thm:main1-nd} 
Assume $99<q=nd+\epsilon$ with $n,d$ as in case $(*q)$. Let $K=K(q^2,d(q-\epsilon),\kappa)$ be a simple knot that admits a longitudinal surgery to $S^1\times S^2$. Let $[K]$ denote its class in $H_1(L(q^2,d(q-\epsilon)))$. Let $\mu, \,\mu'$ be the standard generators of $H_1(L(q^2,d(q-\epsilon)))$.
Then, 
\[\pm[K] = q \mu, dq \mu, 2q \mu.\]
% $\spadesuit$ {\em (For reference, by looking more closely at Corollary~\ref{cor:c+d:q| d-1<4/3 q}$(1)$, I think we can show that if $[K]=\pm 2q \mu$ then $K$ is the union of two arcs on either side of a level surface that writes $[K]=\pm (\mu +n \mu')$. Also for $n \geq 3, |d|>7, y=2, y'=2n$.)}
\end{thm}

\begin{proof}
Note that if $K$ admits longitudinal $S^1\times S^2$ surgery, then so does $K$ with its orientation reversed. So we may work with $[K]$ up to sign. If $y=1$, then $[K]=\pm yq \mu=\pm q \mu$. If $y'=1$, then $[K]=\pm y'q \mu'=\mp \epsilon dq \mu$. We hereafter assume $y, y'>1$. 

{\bf Assume $y \leq y'$.} 
As $y'>1$ we may apply
Theorem~\ref{thm:q+y,m_b=1} to $\Gamma'$ to conclude that there is $b$ such that one of the following holds:

\begin{enumerate}
\item  $b \geq 2q-y'-2$ and $m_b' \leq 3$;
\item  $b \geq q+y'-1$ and $M_b' < 2q$;
\item  $b \geq q+y'-1$ and $m_b' \leq 1$;
\item  $q<4y'$.
\end{enumerate}

If $(4)$ holds, then Theorem~\ref{thm:q<4y':1-ndprime} proves the Theorem. We hereafter assume $q >4y'$ ($y',q$ are relatively prime). If $(1)$ holds, Theorem~\ref{thm:m_b'leq3} proves the Theorem.

Assume that $(2)$ holds for some $b \geq q+y'-1$ where $M_b' < 2q$. Then $b'=(2q-2)-b \leq q-y'-1 \leq q-y-1 =x-1$, and Theorem~\ref{thm:m:b<x} shows that $m_{b'}=0$.
By Corollary~\ref{cor:Mmdecom}, $\pm[K]=c \mu \pm f \mu'$ where $0<c < 2q<q^2/2$ and $f = 0$. Theorem~\ref{thm:c+d}$(4)$  shows that $2q>c=yq$, contradicting our assumption that $y > 1$.

Assume that $(3)$ holds with $m_b' \leq 1$. As $b \geq q+y'-1$, we see again that $b'=(2q-2)-b<x$.
Thus by Theorem~\ref{thm:M:b<x+1}$(1)$, $M_{b'} < q/y$. By Corollary~\ref{cor:Mmdecom}, $\pm[K]=c \mu \pm f \mu'$ where $c \leq 1$ and $0 \leq f < q/2$. Corollary~\ref{cor:c+d:q| d-1<4/3 q}$(1)$ now implies the Theorem. 

{\bf Assume $y > y'$.} The proof is identical to the proof for $y \le y'$ with the roles of \break $\Gamma, \,x, \,y, \,m, \,M, \,c, \,f, \,\mu$, $\alpha$-weight exchanged with $\Gamma', \,x', \,y', \,m', \,M', \,f, \,c, \,\mu'$, $\beta$-weight. Additionally, Theorem~\ref{thm:M:b<x+1}$(1)$, \ref{thm:c+d}$(4)$, \ref{thm:m_b'leq3}, \ref{thm:q<4y':1-ndprime}, and Corollary~\ref{cor:c+d:q| d-1<4/3 q}$(1)$ are replaced by Theorem~\ref{thm:M:b<x+1}$(2)$, \ref{thm:c+d}$(1)$, \ref{thm:m_bleq3}, \ref{thm:q<4y:1-ndprime}, and Corollary~\ref{cor:c+d:q| d-1<4/3 q}$(2)$, respectively.
\end{proof}

\subsubsection{Spherical Theorem in case $(*2q)$}\label{subsubsec:mainthm2q}

The following proves the Spherical Theorem, Theorem~\ref{thm:Gr-Li-Lu}, for $L(q^2,m)$ in the case $(*2q)$ where
\[m=d(q-\epsilon),  \qquad \epsilon=\pm1, \qquad 2q= nd-\epsilon, \qquad n \ge 3,\,\, d\geq 3 \,\, \hbox{both odd}.\]

\begin{thm}\label{thm:main1-nd:2q+1}
Assume $228 < q = \pth{nd-\epsilon}/{2}$ with $n,d$ as in case $(*2q)$. Let $K=K(q^2,d(q-\epsilon),\kappa)$ be a simple knot that admits longitudinal $S^1\times S^2$ surgery. Let $[K]$ denote its class in $H_1(L(q^2,d(q-\epsilon))$. Let $\mu, \mu'$ be the standard generators of $H_1(L(q^2,d(q-\epsilon))$.
Then,
\[\pm[K]=q \mu, dq \mu.\]
\end{thm}

%$\spadesuit$ {\em (When $m_b'=1$, any way we can get $b \geq 2q-y$?)}
%
%$\spadesuit$ {\em (When $q<4y'$. It seems that $K=dq \mu$ can't happen here, for then $y'=1$ and $q<4y'=4$. When $K=q, 2q \mu$ we get $y=1,2$ and then $y'=n,2n$. Then $q<4y'=4n, 8n$. In particular, $|d| \leq 7$. Maybe more can be done here.)}
%
%$\spadesuit${\em (When $q<4y$. How can we get $q\mu, 2q \mu$? In that case we have $y \leq 2$ and $q<4y \leq 8$. Note that if $y=d$, then we conclude $q<4d$ and then $n \leq 3$.  If $q<3y$, then we conclude $n \leq 2$.)}

\begin{proof}
%$K=K(q^2,d(q-1),\kappa)$ and $[K]=\kappa \mu = l \mu'$ where $\kappa=t'q$ and $l=tq$. 
%Recall that $y = \min \{t',q-t'\}, y'= \min \{t, q-t\}$. %By Lemma ~\ref{lem:y=1}, we may assume $y >2, y'>1$. 
The proof is identical to the proof of Theorem~\ref{thm:main1-nd}, with 
Theorem~\ref{thm:c+d}, \ref{thm:m_bleq3}, \ref{thm:m_b'leq3}, \ref{thm:q<4y:1-ndprime}, \ref{thm:q<4y':1-ndprime},  and Corollary~\ref{cor:c+d:q| d-1<4/3 q}
replaced by
Theorem~\ref{thm:c+d:2q+1}, \ref{thm:m_bleq3:2q+1}, \ref{thm:m_b'leq3:2q+1}, \ref{thm:q<4y:2q+1}, \ref{thm:q<4y':2q+1}, and Corollary~\ref{cor:c+d:d|2q+1<q/2}, respectively.
\end{proof}

\subsection{When \pdfmath{b \geq 2q-y-2}{b>=2q-y-2} and \pdfmath{m_b \leq 3}{mb <= 3}}\label{subsec:b>=2q-y-2}
\hfill

In this subsection we consider conclusion $(1)$ (the $\Gamma$ and $\Gamma'$ versions) of the Engine Theorem, Theorem~\ref{thm:q+y,m_b=1}, in the context of $L(q^2,m)$ with $m=q(d-\epsilon)$.
Recall that $\mu'=m \mu$.  Now $m \equiv -\epsilon d \,\,\, (q)$ and $m^{-1} \equiv -\epsilon d^{-1} \,\,\, (q)$. When $q=nd+\epsilon$, $d^{-1} \equiv -\epsilon n \,\,\, (q)$ and $m^{-1} \equiv n \,\,\, (q)$; when $2q = nd-\epsilon$, $d^{-1} \equiv \epsilon n \,\,\, (q)$ and $m^{-1} \equiv -n \,\,\, (q)$.

\begin{lem}\label{lem:y<dy'}
$y \leq dy', \,y' \leq ny$.
\end{lem}

\begin{proof}
Up to sign $[K]=yq \mu=\pm y'q \mu'=\mp \epsilon dy'q \mu$. Thus, $y \equiv \mp \epsilon dy'  \,\,\, (q)$. Let $\bk{\cdot}_p$ denote the least positive residue. As $0<y<q/2$, either $y = \bk{dy'}_q\le dy'$, or $y < q/2 < \bk{dy'}_q \le dy'$. The second inequality follows similarly from the facts that $y' \equiv \pm ny \,\,\, (q)$ and $y'< q/2$.
\end{proof}

\subsubsection{$q=nd+\epsilon$} In this subsubsection of subsection~\ref{subsec:b>=2q-y-2} we consider $L(q^2,m)$ when $(*q)$ holds. 

For $(*q)$, we make the following observations:
\begin{enumerate}[label=(\Alph*)]
\item $\mu'=m \mu$, $m \equiv -\epsilon d \pmod{q}$, and $m^{-1} \equiv n \pmod{q}$. Thus in $H_1(L(q^2,m))$, for any integer $x$, we have $xq\mu'=xq(-\epsilon d) \mu, \,xq \mu= nxq \mu'$ in $H_1(L(q^2,m))$. % We have just recall this before Lemma~\ref{lem:y<dy'}
\item $[K]=\pm y q \mu$ and $0<y<q/2$. Thus if $[K]=xq \mu$ then $y$ is the least positive residue modulo $q$ of either $x$ or $-x$---whichever has residue less than $q/2$.
\item $[K]=\pm y' q \mu'$ and $0<y'<q/2$. Thus if $[K]=xq \mu=nx q \mu'$ then $y'$ is the least positive residue modulo $q$ of either $nx$ or $-nx$---whichever has residue less than $q/2$.
\end{enumerate}

\begin{thm}\label{thm:m_bleq3}
Assume $L(q^2,m)$ as in case $(*q)$. Assume further that  $y \geq y'$, $q>10$, and there is a $b \geq 2q-y-2$ with $m_b \leq 3$. Then 
\[\pm [K] = q \mu, dq \mu.\]
\end{thm}

%$\spadesuit$ {\em (Is the assumption $y\ge y'$ used? This assumption is not in the analog following.)}

%$\blacklozenge$ {\em (I don't think so if we are willing to include $\pm[K]=2q\mu$, which is fine as we include in the next theorem anyway. I think that you separated $2q\mu$ to try to argue $2q\mu$ happens only if $n=2$, which is proved using the Alexander polynomial. In fact, I am not sure if the order of $y$ and $y'$ has ever been used in any of 6.16-6.19. However, we are free to assume the order, by the way how we proved the main theorems 6.13, 6.14, in particular, (2) and (3).
%)}

\begin{proof}
We may assume $y,y'>1$, else the Theorem holds. Since $y \geq y'$, if $y=2$, then $y'=2$. But then up to sign $[K]$ is $2q \mu$ and up to sign $2q \mu'$. That is $2q \mu = \pm 2 q \mu'=\mp 2\epsilon d q \mu$. We then conclude that $2 \pm 2\epsilon d \equiv 0 \ (q)$, contradicting that $q=nd + \epsilon$ with $n \geq 2, d\ge 3$. Thus, we hereafter assume $y>2$. As $q/2>2$, this means that we may further assume that $\pm [K] \neq 2q \mu$.

Let $b'=2q-2-b \leq y<q-y'=x'$ as $y, y'<q/2$. By Lemmas~\ref{thm:m:b<x} and ~\ref{lem:y<dy'}, $M_{b'}' < b'/y' +1 \leq y/y' +1\leq d +1$. As $M_{b'}'$ is an integer, $M_{b'}' \leq d$. By Corollary~\ref{cor:Mmdecom}, we can write $\pm[K]=c \mu \pm f \mu'$ with $0 \leq f=m_b \leq 3, \,0 \leq c \leq d, \,c < y/y' +1 < y$. By Lemma~\ref{thm:c+d}(4), $f > 0$. By Theorem~\ref{thm:c+d}(3), we have three possibilities $(a)$--$(d)$ for $[K]$, which we treat in turn.

\noindent {\bf Assume Theorem~\ref{thm:c+d}$(3)(a)$ holds.}

Then $[K]=(\eta -\epsilon f d)q \mu$ with $\eta \in \{1,2\}$ and $c=\eta q-fd$.
% We further assume either $n \not= 2$ or $f \not= 3$, otherwise we will treat in (b).
Then, $d \geq c= \eta (nd+\epsilon) - f d$ implies that $4d\ge (f+1)d \ge \eta(nd+\epsilon)\ge \eta nd - 2$, so $(\eta n - 4)d\le 2$. Thus, $\eta n - 4 < 0$, as $d\ge 3$. Since $n \ge 2$, we conclude that $\eta = 1$ and $n = 2, 3$. Thus, $c = q- fd$ and $(f+1)d \ge nd+\epsilon$, which implies one of the following four cases: $(1)$ $f = n = 2, 3$; $(2)$ $f = 2$, $n = 3$, and $\epsilon = -1$; $(3)$ $f = 3$ and $n = 2$; or $(4)$ $f = 1$ and $n = 2$, and $\epsilon = -1$.

\begin{enumerate}
    \item Assume $f = n = 2, 3$. $\bk{K} = (1 - \epsilon nd)q\mu = 2q\mu$, which we assume does not occur.
    \item Assume $f = 2$, $n = 3$, and $\epsilon = -1$. Then, $q = 3d-1, \,[K] = (1+2d)q\mu, \,y = q-(1+2d) = d-2$ as $q > 10$. Thus, $c = q-2d = d-1 > y$, contradicting $c<y$.
    \item Assume $f=3, \,n=2$. Then $10< q=2d+\epsilon, \,d\ge c=2q-3d=d+2\epsilon$, so $\epsilon = -1, \,c=d-2, \,d\ge 5$. Thus, $[K]=(1+3d)q\mu = (3-d)q\mu$. Therefore, $y=d-3$ and $y'=5$ using that $q>10$ and facts $(B),(C)$ above. Thus, $d-2 = c < y/y'+1 = (d-3)/5+1$, contradicting $d \ge 5$.
    \item Assume $f = 1, \,n = 2$, and $\epsilon = -1$. Then, $q = 2d-1, \,[K] = (1+d)q\mu, \,y = q-(1+d) = d-2$ using $q > 6$ and fact $(B)$. Thus, $c = q-d = d-1 > y$, contradicting $c<y$.
\end{enumerate}

\noindent {\bf Assume Theorem~\ref{thm:c+d}$(3)(b)$ holds.}

As $c \leq d<q$ the only possibility is that $f=1$ and $[K]=\epsilon dq \mu$, verifying the Theorem.

\noindent {\bf Assume Theorem~\ref{thm:c+d}$(3)(c)$ holds.}

Then $n=2, \,f=3$. As $c \leq d<q=2d+\epsilon$, only conclusion $c=d-\epsilon$ is possible.

Assume $c=d-\epsilon$ and $[K]= (\epsilon d-2)q\mu$. As $[K]=\pm yq\mu$, $0<y<q/2$, and $0<d-2 < (2d+1)/2$, $0\le d-3<(2d-1)/2$, we conclude that $y=d-2$ if $\epsilon = 1$, or $y=d-3$ if $\epsilon = -1$ (fact $(B)$ above). Applying fact $(C)$, looking at the least positive residues mod $q$ of $\pm n(\epsilon d-2)$,  we conclude $y'=5$. Then $d-\epsilon =c<y/y' +1 \leq {(d-2)}/{5}+1$, contradicting that $d \geq 3$.

\noindent {\bf Assume Theorem~\ref{thm:c+d}$(3)(d)$ holds.}

Then $\epsilon = -1$, $[K]=4q\mu$, and $c=2q-1$ as we assumed $\pm[K]\neq 2q\mu$. Thus, $y=4$. But 
$c=2q-1>20$ contradicts $c<y=4$.
\end{proof}

\begin{thm}\label{thm:m_b'leq3}
Assume $L(q^2,m)$ as in case $(*q)$.
Assume further that $q>4$ and there is a $b \geq 2q-y'-2$ with $m_b' \leq 3$. Then, 
\[\pm [K] = 2q \mu, dq \mu.\] %(decomposing as $\mu \pm n \mu'$).
\end{thm}

\begin{proof}
Let $b'=2q-2-b \leq y' < q-y=x$ ($y, y'<q/2$). Lemmas~\ref{thm:m:b<x} and ~\ref{lem:y<dy'} then say that $M_{b'} < b'/y +1 \leq y'/y + 1 \leq n+1$. As an integer, $M_{b'} \leq n$. By Corollary~\ref{cor:Mmdecom}, we can write $\pm[K]=c \mu \pm f \mu'$ with $0 \leq c=m_b' \leq 3, \,0 \leq f \leq n<q$. By Theorem~\ref{thm:c+d}(1), $c>0$. By Theorem~\ref{thm:c+d}(2), we have possibilities $(a)$--$(e)$ for $[K]$,  which we treat in turn.

\noindent {\bf Assume Theorem~\ref{thm:c+d}$(2)(a)$ holds.}
Then $2q - 3n \le 2q-cn=f \leq n$ contradicts $d\geq 3$. 

\noindent {\bf Assume Theorem~\ref{thm:c+d}$(2)(b)$ holds.}
Then $f=q+cn>q$, contradicting that $f <q$. 

\noindent {\bf Assume Theorem~\ref{thm:c+d}$(2)(c)$ holds.}
Then, $\epsilon [K]= (2c-d)q \mu$ and $q-cn=f \leq n$. Thus, $nd + \epsilon = q \leq (1+c)n$ implies  $\epsilon \leq (1+c-d)n$, giving $1+c\ge d$.  Then $d=3$, $q=3n+\epsilon$, and $c=3$. We conclude that $\epsilon [K]= 3q \mu=dq \mu$ as claimed by the Theorem.

\noindent {\bf Assume Theorem~\ref{thm:c+d}$(2)(d)$ holds.}
As $cn=f\leq n$, we conclude that $c=1$, and $\pm [K]=2q \mu$ as claimed by the Theorem.

\noindent {\bf Assume Theorem~\ref{thm:c+d}$(2)(e)$ holds.}
Then, $\pm [K] = dq\mu$ as claimed by the Theorem.
\end{proof}

\subsubsection{$2q = nd-\epsilon$} In this subsubsection of subsection~\ref{subsec:b>=2q-y-2} we consider $L(q^2,m)$ when $(*2q)$ holds.

\begin{thm}\label{thm:m_bleq3:2q+1}
Let $L(q^2,m)$ be as in case $(*2q)$. Assume $q>18$, $y' \leq y$, and there is a $b \geq 2q-y-2$ with $m_b \leq 3$. Then 
\[\pm [K] =q \mu, dq \mu.\]
\end{thm}

\begin{proof} 
Let $b'=2q-2-b \leq y$. As $y,y'<q/2$, $y < q-y'=x'$. Lemmas~\ref{thm:m:b<x} and ~\ref{lem:y<dy'} say that $M_{b'}' < b'/y' +1 \leq y/y' +1\leq d+1$. As $M_{b'}'$ is integral, we conclude $M_{b'}' \leq d$. By Corollary~\ref{cor:Mmdecom}, we can write $\pm[K]=c \mu \pm f \mu'$ with $0 \leq c <y/y' +1, \,c \leq d$, and $0\le f\le 3$. By Theorem~\ref{thm:c+d:2q+1}(4), we may assume that $f>0$. By Theorem~\ref{thm:c+d:2q+1}(3), and there are five possibilities $(a)$--$(e)$ for $[K]$. Conclusion $(d)$, $(e)$ immediately verifies the Lemma. We treat $(a)$--$(c)$ in turn.

\noindent {\bf Assume Theorem$~\ref{thm:c+d:2q+1}(3)(a)$ holds.} Then, $\pm[K]=(\eta + \epsilon fd)q \mu=-(\eta n+f)q \mu'$ and $c=\eta q + fd$ where $\eta  \in \{0,1\}$. As $c\le d$, we conclude that $\eta =0, \,f=1$. Then $\pm[K]=\epsilon dq \mu$ as claimed.

\noindent {\bf Assume Theorem$~\ref{thm:c+d:2q+1}(3)(b)$ holds.} Then, $\pm[K]=(\eta- \epsilon fd)q \mu=-(\eta n-f)q \mu'$ and $c=\eta q-fd$ where $\eta \in \{1,2,3,4\}$. Thus $d \geq c=\eta q-fd \geq 0$ implies $(f+1)d \geq \eta \pth{nd-\epsilon}/{2} \geq f d$, which can be rewritten as 
\[\frac{2(f+1)}{\eta}+\frac{\epsilon}{d}\geq n \geq \frac{2f}{\eta}+\frac{\epsilon}{d}.\]

Assume $\eta \ge 3$. Then $3 \ge 2\pth{f+1}/{\eta}+{\epsilon}/{d} \ge n\ge 3$, we conclude that $n=d=3, \,\epsilon = 1$, but then $q=4$, a contradiction.

Assume $\eta =2$. Then $5 >f+1+{\epsilon}/{d} \geq n \geq f +{\epsilon}/{d}$. As $n$ is odd, we conclude $n=3$ and $f = 1, 2$, or $\epsilon=-1$ and $f=3$. Then $q=\pth{3d-\epsilon}/{2} > 18$, which implies $d\ge 13$, and $[K]=(2-\epsilon fd)q \mu=-(6-f)q \mu'$ and $c=2q-fd=3d-\epsilon-fd=(3-f)d-\epsilon$. As $0<y'<q/2$ and $[K]=\pm y'q \mu'$, we conclude that $y'=6-f$ as $q>18$. Assume $f=1$. Then $0<y<q/2$ implies that $y =\pth{3d-\epsilon}/{2}-(d-2\epsilon) = (d+3\epsilon)/2$. We get that $2d-\epsilon=c<y/y' +1=\pth{d+3\epsilon}/{10}+1$, contradicting that $d\ge 13$. Similarly, when $f=2$, $y =(2d-2\epsilon) - \pth{3d-\epsilon}/{2} = (d-3\epsilon)/2$, we get that $d-\epsilon = c<y/y' +1 = \pth{d-3\epsilon}/{8}+1$, contradicting that $d\ge 13$. Finally, when $f=3$ and $\epsilon=-1$, $\pm[K]=(3d+2)q\mu = (2q+1)q\mu = q\mu$, as desired.

Finally, assume $\eta =1$. Then $2(f+1) +{\epsilon}/{d} \geq n \geq 2f+ {\epsilon}/{d}$. Thus $n= 2f +1$ as $n$ is odd.  Then $\pm[K]=(1-\epsilon fd)q \mu=-((2f+1)-f)q \mu'=-(f+1)q \mu'$ and $c=q-fd=\pth{(2f+1)d-\epsilon}/{2}-fd=\pth{d-\epsilon}/{2}$. As $0<y'<q/2$, $[K]= \pm y'q\mu'$, and $\pm[K]=-(f+1)q \mu'$, we conclude that $y'=f+1$. As $0<y<q/2$, $[K]= \pm yq \mu$, and $\pm \epsilon[K]=(\epsilon(1-\epsilon fd)+\pth{(2f+1)d-\epsilon}/{2})q \mu=(\pth{d+\epsilon}/{2})q \mu$, we conclude that $y=\pth{d+\epsilon}/{2}$ ($q>18$). Then we have $\pth{d-\epsilon}/{2}=c<y/y' +1=\pth{d+\epsilon}/\pth{2(f+1)}+1\le (d+\epsilon)/4+1$, so $d<4+3\epsilon \le 7$. Thus, $d\le 5$ as $d$ is odd and $n=2f+1 \le 7$, contradicting $q>18$.
% In the case when $f = 3$, we derive $n = 7$, $\epsilon = 1$ and $q\le 10$, so $q<10$ condition is sharp

\noindent {\bf Assume Theorem$~\ref{thm:c+d:2q+1}(3)(c)$ holds.} $\pm [K] = -\epsilon 2dq\mu = 2q\mu'$ and $f = 3, \,n = 5, \,c = (d + \epsilon)/2$. As $q>4$, $y'=2$. As $\pm [K]=(-\epsilon 2d + \epsilon \pth{5d-\epsilon}/{2})q \mu=(\epsilon \pth{d-\epsilon}/{2})q \mu$, we conclude $y=\pth{d-\epsilon}/{2}$, but then $\pth{d+\epsilon}/{2}=c <y/y' +1=\pth{d-\epsilon}/{4}+1$, which implies $\epsilon=-1$ and $d \le 5$, contradicting $q>18$.
\end{proof}

\begin{thm}\label{thm:m_b'leq3:2q+1}
Let $L(q^2,m)$ be as in case $(*2q)$. Assume $q>18$, $y \leq y'$, and there is a $b \geq 2q-y'-2$ with $m_b' \leq 3$. Then,
\[\pm [K] =q \mu, dq \mu.\]
\end{thm}

\begin{proof}
The proof is identical to the proof for Theorem~\ref{thm:m_bleq3:2q+1} working with $[K]$ up to a sign and with the roles of $n, \,x, \,y, \,m, \,M, \,c, \,\mu$, $\alpha$-weight exchanged with $d, \,x', \,y', \,m', \,M', \,f, \,\mu'$, $\beta$-weight. Additionally, Theorem~\ref{thm:c+d:2q+1}$(3)$$(4)$ are replaced by Theorem~\ref{thm:c+d:2q+1}$(1)$$(2)$. 
\end{proof}

\subsection{\pdfmath{q<4y}{q<4y} or \pdfmath{q<4y'}{q<4y'}}\label{subsec:q<4yor4y'}

In this subsection we consider conclusion $(4)$ (the $\Gamma$ and $\Gamma'$ versions) of the Engine Theorem, Theorem~\ref{thm:q+y,m_b=1}, in the context of $L(q^2,m)$ with $m=q(d-\epsilon)$.

\subsubsection{$q = nd+\epsilon$} In this subsubsection of subsection~\ref{subsec:q<4yor4y'} we consider $L(q^2,m)$ when $(*q)$ holds. 

\begin{lem}\label{lem:d+2,q<4y'}
Assume $86\le q=nd+\epsilon < 4y'$ and $[K]=(d- \zeta a) q \mu$ for $a \in \{1,2\}$ with $\zeta=\pm 1$. Then,
\[[q]_{y'}  \geq y.\]
\end{lem}

\begin{proof} From $86 \le q < 4y'$ we know $y' \ge 22$. If $n=2$, then Lemma~\ref{lem:dne2an>2}(I)(b) says that $y' \leq 2a+1 \leq 5$, contradicting $y' \geq 22$. We hereafter assume $n>2$.

First assume $d>2a$. By Lemma~\ref{lem:dne2an>2}(I)(c), we have 
\[y=d - \zeta a, \ y'= an + \zeta\epsilon.\]
Then $q<4y'$ becomes $nd+\epsilon<4an+4\zeta\epsilon$. That is, $(d-4a)n < (4\zeta -1)\epsilon \leq 5$, so $d-4a = 1$ and $n=3$, or $d-4a < 1$. If $d-4a=1$ and $n=3$, then $d=5, 9$, contradicting $q\ge 86$, so assume $d-4a < 1$. As $d$ is odd, $d-4a \neq 0$, so $d-4a \leq -1$. Thus, $2a<d \leq 4a-1$. That is, when $a=1$, $d=3$, and when $a=2$, $d=5, 7$. Thus, $q = nd+\epsilon \ge 86$ implies $n\ge 13$, and $y=d-\zeta a\le 9$.

Assume $a=1$ and $d=3$. Then, $q = 3n+\epsilon$, $y=3-\zeta$ and $y'=n+\zeta\epsilon$. Note that $q-2y'=n+\epsilon(1-2\zeta)>0$ and $q-4y'=-n+\epsilon(1-4\zeta)<0$, as $n\ge 13$. Thus, $[q]_{y'} = q-2y'$ or $q-3y'$, depending on if $q-3y'=\epsilon(1-3\zeta)=-\zeta\epsilon(3-\zeta)>0$, that is, if $\zeta\epsilon<0$. In either case, $[q]_{y'}\ge 3-\zeta = y$.

Now assume $a = 2$. Then $d=5$ or $d=7$, $y'=2n+\zeta\epsilon, \,y=d-\zeta a\le 9$. Assume $d=5$. Then $q=5n+\epsilon$. Thus, $[q]_{y'}=q-2y'=n+\epsilon(1-2\zeta)\ge n-3\ge y$, as $n\ge 13$, $y\le 9$.
Assume $d=7$. Then $q=7n+\epsilon$. Thus, $[q]_{y'}=q-3y'=n+\epsilon(1-3\zeta)\ge n-4\ge y$, as $n\ge 13$, $y\le 9$.

Finally assume $d \leq 2a$. Then $a=2, \,d=3, \,q = 3n+\epsilon, \,n\ge 29$. By Lemma~\ref{lem:dne2an>2}(I)(c)(iii), we have $y=3 - 2\zeta$ and $y'=n+\epsilon(1-\zeta) > 10 \ge 2y$. Note that $q-3y'=\epsilon\zeta(3-2\zeta) = \epsilon\zeta y$. Thus, $[q]_{y'}=y$ or $y'-y$. In either case, $[q]_{y'}\ge y$.
\end{proof}

\begin{lem}\label{lem:d+2,q<4y}
Assume $68\le q = nd+\epsilon < 4y$ and $[K]=(d- \zeta a) q \mu$ for $a \in \{1,2\}$ with $\zeta=\pm 1$. 
Then, %either $n=2$ and $[K]=(d-1)q \mu=-dq \mu$ or
\[[q]_{y} \geq y'.\] 
\end{lem}

\begin{proof}
We assume $4y>q \geq 68$. In particular, $y \geq 18$. First we claim that either $n = 2$ or both $n > 2$ and $d > 2a$. Otherwise, $n > 2$, and $d = 3$ and $a = 2$. By Lemma~\ref{lem:dne2an>2}, $y = d - \zeta a \le 5$, contradicting that $y \geq 18$. We now consider the alternatives.

{\bf Assume either $(i)$ $n > 2$ and $d >2a$ or $(ii)$ $n=2$ with $\zeta=1$.} By Lemma~\ref{lem:dne2an>2}, we have $y=d - \zeta a$ and $y'=an + \zeta\epsilon$. Then $q<4y$ becomes $nd+\epsilon < 4d - \zeta 4a$. That is, $d(n-4)<- \zeta 4a-\epsilon \leq 9$. Then $n \leq 6$. If $n>4$, then $d \leq 9$ and $q=nd+\epsilon \leq 55$, a contradiction. Thus $n \leq 4$.
% Then $68 \leq q = nd+\epsilon \leq 4d+\epsilon$ implies $d \geq 17$.  

From $y=d - \zeta a$ and $y'=an + \zeta\epsilon$, we obtain $ny + \zeta y' = (nd - \zeta an) + (\zeta an + \epsilon) = nd+\epsilon = q$. We also know that $y-y' > y-2y' = y - 2(an+\zeta\epsilon) \ge y - 18 \ge 0$. Therefore, $[q]_y$ is either $y'$ if $\zeta = 1$, or $y - y'$ if $\zeta = -1$ --- both at least $y'$ as claimed.

{\bf Assume $n=2,\, \zeta=-1$.} Then, $q=2d+\epsilon \ge 68$, so $d\ge 34$. By Lemma~\ref{lem:dne2an>2}, $y=d-(a-\epsilon), \,y'=2a-\epsilon$, so $y-y' = d-(3a-2\epsilon) > 0$ as $d>8$. As in the previous case, we compute $2y + y' = 2d+\epsilon = q$. Thus, $[q]_y = y'$, as $d>8$.
\end{proof}

\begin{lem}\label{lem:c<1,f<1,q<4y:1-nd}
Assume $68 < q < 4y$ and $\pm [K]=c \mu \pm f \mu'$ with either
\begin{enumerate}
\item $0 \leq c \leq 1$ and $0 \leq f < (1/e + 2/y)(q-e)$; or 
\item $0 \leq f \leq 1$ and $0 \leq c < 2q$.
\end{enumerate}

Then, $\pm[K]=q \mu, dq \mu$. 
\end{lem}

\begin{proof}
We assume $68 < q < 4y$, and then $y \geq 18$. By Theorem~\ref{thm:q<4y,e+1:1-nd}, we may assume that
\[y'>e+1.\]

\noindent {\bf Assume $0 \leq c \leq 1$ and $0 \leq f < (1/e + 2/y)(q-e)$}.

Note then that $f<2q$. By Theorem~\ref{thm:c+d}(1), $c = 0$ implies $[K] = \pm dq\mu$, so we may assume that $c=1$. Applying Theorem~\ref{thm:c+d}(2) we may write $[K]$ in the one of the following forms:
\begin{enumerate}[label=({\roman*})]
    \item $\pm [K]=(d- \zeta 2) q \mu$ for $\zeta=\pm 1$. 
    \item $\pm [K]=2q \mu$.
    \item $\pm [K]=(2d-2) q \mu$ and $f=2q-n$.
\end{enumerate}

Assume (i) holds. By Lemma~\ref{lem:d+2,q<4y}, $[q]_{y} \geq y'$, but $e \geq [q]_{y}$, contradicting that $y'>e+1$.

Assume (ii) holds. By Lemma~\ref{lem:dne2an>2}(II), we conclude that $y=2$, which contradicts $y\ge 18$.

Finally, assume that (iii) holds. By assumption we then have that $2q-n=f< (1/e + 2/y)(q-e)< 10q/9$ as $y\ge 18$, but this contradicts that $q=nd+\epsilon$ and $d \geq 3$.

\noindent {\bf Assume $0 \leq f \leq 1$ and $0 \leq c < 2q$}.

By Theorem~\ref{thm:c+d}(4), $f = 0$ implies $[K] = \pm q\mu$, so we may assume that $f=1$. Applying Theorem~\ref{thm:c+d}(3), noting that $n \geq 2$, we may write  $\pm [K]=(d- \zeta a) q \mu$ for $a \in \{1,2\}$ and $\zeta=\pm 1$. By Lemma~\ref{lem:d+2,q<4y}, $[q]_{y} \geq y'$, but $e \geq [q]_{y}$, contradicting that $y'>e+1$.
\end{proof}

\begin{lem}\label{lem:c<1,f<1,q<4y':1-nd}
Assume $99<q<4y'$ and $\pm [K]=c \mu \pm f \mu'$ with either
\begin{enumerate}
\item $0 \leq  c \leq 1$ and $0 \leq f < 2q$; or   
\item $0 \leq  f \leq 1$ and $0 \leq c < 2q$.
\end{enumerate}
Then, $\pm[K]=q \mu, 2q \mu, dq \mu$. 
\end{lem}

\begin{proof} 
Since we assume $99<q<4y'$, $y' \geq 25$. By Theorem~\ref{thm:q<4y',e'+1:1-nd}, we may assume that 
\[y>e'+1.\]

\noindent {\bf (1). Assume $0 \le c \leq 1$ and $0 \leq f < 2q$}.

By Theorem~\ref{thm:c+d}(1), we may assume that $c=1$. Applying Theorem~\ref{thm:c+d}(2) we may write $[K]$ in the one of the following forms:

\begin{enumerate}[label=({\roman*})]
    \item $\pm [K]=(d- 2 \zeta) q \mu$ for $\zeta=\pm 1$. 
    \item $\pm [K]=2q \mu$. 
    \item $\pm [K]=(2d-2) q \mu$ and $f=2q-n$.
\end{enumerate}

{\bf Assume $(i)$ holds.} By Lemma~\ref{lem:d+2,q<4y'}, $[q]_{y'} \geq y$, but $e' \geq [q]_{y'}$, contradicting that $y>e'+1$.

{\bf Assume $(ii)$ holds.} We obtain one of the conclusions.

{\bf Assume $(iii)$ holds.} As $[K]$ has order $q$, $q$ must be odd; as $d$ is odd, $n$ must be even. 

First assume $d \neq 3$. By Lemma~\ref{lem:2dpmq:1-nd}, $y' =2n+2\epsilon$. Then, $99 < 4y' = 4(2n+2\epsilon)$ implies $n\ge 12$.
As $99<q=nd+\epsilon < 4(2n+2\epsilon)$, we conclude that $d = 5, 7$ and $n \ge 16$.
Then again by Lemma~\ref{lem:2dpmq:1-nd}, we conclude that $y=2d-2$.
Then, $e'\ge [q]_{y'} = (nd+\epsilon)-(2n+2\epsilon)(d-1)/2 = n - \epsilon(d-2) \ge n-d+2$.
Thus, $2d-2 = y > e'+1 = n-d+3 \ge 19-d$, which implies $d > 7$, contradicting $d=5, 7$.

Thus, assume $d=3$. As $99 < q = 3n+\epsilon$, $n \ge 33$. By Lemma~\ref{lem:2dpmq:1-nd}, we conclude that $y' =n-\epsilon, \,y=2d-2=4$. Then $e'\ge [q]_{y'}=4$ or $n-3$, both contradicting that $e' < y - 1 = 3$.

\noindent {\bf (2). Assume $0 < f \leq 1$ and $0 \leq c < 2q$}.

By Theorem~\ref{thm:c+d}(4), we may assume that $f=1$. Applying Theorem~\ref{thm:c+d}(3), noting that $n \geq 2$, we may write  $\pm [K]=(d- \zeta a) q \mu$ for $a \in \{1,2\}$ and $\zeta=\pm 1$. By Lemma~\ref{lem:d+2,q<4y'}, $[q]_{y'} \geq y$, but $e' \geq [q]_{y'}$, contradicting that $y>e'+1$.
\end{proof}

\begin{lem}\label{lem:q<4y,q<4y':1-nd}
Assume $L(q^2,m)$ as in case $(*q)$. If $q <4y,\, q<4y'$, and $q > 75$, then $\pm [K]= q \mu, dq \mu$.
\end{lem}

\begin{proof} 
We assume that $q > 75$. Then $y,y' > 19$.

I. Assume $e \geq e'$.

By Lemma~\ref{lem:q<4y}, if we let $b=q+e-1$, then $m_b \leq 2$. Then $b'=2q-2-b = q-e-1 \leq q-e'-1 \leq q-[q]_{y'}-1$. By Theorem~\ref{thm:M:b<x+1}(4), $M'_{b'} < ({1}/{e'} + {2}/{y'})(q-e') < 2q$. Thus we may write $\pm [K]=c \mu \pm f \mu'$ with $0\leq f \leq 2, 0 \leq c <({1}/{e'} + {2}/{y'})(q-e')<2q$.  By Lemma~\ref{lem:c<1,f<1,q<4y:1-nd}, we may assume that $f=2$. Applying Theorem~\ref{thm:c+d}(3), we consider the alternatives $(a)$--$(d)$. 
First note that $(c)$ does not hold as $f=2$. Considering $(b)$ with $\eta = 0$, note that $\pm [K] = \epsilon 2dq \mu = -2q \mu'$ cannot occur, as then $y'=2$, contradicting that $y'>19$. Similarly, conclusion $(d)$ cannot hold as $y>19$. 

Then $(a),(b)$ leave us with the possibilities that $\pm [K]=(-2d + a)q\mu = (2\epsilon+an)q\mu'$ where $a \in \{-2,-1,1,2\}$. As $y \equiv \pm (2d-a) \, \mod \, (q)$ and $y' \equiv \pm (an+2\epsilon) \, \mod \, (q)$, we have that $y \leq 2d-a$ and $y' \leq an+2\epsilon \leq 2n+2$. From $y,y' > 19$, we further know that $n, d>8$. Then $nd+\epsilon=q<4y' \leq 8n+8$ implies that $d \leq 8$, a contradiction.

II. Assume $e <e'$.

By Lemma~\ref{lem:q<4y}, if we let $b=q+e'-1$, then $m_b' \leq 2$. Then $b'=2q-2-b = q-e'-1 < q-e-1 \leq q-[q]_{y}-1$. By Theorem~\ref{thm:M:b<x+1}(3), $M_{b'} < ({1}/{e}+{2}/{y})(q-e)<2q$. Thus we may write $\pm [K]=c \mu \pm f \mu'$ with $0\leq c \leq 2, 0 \leq f <({1}/{e}+{2}/{y})(q-e)<2q$.  By Lemma~\ref{lem:c<1,f<1,q<4y:1-nd}, we may assume that $c=2$. Applying Theorem~\ref{thm:c+d}(2), we consider the alternatives $(a)$--$(d)$, as $(3)$ verifies the Lemma.

Assume $(a)$ occurs. Then $2q-2n=f<({1}/{e}+{2}/{y})(q-e)<10q/9$ as $y > 19$. That is,  $8(nd+\epsilon)/9<2n$ contradicting that $n\geq 2, \,d \geq 3$.

Assume $(b)$, $(c)$ occur. Then $\pm [K]=(d -4\zeta)q\mu$ where $\zeta=\pm 1$. By Lemma~\ref{lem:dne2an>2}, we have $y \leq d+4$ and $y' \leq 4n+1$. From $y, y' > 19$, we further know that $d>15, \,n>4$.
Thus, as $d>15$ is odd, $nd +\epsilon = q<4y \leq 4d+16$ implies that $n \leq 4$, contradicting $n>4$.

If $(d)$ occurs, then  $[K]=4q \mu$, so $y = 4$ as $q > 10$, which contradicts that $y > 18$.
\end{proof}

\begin{thm}\label{thm:q<4y:1-ndprime} Assume $L(q^2,m)$ as in case $(*q)$. If $75 < q <4y$, then
\[\pm [K]=q \mu, dq \mu.\]
\end{thm}

\begin{proof}
We apply Theorem~\ref{thm:q+y,m_b=1} to $\Gamma'$ and consider the alternatives. By Lemma~\ref{lem:q<4y,q<4y':1-nd}, we may assume that $q \geq 4y'$. This rules out Theorem~\ref{thm:q+y,m_b=1}(4) and implies that $y>y'$. As $y'=1$ implies that $\pm [K]=dq \mu$, we may assume $y>y'>1$. 

%$\spadesuit$ {\em (Doesn't the Engine Theorem give $\mathcal{T}, \mathcal{T'}$-levels for $b$? If so, we should probably state that---or decide it is not necessary for Theorem~\ref{thm:M:b<x+1}(3)(4).}

%$\blacklozenge$ {\em (In the statement of Initial Growth Theorem~\ref{thm:M:b<x+1}(3)(4), the $b$ was given by $\T$-levels. However, in the addenda and The Engine Theorem(s), $\Gamma_b$ are usually not necessarily $\T$-levels, as for any statements for $m_b$, we usually reach the maximum level before a $\T$-level so that $m_b$ wouldn't change. Thus we could specify that $\Gamma_b$ as a $\T$-level in the Initial Growth Theorem to make the statement stronger if we want, but I don't think we used this fact.)}

\noindent {\bf (1) There is a $b \geq 2q-y'-2$ with $m_b' \leq 3$.}

Then $b'=2q-2-b \leq y' <y<q-y=x$. By Lemma~\ref{thm:m:b<x}, as $b'<x$, $M_{b'} < {b'}/{y}+1<2$. That is $M_{b'} \leq 1$. Then we may write $\pm[K]=c \mu \pm f \mu'$ where $0 \leq c \leq 3<q/2$ and $0 \leq f \leq 1$. Corollary~\ref{cor:c+d:q| d-1<4/3 q}(2) now proves the Theorem.

\noindent {\bf (2) There is a $b \geq q + y'-1$ such that $M_b'<2q$.}

Then $b'=2q-2-b \leq q-y'-1 < q$. By Lemma~\ref{lem:q<4y}, $m_{b'} \leq 1$. Thus we can write $\pm[K]=c \mu \pm f \mu'$ with $0 \leq f \leq 1, \,0 \leq c < 2q$. Now Lemma ~\ref{lem:c<1,f<1,q<4y:1-nd}(2) proves the Theorem. 

\noindent {\bf (3) There is a $b \geq q +y'-1$ with $m_b' \leq 1$.}

By Theorem~\ref{thm:q<4y,e+1:1-nd}, we may assume $y'>e+1$. Then $b'=2q-2-b \leq q-y'-1 < q-e-2 \leq q - [q]_y -2$. By Theorem~\ref{thm:M:b<x+1}(3), $M_{b'} < ({1}/{e}+{2}/{y}) (q-e)$. Then we may write $\pm[K]=c \mu \pm f \mu'$ where $0 \leq c \leq 1, 0 \leq f < ({1}/{e}+{2}/{y}) (q-e)$. Now Lemma~\ref{lem:c<1,f<1,q<4y:1-nd}(1) proves the Theorem.
\end{proof}

% $\spadesuit$ {\em (The proof below is just the proof above using $\Gamma$ in place of $\Gamma'$ (unprimed in place of primed), using the first half of Theorem~\ref{thm:q+y,m_b=1} instead of the second, and Lemma~\ref{lem:c<1,f<1,q<4y':1-nd} in place of Lemma~\ref{lem:c<1,f<1,q<4y:1-nd}. We could just say that---though it is a pretty short proof and the bound on $q$ changes.)}

\begin{thm}\label{thm:q<4y':1-ndprime} Assume $L(q^2,m)$ as in case $(*q)$. If $99 < q <4y'$, then
\[\pm [K]=q \mu, 2q \mu, dq \mu.\]
\end{thm}

\begin{proof}
The proof is identical to the proof of Theorem~\ref{thm:q<4y:1-ndprime} with the roles of \break $\Gamma', \,x, \,y', \,m', \,M, \,e, \,c, \,dq\mu$ exchanged with $\Gamma, \,x', \,y, \,m, \,M', \,e', \,f, \,q\mu$. Also, Corollary~\ref{cor:c+d:q| d-1<4/3 q}(2), Lemma ~\ref{lem:c<1,f<1,q<4y:1-nd}(1),(2), Theorem~\ref{thm:M:b<x+1}(3) are replaced by Corollary~\ref{cor:c+d:q| d-1<4/3 q}(1), Lemma ~\ref{lem:c<1,f<1,q<4y':1-nd}(2),(1), Theorem~\ref{thm:M:b<x+1}(4), respectively.
% We apply Theorem~\ref{thm:q+y,m_b=1} to $\Gamma$ and consider the alternatives. By Lemma~\ref{lem:q<4y,q<4y':1-nd}, we may assume that $q \geq 4y$. This rules out Theorem~\ref{thm:q+y,m_b=1}(4) and implies that $y'>y$. As $y=1$ implies that $\pm [K]=q \mu$, we hereafter assume $y'>y>1$. 
%  $\spadesuit$ {\em (Doesn't the Engine Theorem give $\mathcal{T}, \mathcal{T'}$-levels for $b$? If so, we should probably state that---or decide it is not necessary for Theorem~\ref{thm:M:b<x+1}(3)(4).}
% \noindent {\bf I. There is a $b \geq 2q-y-2$ with $m_b \leq 3$.}
% Then $b'=2q-2-b \leq y <y'<q-y'=x'$. By Lemma~\ref{thm:m:b<x}, as $b'<x'$, $M_{b'}' < {b'}/{y'}+1<2$. That is $M_{b'}' \leq 1$. Then we may write $\pm[K]=c \mu \pm f \mu'$ where $0 \leq f \leq 3<q/2$ and $0 \leq c \leq 1$. Corollary~\ref{cor:c+d:q| d-1<4/3 q}(1) now proves the Theorem.
% \noindent {\bf II. There is a $b \geq q + y-1$ such that $M_b<2q$.}
% Then $b'=2q-2-b \leq q-y-1 < q$. By Lemma~\ref{lem:q<4y}, $m_{b'}' \leq 1$. Thus we can write $\pm[K]=c \mu \pm f \mu'$ with $0 \leq c \leq 1, 0 \leq f < 2q$. Now Lemma ~\ref{lem:c<1,f<1,q<4y':1-nd}(1) proves the Theorem. 
% \noindent {\bf IV. There is a $b \geq q +y-1$ with $m_b \leq 1$.}
% By Theorem~\ref{thm:q<4y',e'+1:1-nd}, we may assume $y>e'+1$. Then $b'=2q-2-b \leq q-y-1 < q-e'-2 \leq q - [q]_{y'} -2$. By Theorem~\ref{thm:M:b<x+1}(4), $M_b' < ({1}/{e'} + {2}/{y'}) (q-e')<2q$. Then we may write $\pm[K]=c \mu \pm f \mu$ where $0 \leq c <2q, 0 \leq f \leq 1$. Now Lemma~\ref{lem:c<1,f<1,q<4y':1-nd}(2) proves the Theorem.
\end{proof}

\subsubsection{$2q+\epsilon=nd$}
In this subsubsection of subsection~\ref{subsec:q<4yor4y'} we consider $L(q^2,m)$ when $(*2q)$ holds, where $2q =nd- \epsilon$. 

\begin{lem}\label{lem:d+2,q<4y',2q=nd-1}
Assume $228 < q=(nd-\epsilon)/2 < 4y'$ and $[K]$ is of one the following forms:
\begin{enumerate}
    \item $\pm [K]=(d- \zeta a) q \mu = \zeta(an -\zeta\epsilon) q \mu'$ for $a \in \{1,2\}$ with $\zeta=\pm 1$;
    \item $\pm [K]=(2d- \zeta) q \mu = \zeta(n -2\zeta\epsilon) q \mu'$  with $\zeta=\pm 1$.
\end{enumerate}
Then
\[[q]_{y'}  \geq y.\]
\end{lem}

\begin{proof}
From $228 < q < 4y'$ we know that $y' > 57$.

\textbf{Assume $(1)$ holds.} We apply Lemma~\ref{lem:nupmlambda} with $\nu=1,\,h=-\zeta a,\,\rho=-\zeta$.

If $n\le 9$, then $d>9\ge 4a+1$ as $q>228$, and Lemma~\ref{lem:nupmlambda}$(2)(a)$ says that $y'= 3a - \zeta\epsilon\le 7$, contradicting $q < 4y'$. We hereafter assume $n > 9$ odd.

First assume $d\ge 4a+1$ and $n\ge 5$. By Lemma~\ref{lem:nupmlambda}$(1)(a)$ and $(2)(a)$, we have 
\[y=d - \zeta a, \quad y'= an - \zeta\epsilon.\]
Then $q<4y'$ becomes $nd-\epsilon<8an-8\zeta\epsilon$. That is, $(d-8a)n < (-8\zeta +1)\epsilon \leq 9$, so $d-8a \le -1$ as $n>9$ and $d$ is odd. Thus, $4a<d \leq 8a-1$. That is, when $a=1$, $d \in \{5, 7\}$ and when $a=2$, $d \in \{9, 11, 13, 15\}$. Note that $q = (nd-\epsilon)/2 > 228$. We have the following cases:
\begin{enumerate}
    \item $a=1, \,d=5$: We conclude $n>91, \,q=5n/2-\epsilon/2, \,y=5-\zeta, \,y'=n-\zeta\epsilon$, so $[q]_{y'} = q-2y'=(n+\epsilon(4\zeta-1))/2> 6\ge y$. % n >=17, q >= 43
    \item $a=1, \,d=7$: We conclude $n>65, \,q=7n/2-\epsilon/2, \,y=7-\zeta, \,y'=n-\zeta\epsilon$, so $[q]_{y'} = q-3y'=(n+\epsilon(6\zeta-1))/2> 8\ge y$. % n >=23, q >= 81
    \item $a=2, \,d=9$: We conclude $n>50, \,q=9n/2-\epsilon/2, \,y=9-2\zeta, \,y'=2n-\zeta\epsilon$, so $[q]_{y'} = q-2y'=(n+\epsilon(4\zeta-1))/2> 11\ge y$. % n >=27, q >= 122
    \item $a=2, \,d=11$: We conclude $n>41, \,q=11n/2-\epsilon/2, \,y=11-2\zeta, \,y'=2n-\zeta\epsilon$, so $[q]_{y'} = q-2y'=(3n+\epsilon(4\zeta-1))/2> 13\ge y$. % n > 9, q > 50
    \item $a=2, \,d=13$: We conclude $n\ge 37, \,q=13n/2-\epsilon/2, \,y=13-2\zeta, \,y'=2n-\zeta\epsilon$, so $[q]_{y'} = q-3y'=(n+\epsilon(6\zeta-1))/2\ge 15\ge y$. % n >35, q > 228
    \item $a=2, \,d=15$: We conclude $n>30, \,q=15n/2-\epsilon/2, \,y=15-2\zeta, \,y'=2n-\zeta\epsilon$, so $[q]_{y'} = q-3y'=(3n+\epsilon(6\zeta-1))/2> 17\ge y$. % n >13, q > 98
\end{enumerate}

Assume $4a-1\ge d\ge 2a+1$ and $n\ge 5$. By Lemma~\ref{lem:nupmlambda}$(1)(a)$ and $(2)(b)$, we have 
\[y=d - \zeta a, \ y'= ((d - 2a)n - \epsilon(1 - 2\zeta))/2.\]

Then, when $a=1$, $d=3$ and when $a=2$, $d \in \{5, 7\}$. Note that $nd=2q+\epsilon >455$. We have the following cases:

\begin{enumerate}
    \item $a=1, \,d=3$: We conclude $n>151, \,q=3n/2-\epsilon/2, \,y=3-\zeta, \,y'= n/2-(\epsilon(1-2\zeta))/2\ge (n-3)/2\ge 2y$, and $q=3y'-\epsilon\zeta y$. Thus, either $[q]_{y'}=y$ or $y'-y\ge y$. % n >=19, q >= 29
    \item $a=2, \,d=5$: We conclude $n>91, \,q=5n/2-\epsilon/2, \,y=5-2\zeta, \,y'= n/2-(\epsilon(1-2\zeta))/2\ge (n-3)/2$. Thus, $q-4y' = (n+\epsilon(3-8\zeta))/2 > 0$, so $q>4y'$, a contradiction. % n >=11, q >= 28
    \item $a=2, \,d=7$: We conclude $n>65, \,q=7n/2-\epsilon/2, \,y=7-2\zeta, \,y'=3n/2-(\epsilon(1-2\zeta))/2$, so $[q]_{y'} = q-2y'=(n-\epsilon(4\zeta-1))/2\ge 9\ge y$. % n >=23, q >= 81
\end{enumerate}

Finally, assume $d \leq 2a - 1$ and $n\ge 5$. Then $a=2, \,d=3, \,q = 3n/2-\epsilon/2, \,n> 151$. By Lemma~\ref{lem:nupmlambda}$(1)(a)$ and $(2)(c)$, we have $y=3 - 2\zeta$ and $y'=(n + \epsilon - 2\epsilon\zeta)/2> 10\ge 2y$. Note that $q=3y'+\epsilon\zeta y$. Thus, $[q]_{y'}=y$ or $y'-y$. In both cases, $[q]_{y'} \ge y$.

\textbf{Assume $(2)$ holds.} We apply Lemma~\ref{lem:nupmlambda} with $\nu=2,\,h=\rho=-\zeta$.

If $n\le 17$, then $d>5$ as $q>228$, and Lemma~\ref{lem:nupmlambda}$(2)(a)$ says that $y'= n-2\zeta\epsilon$, contradicting $q < 4y'$. We hereafter assume $n>17$ odd.

First assume $d=3$. Then $q=3n/2-\epsilon/2$ and $n>152$. From Lemma~\ref{lem:nupmlambda}$(1)(a)$ and $(2)(a)$ we conclude $y=6-\zeta\le 7, \,y'=(n+(4\zeta-1)\epsilon)/2> 73$. Note that $q=3y'-\epsilon\zeta y$, so $[q]_{y'} = y$ or $y'-y$. In both cases, $[q]_{y'}\ge y$.

Therefore, assume $d\ge 5$ and $n>17$. By Lemma~\ref{lem:nupmlambda}$(1)(a)$ and $(2)(a)$, we have 
\[y=2d - \zeta, \quad y'= n - 2\zeta\epsilon.\]
Then $q<4y'$ becomes $nd-\epsilon<8n-16\zeta\epsilon$. That is, $(d-8)n < (-16\zeta +1)\epsilon \leq 17$, so $d-8 \le -1$ as $n > 17$ and $d$ is odd. This together with $d\ge 5$ implies that $d = 5, 7$. Note that $q = (nd-\epsilon)/2 > 228$. We have the following cases:

\begin{enumerate}
    \item $d=5$: We conclude $n>91, \,q=5n/2-\epsilon/2, \,y=10-\zeta, \,y'=n-2\zeta\epsilon$, so $[q]_{y'} = q-2y'=(n+\epsilon(8\zeta-1))/2> 11\ge y$. % n >=31, q >= 78
    \item $d=7$: We conclude $n>65, \,q=7n/2-\epsilon/2, \,y=14-\zeta, \,y'=n-2\zeta\epsilon$, so $[q]_{y'} = q-3y'=(n+\epsilon(12\zeta-1))/2> 15\ge y$. % n >=43, q >= 156
\qedhere
\end{enumerate}
\end{proof}

\begin{lem}\label{lem:d+2,q<4y,2q=nd-1}
Assume $228 < q=(nd-\epsilon)/2 < 4y$ and $[K]$ is of one of the following forms:
\begin{enumerate}
    \item $\pm [K]=(n- \zeta a) q \mu' =  \epsilon\zeta(ad-\epsilon\zeta) q\mu$ for $a \in \{1,2\}$ with $\zeta=\pm 1$;
    \item $\pm [K]=(2n- \zeta) q \mu' =  \epsilon\zeta(d-2\epsilon\zeta) q\mu$ with $\zeta=\pm 1$.
\end{enumerate}
Then
\[[q]_{y} \geq y'.\] 
\end{lem}

\begin{proof}
We apply Lemma~\ref{lem:nupmlambda} with $\nu=a,\,h=\rho=-\epsilon \zeta$ in $(1)$ and $\nu=1,\,h=-2\epsilon\zeta, \rho=-\epsilon \zeta$ in $(2)$. Then we follow the proof of Lemma~\ref{lem:d+2,q<4y',2q=nd-1}, arguing analogously on the range of $n$, working with $[K]$ up to sign, with the roles of $y, \,d, \,[q]_{y'}, \,q\mu$ exchanged respectively with $y', \,n, \,[q]_{y}, \,q\mu'$,
and using Lemma~\ref{lem:nupmlambda}(1) in place of Lemma~\ref{lem:nupmlambda}(2) and vice versa.
\end{proof}

% \begin{lem}\label{lem:c<1,f<1,q<4y:2q+1}
% Assume $228 < q < 4y$ and $\pm [K]=c \mu \pm f \mu'$ with either
% \begin{enumerate}
% \item $0 \leq f \leq 1$ and $0 \leq c < 2q$; or
% \item $0 \leq c \leq 1$ and $0 \leq f < (1/e + 2/y)(q-e)$.
% \end{enumerate}

% Then, $\pm[K]=q \mu, dq \mu$. 
% \end{lem}

% \begin{proof}
% Assume $\pm[K]\neq q\mu, dq\mu$. Note that in Theorem~\ref{thm:c+d:2q+1}$(2)(b)$, $\eta=3, 4$ does not occur in the case of $(2)$ as $c\le 1, f<2q$. Similarly, in Theorem~\ref{thm:c+d:2q+1}$(3)(b)$, $\eta=3, 4$ does not occur in the case of $(1)$ as $f\le 1, c<2q$. Note that since we are working with $[K]$ up to sign, it suffices to prove the conclusion with a global factor of $\epsilon$ or $\zeta$. The proof is then identical to that of Lemma~\ref{lem:c<1,f<1,q<4y:1-nd}, replacing $(d- \zeta 2) q \mu$ with $(d- \zeta) q \mu$ (in case (i)), $(2d-2) q \mu$ with $(2d-1) q \mu$ (in case (iii)), $q=nd+\epsilon$ with $2q=nd-\epsilon$, and replacing Theorem~\ref{thm:q<4y,e+1:1-nd},~\ref{thm:c+d}, Lemma~\ref{lem:d+2,q<4y},~\ref{lem:dne2an>2}(II) with Theorem~\ref{thm:q<4y,e+1:2q+1},~\ref{thm:c+d:2q+1}, Lemma~\ref{lem:d+2,q<4y,2q=nd-1},~\ref{lem:nupmlambda}$(1)(a)$, respectively.
% \end{proof}

\begin{lem}\label{lem:c<1,f<1,q<4y:2q+1} Assume $228 < q < 4y$ and $\pm [K]=c \mu \pm f \mu'$ with either
\begin{enumerate}
\item $0 \leq c \leq 1$ and $0 \leq f < 2q$; or
\item $0 \leq f \leq 1$ and $0 \leq c < 2q$.
\end{enumerate}

Then, $\pm[K]=q \mu, dq \mu$. 
\end{lem}

\begin{proof} We follow the proof of Lemma~\ref{lem:c<1,f<1,q<4y:1-nd} with appropriate substitutions, and working with $[K]$ up to sign. By Theorem~\ref{thm:q<4y,e+1:2q+1}, we assume $y'>e+1$. We apply Theorem~\ref{thm:c+d:2q+1} so we may  assume either $(i)$ $[K]=(d \pm 1)q\mu$, $(ii)$ $[K]=(2d-1)q \mu$, or $(iii)$ $[K]=(d \pm 2)q \mu$.  Note $\eta=3, 4$ does not occur in Theorem~\ref{thm:c+d:2q+1}$(2)(b)$ as $c= 1, f<2q$. Similarly, in Theorem~\ref{thm:c+d:2q+1}$(3)(b)$, $\eta=3, 4$ does not occur as $f= 1, c<2q$.

In any of $(i)$, $(ii)$, $(iii)$, Lemma~\ref{lem:d+2,q<4y,2q=nd-1} gives the contradiction that $e \geq [q]_{y} \geq y' > e+1$. 
\end{proof}

\begin{lem}\label{lem:c<1,f<1,q<4y':2q+1}
Assume $228<q<4y'$ and $\pm [K]=c \mu \pm f \mu'$ with either
\begin{enumerate}
\item $0 \leq  c \leq 1$ and $0 \leq f < 2q$; or   
\item $0 \leq  f \leq 1$ and $0 \leq c < 2q$.
\end{enumerate}

Then, $\pm[K]=q \mu, dq \mu$.
\end{lem}

\begin{proof}
Notice the symmetry in the statements of Theorem~\ref{thm:c+d:2q+1}(2) and (3). The proof is then identical to that of Lemma~\ref{lem:c<1,f<1,q<4y:2q+1} working with $[K]$ up to sign, with the roles of $y, \,e, \,d, \,f, \,[q]_{y}, \,q\mu$ exchanged respectively with $y', \,e', \,n, \,c, \,[q]_{y'}, \,q\mu'$ 
and replacing Theorem~\ref{thm:q<4y,e+1:2q+1}, Lemma~\ref{lem:d+2,q<4y,2q=nd-1} respectively with Theorem~\ref{thm:q<4y',e'+1:2q+1}, Lemma~\ref{lem:d+2,q<4y',2q=nd-1}.
\end{proof}

\begin{lem}\label{lem:q<4y,q<4y':2q+1}
Assume $L(q^2,m)$ as in $(*2q)$. If $q <4y,4y'$, and $q > 228$, then $\pm [K]= q \mu, dq \mu$.
\end{lem}

\begin{proof} 
Assume $\pm [K]\neq q \mu, dq \mu$. We assume that $q > 228$. Then $y,y' > 57$.

With the same arguments as in the proof of Lemma~\ref{lem:q<4y,q<4y':1-nd} replacing Lemma~\ref{lem:c<1,f<1,q<4y:1-nd} with Lemma~\ref{lem:c<1,f<1,q<4y:2q+1}, we write $\pm[K]=c\mu\pm f\mu'$, with either $f=2$ and $0\le c<(1/e'+2/y')(q-e')<2q$ when $e\ge e'$, or $c=2$ and $0\le f<(1/e+2/y)(q-e)<2q$ when $e<e'$.

Assume $e\ge e'$. We apply Theorem~\ref{thm:c+d:2q+1}(3) and only consider (a) and (b) as $f=2$. Note that $\eta=4$ in (3)(b) implies $c=4q-2d > 2q$, a contradiction. Then, from Theorem~\ref{thm:c+d:2q+1}(3)(a)(b) we conclude that $y \equiv \pm(\eta \pm \epsilon 2d)\pmod{q}$ and $y' \equiv \pm (\eta n \pm 2) \pmod{q}$. Thus $y\le 2d+3$ and $y'\le 3n+2$. From $y,y' > 57$, we further know that $d> 27, \,n>18$. But this contradicts that $nd-\epsilon=2q<8y \leq 16d+24$.

The proof for $e<e'$ is identical to the above case, with the roles of $c, \,n, \,y$ exchanged with $f, \,d, \,y'$, replacing Theorem~\ref{thm:c+d:2q+1}(3) with Theorem~\ref{thm:c+d:2q+1}(2).
\end{proof}

\begin{thm}\label{thm:q<4y:2q+1} Assume $L(q^2,m)$ as in case $(*2q)$. If $228 < q <4y$, then
\[\pm [K]=q \mu, dq \mu.\]
\end{thm}

\begin{proof}
The proof is identical to that of Theorem~\ref{thm:q<4y:1-ndprime}, replacing 
Lemma~\ref{lem:q<4y,q<4y':1-nd},~\ref{lem:c<1,f<1,q<4y:1-nd}, Corollary~\ref{cor:c+d:q| d-1<4/3 q}(2),  Theorem~\ref{thm:q<4y,e+1:1-nd} 
with
Lemma~\ref{lem:q<4y,q<4y':2q+1},~\ref{lem:c<1,f<1,q<4y:2q+1}, Corollary~\ref{cor:c+d:d|2q+1<q/2}(2),  Theorem~\ref{thm:q<4y,e+1:2q+1}.
\end{proof}

\begin{thm}\label{thm:q<4y':2q+1} Assume $L(q^2,m)$ as in case $(*2q)$. If $228 < q <4y'$, then
\[\pm [K]=q \mu, dq \mu.\]
\end{thm}

\begin{proof}
The proof is identical to that of Theorem~\ref{thm:q<4y:2q+1}, working with $[K]$ up to sign, with the roles of $\Gamma', \,y, \,e, \,f, \,M, \,m, \,[q]_{y}$ exchanged respectively with $\Gamma, \,y', \,e', \,c, \,M', \,m', \,[q]_{y'}$, and replacing 
Lemma~\ref{lem:c<1,f<1,q<4y:2q+1}, Corollary~\ref{cor:c+d:d|2q+1<q/2}(2), Theorem~\ref{thm:q<4y,e+1:2q+1} 
with
Lemma~\ref{lem:c<1,f<1,q<4y':2q+1}, Corollary~\ref{cor:c+d:d|2q+1<q/2}(1), Theorem~\ref{thm:q<4y',e'+1:2q+1}.
\end{proof}

\section{Alexander polynomial width}\label{sec:alexwidth}

Recall the class of lens spaces given by part $(3)$ of Theorem~\ref{thm:lisca}: $L(p,m)$ where $p = q^2$, $m = d(q - \epsilon)$ for some $\epsilon = \pm 1$,  $q = nd + \epsilon$ with $d \ge 3$ odd and $n > 1$. This is the context of part $(4)$ of Theorem~\ref{thm:BBL} where it is shown that the simple knots $\pm K(p, m, \kappa)$ admit surgeries to $S^1 \times S^2$ when $\kappa=2q$ and $n=2$. This section is to prove a converse, Theorem~\ref{thm:q=1+nd,nmustbe2}, that only in these lens spaces where $n=2$ will $\pm K(p, m, \pm 2q)$ admit a surgery to $S^1 \times S^2$.
The argument of this section is quite different in flavor from those in the rest of the paper---more like those of \cite{cebanu2012propertyR} and \cite{greene2022fiberedSimpleKnots}. Let $K=K(p, m, 2q)$ admit a surgery to $S^1 \times S^2$. We define a function $g$ that generates a sequence whose generating function is the Alexander polynomial of $K$, multiplied by a monic polynomial. To show that $n=2$ we compare the Alexander polynomial's width obtained by two methods: analyzing the function $g$ and examining the Seifert surface bounded by $K$. 
Finding the width of the Alexander polynomial from $g$ is an issue of finding the minimum and maximum  values of $g$---to which subsections \ref{subsec:gresults} to \ref{subsec:gvalues} are devoted. In these subsections, we establish various symmetries of $g$. These allows us to reduce the problem of finding the extrema of $g$ to those over $[1,2q]$, which can then be done explicitly.

The above gives, more generally, a procedure for deciding if a simple knot in a lens space admits a surgery to $S^1 \times S^2$. In subsection~\ref{subsec:q<600}, we use this to verify computationally that Theorem~\ref{thm:Gr-Li-Lu} holds for small lens spaces ($q \le 600$). 

\subsection{Settings and the main results}\label{subsec:gresults}\hfill

Consider the simple knot $K=K(p, m, \kappa)$ with $p=q^2$. As in \cite{greene2022fiberedSimpleKnots}, let
\[l = [m^\ast \kappa]_p,\quad \theta = {p}/{\gcd(p, \kappa)},\quad t = {l}/{\gcd(p, l)} = {l}/{\gcd(p, \kappa)} = {l\theta}/{p},\]
where $\bk{\cdot}_p$ denotes the least positive residue mod $p$, and $m^\ast$ is defined such that $[mm^\ast ]_p = 1$.
Then, $\theta$ is the order of the simple knot $K$ in $H_1(L(p, m))$, $0 < t < \theta$, and $\gcd(t, \theta) = 1$. Let $X = \bc{1, 2, \dots, p} = \bc{1, 2, \dots, q^2}$ and $S \sse X$ denote the values
\[\bc{[m]_p, [2m]_p, \dots, [lm]_p}.\]
Denote $X\sm S$ by $S^c$. Let $\delta: X\to \Z$ be
\[\delta(y) =
    \begin{cases}
      t - \theta, & y\in S,\\
      t, & y\not\in S.
    \end{cases}\]
Then, define the unique function $g: \bc{0, 1, \dots, q^2}\to \Z$ obtained by $g(0) = 0$ and
\[g(x) = g(x - 1) + \delta(x),\ 0 < x \le q^2,\]
which we extend as a periodic function on $\Z$, noting that $g(q^2) = t\abs{S^c} + (t - \theta)\abs{S} = t(p - l) + (t - \theta)l = 0 = g(0)$. Note also that the $g$ we define here differs by an additive constant from the function $g$ in \cite{greene2022fiberedSimpleKnots}. We will establish the following properties of $g$ in subsection~\ref{subsec:gvalues}.

\newcounter{extremathmq}
\setcounter{extremathmq}{\value{thm}}
\begin{thm}\label{thm:extrema,q=1+nd}
    Let $K=K(q^2, m, 2q)$, where $q = nd + \epsilon$, $m = d(q-\epsilon)$, $\epsilon = \pm 1$, $d \ge 3$ is odd, and $n >1$ is even. %Assume a longitudinal surgery on $K$ yields $S^1\times S^2$. 
    %Assume $m = d(q-\epsilon)$ and $q = \epsilon+nd$ with $\epsilon = \pm 1$, $n>0$ even, and $d\ge 3$ odd. 
    The associated function $g$ has extrema given by:
    \[g_\MM = q-2,\quad g_\mm = -q-2n+4,\quad g_\MM - g_\mm = 2q + 2n - 6.\]
\end{thm}

At the same time, the Alexander polynomial width of simple knots with longitudinal $S^1\times S^2$ surgery is highly constrained.

\begin{lem}\label{lem:widthbyhom}
Let $K$ be a simple knot in $L(q^2,m)$ that admits a longitudinal surgery to $S^1\times S^2$. 
%The width of the generating function for $(g(i))_{i=0}^{q^2-1}$ is $2q - 2$. That is, 
Let $g$ be the function associated to $K$ at the beginning of subsection~\ref{subsec:gresults}. Then
$g_\MM-g_\mm = 2q - 2$.
\end{lem}

\begin{proof}
Assume $K$ admits a longitudinal surgery to $S^1\times S^2$. By Corollary~\ref{cor:Korder}, the exterior of $K$, $X_K$, fibers over the circle with fiber, $Q$, a planar surface with $q$ boundary components; furthermore, $[K] \in H_1(L(q^2,m))$ has torsion $\theta=q$. If $M$ is the matrix corresponding to the monodromy acting on $H_1(Q)$, then the Alexander polynomial is the characteristic polynomial of $M$. Since the rank of $H_1(Q)$ is $q - 1$, the width of the Alexander polynomial $\Delta_K$ is $q - 1$. Assume in addition $K$ is simple in $L(q^2,m)$. By \cite{greene2022fiberedSimpleKnots}, the generating function of the sequence $(g(i))_{i=0}^{q^2-1}$ equals $\Delta_K$, multiplied by a monic polynomial
    \[\sum_{i = 0}^{q^2 - 1}T^{g(i)} \sim \frac{T^{\theta} - 1}{T - 1}\cdot \Delta_K,\]
    where $\theta = q$ is the torsion of $K$, and the symbol $\sim$ denotes the equality in $\Z[T, T^{-1}]$ up to multiplication by a unit. Therefore,
    \[
    g_\MM - g_\mm = \theta - 1 + \text{width}(\Delta_K) = 2q - 2. \qedhere
    \]
   % In particular, a knot in a lens space with $S^1\times S^2$ surgery, it is fibered with fiber a $2$-sphere with $\abs{q}$ punctures by Corollary~\ref{cor:braid}, so $g_\MM-g_\mm$ is $2\abs{q} - 2$.
\end{proof}

We are now in a position to prove Theorem~\ref{thm:q=1+nd,nmustbe2}, the main theorem of this section.

\begin{thm}\label{thm:q=1+nd,nmustbe2}
    Let $K=K(q^2, m, 2q)$, where $q = \epsilon + nd$, $m = d(q-\epsilon)$, $\epsilon = \pm 1$, $d \ge 3$ odd, and $n > 1$. If a longitudinal surgery on $K$ yields $S^1\times S^2$, $n$ must be $2$. In particular, $K$ is listed in Theorem~\ref{thm:BBL}.
\end{thm}

\begin{proof}
As $K$ admits a longitudinal surgery to $S^1\times S^2$, Corollary~\ref{cor:Korder} shows that $[K] \in H_1(L(q^2,m))$ has torsion $q$. As $K=K(q^2, m, 2q)$, we conclude that $\gcd(q,2)=1$ and consequently that $n$ is even. 
By Theorem~\ref{thm:extrema,q=1+nd} and Lemma~\ref{lem:widthbyhom}, $g_\MM - g_\mm = 2q + 2n - 6 = 2q - 2$, which yields $n = 2$.
\end{proof}

\subsection{General properties of \pdfmath{g}{g}}\label{subsec:genpropofg}
To prove Theorem~\ref{thm:extrema,q=1+nd}, we first investigate several properties of the function $g$.
\subsubsection{$K(q^2, m, cq)$}

We begin by considering a general simple knot of the form $K(q^2, m, cq)$ with $0 < c < q$. Note that in this subsubsection, we \emph{do not} assume that $m$ has the form $d(q-\epsilon)$, but we do assume, as we may, $m<q^2$ and $m$ is relatively prime to $q$.

First, we establish a criterion for verifying whether an element in $X$ belongs to $S$. Recall that $\bk{\cdot}_{q^2}$ denotes the least positive residue mod $q^2$. For simplicity, {\bf we denote} $[\cdot]_{q^2}$ by $[\cdot]$. Note that $[\cdot]$ takes values in $X = \bc{1, 2,\cdots, q^2}$ and that $0=g(0)=g(q^2)=g([0])$. 

\begin{ctr}\label{ctr:XinS}
    For any element $x\in X = \bc{1, 2,\cdots, q^2}$, $x\in S = \bc{[m], [2m], \cdots, [lm]}$ if and only if $[xm^\ast ]\le l$.
\end{ctr}

\begin{proof}
     $m^\ast $ is the multiplicative inverse of $m$ mod $q^2$.
\end{proof}

We denote $\overline{x} = [cq-x]$ for all $x\in X$ and observe the following properties of $S$. 

\begin{lem}~\label{lem:symmetry}
    The set $S\sm\bc{cq}$ is {\em symmetric about point $cq/2$}, i.e., $x\in S$ if and only if $\overline{x}\in S$, unless $x = cq$ or $\overline{x} = cq$.
\end{lem}

\begin{proof}
    Assume $x\in S\sm\bc{cq}$. Note that $[cqm^\ast] = l$, so $1\le [xm^\ast]\le l - 1$ by Criterion~\ref{ctr:XinS}. Moreover, $[\overline{x}m^\ast] = [cqm^\ast  - xm^\ast] = l - [xm^\ast] < l$. Thus, $\overline{x}\in S$. Conversely, assume $\overline{x}\in S\sm\bc{cq}$. We then apply the same arguments as above to $\overline{x}$.
\end{proof}

\begin{lem}\label{lem:plusm}
    Assume $x\in X$ with $x\neq q^2, cq$. Then, $x \in S$ if and only if $[x + m] \in S$. In the exceptional cases, $q^2\not\in S$ while $[q^2 + m]\in S$, $cq\in S$ while $[cq + m]\not\in S$.
\end{lem}

\begin{proof}
    If $x=q^2$, $[xm^\ast]=[0]=q^2>l$ and $[mm^\ast]=1<l$, so by Criterion~\ref{ctr:XinS}, $q^2\not\in S$ but $[q^2+m]\in S$. If $x=cq$, $[xm^\ast]=[cqm^\ast]=l$ and $[(cq+m)m^\ast]=l+1 > l$, so by Criterion~\ref{ctr:XinS}, $cq\in S$ but $[cq+m]\not\in S$.
    
    Now assume $x \neq q^2, cq$, so $[xm^\ast]\neq q^2, l$. Thus, $[x+m] \in  S$, iff $[xm^\ast+1]=[(x+m)m^\ast] \leq l$, iff $[xm^\ast] \leq l-1$, iff $x \in S$.
\end{proof}

\begin{lem}\label{lem:periodicity}
    The function $g$ follows the periodic pattern that for all $cq \le x < q^2$,
    \[g([x + m]) - g([x]) = \xi,\]
    where $\xi$ is a fixed integer.
\end{lem}

\begin{proof}
    By Lemma~\ref{lem:plusm}, $\delta([y]) = \delta([y + m])$, for all $y\not\in\bc{cq, q^2}$. In particular, $\delta([cq + 1]) = \delta([cq + m + 1])$. Let $\xi = g([cq + m]) - g([cq])$. Then,
    \begin{align*}
        \xi &= g([cq + m]) - g([cq]) = g([cq + m]) + \pth*{\delta([cq + 1]) - g([cq + 1])}\\
        &= \pth*{g([cq + m]) + \delta([cq + 1 + m])} - g([cq + 1]) = g([cq + 1 + m]) - g([cq + 1]).
    \end{align*}
    By analogous computations, $g([x + m]) - g([x]) = \xi$, for all $cq \le x < q^2$.
\end{proof}

\begin{lem}\label{lem:nonperiodicity}
    For all $0\le x < cq$,
    \[g([x + m]) - g([x]) = \xi - \theta.\]
\end{lem}

\begin{proof}
    For any $0 < x < cq$, by Lemma~\ref{lem:plusm}, $\delta([x+m]) = \delta([x])$. Applying the same arguments as in Lemma~\ref{lem:periodicity},
    \[g([x+m]) - g([x]) = g([x+m-1]) - g([x-1]) = \dots = g([m]) - g([0])=g([m])-g([q^2]).\]
    As $m \in S$ and $q^2 \notin S$, $\delta([m]) = t-\theta$ and $\delta([q^2])=t$, so $\delta([m]) - \delta([q^2]) = -\theta$. Thus,
    \begin{align*}
    g([m]) - g([q^2]) &= \pth*{g([m-1])+\delta([m])} - \pth*{g([q^2-1])+\delta([q^2])} \\ &= \pth*{g([q^2+m-1])-g([q^2-1])}-\theta \\ &= \xi - \theta,
    \end{align*}
    the last equality by Lemma~\ref{lem:periodicity}.
\end{proof}

\begin{lem}~\label{lem:epsilon}
    The function $g$ follows the periodic pattern that for $x$ in $X$,
    \[g([x + m]) - g([x]) =
        \begin{cases}
            (c-q)/\gcd\pth{q, c}, & x = q^2\text{ or }x< cq,\\
            c/\gcd\pth{q, c}, & cq\le x < q^2.
        \end{cases}\]
\end{lem}

\begin{proof}
    We determine the value of $\xi$ in Lemmas~\ref{lem:periodicity} and~\ref{lem:nonperiodicity}. 
    \begin{align*}
        0 = \sum_{x = 0}^{q^2 - 1}\pth*{g([x+m]) - g([x])} &= \sum_{x = 0}^{cq - 1}\pth*{g([x+m]) - g([x])} + \sum_{x = cq}^{q^2 - 1}\pth*{g([x+m]) - g([x])}\\
        &= cq(\xi - \theta) + (q^2-cq)\xi = -cq\theta + \xi q^2
    \end{align*}

     The first equality follows as every value of $g$ appears exactly twice with opposite signs; the third equality follows from Lemmas~\ref{lem:periodicity} and~\ref{lem:nonperiodicity}.
    Thus, $\xi = c\theta/q = c/\gcd\pth{q, c}$, as $\theta = q^2/\gcd\pth{q^2, cq} = q/\gcd\pth{q, c}$. The desired result now follows from Lemmas~\ref{lem:periodicity} and~\ref{lem:nonperiodicity}.
\end{proof}

We then obtain the following corollary.
\begin{cor}\label{cor:extremarange}
    $g$ always attains its maximum at $[x]$ for some $0\le x < cq$ and its minimum at $[x]$ for some $m \le x < m+cq$.
\end{cor}

% $\spadesuit$ {\em (This Corollary would be easier to state if we thought of $g$ as defined over Z.)}

\subsubsection{Specializing to $K(q^2, d(q-\epsilon), 2q)$.}

We hereafter assume we are in the context of Theorem~\ref{thm:lisca}(3) where $m=d(q-\epsilon)$ with $\epsilon=\pm 1, \,q=nd+\epsilon, \,d \geq 3$ odd, $n>1$. We specialize to when $c=2$, that is to the simple knots $K=K(q^2, d(q-\epsilon), 2q)$. We furthermore specialize to the case where {\bf $n$ is even} (which follows if $K$ admits a longitudinal surgery to $S^1\times S^2$ since it will have torsion $q$). %. By Corollary~\ref{cor:Korder}, $K$ has torsion $q$. In the terminology above, $\theta=q$. We conclude that $\gcd(q,2)=1$ and consequently {\bf that $n$ is even}. 

\begin{lem}\label{lem:values,q=1+nd}
    For $K=K(q^2, d(q-\epsilon), 2q)$, 
    \[m^\ast = [n(2\epsilon q + 1)],\ \theta = q.\]
    \[l = \abs{S} = 2nq,\ \abs{S^c} = (q-2n)q,\ t = 2n,\ t - \theta = 2n-q.\]
\end{lem}

\begin{proof}
    We have $n(2\epsilon q + 1)\cdot m = n(2\epsilon q + 1)d(q - \epsilon) \equiv -nd(q + \epsilon) = (\epsilon-q)(\epsilon+q) \equiv 1\Mod{q^2}$. Thus, $m^\ast = [n(2\epsilon q + 1)]$.
    
    Recall that $\theta=q/\gcd(q, 2)=q$ and $d\ge 3$. Then $\abs{S} = l = [m^\ast \kappa]= [n(2\epsilon q + 1)\cdot 2q] = 2nq$, so $\abs{S^c} = (q-2n)q$. Moreover, $t = {l\theta}/{q^2} = {2nq^2}/{q^2} = 2n$, $t - \theta = 2n - q$.
\end{proof}

The following theorem will be useful in subsection~\ref{subsec:gvalues}.

\begin{lem}\label{lem:gproperties}
The function $g$ has the following properties.
    \begin{enumerate}[label={(\arabic*)}]
        \item $g(x) = -g(2q - 1 - x) + g(q - 1) + g(q)$ for all $q+1\le x\le 2q-1$.
        \item If $\abs{\bc{[x + 1], [x + 2], \cdots, [x + d]}\cap S} = 2$, then $g([x + d]) - g([x]) = -2\epsilon$.
    \end{enumerate}
\end{lem}

\begin{proof}
    \begin{enumerate}[label={(\arabic*)}]
        \item By Lemma~\ref{lem:symmetry}, for $1\le j\le q-1$, $g(q - 1) - g(q - j - 1) = \delta(q-j) + \dots + \delta(q-1) = \delta(q+1) + \dots + \delta(q+j) = g(q + j) - g(q)$, so $g(q + j) = -g(q - j - 1) + g(q - 1) + g(q)$. In other words, $g(x) = -g(2q - 1 - x) + g(q - 1) + g(q)$ for all $q+1\le x\le 2q-1$.
        \item By direct computation, $g([x + d]) - g([x]) = \delta([x+1]) + \dots + \delta([x + d]) = (d-2)t + 2(t-\theta) = (d-2)2n + 2(2n-q) = 2nd - 2q = -2\epsilon$.
        \qedhere
    \end{enumerate}
\end{proof}

\subsection{Extrema of \pdfmath{g}{g}}\label{subsec:gvalues}\hfill

Corollary~\ref{cor:extremarange} shows that $g$ takes its maximum when $0\le x < 2q$ and its minimum when $m\le x < m + 2q$. We find these extrema by analyzing $g(x)$ with $0\le x< 2q$ using Criterion~\ref{ctr:XinS}. 

\subsubsection{$q=nd+\epsilon$}

To begin, we determine whether certain elements are contained in $S$.

\begin{lem}\label{lem:eltsinS,q=1+nd}
    For all $1\le j\le n$,
    \[\pth{(2j - 1)d + \epsilon}/{2}\in S,\ jd + \epsilon\in S.\]
\end{lem}

\begin{proof}
Recalling that $nd=q-\epsilon$, for $1\le j\le n$, we have
    $[m^\ast\pth{(2j - 1)d + \epsilon}/{2}] = [n(2\epsilon q + 1)\pth{(2j - 1)d + \epsilon}/{2}]
    = [((2j - 1)(q - \epsilon) + \epsilon n)(\epsilon q + 1/2)]
    = [(n-j)q + (q+\epsilon(n-2j+1))/2] \le l$,
as $0 < q-n+1 \le q+\epsilon(n-2j+1) \le q+n-1 < 2q$ and $(n-j+1)q < 2nq = l$. Thus, by Criterion~\ref{ctr:XinS}, $\pth{(2j - 1)d + \epsilon}/{2}\in S$, $1\le j\le n$.

Furthermore, for $1\le j\le n$,
    $\bk*{(jd + \epsilon)m^\ast} = \bk*{(jd + \epsilon)n(2\epsilon q + 1)} = [(j(q-\epsilon) + \epsilon n)(2\epsilon q + 1)] = [2nq - (j(q+\epsilon)-\epsilon n)] < l$,
as $0\le q-n-1 \le j(q+\epsilon)-\epsilon n \le n(q+1)+n < 2nq = l$. Then, by Criterion~\ref{ctr:XinS} again, $jd + \epsilon\in S$.
\end{proof}

To bridge the gap between these elements, it turns out that the elements in $S$ are sparse:

\begin{lem}\label{lem:Xsparse,q=1+nd}
     Assume $d\ge 5$. 
     \begin{enumerate}
     \item Assume $\epsilon=1$. If $x \in S$, then $x$ is isolated on either side by $(d-3)/2$ elements of $X$ in $S^c$. That is, if $x\in S$, then for $1 \le j \le (d-3)/{2}$, $[x \pm j] \in S^c$.
     
    %\[[x - \pth{d-3}/{2}],\dots, [x - 1], [x + 1],\dots, [x + \pth{d-3}/{2}]\in S^c\cup [2q+1, q^2]\]
    \item Assume $\epsilon=-1$ and $x \in S \cap [1,2q]$. Then within the interval $[1,2q]$, $x$ is similarly isolated by elements of $S^c$. Specifically, for $1 \le j \le (d-3)/{2}$, 
        \begin{itemize}
            \item $[x-j] \in S^c$;
            \item $[x+j] \in S^c$ when $[x+j] \in [1,2q]$.
        \end{itemize}
    \end{enumerate}
    
\end{lem}

\begin{proof}
    % Assume $d\ge 5$, for otherwise the lemma holds.
    
    Recall that $l=2nq$ and $m^\ast=[n(2\epsilon q + 1)] = [\epsilon n(2q + \epsilon)]$. For $1\le j\le (d-3)/2$, to estimate $jm^\ast$, we have $l + n\epsilon \le jn(2q+\epsilon) \le ((d-3)/2)n(2q+\epsilon) = (q-3n-\epsilon)(q+\epsilon/2) = q^2 - ((3n+\epsilon/2)q+(3\epsilon n+1)/2) < q^2-2nq = q^2 - l$. That is
    \[l + n\epsilon \le jn(2q+\epsilon) < q^2-l.\]

    If $x\in S$, by Criterion~\ref{ctr:XinS}, $[xm^\ast] \le l$. We argue the Lemma depending on the sign $\epsilon$.

    \begin{enumerate}
        \item Assume $\epsilon = 1$. Then, for all $1\le j\le (d-3)/2$, using the above inequalities, we have
        \[[(x+j)m^\ast] = [xm^\ast] + jn(2q+1) > l,\]
        \[[(x-j)m^\ast] = [xm^\ast] + q^2 - jn(2q+1) > l.\]
        Therefore, $[x+j], [x-j]\in S^c$ by Criterion~\ref{ctr:XinS} as claimed in $(1)$.

        \item Assume $\epsilon = -1$. Then $\bc{[m],\dots, [nm]} = \bc{[id(q+1)] \,| \,1\le i\le n} = \{d(q+1),\dots, (n-1)d(q+1), 2q+1\}$, which does not contain any point in $[1, 2q]$. Thus, for any $z\in [1, 2q]$, $[z m^\ast] > n$. Assume $x \in S \cap [1,2q]$. Then $n<[x m^\ast] \le l$. With $1\le j\le (d-3)/2$, recall from above that
        \[l-n \le jn(2q-1) <q^2-l.\]  
        Thus,
        \[[(x-j)m^\ast]= [xm^\ast] + jn(2q-1) > l.\]
    We conclude $x-j \in S^c$ as claimed in $(2)$.

    Similarly, $q^2+n \ge [xm^\ast] +q^2-jn(2q-1) > l+n$. If furthermore, $[x+j] \in [1,2q]$, hence $[(x+j)m^\ast]>n$, then
    \[[(x+j)m^\ast] = [xm^\ast] + q^2 - jn(2q-1) > l.\]
    We conclude $x+j \in S^c$, also claimed in $(2)$. \qedhere
    \end{enumerate}
\end{proof}

\begin{lem}\label{lem:eltsinSc,q=1+nd,d=3}
    Assume $d = 3$. Then, for all $1\le j\le n-1$, $3j\in S^c$. If $\epsilon = 1$, $3n\in S^c$.
\end{lem}

\begin{proof}
Since $d = 3$, $q = 3n + \epsilon$. Therefore, for $1\le j\le n-1$, $l = 2nq < 2nq + (1+\epsilon)q-\epsilon(n-1) = 2nq+(q-3n+1)q-\epsilon(n-1) = q^2-(n-1)(q+\epsilon) \le q^2-j(q+\epsilon) < q^2$. Thus, $[3jm^\ast] = [3jn(2\epsilon q+1)] = [j(q-\epsilon)(2\epsilon q+1)] = [-j(q + \epsilon)] = q^2 - j(q + \epsilon) > l$. By Criterion~\ref{ctr:XinS}, $3j\in S^c$. Moreover, if $\epsilon = 1$, then $[3nm^\ast] = [(q-1)n(2q+1)] = q^2-nq-n = 2nq+q-n > 2nq = l$, so $3n\in S^c$.
\end{proof}

Combining Lemma~\ref{lem:eltsinS,q=1+nd},~\ref{lem:Xsparse,q=1+nd}, and~\ref{lem:eltsinSc,q=1+nd,d=3}, we obtain

\begin{lem}\label{lem:gatfirstq,q=1+nd}
    Lemma~\ref{lem:eltsinS,q=1+nd} yields all elements in $S$ among the first $q$ elements in $X$, i.e.,
    \[
    S\cap[1, q] = \bc{\pth{(2j - 1)d + \epsilon}/{2}, \quad 1\le j\le n; \quad jd + \epsilon, \quad1\le j\le n}.
    \]
\end{lem}

\begin{proof}
There are either $(d-3)/2$ or $(d-1)/2$ elements of $X$ between the elements of $S$ listed in Lemma~\ref{lem:eltsinS,q=1+nd}. Assume $d \ge 5$. By Lemma~\ref{lem:Xsparse,q=1+nd}, for any of the listed elements in $S$, both the previous $(d - 3)/2 \geq 1$ elements and the next $(d - 3)/2$ elements lie in $S^c$. This shows that all elements between the values listed in Lemma~\ref{lem:eltsinS,q=1+nd} are in $S^c$ except possibly the value $1$ when $\epsilon=1$. In this case, we see that $1 \in S^c$ since $[1m^\ast] = n(2q+1) > 2nq = l$. This proves the result when $d \geq 5$. When $d=3$, Lemma~\ref{lem:eltsinS,q=1+nd} lists all values that are not multiples of $3$ as in $S$.  Lemma~\ref{lem:eltsinSc,q=1+nd,d=3} shows that multiples of $3$ are not in $S$, showing that the list of elements of $S$ is complete.
\end{proof}

Let $\id_S: X\to \bc{0, 1}$ be the indicator function
\[\id_S(x) =
    \begin{cases}
      1, & x\in S,\\
      0, & x\not\in S.
    \end{cases}\]
The values of $\id_S$ for the first $q$ points in $X$ are given in Table~\ref{tab:q=1+nd,d>=3}.
\begin{table}[ht!]
  \centering
  \resizebox{1\textwidth}{!}{
  \begin{tabular}{ccccccccccc}
    \hline\hline
    \multicolumn{10}{c}{$\id_S$ at the first $q$ elements ($\epsilon = 1$)} \\
    \hline
    $x$   & & & & & & & & & $1$\\
    $\id_S$& & & & & & & & & $0$\\
    \hline
    $x$    & $2$ & $\cdots$ & $\frac{d-1}{2}$ & $\frac{d+1}{2}$ & $\frac{d+3}{2}$ & $\cdots$ & $\cdots$ & $d$ & $d + 1$ \\
    $\id_S$ & $0$ & $\cdots$ & $0$ & $1$ & $0$ & $\cdots$ & $\cdots$ & $0$ & $1$\\
    \hline
    $x$    & $d+2$ & $\cdots$ & $\frac{3d-1}{2}$ & $\frac{3d+1}{2}$ & $\frac{3d+3}{2}$ & $\cdots$ & $\cdots$ & $2d$ & $2d + 1$ \\
    $\id_S$ & $0$ & $\cdots$ & $0$ & $1$ & $0$ & $\cdots$ & $\cdots$ & $0$ & $1$\\
    \hline
    & $\vdots$ & $\vdots$ & $\vdots$ & $\vdots$ & $\vdots$ & $\vdots$ & $\vdots$ & $\vdots$ & $\vdots$\\
    \hline
    $x$    & $(n-1)d+2$ & $\cdots$ & $\frac{(2n-1)d-1}{2}$ & $\frac{(2n-1)d+1}{2}$ & $\frac{(2n-1)d+ 3}{2}$ & $\cdots$ & $\cdots$ & $nd$ & $nd + 1 = q$ \\
    $\id_S$ & $0$ & $\cdots$ & $0$ & $1$ & $0$ & $\cdots$ & $\cdots$ & $0$ & $1$\\
    \hline\hline
    \multicolumn{10}{c}{$\id_S$ at the first $q$ elements ($\epsilon = -1$)} \\
    \hline
    $x$    & & $1$ & $\cdots$ & $\frac{d-3}{2}$ & $\frac{d-1}{2}$ & $\frac{d+1}{2}$ & $\cdots$ & $d-2$ & $d-1$ \\
    $\id_S$ & & $0$ & $\cdots$ & $0$ & $1$ & $0$ & $\cdots$ & $0$ & $1$\\
    \hline
    $x$    & $d$ & $d+1$ & $\cdots$ & $\frac{3d-3}{2}$ & $\frac{3d-1}{2}$ & $\frac{3d+1}{2}$ & $\cdots$ & $2d - 2$ & $2d - 1$ \\
    $\id_S$ & $0$ & $0$ & $\cdots$ & $0$ & $1$ & $0$ & $\cdots$ & $0$ & $1$\\
    \hline
    & $\vdots$ & $\vdots$ & $\vdots$ & $\vdots$ & $\vdots$ & $\vdots$ & $\vdots$ & $\vdots$ & $\vdots$\\
    \hline
    $x$    & $(n-1)d$ & $(n-1)d + 1$ & $\cdots$ & $\frac{(2n-1)d-3}{2}$ & $\frac{(2n-1)d-1}{2}$ & $\frac{(2n-1)d+ 1}{2}$ & $\cdots$ & $nd - 2$ & $nd - 1 = q$ \\
    $\id_S$ & $0$ & $0$ & $\cdots$ & $0$ & $1$ & $0$ & $\cdots$ & $0$ & $1$\\
    \hline\hline
  \end{tabular}}
  \caption{$\mathbf{1}_S$ values for the first $q$ elements in $X$}
  \label{tab:q=1+nd,d>=3}
\end{table}

\subsubsection{$q=nd+1$} We compute the extrema of $g$ for $\epsilon = 1$.

\begin{lem}\label{lem:localextrema,q=1+nd}
    Among $\bc{1, 2, \cdots, 2q - 1}\cup\bc{q^2}\sse X$, the maximum value attained by $g$ is $g(d) = q - 2$, and the minimum value is $g(2q - d - 1) = 2-2n$.
\end{lem}

\begin{proof}
    As shown in the upper half of Table~\ref{tab:q=1+nd,d>=3}, $\id_S(x+d)=\id_S(x)$ for any $2 \le x \le q-d=(n-1)d+1$. This fact along with Lemma~\ref{lem:gproperties}(2) implies that $g(x) = g(x+d) + 2 > g(x+d)$ for  $1 \le x \le (n-1)d+1$. Note also that $g(1)=g(0)+t=g(q^2)+t=t$. 
    Therefore, for $[x] \in \bc{q^2, 1, 2, \cdots, q}$, $g$ attains its maximum within $\bc{1, 2, \dots, d}$. On the other hand, $g(x+1)=g(x) + \delta(x+1)$ implies that if $x$ is a maximum for $g$, then $x+1$ must be contained in $S$. From Table~\ref{tab:q=1+nd,d>=3}, we conclude that a maximum occurs either when  $x = (d-1)/2$ or $x = d$. We know
    \[
    g\pth*{\frac{d - 1}{2}} = \sum_{y = 1}^{(d - 1)/2}\delta(y) = \frac{d - 1}{2}t = \frac{d - 1}{2}2n = nd - n = q - n - 1;
    \]
    \[
    g\pth*{d} = \sum_{y = 1}^{d}\delta(y) = 1(t - \theta) + (d - 1)t = (2n-q) + (d - 1)2n = 2nd - q = q - 2.
    \]
    Since $n\ge 2$, the maximum is $g(d) = q - 2$ among $\bc{q^2, 1, 2, \cdots, q}$.
    
    Similarly, among $\bc{q^2, 1, 2, \cdots, q}$, $g$ must take its minimum within $\bc{(n-1)d+2, \dots , nd+1=q}$ (the calculation below shows that $0=g(q^2)$ is not a minimum for $g$). Also, $g(x) = g(x - 1) + \delta(x)$ implies that if $g$ attains a minimum at $x$ then
    $x$ must be contained in $S$. From Table~\ref{tab:q=1+nd,d>=3}, we conclude $g$ takes its minimum at $x = ((2n-1)d+1)/2$ or $x = q$.  From Table~\ref{tab:q=1+nd,d>=3} and the fact that $g(x+d) -g(x) = -2$ for all $1 \le x \le (n-1)d+1$,
    \[g(q) = \sum_{y = 1}^{q}\delta(y) = 2n(t-\theta) + (q - 2n)t = (2n)(2n-q) + (q - 2n)2n = 0;\]
    \[
    \begin{aligned}
    g\pth*{\frac{(2n - 1)d + 1}{2}} &= g\pth*{\frac{d - 1}{2}} + (t-\theta) + (n-1)(-2)\\
    &= (q-n-1) + (2n-q) + (2-2n) = 1 - n.
    \end{aligned}
    \]
    As $n \ge 2$, $1 - n < 0$. Thus, among $\bc{q^2, 1, 2, \cdots, q}$, the minimum is $ g(\pth{(2n - 1)d + 1}/{2}) = 1 - n$.
    
    Finally, consider the extrema of $g$ on $\bc{q^2, 1, 2, \cdots, 2q - 1}$. Note that  $g(q - 1) = g(q) - (t-\theta) = q-2n$.
    By Lemma~\ref{lem:gproperties}(1) then, $g(x) = -g(2q - 1 - x) + q - 2n$ for all $q+1\le x\le 2q-1$. 
    Thus, by comparing the maximum of $g$ on the interval $[0,q]$,  $g(d) = q-2$,  and that on $[q+1,2q-1]$, $g(q+(d-1)/2)= -g(((2n - 1)d + 1)/2)+q-2n= n-1+q-2n = q-n-1$, we conclude that the maximum value of $g$ on $[0,2q-1]$ is $g(d) = q-2$. Similarly, comparing the minimum value of $g$ on $[0,q]$, $g(((2n - 1)d + 1)/2) = 1-n$, with the minimum on $[q+1,2q-1]$, $g(2q-d-1) = -g(d) + q-2n = 2-q+q-2n = 2-2n$, we conclude that the minimum value of $g$ on  $[0,2q-1]$ is $g(2q-d-1) = 2-2n$.
\end{proof}

\newcounter{savedthm2}
\setcounter{savedthm2}{\value{thm}}
\setcounter{thm}{\value{extremathmq}}
We now proceed to prove
\begin{thm}[$\epsilon = 1$]
    Assume $m = d(q-1)$, $q = nd+1$ with $n>1$ even and $d\ge 3$ odd. The function $g$ has extrema given by:
        \begin{enumerate}[label={(\arabic*)}]
          \item $g_\MM = g(d) = q-2$.
          \item $g_\mm = g(m+2q-d-1) = -q-2n+4$.
          \item $g_\MM - g_\mm = 2q + 2n - 6$.
        \end{enumerate}
\end{thm}

\begin{proof}
    Corollary~\ref{cor:extremarange} states that $g(x)$ must achieve its maximum at $[x]$ when $0\le x < 2q$ and its minimum when $m \le x < m + 2q$. Lemma~\ref{lem:localextrema,q=1+nd} shows that
    \[g_\MM = (g|_{[0, 2q)})_\MM = g(d) = q - 2.\]
    Lemma~\ref{lem:epsilon} says that for $x = q^2\text{ or }x< 2q$, $g([x + m]) - g([x]) = 2-q$. This implies that $(g|_{[2q,m+2q)})_\mm=(g|_{[0,2q)})_\mm+2-q$. Thus
    \[g_\mm = (g|_{[m, m+2q)})_\mm = (g|_{[0, 2q)})_\mm + 2-q = g(2q-d-1) + 2 - q =  -q - 2n + 4.
    \qedhere
    \]
\end{proof}
\setcounter{thm}{\value{savedthm2}}

\subsubsection{$q=nd-1$} In this subsubsection, we consider the complementary case where $\epsilon = -1$.

\newcounter{savedthm3}
\setcounter{savedthm3}{\value{thm}}
\setcounter{thm}{\value{extremathmq}}
\begin{thm}[$\epsilon = -1$]
    Assume $m = d(q+1)$, $q = nd - 1$, with $n>1$ even and $d\ge 3$ odd. The function $g$ has extrema given by:
    \begin{enumerate}[label={(\arabic*)}]
        \item $g_\MM = g(2q-d) = q-2$.
        \item $g_\mm = g(m+d-1) = -q-2n+4$.
        \item $g_\MM - g_\mm = 2q + 2n - 6$.
    \end{enumerate}
\end{thm}

\begin{proof}
We begin following the proof of Lemma~\ref{lem:localextrema,q=1+nd}, modifying as necessary. We first find the extrema of $g$ among the first $q$ elements. Using the lower half of Table~\ref{tab:q=1+nd,d>=3}, Lemma~\ref{lem:gproperties}  implies that $g(x+d)=g(x)+2$ for $0 \le x \le (n-1)d-1$. As argued in  Lemma~\ref{lem:localextrema,q=1+nd},
$g$ attains its maximum either at $q-1$ or at $((2n-1)d-3)/2$, and attains its minimum at either 
$d-1$ or $(d-1)/2$. We similarly compute that 
$g(q) = 0$; $g(q-1) = g(q)-\delta(q)=q-2n$; $g((d-3)/2)=(d-3)(2n)/2=(d-3)n$; $g((d-1)/2)=g((d-3)/2)+\delta((d-1)/2)=1-n$; $g(d-1) = (d-3)t + 2(t-\theta) = 2-2n$. Finally, $g(\pth{(2n-1)d-3}/{2})=g(\pth{d-3}/{2})+2(n-1)=q-n-1$.

We conclude that on $[0,q]$, the maximum value is $g(q-(d+1)/2)=g(((2n-1)d-3)/2)=q-n-1$ and the minimum is $g(d-1)=2-2n$.

Now consider the range $0\le x < 2q$. By Lemma~\ref{lem:gproperties}(1), $g(x) = -g(2q - 1 - x) + q - 2n$ for all $q+1\le x\le 2q-1$, so $g(q+(d-1)/2)) = -g(q-(d+1)/2) + q - 2n = 1-n > g(d-1)$, and $g(2q-d) = -g(d-1) + q - 2n = q-2 > g(q-(d+1)/2)$. Therefore, $(g|_{[0, 2q)})_\MM = g(2q-d) = q- 2$ and $(g|_{[0, 2q)})_\mm = g(d-1) = 2-2n$.
    
Finally, for the entire domain $X$, we follow the proof of Theorem~\ref{thm:extrema,q=1+nd} when $\epsilon=1$. Applying  Lemma~\ref{lem:epsilon} and Corollary~\ref{cor:extremarange}, we have 
\[g_\MM  = (g|_{[0, 2q)})_\MM = g(2q-d) = q-2,\]
\[
g_\mm = (g|_{[m, m+2q)})_\mm = g(m+d-1) = g(d-1) + 2-q = -q - 2n + 4.
\qedhere\]
\end{proof}
\setcounter{thm}{\value{savedthm3}}

\subsection{Verification for simple knots in small lens spaces}\label{subsec:q<600}
Lemma~\ref{lem:widthbyhom} gives a criterion that allows us to verify computationally that Theorem~\ref{thm:Gr-Li-Lu} holds for small lens spaces. 

%Let $K = K(q^2, m, cq)$ be a spherical simple knot in Lisca's form. By Corollary~\ref{cor:Korder}, $K$ has torsion order $\abs{q}$, which implies $\gcd\pth{c, q} = 1$. Furthermore, by Corollary~\ref{cor:braid}, $K$ is a spherical braid. This allows us to check whether $g_\MM-g_\MM = 2\abs{q} - 2$ holds for any given tuple $(q, m, c)$. Using a computer program, we verify that the knots listed in Theorem~\ref{thm:BBL} are complete for $q \le 600$.
\begin{thm}\label{thm:q<600}
Let $L(q^2,m)$ be in one of the forms stated in Theorem~\ref{thm:lisca}. Assume $q \leq 600$. If $K$ is a simple knot in $L(q^2,m)$ that admits a longitudinal surgery to $S^1\times S^2$, then it is listed in Theorem~\ref{thm:BBL}.
\end{thm}

\begin{proof} Let $K=K(q^2,m,\kappa)$ be a simple knot that admits a longitudinal surgery to $S^1\times S^2$. By Corollary~\ref{cor:Korder}, $[K]$ has order $q$ in first homology. By Lemma~\ref{lem:widthbyhom}, $g_\MM-g_\mm=2q-2$. For $q \le 600$, a MATLAB program verifies that both these conditions only hold for the knots listed in Theorem~\ref{thm:BBL}.
\end{proof}

\begin{remark}
Related computational verifications were done by J. Rasmussen for $q^2\le 100{,}000$ in \cite[Section~6.3]{rasmussen2007lensSpaceSurgeries} and by Cebanu for $q\le 500$ in \cite{cebanu2012propertyR}. A Python implementation and the complete dataset of computed extrema for $q\le600$ are available in \cite{Liang2026Verifier}.
\end{remark}

\section{The General Setting. Homology class detection}\label{sec:homologydetection}

In this section we gather for reference some results about knots in lens spaces with planar Seifert surfaces. Note that Dehn surgery along the boundary slope of the Seifert surface of such a knot, the knot's rational longitude, is $(S^1 \times S^2) \# N$ where $N$ is a rational homology sphere. 

A {\em Heegaard Floer homology solid torus} is a rational homology solid torus for which Dehn twisting along the rational longitude does not change the bordered Heegaard Floer invariants (e.g.  \cite[Definition 1.10]{rasmussenrasmussen2017floerSimpleManifolds} and \cite[Definition 26]{hanselmanrasmussenwatson2022torusBoundary}). 

We thank Jake Rasmussen and Tye Lidman for pointing out the results from \cite{rasmussenrasmussen2017floerSimpleManifolds} and \cite{alisha2025detectingHeegaardFloerSolidTori} that establish the following.

\newcounter{savedthm5}
\newcounter{savedsec3}
\setcounter{savedthm5}{\value{thm}}
\setcounter{savedsec3}{\value{section}}
\setcounter{thm}{\value{equivalence}}
\setcounter{section}{\value{secr>1}}

\begin{thm} Let $K$ be a knot in a lens space with irreducible exterior $X_K$. 
The following are equivalent:
\begin{enumerate}
    \item $K$ admits a surgery to $(S^1 \times S^2) \# N$ where $N$ is a rational homology sphere.
    \item $X_K$ is a generalized solid torus \cite[Definition 7.2]{rasmussenrasmussen2017floerSimpleManifolds}.
    \item $K$ has a planar Seifert surface.
    \item $X_K$ is a Heegaard Floer homology solid torus \cite[Definition 26]{hanselmanrasmussenwatson2022torusBoundary}.
\end{enumerate}
\end{thm}

\begin{proof} As $K$ lies in a lens space, $X_K$ is a rational homology solid torus. We may assume that $X_K$ has incompressible boundary, as all four statements hold if $X_K$ is a solid torus.

\noindent
$(1)\implies(2)$: Proposition 7.8 of \cite{rasmussenrasmussen2017floerSimpleManifolds} applied to the dual of the surgery, lying in $(S^1 \times S^2) \# N$, shows that $X_K$ is generalized solid torus.

\noindent
$(2)\implies(3)$: Two sentences of \cite{rasmussenrasmussen2017floerSimpleManifolds} following Definition 7.2, using Corollary 2.3, show that $K$ has a planar Seifert surface.
%Let $F$ be a minimal genus Seifert surface for $K$. We may assume $F$ is not a disk, hence $\partial X_K$ is incompressible. Definition 7.2 (of generalized solid torus) along with Corollary 2.3 then imply that $-\chi(F) \leq |F|-2$ and hence $F$ is planar. 

\noindent
$(3)\implies(1)$: Surgery along the rational longitude of $K$ produces a manifold with a non-separating $2$-sphere, hence of the form $(S^1 \times S^2) \# N$. $N$ must be a rational homology sphere as $K$ is a knot in a lens space. %with irreducible exterior. 

Finally, Theorem~1.1 of \cite{alisha2025detectingHeegaardFloerSolidTori} establishes $(1) \iff (4)$.
\end{proof}

\setcounter{thm}{\value{savedthm5}}
\setcounter{section}{\value{savedsec3}}

The following shows that the homology classes of knots in lens spaces that satisfy the conditions of Theorem~\ref{thm:generalizedsolidtorus} are determined by the simple knots satisfying these conditions. A knot $K$ in a rational homology sphere $M$ is called {\em Floer simple} if the knot Floer homology $\widehat{HFK}(K) = \mathbb{Z}^{|H_1(M)|}$. A simple knot in a lens space is Floer simple. 

%that admit an $S^1\times S^2$ surgery are determined by the simple knots that do so. This is Theorem~\ref{thm:homdetection}, the goal of the section. It is needed for the proof of Theorem~\ref{thm:homologyclassification}.

\newcounter{savedthm7}
\newcounter{savedsec5}
\setcounter{savedthm7}{\value{thm}}
\setcounter{savedsec5}{\value{section}}
\setcounter{thm}{\value{KmeansK'also}}
\setcounter{section}{\value{secKmeansK'also}}

\begin{thm}%\ref{thm:r>0reduciblemeansreductiblesimple}
Let $K$ be a knot in $L(p,m)$ with irreducible exterior. Assume that, for some slope $\lambda$, $K(\lambda)=(S^1 \times S^2) \# N$ where $N$ is a rational homology sphere. Then $K$ is a Floer simple knot and has a planar Seifert surface, $Q$, of slope $\lambda$. Furthermore, $p=q^2 |H_1(N)| \Delta(\lambda,\mu_K)$, where $\mu_K$ is the meridian of $K$ and where $q=|\partial Q|$. 

Let $K'$ be the simple knot in the same homology class as $K$.
Then the following hold:
\begin{enumerate}
    \item $K'$ has a planar Seifert surface, $Q'$ of some slope $\lambda'$;
    \item $|\partial Q'|=q$, and $\Delta(\lambda, \mu_{K})=\Delta(\lambda', \mu_{K'})$; and
    \item $K(\lambda')=(S^1 \times S^2) \# N'$ for a lens space $N'$ with $|H_1(N')|=|H_1(N)|$.
\end{enumerate}
\end{thm}

\begin{proof}
Let $K$ be as stated. By Theorem~\ref{thm:generalizedsolidtorus}, $X_K$ is a Heegaard Floer homology solid torus. By Theorem 27 of  \cite{hanselmanrasmussenwatson2022torusBoundary}, all fillings except along the rational longitude are L-spaces. Corollary 3.6 of \cite{rasmussenrasmussen2017floerSimpleManifolds} then shows that $K$ is a Floer simple knot in $L(p,m)$. 
%Again, by Theorem~\ref{thm:generalizedsolidtorus}, $X_K$ is a generalized solid torus.  $\spadesuit$ {\em (Why didn't Ni-Vafee use this argument, were they just worried about being fibered? Do we need to use an argument like they did instead of Corollary 3.6?)}
%By Theorem~7.8 of \cite{rasmussenrasmussen2017floerSimpleManifolds}, the rational homology solid torus $X_K$ is a generalized solid torus in the sense of that paper, in particular it is Floer simple. 
%Corollary 3.6 of the same paper shows that $K$ is a Floer simple knot. 
%The lines after Definition~7.2 in that paper show that
%Again, by Theorem~\ref{thm:generalizedsolidtorus}, $X_K$ is a generalized solid torus. 
By Theorem~\ref{thm:generalizedsolidtorus}, $K$ has a planar Seifert surface, $Q$, that minimizes the Thurston norm. Let $q=|\partial Q|$ and let $\lambda$ be the slope of $Q$.  We may write $H_1(X_{K})=\mathbb{Z} \oplus H$ where $H$ is a finite abelian group. 
%As filling $X_{K'}$ along $\mu$ has first homology $Z$, we conclude that $H$ is generated by $[\mu]$. 
%As $K'$ is a spherical braid,
%H_1(X_{K'})$ is generated by $[\mu]$ and $[\nu]$, where $\nu$ is a simple closed curve intersecting $Q$ once. 
$[\lambda]$ has order $q$ in $H_1(X_{K})$
($Q$ shows that the order of $K$ is a factor of $q$, and in a filling $X_K(\alpha)$ where $\Delta(\alpha, \lambda)=1$ the dual knot will have order $q$ by Lemma 2.13 of \cite{calegarigordon2013smallRationalGenus}). Note that $[\lambda]$ projects to $0$ in the $\mathbb{Z}$ summand. Thus $|H_1(N)|q=|H|$.
%with free factor generated by $[\nu]$ and torsion subgroup generated by $\mu$. 
Apply Lemma 3.2 of \cite{watson2012khovanovObstructions}, where $M=X_{K}$, $\lambda_M=\lambda$, $|l|=q$, and $M(\alpha)=X_{K}(\mu_K)=L(p,m)$. We conclude that  $p=c_M \Delta(\lambda,\mu_K)$ where $c_M$ is a product of $|l|=q$ and $|H|$. That is, $p= q |H| \Delta(\lambda,\mu_K)=q^2|H_1(N)|\Delta(\lambda,\mu_K)$. 

%Applying the same argument to $K$ gives that $p=q |H| |\lambda \cdot \mu|$. Thus $|H|=|H'|$ and we have that $|H_1(N')|=|H'|/q=|H_1(N)|$. Finally, note that $p=q |H| |\lambda' \cdot \mu'|=q^2 |H_1(N')||\lambda' \cdot \mu'|$.

Let $K'$ be the simple knot in the same homology class as $K$. Let $Q'$ be a Seifert surface for $K'$ that minimizes the Thurston norm.  Let $q'=|\partial Q'|$ and let $\lambda'$ be the slope of $Q'$. Note that by Lemma 2.3 of \cite{calegarigordon2013smallRationalGenus}, $q \Delta(\lambda,\mu_K)=q' \Delta(\lambda',\mu_{K'})$ is the order of $[K]=[K'] \in L(p,m)$ where $\mu_K, \mu_{K'}$ are the meridional slopes of $K,K'$. As $[K]=[K']$, we have \[ \frac{a}{\Delta(\lambda,\mu_K)}=[K]\cdot [K]=[K'] \cdot [K']=\frac{a'}{\Delta(\lambda',\mu_{K'})} \qquad (\mod 1) \] 
given by the linking pairing on $H_1(L(p,m))$ (e.g. see the discussion after Definition~2.1 in \cite{rasmussen2007lensSpaceSurgeries}) where $a$ is relatively prime to $\Delta(\lambda,\mu_K)$ and $a'$ relatively prime to $\Delta(\lambda',\mu_{K'})$. We conclude that  $\Delta(\lambda',\mu_{K'})=\Delta(\lambda,\mu_{K})$. By the above then, $q=q'$.

By Theorem 1.2 of \cite{niwu2014correctionTerms}, both $K$ and $K'$, as Floer simple knots, minimize the rational genus among all knots in their common homology class. In particular,
         \[  \frac{\max \{0,-\chi(Q')\} }{2q' \Delta(\lambda',\mu')}=\frac{\max \{0,-\chi(Q')\}}{2 \Delta(\partial Q',\mu_{K'})}=\norm{K'} = \norm{K}= \frac{\max \{0,-\chi(Q)\}}{2 \Delta(\partial Q,\mu_K)}=\frac{\max \{0,q-2\}}{2 q \Delta(\lambda,\mu_K)} \]
As $\Delta(\lambda',\mu_{K'})=\Delta(\lambda,\mu_{K})$ and $q=q'$, we conclude that $\max \{0,-\chi(Q')\}=\max \{0,q-2\}$. When $q>1$ we conclude that $\chi(Q')=2-q'$ and $Q'$ is planar. When $1=q=q'$, we conclude that $\chi(Q') \geq 0$ and that $Q'$ is also a disk.

Capping off $Q'$ with meridional disks, we see that $K(\lambda')$ contains a non-separating $2$-sphere. Thus we may write $K'(\lambda')=(S^1 \times S^2) \# N'$. As a simple knot is a tunnel number one knot, the additivity of Heegaard genus implies that $N$ is a lens space. The first part of the Theorem applied to $K'$ and $Q'$ shows that $(q')^2|H_1(N')|\Delta(\lambda',\mu_{K'})=p=q^2|H_1(N)|\Delta(\lambda,\mu_K)$. Thus $|H_1(N')|=|H_1(N')|$.
\end{proof}

\setcounter{thm}{\value{savedthm7}}
\setcounter{section}{\value{savedsec5}}

\section{Planar rational Seifert surfaces which are not longitudinal}\label{sec:nonlongitudinal}

In this section we give the classification of simple knots in lens spaces that have planar rational Seifert surfaces which are not longitudinal. 

%\begin{defin}
%Let $X$ be a $3$-manifold with $\partial X$ a torus. A {\it slope} on $X$ is the isotopy class of an essential simple closed on $\partial X$. The {\it distance}, $\Delta(\alpha,\beta)$, between two slopes $\alpha, \beta$ is defined to the absolute value of their intersection number on $\partial X$.
%\end{defin}

% $\spadesuit$ {\it (Maybe we should say that our Seifert surfaces are connected. Did this adding it to the beginning of section 2. )}

% $\spadesuit$ {\it (Should we define the winding number to be the absolute value of the winding number, so we don't have to keep saying it? Do we need to say that a torus knot is one that has winding number $>1$ in each Heegaard solid torus. Need to look up the terminology  Gordon. Did this and added it to the beginning of section 2.)}

\newcounter{savedthm6}
\newcounter{savedsec4}
\setcounter{savedthm6}{\value{thm}}
\setcounter{savedsec4}{\value{section}}
\setcounter{thm}{\value{maindelta>1}}
\setcounter{section}{\value{secdelta>1}}

\begin{thm}
[Non-longitudinal theorem]
%\label{thm:main:nonlongitudinal}
Let $K$ be a knot in $L(p,m)$ that admits a planar Seifert surface with $q>1$ boundary components and non-longitudinal slope $\lambda$.
%, with $\Delta(\lambda, \mu_K)>1$.
Then 
\begin{enumerate}
    \item $K(\lambda)=S^1 \times S^2$;
    \item $L(p,m)$ is diffeomorphic, possibly orientation-reversing, to $L(k q^2, 1+kdq)$, where $\gcd (d,q)=1$, $q>d>0$, and $k=\Delta(\lambda, \mu_K)>1$;
    \item $K$ is the $(q,q)$-torus knot in $L(k q^2, 1+kdq)$;
    \item $\pm [K]=q \mu \in H_1(L(k q^2, 1+kdq))$; and
    \item $K$ is simple.
    \end{enumerate}
\end{thm}

%\begin{thm}
%Let $K$ be a knot in $L(p,m)$ that admits a planar Seifert surface with $q>1$ boundary components and slope $\lambda$ which is not longitudinal. Then 
%\begin{enumerate}
%    \item $K(\lambda)=S^1 \times S^2$.
%    \item $L(p,m)$ is  diffeomorphic  (possibly orientation-reversing) to $L(k q^2, 1+kdq)$, where $\gcd (d,q)=1$, $q>d>0$, and $k=\Delta(\lambda, \mu_K)>1$ for the meridian, $\mu_K$, of $K$. 
%    \item $K$ is the torus knot in $L(k q^2, 1+kdq)$ where the winding number of $K$ is $q$ in each Heegaard solid torus.
%    \item $\pm [K]=q \mu \in H_1(L(k q^2, 1+kdq))$. $K$ is the simple knot in  this homology class.
%    \end{enumerate}
%\end{thm}

\setcounter{thm}{\value{savedthm6}}
\setcounter{section}{\value{savedsec4}}

\begin{remark}\label{simplicity}\hfill
\begin{enumerate}
\item Theorem~\ref{thm:main:nonlongitudinal} does not assume $K$ is simple.
\item Theorem~\ref{thm:realizationNlongitudinal}(1) shows that, conversely, the simple knot $\pm [K]=q \mu \in H_1(L(k q^2, 1+kdq))$ of Theorem~\ref{thm:main:nonlongitudinal} has a planar Seifert surface with $q$ boundary components.
\end{enumerate}
\end{remark}

%$\spadesuit$ {\em (Should just use Theorem 2.0.3 of [CGLS] here to rule out non-integral reducible surgeries below rather than Theorem~\ref{thm:r>1onlyS1xS2}. Though we might want to point that an argument of that form, using techniques of this paper, also works. Should also point out Boyer-Zhang. Also for the $S^1 \times S^2$ case, need to consider the case when the knot exterior is Seifert fibered with orbit surface a Mobius band; i.e. include that in Lemma~\ref{lem:r>1torusknot}.}

\begin{proof}
Let $K$ be as given. Note that if the exterior of $K$, $X_K$, is reducible then $K$ lies in a ball, $B^3$, in $L(p,m)$. Furthermore, $K$ bounds a planar Seifert surface in $B^3$ of slope $\lambda$---contradicting that $\lambda$ is non-longitudinal. We conclude that $X_K$ is irreducible. As $q>1$, $X_K$ is also boundary irreducible.

Now $K(\lambda)=(S^1 \times S^2) \# N$, where $N$ is a rational homology sphere. In particular $K(\lambda)$ contains a non-separating $2$-sphere which cannot be isotoped to be disjoint from the surgery-dual knot. Applying Theorem~2.0.3 of \cite{cullergordonlueckeshalen1987dehnSurgery}, with $r=\lambda$ and $s$ the meridian of $K$ ($\Delta(r,s)>1$), we see that conclusion $(iv)$ of that Theorem must apply. In particular,  $K(\lambda)=S^1 \times S^2$, proving $(1)$ of Theorem~\ref{thm:main:nonlongitudinal}. The Cyclic Surgery Theorem  of \cite{cullergordonlueckeshalen1987dehnSurgery} says then that $X_K$ is a Seifert fiber space. Let $K' \subset S^1 \times S^2$ be the dual knot of the surgery. Applying Theorem~6.9 of \cite{bakerbuck2013rationalSubtangleReplacements} (where $S^1 \times S^2$ is considered the lens space $L(0,1)$), we that $K'$ is the $(q,d)$-torus knot in $S^1 \times S^2$ and $L(p,m)$ is $L(0+nq^2,1+ndq)$ where $|n|=k=\Delta(\lambda, \mu_K)>1$ for the meridian, $\mu_K$, of $K$. This gives conclusion $(2)$ of Theorem~~\ref{thm:main:nonlongitudinal} (as a torus knot $\gcd (d,q)=1$, and up to diffeomorphism of the lens space, $n=k$, and $d$ can be taken positive and less than $q$). From the surgery description of $L(nq^2, 1+ndq)$ given in Theorem 6.9, we see that $K$ is the $(q,q)$-torus knot in $L(kq^2,1+kdq)$; furthermore, from the resulting diagram of $L(nq^2,1+ndq)$, we see that $K$ is simple.  Hence $(3),(4)$ follow as well.
For $(2)$--$(4)$ also see \cite{darcysumners2000rationalTangleDistances}.
%Lemma~\ref{lem:r>1torusknot} where $\Delta(\lambda,\mu)>1$. Note that by Lemma~\ref{lem:r>1torusknot}(3), the torsion of $[K]$ is $\Delta(\lambda,\mu)q$, hence the $q$ there is the number of boundary components of $Q$.
\end{proof}

\section{Realizing the surgeries described in  
Theorems~\ref{thm:main:nonlongitudinal} and \ref{thm:mainNextended}}\label{sec:Appendix}

%\begin{thm}\label{thm:realizationNlongitudinal}
%For positive integers $q,k,d$, with $d<q,(d,q)=1$, there is a simple knot $K \subset L(kq^2,kdq + 1)$ that admits a longitudinal surgery to $(S^1 \times S^2) \# L(k,1)$. Furthermore, these occur in each of the following classes.
%\begin{enumerate}
 %   \item $\pm [K]=kdq \mu \in H_1(kq^2,kdq + 1)$; or
 %   \item $\pm [K]=kq \mu \in H_1(kq^2,kdq + 1)$. $K$ is an $k$-fold cable on a companion, $C$, that lies on the Heegaard 2-torus of $L(p,m)$; $C$ has winding number $q$ in both Heegaard solid tori. The slope $\lambda$ on $\partial X_K$ is the Seifert fiber of the cabling.
%\end{enumerate}
%\end{thm}

In this section, we describe those simple knots upon which a surgery has a non-separating $2$-sphere but are not described in \cite{baker2016knotsS1xS2}---because the surgeries are non-prime or non-longitudinal. As simple knots are strongly invertible, we describe these surgeries via tangle replacement (similar to \cite{baker2016knotsS1xS2}) in a double branched cover description of the lens space. These examples realize the surgeries described in  
Theorems~\ref{thm:main:nonlongitudinal} and \ref{thm:mainNextended}. We are grateful to Ken Baker for bringing these examples to our attention and allowing us to include them here.

\begin{thm}\label{thm:realizationNlongitudinal}
Fix positive integers $q,k,d$, with $d<q,(d,q)=1$. Let $K \subset L(k q^2, kdq+1)$ be a simple knot with rational longitude $\lambda$.

\begin{enumerate} 
\item
If $\pm [K]=q \mu \in H_1(L(k q^2, 1+kdq))$, then 
\begin{itemize}
    \item 
    $K(\lambda)=S^1 \times S^2$,
    \item
    $\Delta(\lambda,\mu_K)=k$, and
    \item 
    $K$ is the $(q,q)$-torus knot in $L(k q^2, kdq+1)$.
\end{itemize}
\item
If $\pm [K]=kq \mu \in H_1(L(kq^2,kdq + 1))$, then 
\begin{itemize}
    \item 
    $K(\lambda)=(S^1 \times S^2) \# L(k,1)$,
    \item
    $\Delta(\lambda,\mu_K)=1$,
    \item 
    the slope $\lambda$ on $\partial X_K$ is the Seifert fiber of the cabling, and
    \item 
    $K$ is a $k$-fold cable on the $(q,q)$-torus knot.
\end{itemize}
\item
If $\pm [K]=kdq \mu \in H_1(kq^2,kdq + 1)$, then
\begin{itemize}
    \item 
    $K(\lambda)=(S^1 \times S^2) \# L(k,1)$ and
    \item 
     $\Delta(\lambda,\mu_K)=1$.
\end{itemize}

\end{enumerate}
\noindent
In each case, $K$ has a planar rational Seifert with $q$ boundary components, each of slope $\lambda$.
\end{thm}

\begin{remark}\label{rmk:hyperbolicity}
The knots of Theorem~\ref{thm:realizationNlongitudinal}$(1)$, $(2)$ are not hyperbolic. However, Theorem 2 of \cite{bakercabling} shows that for infinitely many values of $q,d$ the knot $K$ of Theorem~\ref{thm:realizationNlongitudinal}$(3)$ (which is $-K^{q,d}_{0,k}$ of \cite{bakercabling}) is hyperbolic.
\end{remark}

\begin{proof}
Find positive integers $b_1, \dots, b_l$ for some odd, positive $l$ such that $\frac{q}{d}=[b_1, \dots, b_l]^-$.
By Lemma 3.3 of \cite{baker2016knotsS1xS2}, 
\[ \frac{kq^2}{kdq+1}=[b_1, \dots, b_l, k,-b_l, \dots, -b_1]^{-}.\]
Let $K(kq^2,kdq+1)$ be the $2$-bridge knot represented by the sequence $L(-b_1, \dots, -b_l, -k, b_l,\allowbreak \dots,\allowbreak b_1)$ in the notation of section~1.6.3 of \cite{baker2016knotsS1xS2}, and in particular as described in Figure~5 (right-hand side) of that section. As described there, the double-branched cover of $K(kq^2,kdq+1)$ is $L(kq^2,kdq+1)$. See the far right-hand side of Figure~7 in section~2.1 of \cite{baker2016knotsS1xS2} for the case that $k=4$.

\begin{figure}
    {
    \includegraphics[width=1in]{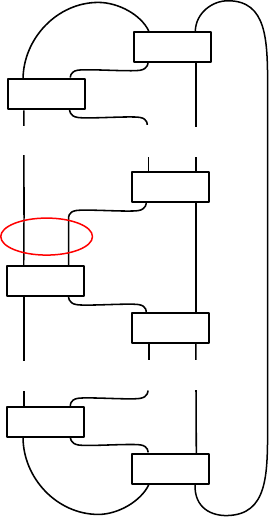}
    \put(-31,125){\tiny $-b_1$}
    \put(-66,112){\tiny $-b_2$}
    \put(-27,99){\tiny $\vdots$}
    \put(-62,99){\tiny $\vdots$}
    \put(-33,87){\tiny $-b_\ell$}
    \put(-65,61){\tiny $-k$}
    \put(-30,49){\tiny $b_\ell$}
    \put(-27,35){\tiny $\vdots$}
    \put(-62,35){\tiny $\vdots$}
    \put(-65,24){\tiny $b_{2}$}
    \put(-30,11){\tiny $b_{1}$}
    }
    \qquad
    {
    \includegraphics[width=1in]{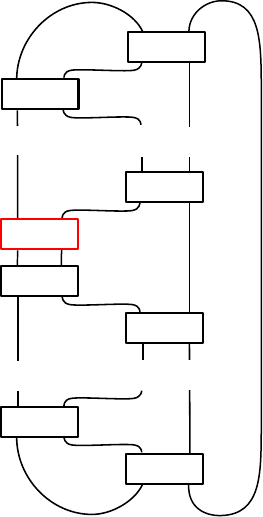}
    \put(-32,128){\tiny $-b_1$}
    \put(-67,114){\tiny $-b_2$}
    \put(-27,99){\tiny $\vdots$}
    \put(-62,99){\tiny $\vdots$}
    \put(-33,89){\tiny $-b_\ell$}
    \put(-65,76){\tiny\color{red}$k$}
    \put(-66,63){\tiny $-k$}
    \put(-30,50){\tiny $b_\ell$}
    \put(-27,35){\tiny $\vdots$}
    \put(-62,35){\tiny $\vdots$}
    \put(-65,24){\tiny $b_{2}$}
    \put(-30,11){\tiny $b_{1}$}
    }
    \caption{A tangle replacement corresponding to a non-longitudinal surgery of $K \subset L(kq^2,kdq+1)$ to $S^1 \times S^2$.}
    \label{fig:nonlong}
\end{figure}

Consider the tangle replacement pictured in Figure~\ref{fig:nonlong} (see also \cite[Theorem 1(3)]{darcysumners2000rationalTangleDistances} or \cite[Theorem 5.2]{bakerbuck2013rationalSubtangleReplacements}). The result is the unlink. Via the double branched cover, the tangle ball becomes the neighborhood of a knot $K$ in $L(kq^2,kdq+1)$ that admits a surgery $K(\lambda) = S^1 \times S^2$ with $\Delta(\lambda, \mu_K)=k$ as claimed in $(1)$. A Conway sphere lying just above the tangle ball corresponds to a Heegaard torus in $L(kq^2,kdq+1)$. From this, and the definition of the $b_i$, one sees that $K$ is a torus knot with winding number $q$ on each side of the Heegaard torus. From this we conclude that $\pm [K]=q \mu$. As a torus knot, $K$ is a $1$-bridge (a $(1,1)-$) knot in the lens space. By Theorem~\ref{thm:r>0reduciblemeansreductiblesimple}, $K$ is Floer simple.
%($\spadesuit$ {\em (Or maybe just use Proposition 1.2 of Ni-Vafee.)})
 Theorem~\ref{thm:hedden} verifies that $K$ is simple. This is $(1)$ of Theorem~\ref{thm:realizationNlongitudinal}. 

%$\spadesuit$ {\em (Maybe an alternative argument for this case. Let $K_1$ be the knot in $L(q^2,dq+1)$ gotten by taking $k=1$ in the above picture. Let $K$ be knot as above. Then $K$ is the dual knot to doing $(k-1)/1$-surgery on $K_1$ (using the framing coming from the Heegaard torus). By \cite{baker2016knotsS1xS2}, $K_1$ is the simple knot in class $q \mu$ in $L(q^2, dq+1)$. One checks that the simple knot in this class is a torus knot (i.e. the two constituent arcs intersection only in their endpoints). Then we get a standard Heegaard presentation of $L(kq^2,kdq+1)$ by twisting $k-1$ times by parallel copy of $K_1$ on the Heegaard torus. The two arcs forming $K_1$ remain in this Heegaard presentation of $L(kq^2,kdq+1)$, which is isotopic to $K$ as the dual of the surgery. That is, $K$ is the union of two arcs in a standard presentation of $L(kq^2,kdq+1)$; hence simple.)} 

The banding shown on the left side of Figure~\ref{fig:banding} changes $K(kq^2,kdq+1)$ to the split link, $S$, consisting of the unknot and the $(2,-k)$-torus knot. In the double cover of $S^3$ branched over $K(kq^2,kdq+1)$, a ball neighborhood of this band lifts to a neighborhood of a knot $K$. The banding operation then corresponds to a longitudinal Dehn surgery along $K \subset L(kq^2,kdq+1)$, yielding the double cover of $S^3$ branched over $S$, namely, to $(S^1 \times S^2) \# L(k,1)$. Following the notation in section~3.2 of \cite{baker2016knotsS1xS2}, replacing $L(p,q)$ with $L(kq^2, kdq+1)$, we obtain $\pm[K] = \mu_1 - \mu_{2l+1} = (kdq+1)\mu_{2l+1} - \mu_{2l+1} = kdq\mu_{2l+1} = kdq\mu$. In  \cite{bakercabling}, Baker described a 4-parameter family of knots $K^{p,q}_{r,s}$ that surger from a lens spaces to the connected sum of lens spaces $L(r,1) \# L(s,1)$. $K$ is in that family (with $r=0,s=-k$) and the surgery, corresponding to the banding, is the one described there (compare Figure 8 of \cite{bakercabling} with the left-hand side of Figure~\ref{fig:banding}). In \cite{bakercabling} (see Lemma 17(1) and Remark 18), as part of this family with $r=0$, $K$ is shown to be simple in $L(kq^2,kdq+1)$. Thus $K$ is the knot described in Theorem~\ref{thm:realizationNlongitudinal}$(3)$.
%$\spadesuit$ {\em Alternatively, to see that $K$ is simple, we need only show that $K$ is a $(1,1)$-knot and then argue as in the preceding paragraph or use Lemma 17(1) in \cite{bakercabling}. Consider a Conway sphere given by a horizontal sphere that intersects the branch locus in four points and the band in a single horizontal arc. On each side of the Conway sphere we can undo each of the twist boxes by an isotopy disjoint from the band. We then see a half disk cobounded by the band and the branch locus. This half disk lifts in the double branch cover to a trivializing disk of the part of $K$ on that side of the Heegaard torus. That is, $K$ can be written as the union of two trivial arcs, one on each side of the Heegaard torus, and $K$ is a  $(1,1)$-knot in the lens space.}
% The proof of Theorem 1.7(1) of section~3.2 for $[K_{GOFK}]$, applied to the above $K_1$ shows that $\pm [K]=kdq \mu \in H_1(L(kq^2,kdq+1))$ as claimed.

\begin{figure}
    {
    \includegraphics[width=1.2in]{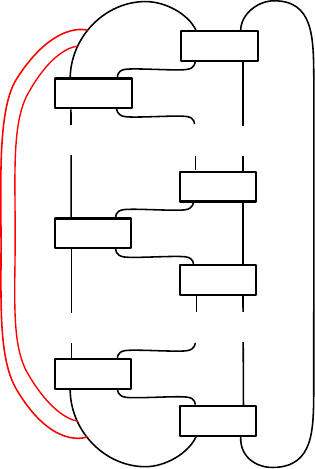}
    \put(-32,114){\tiny $-b_1$}
    \put(-68,101){\tiny $-b_2$}
    \put(-27,87){\tiny $\vdots$}
    \put(-62,87){\tiny $\vdots$}
    \put(-33,76){\tiny $-b_\ell$}
    \put(-66,63){\tiny $-k$}
    \put(-30,50){\tiny $b_\ell$}
    \put(-27,35){\tiny $\vdots$}
    \put(-62,35){\tiny $\vdots$}
    \put(-65,24){\tiny $b_{2}$}
    \put(-30,11){\tiny $b_{1}$}
    }
    \qquad
    {
    \includegraphics[width=1in]{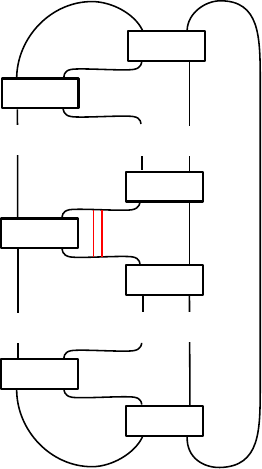}
    \put(-32,115){\tiny $-b_1$}
    \put(-68,102){\tiny $-b_2$}
    \put(-27,87){\tiny $\vdots$}
    \put(-62,87){\tiny $\vdots$}
    \put(-33,76){\tiny $-b_\ell$}
    \put(-66,63){\tiny $-k$}
    \put(-30,50){\tiny $b_\ell$}
    \put(-27,35){\tiny $\vdots$}
    \put(-62,35){\tiny $\vdots$}
    \put(-65,24){\tiny $b_{2}$}
    \put(-30,11){\tiny $b_{1}$}
    }
    \caption{Bandings corresponding to longitudinal surgeries from $L(kq^2,kdq + 1)$ to $(S^1 \times S^2) \# L(k,1)$.}
    \label{fig:banding}
\end{figure}

The banding pictured on the right side of Figure~\ref{fig:banding} changes $K(kq^2,kdq+1)$ to the split link $S$ described in the preceding paragraph. Again, taking double-branched covers, we see that the banding corresponds to a longitudinal surgery on a knot $K \subset L(kq^2,kdq+1)$ yielding $(S^1 \times S^2) \# L(k,1)$. Furthermore, a ball engulfing the $-k$ twists plus the pictured band lifts to the neighborhood of the knot $K'$ in $L(kq^2,kdq+1)$ described in Figure~\ref{fig:nonlong}. One checks that $K$ is the $(1,\pm k)$-cable of $K'$ where the longitude for the neighborhood of $K'$ in this cabling is given by the framing $K'$ from the Heegaard torus. This shows that $K$ is a $(1,1)$-knot in the lens space. By the preceding arguments (e.g. Lemma 17(1) in \cite{bakercabling}), this shows $K$ is simple. Furthermore, $\pm [K]=k[K']=kq \mu$, identifying this as the knot of Theorem~\ref{thm:realizationNlongitudinal}$(2)$. Note that by taking the ball engulfing the $-k$ twists so that the banding lies on its boundary, one sees that the slope of the surgery, $\lambda$, is the Seifert fiber of the cabling as claimed. 
%which lies on a a Heegaard torus. One checks that $C$ has winding number $q$ in each Heegaard solid torus. $K_2$ is then seen as the $(1,k)$-cable of $C$, with longitude on $C$ for the cabling given by the framing from the Heegaard torus. This shows that $K_2$ is simple. Thus $K_2$ is of the form claimed in part $(2)$ of the Theorem.

%We quoted at the end of the introduction Hedden's result that (1,1)-knots that are Floer simple are simple, so 
%by Theorem~\ref{thm:generalizedsolidtorus} 
%if we show $K(\lambda)=(S^1 \times S^2) \# N$ and $K$ is (1,1) then Theorem~\ref{thm:r>0reduciblemeansreductiblesimple} shows $K$ is Floer simple hence simple. So this can be used to show simplicity. 
% $\blacklozenge$ {\em (What about this: ``In Figure~\ref{fig:banding}, consider a Conway sphere given by a horizontal sphere that intersects the branch locus in four points and the knot once. In the quotient, the knot $K_1$ splits into two arcs $\delta_1\cup \delta_2$, one lying on each side of the Conway sphere. On each side, after resolving all twist boxes fixing the arc $\delta_i$, this arc cobounds a half disk with part of the branch locus. Thus, each $\delta_i$ lifts to a trivial arc on its side of the Heegaard torus corresponding to the Conway sphere, showing $K$ as a $(1,1)$-knot in the lens space. }

By Theorem~\ref{thm:generalizedsolidtorus}, each of the knots $K$ above have planar Seifert surfaces of the surgery slope $\lambda$ with, say, $c_K$ boundary components. In Case $(1)$ of Theorem~\ref{thm:realizationNlongitudinal}, $[K]$ has order $kq$. By Lemma 2.13 of \cite{calegarigordon2013smallRationalGenus}, and $\Delta(\lambda, \mu_K)=k$, we conclude that $c_K=q$ as claimed. In Cases $(2),(3)$, $[K]$ has order $q$ and $\lambda$ is longitudinal. By Lemma 2.13 of \cite{calegarigordon2013smallRationalGenus}, we conclude that $c_K=q$ as claimed. 
\end{proof}

\bibliographystyle{amsalpha2}

\bibliography{References}

\end{document}